\documentclass[english]{article}

\usepackage{geometry}
\usepackage[T1]{fontenc}
\usepackage{amsmath,amssymb,amsthm,graphicx,bm,dsfont}
\usepackage{epsfig, fontawesome}
\usepackage[authoryear]{natbib} 
\usepackage{psfrag}
\usepackage{enumerate} \usepackage{setspace}
\usepackage{float} 
\usepackage{color}
\usepackage{array}
\usepackage{microtype, wasysym, comment,mathtools}
\usepackage{subfigure}
\usepackage{booktabs} 
\usepackage{enumitem}
\usepackage[colorlinks=true,allcolors=blue]{hyperref}%
\usepackage{graphicx} 
\usepackage{epstopdf}
\usepackage[ruled,vlined]{algorithm2e}
\usepackage{algorithmic}
\definecolor{yuting}{RGB}{255,69,0}
\definecolor{gen}{RGB}{199,21,133}

\DeclareMathOperator{\ind}{\mathds{1}}

\newcommand{\vstar}{v^\star}
\newcommand{\lt}{\left}
\newcommand{\rt}{\right}
\newcommand{\dx}{\mathrm{d} x}

\newcommand{\mymid}{\,|\,}
\newcommand{\myE}{\mathbb{E}}

\newcommand{\hatxi}{\widehat{\xi}}
\newcommand{\hatzeta}{\widehat{\zeta}}
\newcommand{\hatalpha}{\widehat{\alpha}}
\newcommand{\hatgamma}{\widehat{\gamma}}
\newcommand{\hatbeta}{\widehat{\beta}}

\newcommand{\hats}{\widehat{s}}

\newcommand{\overalpha}{\overline{\alpha}}
\newcommand{\overgamma}{\overline{\gamma}}
\newcommand{\poly}{\mathsf{poly}}
\newcommand{\thetatilde}{\widetilde{\theta}}
\newcommand{\lone}[1]{\|#1\|_1}

\newtheorem{assumption}{\textbf{Assumption}}
\newtheorem{claim}{\textbf{Claim}}
\newtheorem{remark}{\textbf{Remark}}

\theoremstyle{plain}

\newtheorem{theo}{Theorem}[section]

\newtheorem{lem}{Lemma}[section]
\newtheorem{prop}{Proposition}[section]
\newtheorem{cor}{Corollary}[section]

\theoremstyle{definition} 

\newtheorem{nota}{Notation}[section]
\newtheorem{de}{Definition}[section]
\newtheorem{exa}{Example}[section]
\newtheorem{as}{Assumption}[section]
\newtheorem{alg}{Algorithm}[section]

\newcommand{\btheo}{\begin{theo}}
\newcommand{\bde}{\begin{de}}
\newcommand{\ble}{\begin{lem}}
\newcommand{\bpr}{\begin{prop}}
\newcommand{\bno}{\begin{nota}}
\newcommand{\bex}{\begin{exa}}
\newcommand{\bcor}{\begin{cor}}
\newcommand{\spro}{\begin{proof}}
\newcommand{\bas}{\begin{as}}
\newcommand{\balg}{\begin{alg}}

\newcommand{\etheo}{\end{theo}}
\newcommand{\ede}{\end{de}}
\newcommand{\ele}{\end{lem}}
\newcommand{\epr}{\end{prop}}
\newcommand{\eno}{\end{nota}}
\newcommand{\eex}{\end{exa}}
\newcommand{\ecor}{\end{cor}}
\newcommand{\fpro}{\end{proof}}
\newcommand{\eas}{\end{as}}
\newcommand{\ealg}{\end{alg}}

\theoremstyle{plain}

\newtheorem{theos}{Theorem}
\newtheorem{props}{Proposition}
\newtheorem{lems}{Lemma}
\newtheorem{cors}{Corollary}

\theoremstyle{definition}
\newtheorem{exas}{Example}
\newtheorem{algs}{Algorithm}
\newtheorem{asss}{Assumption}
\newtheorem{defns}{Definition}

\newcommand{\btheos}{\begin{theos}}
\newcommand{\etheos}{\end{theos}}
\newcommand{\bprops}{\begin{props}}
\newcommand{\eprops}{\end{props}}
\newcommand{\bdes}{\begin{defns}}
\newcommand{\edes}{\end{defns}}
\newcommand{\blems}{\begin{lems}}
\newcommand{\elems}{\end{lems}}
\newcommand{\bcors}{\begin{cors}}
\newcommand{\ecors}{\end{cors}}
\newcommand{\bexs}{\begin{exas}}
\newcommand{\eexs}{\end{exas}}
\newcommand{\balgs}{\begin{algs}}
\newcommand{\ealgs}{\end{algs}}
\newcommand{\bass}{\begin{asss}}
\newcommand{\eass}{\end{asss}}

\newcommand{\ltwo}[1]{\|#1\|_2}

\newcommand{\real}{\ensuremath{\mathbb{R}}}

\newcommand{\thetastar}{\ensuremath{{\theta}^\star}}

\newcommand{\inprod}[2]{\ensuremath{\langle #1 , \, #2 \rangle}}

\newcommand{\thetahat}{\ensuremath{\widehat{\theta}}}

\newcommand{\mprob}{\ensuremath{\mathbb{P}}}

\newcommand{\defn}{\coloneqq}
\newcommand{\argmin}{\arg\!\min}

\newcommand{\Exs}{\ensuremath{\mathbb{E}}}

\newcommand{\sign}{\ensuremath{\mbox{sign}}}

\long\def\comment#1{}

\newcommand{\HACKPROOF}{\begin{proof}}
\newcommand{\HACKENDPROOF}{\end{proof}}

\newcommand{\var}{\ensuremath{\operatorname{var}}}

\newcommand{\Ball}{\ensuremath{\mathbb{B}}}

\newlength{\widebarargwidth}
\newlength{\widebarargheight}
\newlength{\widebarargdepth}

\makeatletter
\long\def\@makecaption#1#2{
        \vskip 0.8ex
        \setbox\@tempboxa\hbox{\small {\bf #1:} #2}
        \parindent 1.5em  
        \dimen0=\hsize
        \advance\dimen0 by -3em
        \ifdim \wd\@tempboxa >\dimen0
                \hbox to \hsize{
                        \parindent 0em
                        \hfil 
                        \parbox{\dimen0}{\def\baselinestretch{0.96}\small
                                {\bf #1.} #2
                                } 
                        \hfil}
        \else \hbox to \hsize{\hfil \box\@tempboxa \hfil}
        \fi
        }
\makeatother

\allowdisplaybreaks

\begin{document}

\title{
A non-asymptotic distributional theory of \\ approximate message passing for sparse and robust regression
}

\author{
Gen Li \thanks{Department of Statistics, The Chinese University of Hong Kong, Hong Kong.} 
\and Yuting Wei\thanks{Department of Statistics and Data Science, the Wharton School, 
	University of Pennsylvania, Philadelphia, PA.}\\[0.2cm]
}

\maketitle

\begin{abstract}

Characterizing the distribution of high-dimensional statistical estimators is a challenging task, due to the breakdown of classical asymptotic theory in high dimension. This paper makes progress towards this by developing non-asymptotic distributional characterizations for approximate message passing (AMP) --- a family of iterative algorithms that prove effective as both fast estimators and powerful theoretical machinery --- for both sparse and robust regression. Prior AMP theory, which focused on high-dimensional asymptotics for the most part, failed to describe the behavior of AMP when the number of iterations exceeds $o\big({\log n}/{\log \log n}\big)$ (with $n$ the sample size). We establish the first finite-sample non-asymptotic distributional theory of AMP for both sparse and robust regression that accommodates a polynomial number of iterations. Our results derive approximate accuracy of Gaussian approximation of the AMP iterates, which improves upon all prior results and implies enhanced distributional characterizations for both optimally tuned Lasso and robust M-estimator. 

\end{abstract}


\setcounter{tocdepth}{2}
\tableofcontents

\section{Introduction}

Determining the distributions of the estimators of interest plays a pivotal role in addressing fundamental questions in uncertainty quantification, hypothesis testing, and risk prediction, among others.  
In classical large-sample theory  \citep{fisher1922mathematical,le2012asymptotic,van2000asymptotic}, this is often achieved by pinning down the limiting distribution, such as asymptotic normality, of the estimators of interest in the limit as the sample size $n$ approaches infinity with the problem dimension held fixed. 
Nevertheless, such large-sample theory often breaks down in modern high-dimensional settings where the ambient dimension $p$ of the unknowns is large as well (e.g., comparable to the sample size), 
due to prevalent issues such as non-negligible bias and inflated variance \citep{el2013robust,donoho2015variance,donoho2016high,sur2019modern,sur2019likelihood}.  
These issues have motivated a recent wave of research activities proposing new paradigms and analyses that enable tractable distributional characterizations in high dimension (see e.g.~\cite{zhang2014confidence,van2014asymptotically,javanmard2014hypothesis,javanmard2018debiasing,ren2015asymptotic,bellec2022biasing,bellec2023debiasing,bellec2022asymptotic,chen2019inference,celentano2020lasso,cai2022uncertainty,xia2021statistical,celentano2021cad,yan2021inference} and the references therein). 
Focusing on linear models, the present paper aims to make progress towards understanding the distribution of a powerful family of statistical estimators, called approximate message passing (AMP) \citep{donoho2009message,feng2021unifying}, 
that are among the most effective when tackling high-dimensional problems.

\subsection{Sparse and robust regression in high dimension}

The current paper is focused on the prototypical problem of estimating a set of unknown parameters in a linear model. 
Given a design matrix $X \in \mathbb{R}^{n \times p}$ (with $X_1,\dots,X_n$ denoting the rows of $X$),  
the classical linear regression model takes the form of
\begin{align}
\label{eqn:linear}
	y = X\theta^{\star} + \varepsilon,
\end{align}
where $y = [y_i]_{1\leq i\leq n} \in \real^{n}$ stands for the observed data vector, 
$\thetastar = [\thetastar_i]_{1\leq i\leq p} \in \real^{p}$ represents some unknown signal of interest, 
and $\varepsilon = [\varepsilon_i]_{1\leq i\leq n} \in \real^{n}$ indicates independent random noise contaminating the observations. 
The aim is to reconstruct the unknown object $\thetastar$ based on $(y,X)$. 
In practice, it is common to encounter situations where either the signal coefficients or the noise distributions exhibit certain structural properties (e.g., sparsity, group sparsity, heavy tails)
 that are known to scientists {\em a priori} (e.g., \citet{tibshirani1996regression,chen2001atomic,donoho2001uncertainty,mitchell1988bayesian,fan2001variable,zou2005regularization,yuan2006model,candes2007dantzig,donoho2015variance,bogdan2015slope,sun2020adaptive,bu2020algorithmic,buhlmann2011statistics,hastie2015statistical,fan2020statistical}).  While numerous instances of linear regression have been studied across various contexts, 
we single out two concrete settings that will serve as a guiding thread throughout this paper. 
\begin{itemize}
	\item \emph{Sparse regression.} 
	Imagine that the signal of interest $\thetastar \in \real^{p}$ in \eqref{eqn:linear} is sparse, namely, 
	\begin{equation}
		\thetastar\text{ is } k\text{-sparse} 
	\end{equation}
	with the sparsity level $k$ much smaller than the ambient dimension $p$. 
	The widespread applicability of sparse linear regression across diverse data science applications, 
		including but not limited to medical imaging, genomics, geophysics, and signal processing,  
		has inspired substantial research activities dedicated to the design and analysis of sparse statistical estimators  (e.g., \citet{tibshirani1996regression,donoho2006compressed,candes2006stable,fan2001variable,zou2005regularization,yuan2006model}).

	\item \emph{Robust regression.}	
	While a dominant fraction of linear regression works operates upon commonly encountered noise assumptions like Gaussians, 
		the prevalence of (sparse) outliers in reality has incentivized research into the ``robustness'' aspect of regression. 
		Originally proposed by  \citet{huber1964robust,huber1973robust}, 
		the gross-errors contamination model assumes that each noise component $\varepsilon_i$ is independently drawn from the following distribution: 
		\begin{equation}
			\varepsilon_i \overset{\text{i.i.d.}}{\sim} (1-\epsilon_{H}) \mathcal{N}(0, \sigma^2) + \epsilon_{H} H,
			\label{eqn:robust-error-intro}
		\end{equation}
		where $H$ denotes some (unknown) contaminating distribution, and $\epsilon_H\in (0,1)$ represents the contamination fraction. 
		In other words, the observed data might contain a fraction $\epsilon_H$ of abnormal data that deviate from situations under Gaussian noise.
		Statistical performances of robust estimators tailored for this model are developed subsequently in  
		\cite{hampel1974influence,bickel1975one,maronna2019robust,fan2014robust,loh2017statistical,sun2020adaptive,el2013robust,donoho2016high,donoho2015variance,el2018impact,lei2018asymptotics}, among others.

\end{itemize}
The current paper concentrates on a particularly challenging scenario called the proportional-growth regime, 
where the number of observations $n$ and the ambient dimension $p$ are on the same order. 
In the sparse regression case, our focus is on the linear sparsity regime, where the sparsity level $k$ is on the same order as $p$ and $n$.  
A family of algorithm that are well-suited for this challenging scenario is called approximate message passing (AMP), 
which we shall elaborate on next.

\subsection{Approximate message passing (AMP)}

Approximate message passing was originally developed in the context of compressed sensing \citep{donoho2009message,bayati2011dynamics} as a family of low-complexity iterative algorithms, and has now been widely recognized as a powerful machinery to assist in understanding the performances of a broad class of statistical procedures, especially in scenarios with low signal-to-noise ratios (SNRs). Its applications span linear and generalized linear models \citep{bayati2011lasso,rangan2011generalized,schniter2014compressive,donoho2016high,sur2019likelihood,barbier2019optimal,mondelli2022approximate,li2021minimum,fan2022approximate,zhang2023spectral,li2023random,celentano2023mean}, low-rank matrix estimation \citep{rangan2012iterative,montanari2021estimation,deshpande2014information,mondelli2021pca,zhong2021approximate,celentano2023local,li2022non}, community detection \citep{deshpande2017asymptotic,ma2021community,wang2022universality}, and more recently, Bayesian sampling from diffusion processes \citep{montanari2023posterior}. The interested reader is refeerred to \citet{feng2021unifying} for a comprehensive overview of AMP and its broad applications.


When applied to the linear model \eqref{eqn:linear} with i.i.d.~Gaussian design, 
the AMP procedure typically maintains running estimates $\{\theta_t\}_{t\geq 0}\subset \mathbb{R}^p$  of the signal $\thetastar$ as well as adjusted residuals $\{r_t\}_{t\geq 0}\subseteq \mathbb{R}^n$. More specifically, 
letting $f_t: \mathbb{R}\rightarrow \mathbb{R}$ and $g_t: \mathbb{R}^n\rightarrow \mathbb{R}^n$ be some properly chosen denoising functions, 
AMP executes the following update rule in the $t$-th iteration:  
\begin{subequations}
\label{eqn:AMP}
\begin{align}
\label{eqn:AMP-rt}
	r_{t} &= y - Xf_t(\theta_{t}) + \big\langle f_t^{\prime}(\theta_{t})\big\rangle \big(\big\langle g_{t-1}^{\prime}(r_{t-1}) \big\rangle\big)^{-1}g_{t-1}(r_{t-1}), \\
\label{eqn:AMP-thetat}
	\theta_{t+1} &=  \big(\big\langle g_t^{\prime}(r_{t})\big\rangle\big)^{-1} X^{\top} g_t(r_{t}) +  f_{t}(\theta_{t}), 
\end{align}
\end{subequations}
where $f_t$, $g_t$, and their derivatives ($f_t'$ and $g_t'$) are applied component-wise to the vector argument, 
and for every integer $m>0$, we adopt a slightly unconventional piece of notation\footnote{Note that here, instead of using a ${1}/{m}$ scaling, we use a ${1}/{n}$ normalization instead. Due to this slight different scaling, there is no need for an additional 
multiplicative factor ${1}/{\delta} \defn p/n$ in the last term of \eqref{eqn:AMP-rt} as in \cite{donoho2009message} or \cite{donoho2016high}.}
\begin{align}
	\langle x \rangle  \coloneqq \frac{1}{n} \sum^m_{i=1} x_{i}, \qquad \text{for } x \in \real^{m}.
\end{align}
The algorithm is initialized at
\begin{align*}
	f_1(\theta_{1}) = 0 \in \real^{p}, \qquad g_0(r_0) = 0, 
\end{align*}
and quantities associated with non-positive iteration numbers are all set to be zero. 
When instantiated to the two concrete settings described above,  
the following denoising functions have been suggested in past works.  
\begin{itemize}
	\item {\em AMP for sparse regression.} When tackling the sparse regression setting, AMP adopts the denoising functions 
\begin{align}
\label{eqn:lasso-func-intro}
	g_t(x) = x \qquad \text{and} \qquad f_t(x) = \sign(x)\big(|x| - \tau_t\big)_{+} =: \mathsf{ST}_{\tau_t}(x)
\end{align}
for some properly selected threshold $\tau_t$ (to be made precisely in Section~\ref{sec:sparse}). 
Notably, $f_t$ is chosen to be the soft-thresholding function in order to promote sparsity. 
As demonstrated in \citet{bayati2011lasso}, the AMP procedure \eqref{eqn:AMP} with the choices \eqref{eqn:lasso-func-intro} 
serves as a fast algorithm to solve, and help assess the risk of, the Lasso estimator in the most sample-starved regime.

\item {\em AMP for robust regression.} When it comes to the robust regression problem, suppose first that we are given a convex loss function 
	$\rho: \mathbb{R} \rightarrow \mathbb{R}_{\geq 0}$. 
	The denoising functions for AMP can then be selected as   
\begin{align}
\label{eqn:robust-func-intro}
	g_t(x) = \frac{n}{p}\Psi(z,b_t) \qquad \text{and}\qquad f_t(x) = x,
\end{align}
where $\Psi$ is defined such that
\begin{align}
	\Psi(z,b) = \rho_b'(z) \qquad \text{with } \rho_b(z) \defn \min_{x} \Big\{\rho(x) + \frac{1}{2b}(x-z)^2 \Big\}.
\end{align}
Here, $\rho_{b}$ (with some $b>0$) can be viewed as a regularized variant of $\rho$, 
		and the regularization parameter $b_{t}$ will be made precisely momentarily (see Section~\ref{sec:robust-Huber}). 
As we shall elaborate on in Section~\ref{sec:robust-Huber}, 
		the AMP procedure \eqref{eqn:AMP} with the choices \eqref{eqn:robust-func-intro} 
		is a rapid method for solving the M-estimator with loss function $\rho$ \citep{donoho2016high}.

\end{itemize}

In addition to their computational efficiency, the aforementioned AMP algorithms often admit exact asymptotic characterizations, in the sense that their risk and dynamics in the high-dimensional asymptotics can  often be characterized in a precise manner. Consequently, AMP has now been widely recognized as both a family of standalone fast statistical estimators and a power machineary for analyzing other optimization-based statistical estimators (e.g., the M-estimator and the Lasso).

\subsection{From asymptotics to non-asymptotics}

\paragraph{Exact asymptotics and state evolution.} 
Recent years have witnessed a flurry of activity in developing theoretical tools towards demystifying the efficacy of AMP. 
More concretely, existing AMP theory reveals that: the behavior of each iteration of AMP, in the high-dimensional asymptotics (i.e., with $n,p$ approaching infinity and $t$ held fixed), 
can often be characterized by a low-dimensional recursive formula dubbed as the \emph{state evolution (SE)}. 
Informally, under i.i.d.~Gaussian design (i.e., $X_{ij} \overset{\mathrm{i.i.d.}}{\sim} \mathcal{N}(0,1/n)$) as well as some other mild conditions, 
past theory introduced the following SE recursion 
\begin{align}
	\alpha^{\star 2}_t = \Exs \big[G_t^2(\gamma^\star_t Z, W)\big], \qquad \gamma^{\star 2}_t = \Exs \big[F^2_t(\alpha^\star_{t-1}Z, V) \big], \qquad t \geq 1,
	\label{eq:SE-informal}
\end{align}
where $G_t$ (resp.~$F_t$) is some function depending on $g_t$ (resp.~$f_t$) to be specified shortly (see Section~\ref{sec:decompose}). 
Here, $Z,V,W$ are independently generated such that (i) $Z\sim\mathcal{N}(0,1)$ and (ii) $V$ (resp.~$W$) is drawn from the empirical distribution of $\{\sqrt{p}\theta_i^{\star}\}$ (resp.~$\{\sqrt{n}\varepsilon_i\}$). 
With this two-dimensional sequence in place, it has been proven that: for any fixed iteration number $t$ and any pseudo-Lipschitz function $\Phi: \mathbb{R} \times \mathbb{R} \to \mathbb{R}$, it holds almost surely that 
\begin{subequations}
		\label{eq:AMP-prediction-SE}
\begin{align}
	\lim_{p\to \infty} \frac{1}{p}\sum_{i=1}^p \Phi\Big(\sqrt{p} \big(\theta_{t+1,i}-\thetastar_i\big), \sqrt{p} \theta^\star_i \Big) 
	&=
	\myE \big[\Phi(\alpha^\star_t Z, V)\big],  
	\\
	\lim_{n\to \infty} \frac{1}{n}\sum_{i=1}^n \Phi\Big(\sqrt{n} (r_{t,i} - \varepsilon_i), \sqrt{n} \varepsilon_i\Big) 
	&=
	\myE \big[\Phi(\gamma^\star_t Z, W)\big],
\end{align} 
\end{subequations}
provided that $p/n$ is a fixed constant. 
For instance, when $\Phi: \mathbb{R}\times \mathbb{R}\rightarrow \mathbb{R}$ is chosen to be $\Phi(a,b)=a^2$, the prediction in \eqref{eq:AMP-prediction-SE} pins down  
the asymptotic squared loss of the AMP iterates as follows
\begin{subequations}
\begin{align}
	\lim_{p\to \infty}  \big\| \theta_{t+1}-\thetastar\big\|_2^2 
	&=
	\myE \big[ (\alpha^\star_t Z)^2\big] = \big(\alpha^\star_t\big)^2,\\
	\lim_{p\to \infty}  \big\| r_{t} - \varepsilon \big\|_2^2 
	&=
	\myE \big[ (\gamma^\star_t Z)^2\big] = \big(\gamma^\star_t\big)^2, 
\end{align}
\end{subequations}
making explicit the operational meanings of  $\alpha^\star_t$ and $\gamma^\star_t$ constructed in the SE recursion \eqref{eq:SE-informal}.

\paragraph{A non-asymptotic theory?} 
While state evolution has played a pivotal role towards understanding AMP in various applications, it is asymptotic in nature, in the sense that it predicts the AMP dynamics in the presence of asymptotically large dimensions with the number of iterations held fixed. 
This limits the prediction power of existing AMP theory in at least two aspects: 
\begin{itemize}
	\item Most prior AMP theory fell short in predicting the convergence rate of AMP below a constant error floor (e.g., it did not predict how many iterations are needed in order to achieve a risk that is $o(1)$ away from that of the fixed point). 
	\item Most prior AMP theory did not provide non-asymptotic rates for the statistical estimation error, nor did it offer non-asymptotic distributional guarantees.   
\end{itemize}
In short, when viewed as a fast iterative algorithm, existing theory for AMP provides an incomplete picture of the convergence rate when compared with that of other optimization algorithms; 
when employed as a theoretical machinery, the AMP theory might sometime lose its benefits as well when compared with other alternative tools (e.g.,  the Gaussian min-max theorem \citep{thrampoulidis2018precise,celentano2020lasso} and the leave-one-out analysis framework \citep{el2018impact,ma2020implicit,chen2019gradient}).

Developing a finite-sample analysis of AMP is instrumental not only in comprehending AMP's efficacy as an optimization algorithm, but also in extending its utility as a fundamental statistical analysis tool. Consequently, it has been an active research direction over the last couple of years. 
The seminal work by \citet{rush2018finite} analyzed AMP for linear models and developed the first result allowing the number of iterations to grow with the problem dimension $n$ --- more precisely, the iteration number $t$ can be as large as $o\big(\log n / \log\log n\big)$; this result is further improved to $O\big(\log n / \log\log n\big)$	for  symmetric AMP in \cite{bao2023leave}. 
Subsequently, \cite{li2022non} presented a general framework for understanding the non-asymptotic performance of AMP in spiked low-rank matrix estimation, allowing the number of iterations to grow as $O\big( n / \poly(\log n)\big)$ and facilitating a more precise non-asymptotic prediction of AMP's behavior.
Of particular interest was the subsequent study by \cite{li2023approximate} concerning the problem of $Z_{2}$ synchronization --- a special case of structured matrix estimation --- revealing fast non-asymptotic convergence of AMP even when initialized randomly. This type of results cannot be derived based on previous SE-based asymptotic analysis.


When it comes to sparse and robust regression, however, the non-asymptotic AMP theory remains highly inadequate. 
On one hand, \citet{rush2018finite} was only able capable of analyzing AMP up to $o\big(\log n / \log\log n\big)$ iterations, which is typically insufficient to uncover the convergence behavior of AMP for higher precision. On the other hand, the non-asymptotic framework in \cite{li2022non} is not readily applicable to linear regression. 
All this gives rise to the following natural questions: 
\begin{center}\emph{Can we develop a non-asymptotic theory for AMP tailored to sparse and robust regression, \\ allowing the number of iterations to grow polynomially in the problem size? }
\end{center}
This question was previously out of reach, and has been posed as an open problem 
in \citet{cademartori2023non}. Addressing this question is crucial in understanding and unleashing the power of AMP across diverse statistical domains. 



\subsection{A peek at our main contributions}

In this paper, we answer the above-mentioned open problem in the affirmative, 
through development of a novel non-asymptotic framework that enables faithful prediction of AMP dynamics even when the number of iterations scales with the problem dimension. 
Based on this framework, we derive finite-sample/finite-time statistical guarantees 
that substantially strengthen the celebrated Gaussian approximation theory of AMP. In what follows, let us highlight several key findings. 

\paragraph{A general analysis recipe.}  
In an attempt to develop non-asymptotic theory for sparse and robust regression, 
we propose a unified recipe that facilitates fine-grained characterizations of the AMP iterates. 

\begin{itemize}
	\item \textit{A fine-grained Gaussian decomposition of AMP iterates.} 
	For any $1\leq t \leq \min\{n,p\}$,  we rigorize a general decomposition of the AMP updates as follows
	\begin{equation}
		\label{eqn:decomp-intro}
	%
		\theta_{t+1} -\thetastar = \sum_{k = 1}^{t} \alpha_{t}^k\psi_k + \zeta_{t} ,
	\end{equation}
		where $\{\psi_k\}_{k=1}^{t}$ are independently generated obeying $\psi_k\sim \mathcal{N}(0,\frac{1}{n} I_p)$, $\zeta_t\in \mathbb{R}^p$ stands for a residual vector, 
		and we denote by $\alpha_t=[\alpha_t^k]_{1\leq k\leq t}\in \mathbb{R}^t$ the coefficient vector. 
		See Theorem~\ref{thm:main-matrix} for details. 
		In particular, for both sparse and robust regression, we can demonstrate that
		\begin{align}
			\|\alpha_t\|_2^2 = \sum_{k=1}^t \big(\alpha_t^k\big)^2 \approx \big(\alpha_t^{\star} \big)^2, 
			\label{eq:approx-alphat-informal}
		\end{align}
		with $\alpha_t^{\star}$ obtained in the asymptotic state evolution sequence \eqref{eq:SE-informal}.

\item \textit{Finite-sample control of the residual terms.}
In light of the above decomposition \eqref{eqn:decomp-intro} of AMP, 
		we further prove in Theorem~\ref{prop:final} that: under certain conditions, the residual terms in \eqref{eqn:decomp-intro} satisfy\footnote{For two functions $f(n)$ and $g(n)$, we write $f(n)\lesssim g(n)$ (or $f(n) = O(g(n))$) if there exists a universal constant $c_{1}>0$ such that $f(n)\leq c_1 g(n)$; similarly, we write $f(n)\gtrsim g(n)$  if $f(n)\geq c_2 g(n)$ for some universal constant $c_{2}>0$. 
If both $f(n)\lesssim g(n)$ and $f(n)\gtrsim g(n)$ hold true, we denote $f(n) \asymp g(n).$ 
}
\begin{align}
\label{eqn:zeta-mag-intro}
	\|\zeta_{t}\|_2 \lesssim \lt(\frac{t\log^2 n}{n}\rt)^{\frac{1}{3}}
\end{align}
for every $t$ obeying $ t \lesssim n/\log^4 n.$ 
This result in conjunction with \eqref{eqn:decomp-intro} and \eqref{eq:approx-alphat-informal} delivers the first finite-sample theory that validates Gaussian approximation of AMP for up to $O(n/\log^4 n)$ iterations.





\end{itemize}

\paragraph{Non-asymptotic AMP theory for sparse and robust regression.} 
The general recipe described above allows us to derive --- in a non-asymptotic fashion --- distributional characterizations of the AMP iterates for both sparse and robust regression, detailed below.

\begin{itemize}
	\item \textit{Sparse regression.} 
	When the unknown signal $\thetastar$ is $k$-sparse, 
	we study the dynamics of AMP designed to promote sparsity, 
	which has intimate connection with the optimally tuned Lasso. 
	We demonstrate that the general framework mentioned above is particularly effective in tackling this setting, 
	with the residual term controlled by \eqref{eqn:zeta-mag-intro}. 
As a concrete consequence, this reveals that the estimation error $\theta_{t+1} -\thetastar$  obeys
\begin{align}
	\theta_{t+1} -\thetastar = v_{t+1} + \zeta_t \qquad \text{with }
	W_1 \left( \mu_{v_{t+1}}, \mathcal{N}\Big(0, \frac{ (\alpha_{t}^{\star})^2 }{n}I_p \Big)  \right)
	\lesssim \frac{\mathsf{poly}(\log n)}{ n^{1/2} } ~\text{ and }~ \|\zeta_t\|_2\lesssim \frac{\log n}{n^{1/3}}
	\label{eq:distribution-guarantee-sparse-informal}
\end{align}
		for any $t\lesssim \mathsf{poly}(\log n)$, 
		where $\mu_{v_t}$ represents the distribution of $v_t$, and 
		$W_1(\cdot,\cdot)$ indicates the Wasserstein distance of order 1 between two distributions (to be defined in \eqref{eqn:wasserstein-p}). 
		Moreover, it can be shown that $\alpha_t^{\star}$ converges to its limiting point (as $t\rightarrow \infty$) exponentially fast. 
	In summary, our result confirms the efficacy of AMP for solving sparse regression, 
		while at the same time improving upon prior theory by providing non-asymptotic distributional guarantees that remain valid up to $n/\mathsf{poly}(\log n)$ iterations. 
	As another implication of our distributional characterization, 
		the distance between the risk of our sparse estimator and the state evolution prediction obeys  
	\begin{align}
		\big\|\mathsf{ST}_{\tau_t}(\theta_t) - \theta^{\star}\big\|_2  
		- \gamma_t^{\star} = O\Big(\frac{\log n}{n^{1/3}}\Big)
	\end{align}
	after a logarithmic number of iterations; this error estimate improves upon the state-of-the-art theory for the Lasso estimator  
		(which was ${O}(\frac{\mathsf{poly}(\log n)}{n^{1/4}})$ as derived in \citet{miolane2021distribution,celentano2020lasso}).
		More details can be found in Section~\ref{sec:sparse}.


	\item \textit{Robust regression.} 
	Another contribution of this paper is to establish non-asymptotic distributional guarantees for AMP tailored to robust regression. 
	More specifically, focusing on the Huber loss, 
		we study the dynamics of the AMP designed to solve the robust M-estimation problem \citep{donoho2016high}. 
		In this case, we demonstrate that the AMP iterates also admit the decomposition~\eqref{eqn:decomp-intro} with the residual term satisfying \eqref{eqn:zeta-mag-intro} for all $t\lesssim n / \mathsf{poly}(\log n)$; 
		as a consequence, the non-asymptotic distributional guarantees \eqref{eq:distribution-guarantee-sparse-informal} continue to be valid in robust regression (albeit with a different state evolution prediction).  
	Another implication of our results is the risk estimate
	\begin{align}
		\ltwo{\theta_{t+1} - \thetastar} - \gamma_t^{\star } = {O}\Big(\frac{\log n}{n^{1/3}}\Big)
	\end{align}
	for all $t = O(\log n)$. When translated to the risk of robust M-estimator, this exhibits a faster rate compared to prior work 
	(note that the previously known bound in \citet{han2022universality} has an error term on the order of ${O}(\frac{\mathsf{poly}(\log n)}{n^{1/500}})$).  
	We refer the readers to Section~\ref{sec:robust-Huber} for detailed discussions.


\end{itemize}

\subsection{Notation} 

In this subsection, we introduce a set of notation that will be useful throughout. 
To begin with, for any integer $n>0$, we denote $[n]=\{1,\dots,n\}$. 
An $\epsilon$-cover of a set $\Theta$ w.r.t.~metric $\rho(\cdot,\cdot)$ refers to a set $\{\theta^1,\theta^2,\ldots,\theta^N\} \subseteq \Theta$ such that, for every $\theta \in \Theta$, there exists some $i \in [N]$ such that $\rho(\theta, \theta^i) \leq \epsilon. $ The $\epsilon$-covering number $N(\epsilon, \rho, \Theta)$ is the cardinality of the smallest $\epsilon$-cover of $\Theta$ w.r.t.~metric $\rho(\cdot,\cdot)$. 
For notational convenience, 
when $\rho$ is taken to be the $\ell_2$-norm, 
we often abbreviate the covering number as $N(\epsilon, \Theta)$.  
In addition, we denote by $\mathbb{B}^d(r)=\{\theta \in \real^d \mid \|\theta\|_2\leq r\}$ the $d$-dimensional Euclidean ball of radius $r$ centered at $0$, and $\mathcal{S}^{d-1}=\{x\in \real^d \mid \|x\|_2=1\}$ the unit sphere in $\real^d$.  
%
%
We often write $0$ (resp. $1$) to denote the all-zero (resp. all-one) vector, and let $I_{n}$ (or simply $I$) indicate the $n\times n$ identity matrix.
When a scalar function is applied to a vector, it should be understood that the function is applied in an entry-wise fashion. 
In addition, for any two functions $f(\cdot)$ and $g(\cdot)$, we write $f(n) \ll g(n)$ or $f(n) = o(g(n))$ if $f(n)/g(n) \to 0$ as $n \to \infty$, and write $f(n) \gg g(n)$ if $g(n)/f(n) \to 0$ as $n \to \infty$.
We use $c_1, c_2, \ldots, C_{1}, C_{2}, \ldots$ to denote universal constants that do not change with salient parameters. Note that these universal constants may change from line to line. 
For any two vectors $a=[a_i]_{1\leq i\leq n}$ and $b=[b_i]_{1\leq i\leq n}$ of the same dimension, 
we denote by $a\circ b=[a_ib_i]_{1\leq i\leq n}$ the Hadamard product.

Given two probability distributions $\mu$ and $\nu$ on $\real^{n}$, the Wasserstein distance of order $q$ between these two distributions is defined as 
\begin{align}
\label{eqn:wasserstein-p}
	W_q(\mu, \nu) \defn \bigg(\inf_{\gamma \in \mathcal{C}(\mu,\nu)} \mathop{\mathbb{E}}\limits_{(x,y)\sim \gamma}\big[ \|x - y\|_2^q \big] \bigg)^{1/q}, 
\end{align}
where $\mathcal{C}(\mu,\nu)$ stands for the set of all couplings of $\mu$ and $\nu$ (i.e., all joint distributions $\gamma(x,y)$ whose marginal distributions are $\mu$ and $\nu$, respectively). 
We often employ $\mu(X)$ or $\mu_{X}$  to denote the distribution of $X$.

\section{Main results}
\label{sec:main-results}

In this section, we present our non-asymptotic decomposition for AMP when applied to both sparse and robust regression, 
following the development of a crucial decomposition of the AMP iterates that make explicit their approximate Gaussianity.

\subsection{A general decomposition for AMP iterates}
\label{sec:decompose}

In this section, we develop a general recipe that helps decompose each AMP iterate into three components: (i) a signal component, (ii) a superposition of Gaussian vectors that captures the main error component, 
and (iii) a residual term (which will be shown to be well-controlled for both sparse and robust regression).


Before proceeding, let us first make note of an equivalent form of the original AMP iterations \eqref{eqn:AMP} as studied in \citet{bayati2011dynamics}. To be precise, 
by setting 
\begin{align}
	\beta_t = \theta_t - \thetastar \qquad &\text{and} \qquad
	s_t = r_t - \varepsilon, \qquad t=0,1,\cdots
\end{align}
(namely, $\beta_t$ (resp.~$s_t$) indicates the error when using $\theta_t$ (resp.~$r_t$) to estimate the true signal $\thetastar$ (resp.~the noise $\varepsilon$)), 
the AMP algorithm \eqref{eqn:AMP} can be equivalently expressed as \citep{bayati2011dynamics}
\begin{subequations}
\label{eqn:AMP-general}
\begin{align}
	s_{t} &= XF_t(\beta_{t}) - \lt\langle F_t^{\prime}(\beta_{t})\rt\rangle G_{t-1}(s_{t-1}), \\
	\beta_{t+1} &= X^{\top}G_t(s_t) - \lt\langle G_t^{\prime}(s_t) \rt\rangle F_{t}(\beta_{t}),
\end{align}
\end{subequations}
where $\{F_t\}_{t\geq 1}$ and $\{G_t\}_{t\geq 0}$ denote two sequences of properly chosen scalar functions  (note that they are applied entrywise to the vector argument here) as
\begin{subequations}
\begin{align}
	G_t(s) = \lt\langle g_{t}^{\prime}(s + \varepsilon)\rt\rangle^{-1}g_{t}(s + \varepsilon)
	\qquad &\text{and} \qquad F_{t}(\beta) = \thetastar - f_{t}(\beta + \thetastar)
\end{align}
initialized to 
\begin{align*}
	G_0(s_0) = 0 \qquad \text{ and } \qquad F_1(\beta_1) &= \thetastar - f_1(\beta_1) = \thetastar. 
\end{align*}
\end{subequations}
%

As a consequence, in order to understand how $\theta_{t}$ evolves during the execution of the algorithm, 
it suffices to focus on the dynamics of $\beta_{t+1}$. 
The following theorem introduces a general decomposition of $(s_t, \beta_{t+1})$, whose proof is postponed to Section~\ref{sec:pf-theorem-main}.





\begin{theos} 
\label{thm:main-matrix}
Consider the linear model~\eqref{eqn:linear} under i.i.d.~Gaussian design (i.e., $X_{ij} \overset{\mathrm{i.i.d.}}{\sim} \mathcal{N}(0,1/n)$).
Suppose the functions $\{G_t\}$ and $\{F_{t}\}$ are differentiable except at a finite number of points. For any $1 \leq t \leq \min\{n,p\},$ the AMP sequence defined in \eqref{eqn:AMP-general} admits the following decomposition:  
\begin{subequations}
\label{defi:dynamics}
\begin{align}
	s_t &= \sum_{k = 1}^{t} \gamma_{t}^k\phi_k + \xi_{t} =: u_t + \xi_t, \label{defi:dynamics-s} \\
	\beta_{t+1} &= \sum_{k = 1}^{t} \alpha_{t}^k\psi_k + \zeta_{t} =: v_{t+1} + \zeta_t, \label{defi:dynamics-beta}
\end{align}
\end{subequations}
where
\begin{itemize}

	\item[(i)] $\{\phi_k\}_{k=1}^{t}$ and $\{\psi_k\}_{k=1}^{t}$ are independent vectors obeying $\phi_{k} \overset{\mathrm{i.i.d.}}{\sim} \mathcal{N}(0, \frac{1}{n} I_n)$ and $\psi_{k} \overset{\mathrm{i.i.d.}}{\sim} \mathcal{N}(0, \frac{1}{n} I_p)$;

	\item[(ii)] the coefficients $\gamma_{t}=[\gamma_t^k]_{1\leq k\leq t} \in \real^{t}$ and $\alpha_{t} =[\alpha_t^k]_{1\leq k\leq t}\in \real^t$ satisfy 
	\begin{align}
	\label{eqn:norm-condition}
		\|\gamma_{t}\|_2 = \|F_t(\beta_t)\|_2\qquad\text{and}\qquad\|\alpha_{t}\|_2 = \|G_t(s_t)\|_2;
	\end{align}

	\item[(iii)] the residual vectors obey $\xi_{t} \in \textsf{span}\{G_1(s_1),\ldots,G_{t-1}(s_{t-1})\}$ and 
	$\zeta_{t} \in \textsf{span}\{F_1(\beta_1),\ldots,F_{t}(\beta_{t})\}$. 
\end{itemize}
\end{theos}
%
\begin{remark}
	Note that the coefficient vectors $\gamma_{t}$ and $\alpha_{t}$ might be statistically dependent on $\{\phi_k\}_{k=1}^{t}$ and $\{\psi_k\}_{k=1}^{t}$. 
\end{remark}

In words, Theorem~\ref{thm:main-matrix} ensures that both $s_{t}$ and $\beta_{t+1}$ can be viewed as weighted superpositions of Gaussian vectors in addition to some residual terms. 
If the residual terms are negligible (which we will demonstrate for both sparse and robust regression), 
$s_{t}$ and $\beta_{t+1}$ are well approximated by $\sum_{k = 1}^{t} \gamma_{t}^k\phi_k$ and $\sum_{k = 1}^{t} \alpha_{t}^k\psi_k$, 
which are both close to spherical Gaussian distributions (in terms of the 1-Wasserstein distance) in the sense that
\begin{align*}
	W_1 \left(\mu\left(\frac{1}{\ltwo{\alpha_t}}\sum_{k = 1}^{t} \alpha_{t}^k\psi_k\right), \mathcal{N}\Big(0,\frac{1}{n}I_p\Big)\right)
	\lesssim 
	\sqrt{\frac{t\log n}{n}}. 
\end{align*}
Here, $\mu(\cdot)$ represents the distribution of a random vector;  see  \citet[Lemma 9]{li2022non} for a proof of this 1-Wasserstein distance result.

In contrast to prior literature, the above decomposition \eqref{defi:dynamics} is deterministic and general in nature, 
requiring very few assumptions (resp.~no assumption) on the denoising functions (resp.~the underlying signal $\thetastar$) 
and making it well-suited for the studies of various models and estimators. 
The most critical challenge for applying Theorem~\ref{thm:main-matrix} then boils down to bounding the magnitudes of the residual terms $\xi_{t}$ and $\zeta_{t}$, 
which often require non-trivial treatments. 
Fortunately, these terms can be very well controlled under both sparse and robust regression, 
which we shall discuss next in  Sections~\ref{sec:sparse} and \ref{sec:robust-Huber}.


\paragraph{Comparisons with prior approaches.} 
Before continuing, we note that the iterative procedure~\eqref{eqn:AMP-general} has also been analyzed in previous works (e.g., \citet[Section 3.2]{bayati2011dynamics}) for understanding the high-dimensional asymptotics 
of AMP for solving various estimators like the Lasso.  
These prior techniques typically rely on the Gaussian conditioning trick \citep{bolthausen2009high,bayati2011dynamics,wu2024sharp} and the state-evolution type analysis, which are drastically different from our proof strategy (as we shall elucidate momentarily). 
As another remark, 
the quantities $\ltwo{\gamma_t}$ and $\ltwo{\alpha_t}$ in Theorem~\ref{thm:main-matrix} are often closely connected to the scalars $\gamma^\star_{t}$ and $\alpha^\star_{t}$ in the asymptotic state evolution \eqref{eq:SE-informal}, which will be made more clear in the next subsections. For this reason, we will sometimes refer to $\ltwo{\gamma_t}$ and $\ltwo{\alpha_t}$ as the finite-sample counterpart of the asymptotic state evolution.



\subsection{Non-asymptotic AMP theory for sparse regression} 
\label{sec:sparse}

With the general decomposition in Theorem~\ref{thm:main-matrix} in place, 
we can readily move forward to investigate concrete models, for which we shall begin with sparse regression. 
Consider the linear model \eqref{eqn:linear} with the underlying signal  $\thetastar \in \real^{p}$ being $k$-sparse.  
The statistical performance of various sparse estimators has been extensively studied, 
with a primary focus on the regime where  $k$ is substantially smaller than $p$ (see, e.g., \cite{donoho2006compressed,candes2006stable,candes2007dantzig,meinshausen2006high,fan2001variable,zou2005regularization,yuan2006model,rudelson2008sparse,wainwright2009sharp,zhao2006model,zhang2010nearly}).
In this work, we focus on the most sample-starved regime with linear sparsity and proportional growth, namely,  
\begin{align}
\label{eqn:lasso-con1}
k \asymp n \asymp p, 
\end{align}
a regime in which AMP proves extremely powerful \citep{bayati2011dynamics,bayati2011lasso,javanmard2013state,maleki2013asymptotic,donoho2013information,su2017false,rush2018finite,fan2022approximate}.



\paragraph{AMP for sparse regression.} 
Let us first remind the readers of the AMP procedure tailored to sparse regression, 
which was introduced both as a fast procedure to find a sparse solution and as a theoretical tool for characterizing the risk of the Lasso estimator \citep{donoho2009message,bayati2011lasso}. 
As mentioned before, the algorithm \eqref{eqn:AMP} adopts the following denoising functions: 
\begin{align}
\label{eqn:lasso-func}
f_t(x) = \sign(x)(|x| - \tau_t)_{+} =: \mathsf{ST}_{\tau_t}(x)
	\qquad \text{and}\qquad 
		g_t(x) = x. 
\end{align}
When it comes to the alternative form \eqref{eqn:AMP-general}, we can simply write 
\begin{align}
\label{eqn:lasso-F-G}
	F_t(\beta) = \theta^{\star} - \mathsf{ST}_{\tau_t}(\theta^{\star} + \beta)\qquad\text{and}\qquad G_t(s) = s + \varepsilon. 
\end{align}

It remains to specify the threshold sequence $\{\tau_t\}$.  
In this work, we concentrate on a specific choice as follows: by augmenting the notation in \eqref{eqn:AMP-rt} and defining
\begin{align}
	r_{t}(\tau) := y-Xf_{t}(\theta_{t};\tau)+\big\langle f_{t}^{\prime}(\theta_{t};\tau)\big\rangle r_{t-1}\qquad\text{with }f_{t}(x;\tau) := \sign(x)(|x|-\tau)_{+}, 	
\end{align}
we select the (adaptive) threshold $\tau_t$ to be  
\begin{align} 
	\label{eq:tau_choice}
	\tau_t \defn \arg\min_{\tau\geq 0} \|r_t(\tau)\|_2.
\end{align}
\begin{remark}
\label{remark-threshold-sparse}
As we will demonstrate later (as in Section~\ref{sec:lasso-se-pf}), one can show that $\tau_{t}$ is very close to 
the quantity $\tau_t^\star$ below:
\begin{align*}
	\tau_t^\star \defn \inf_{\tau: \tau\geq 0} \mathbb{E}\Big[\big\|\theta^{\star} - \mathsf{ST}_{\tau}(\theta^{\star} + \alpha_{t}^{\star} g) \big\|_2^2\Big]
	\qquad \text{with } g\sim \mathcal{N}\Big(0,\frac{1}{n}I_p \Big).
\end{align*}
%
Informally speaking, $\tau_t$ is selected in a data-driven manner aimed at minimizing the mean square estimation error. 
\end{remark}


\paragraph{State evolution for sparse regression.} 
Next, we find it helpful to recall the (limiting version of) state evolution of AMP described in \cite{donoho2009message,bayati2011lasso}. 
Given a fixed sequence of thresholding scalars $\{b_t\}$, for every $t\geq 1$, define a two-dimensional vector $(\alpha^\star_t, \gamma^\star_{t+1})$ recursively as 
\begin{subequations}
\label{eqn:se-sparse-generic}
\begin{align}
	\alpha_t^{\star 2} &= \gamma_t^{\star 2} + \|\varepsilon\|_2^2, \\
	\gamma_{t+1}^{\star 2} &= \mathbb{E}\Big[\big\|\theta^{\star} - \mathsf{ST}_{b_{t}}(\theta^{\star} + \alpha_t^{\star}g)\big\|_2^2\Big],
\end{align}
\end{subequations}
with $g\sim \mathcal{N}\big(0,\frac{1}{n}I_p \big)$ and $\gamma_1^{\star} = \|\theta^{\star}\|_2$. 
In the asymptotic limit (with $p,n\rightarrow \infty$), 
this SE sequence \eqref{eqn:se-sparse-generic} often depends only upon the empirical distribution of $\thetastar$ and the limit of $\|\varepsilon\|_2^2$, 
and is independent from the iterates of the AMP procedure (as long as the threshold sequence is given). 
When the goal is to solve the Lasso estimator for some prescribed regularization parameter $\lambda > 0$ 
\begin{align}
\label{eqn:lasso-lambda}
	\thetahat^{\textsf{Lasso}} \defn \argmin_{\theta} \frac{1}{2}\|y - X\theta\|_2^2 + \lambda \|\theta\|_1, 
\end{align}
the thresholding sequence can be selected to be $b_{t} = a(\lambda) \alpha^\star_{t}$, with $a(\lambda)$ a function of $\lambda$ as specified in \cite{bayati2011lasso}. 

Given our adaptive threshold \eqref{eq:tau_choice}, this subsection focuses on 
\begin{subequations}
\label{eqn:se-sparse}
\begin{align}
	\alpha_t^{\star 2} &= \gamma_t^{\star 2} + \|\varepsilon\|_2^2, \\
	\gamma_{t+1}^{\star 2} &= \mathbb{E}\Big[\big\|\theta^{\star} - \mathsf{ST}_{\tau_t^\star}(\theta^{\star} + \alpha_t^{\star}g)\big\|_2^2\Big],
	\qquad \text{with }
	\tau_t^\star \defn \inf_{\tau: \tau\geq 0} \mathbb{E}\Big[\big\|\theta^{\star} - \mathsf{ST}_{\tau}(\theta^{\star} + \alpha_{t}^{\star} g) \big\|_2^2\Big],
\end{align}
\end{subequations}
where $g\sim \mathcal{N}\big(0,\frac{1}{n}I_p \big)$ and $\gamma_1^{\star} = \|\theta^{\star}\|_2$. 
With this sequence \eqref{eqn:se-sparse} in place, we shall also defining their limiting values (or fixed points): 
\begin{align}
	\alpha^{\star} := \lim_{t\rightarrow \infty} \alpha^{\star}_t
	\qquad \text{and}\qquad
	\gamma^{\star} := \lim_{t\rightarrow \infty} \gamma^{\star}_t. 
\end{align}
Combining this with Remark~\ref{remark-threshold-sparse}, 
we see that the AMP \eqref{eqn:AMP} with the threshold \eqref{eq:tau_choice} attempts to solve the optimally tuned Lasso (namely, picking the choice of $\lambda$ that minimizes the asymptotic estimation risk).


\paragraph{Non-asymptotic analysis for sparse regression.} 
With the above setup and notation in place, we are ready to characterize the non-asymptotic performance of AMP below with the assistance of our general decomposition in Theorem~\ref{thm:main-matrix}.  
The proof of the theorem below is postponed to Section~\ref{sec:pf-thm-lasso}. 
\begin{theos} 
\label{thm:lasso}
Suppose that the $k$-sparse signal $\theta^{\star}$ and the noise vector satisfy 
\begin{align}
\label{eqn:lasso-con2}
\|\theta^{\star}\|_1 \gtrsim \sqrt{k} 
\qquad\text{and}\qquad \|\theta^{\star}\|_2 \asymp \|\varepsilon\|_2 \asymp 1 
\end{align}
with probability exceeding $1 - O(n^{-10})$, 
and assume that  $n > 2k\log ({p}/{k})$ and $p > 2.3k$. 
Then with probability at least $1 - O(n^{-10})$, the AMP iterates \eqref{eqn:AMP} with denoising functions \eqref{eqn:lasso-func}  and threshold \eqref{eq:tau_choice} admit the following decomposition
\begin{align}
\label{eqn:lasso-residual}
\theta_{t+1} - \thetastar = \sum_{j = 1}^{t} \alpha_{t}^j\psi_j + \zeta_{t}
	\qquad \text{and} \qquad
r_t - \varepsilon =  \sum_{j = 1}^{t} \gamma_{t}^j\phi_j + \xi_{t}
\end{align}
for every $t \lesssim \frac{n}{\log^4 n}$, 
where $\{\psi_j\}_{1\leq j\leq t}$ (resp.~$\{\phi_j\}_{1\leq j\leq t}$) are independent Gaussian vectors drawn from $\mathcal{N}\big(0, \frac{1}{n}I_p\big)$ (resp.~$\mathcal{N}\big(0, \frac{1}{n}I_n\big)$), 
the coefficient vectors $\alpha_t=[\alpha_t^j]_{1\leq j\leq t}$ and $\gamma_t=[\gamma_t^j]_{1\leq j\leq t}$ obey 
\begin{subequations}
\begin{align}
\label{eqn:pretty-se-lasso}
\big|\ltwo{\alpha_t}^2 - \alpha^{\star 2}_t\big| \lesssim \Big(\frac{t\log^2 n}{n}\Big)^{1/3}
\qquad \text{and} \qquad
	\big|\ltwo{\gamma_t}^2 - \gamma^{\star 2}_{t}\big| \lesssim \Big(\frac{t\log^2 n}{n}\Big)^{1/3}, 
\end{align}
and the residuals $\{\xi_{t}\}$ and $\{\zeta_{t}\}$ satisfy
\begin{align}
	\|\xi_{t}\|_2, \|\zeta_{t}\|_2 &\lesssim \Big(\frac{t\log^2 n}{n}\Big)^{1/3}.
\end{align}
\end{subequations}
\end{theos}

\begin{remark} 
	In the noiseless case (i.e., $\varepsilon = 0$), 
	the minimum $\ell_{1}$-norm estimator, which corresponds to the $\lambda \to 0$ limit of the Lasso estimator, 
	undergoes a sharp phase transition. 
	As discussed in \cite{amelunxen2014living} and \citet[Page 2201]{celentano2020lasso}, 
	exact recovery by this estimator can only happen when $n \geq 2k(1 + o(k/p))\log(p/k)$, 
	which coincides with the assumption $n > 2k\log ({p}/{k})$ imposed  in Theorem~\ref{thm:lasso}.  
\end{remark}
\begin{remark}
	It is also worth mentioning that Theorem~\ref{thm:lasso} does not restrict the distribution of the noise vector $\varepsilon$, as long as its $\ell_2$-norm is properly controlled to be 
	on the order $1.$ 
\end{remark}
\begin{remark}
The exponent in the probability $1 - O(n^{-10})$ can be replaced with $1-O(n^{-c})$ with an arbitrarily large constant $c>0$. For simplicity, we have made no efforts to obtain the sharpest possible ones. 	
\end{remark}

Let us take a moment to highlight several implications of Theorem~\ref{thm:lasso}. 
\begin{itemize}
	\item {\em Non-asymptotic Gaussian approximation.}
In a nutshell, Theorem~\ref{thm:lasso} demonstrates the proximity of the AMP update $\theta_{t+1}$ and some Gaussian distribution. 
For instance, taking the number of iterations to be $t = c_t \log n$ for some large enough constant $c_t>0$,  
we can guarantee that, 
with probability at least $1 - O(n^{-10})$, 
\begin{align}
\theta_{t+1}=\thetastar+v_{t+1}+\zeta_{t}\qquad\text{with }W_{1}\left(\mu_{v_{t+1}},\mathcal{N}\left(0,\frac{(\alpha_{t}^{\star})^{2}}{n}I_{p}\right)\right)\lesssim\frac{\log n}{n^{1/2}}~\text{ and }~\|\zeta_{t}\|_{2}\lesssim\frac{\log n}{n^{1/3}},	
\end{align}
where $\mu_{v_{t+1}}$ represents the distribution of $v_{t+1}$. 
In fact, given that $\alpha_t^{\star}$ converges to the limiting point $\alpha^{\star}$ exponentially fast (see discussion in Section~\ref{sec:lasso-se-pf}), 
we can further conclude that
\begin{align}
\theta_{t+1}=\thetastar+v_{t+1}+\zeta_{t}\qquad\text{with }W_{1}\left(\mu_{v_{t+1}},\mathcal{N}\left(0,\frac{(\alpha^{\star})^{2}}{n}I_{p}\right)\right)\lesssim\frac{\log n}{n^{1/2}}~\text{ and }~\|\zeta_{t}\|_{2}\lesssim\frac{\log n}{n^{1/3}}	
\end{align}
with probability at least $1 - O(n^{-10})$. 
As far as we know, 
this result offers the first non-asymptotic theory of 
the AMP estimator tailored to sparse regression when $t\gtrsim \frac{\log n}{\log \log n}$, 
which significantly improves upon the best non-asymptotic prior result \citet{rush2018finite} that was only able to accommodate $o\big(\frac{\log n}{\log \log n}\big)$ iterations.



%
\item {\em Improved non-asymptotic risk of the optimally-tuned Lasso.} 
Given the intimate connection between the aforementioned AMP procedure and Lasso --- particularly the one with the regularization parameter carefully chosen to minimize the mean square estimation error --- 
we can immediately see that our result offers a non-asymptotic distributional theory for the optimally-tuned Lasso. 
Note that the best-known distributional theory for the Lasso has been established by \citet{miolane2021distribution,celentano2020lasso} using the Gaussian min-max theorm; 
more concretely, \citet[Theorem 5]{celentano2020lasso} asserts that
\begin{align*}
	\Big|\|\thetahat^{\textsf{Lasso}} - \thetastar\|_2 - \gamma^{\star} \Big| = O\left(\frac{\log n}{n^{1/4}}\right),
\end{align*}
where $\thetahat^{\textsf{Lasso}}$ denotes the solution of the optimally-tuned Lasso. 
Our result indicates that better rates can be obtained with sparse estimators produced by the AMP algorithm. 
In particular, taking $t = c_t \log n$ for some large enough constant $c_t>0$ reveals that 
\begin{align}
	\big\|\mathsf{ST}_{\tau_t}(\theta_t) - \theta^{\star}\big\|_2  = \ltwo{F_t(\beta_t)} = \ltwo{\gamma_t} 
	= 
	\gamma_t^{\star} + O\Big(\frac{\log n}{n^{1/3}}\Big)
	= 
	\gamma^{\star} + O\Big(\frac{\log n}{n^{1/3}}\Big),
\end{align}
where we once again invoke the property that $\gamma_t^{\star}$ converges exponentially fast to $\gamma^{\star}$ (see discussion in Section~\ref{sec:lasso-se-pf}). 
%
It is worth emphasizing that such fine-grained results were unavailable in prior results that used the state-evolution-based analysis of AMP.

\end{itemize}

\subsection{Non-asymptotic AMP theory for robust regression} 
\label{sec:robust-Huber}

%
%

Next, let us turn to robust regression, 
which concerns the linear model \eqref{eqn:linear} with the noise being a mixture of Gaussians and some contamination distribution $H$, i.e., 
\begin{equation}
	\varepsilon_i \overset{\text{i.i.d.}}{\sim} (1-\epsilon_{H}) \mathcal{N}(0, \sigma^2) + \epsilon_{H} H, \qquad 1\leq i\leq n. 
	\label{eqn:robust-error}
\end{equation}
We shall assume throughout that 
\begin{equation}
	\sigma^2 \asymp 1/n. 
\end{equation}
Robust regression was originally proposed by \citet{huber1973robust} and subsequently developed by \citet{bickel1975one} and many others.  
The focus therein was on the case where the signal dimension $p$ is much smaller than the sample size $n$ (see, e.g., \cite{hampel1974influence,maronna2019robust,rousseeuw2005robust,fan2014robust,loh2013regularized,loh2017statistical,sun2020adaptive} and the references therein). 
The modern high-dimensional setting --- where the number of variables $p$ is comparable to the sample size $n$ --- has been recently explored by \cite{el2013robust,donoho2016high,donoho2015variance,el2018impact,lei2018asymptotics,thrampoulidis2018precise,bellec2022asymptotic,bellec2023error,adomaityte2023high}.


\paragraph{AMP for robust regression.} 
A common estimator for robust regression is called the robust M-estimator,  
which selects a non-negative convex loss function $\rho: \real \to \real_{\geq 0}$ and solves the following optimization problem: 
\begin{align}
\label{eqn:robust-reg}
	\thetahat \defn \argmin_{\theta\in \real^p} \mathcal{L}(\theta; y, X), 
	\qquad 
	\text{where } \mathcal{L}(\theta; y, X) \defn \sum_{i=1}^n \rho\big( y_i - \inprod{X_i}{\theta} \big).
\end{align}
In this subsection, we focus on the Huber loss as follows 
\begin{align}
	\label{eq:defn-Huber-loss}
	\rho(z)
	= \rho_{\mathsf{huber}}(z,\lambda) \defn 
	\begin{cases}
		z^2/2, & \text{if } |z| \leq \lambda \\
		\lambda |z| - \lambda^2/2,\quad &\text{otherwise}
	\end{cases}
\end{align}
for some prescribed threshold $\lambda >0$ (chosen such that $\lambda \asymp 1/\sqrt{n}\asymp \sigma$), 
which is arguably the most popular choice to tackle robust regression.

In an attempt to compute and quantify the risk of the above robust M-estimator,  
one can resort to the AMP algorithm \eqref{eqn:AMP} with the following denoising functions \citep{donoho2016high}
\begin{subequations}
\label{eqn:Huber-func}
\begin{align}
	f_t(x) = x \qquad \text{and}\qquad g_t(x) = \frac{n}{p}\Psi(z,b_t),
\end{align}
where we remind the reader that 
\begin{align*}
	\Psi(z,b) = \rho_b'(z) 
	\qquad \text{with }~\rho_b(z) \defn \min_{x} \Big\{\rho(x) + \frac{1}{2b}(x-z)^2 \Big\}
\end{align*}
for some regularization parameter $b>0$. 
When $\rho$ corresponds to the Huber loss \eqref{eq:defn-Huber-loss}, 
it is easily seen that 
\begin{align*}
	\Psi(z,b) = b\psi\Big(\frac{z}{1+b}, \lambda\Big)
	\qquad
	\text{with}
	~~
	\psi(z; \lambda) = \min\big\{\max\{z, -\lambda\}, \lambda \big\}.
\end{align*}
\end{subequations}
Additionally, the algorithm is initialized at $\theta_1 = 0$, $r_0=0$ and $r_{1} = y$, and the parameter $b_{t}$ is chosen such that 
\begin{align}
	\frac{1}{n} \sum_{i=1}^n \Psi'(r_i^{t} , b_{t}) = \frac{p}{n},
	\qquad \text{or equivalently,} \qquad
	\big\langle g_t'(r_t) \big\rangle = 1, 
		\label{eq:bt_choice}
\end{align}
where $\Psi'(\cdot , \cdot)$ denotes differentiation w.r.t.~the first variable.

\begin{remark}
When it comes to the equivalent representation \eqref{eqn:AMP-general}, one can choose 
\begin{align}
\label{eqn:robust-F-G}
F_t(\beta) = -\beta\qquad\text{and}\qquad 
G_{t}(s)=g_{t}(s+\varepsilon)=\frac{n}{p}b_{t}\psi\lt(\frac{s+\varepsilon}{1+b_{t}};\lambda\rt).  
\end{align}
%
%
%
\end{remark}








\paragraph{State evolution for robust regression.} 
%
%
In order to predict the dynamics of the AMP algorithm, 
it is helpful to introduce the following state evolution recursion as introduced in \citet{donoho2016high}. 
Specifically, for any $t \geq 1$, \citet{donoho2016high} defines a sequence of $(\alpha^{\star}_t, \gamma^{\star}_{t+1}) \in \real^{2}$ recursively as follows:
\begin{subequations}
\label{eqn:robust-SE}
\begin{align}
	\alpha_t^{\star 2} &= \Big(\frac{n {b}^\star_t}{p(1 + {b}^\star_t)}\Big)^2\mathbb{E}\Big[\big\|\psi\big(\varepsilon + \gamma_{t}^{\star}g; \lambda(1+ {b}^\star_t)\big)\big\|_2^2\Big], \\
	\gamma_{t+1}^{\star 2} &= \frac{p}{n}\alpha_t^{\star 2},
\end{align}
where $\gamma_1^{\star} = \|\theta^{\star}\|_2$ and ${b}^\star_t$ is chosen to satisfy
\begin{align}
	\frac{1}{n} \sum_{i=1}^n \Exs\big[\Psi'( \varepsilon_i + \gamma_{t}^{\star}g_i, b^\star_t) \big] =	\frac{p}{n}
	\qquad \text{with } g_i \overset{\mathrm{i.i.d.}}{\sim} \mathcal{N}\Big(0, \frac{1}{n}\Big). 
\end{align}
\end{subequations}
%
It is worthwhile to remark that the sequence $(\alpha^{\star}_t, \gamma^{\star}_{t+1}) \in \real^{2}$ does not depend on the actual iterates of the AMP procedure. 
%

\paragraph{Non-asymptotic analysis for robust regression.}
We are now ready to present our non-asymptotic theory for AMP in the context of robust regression. 
With the aid of our general decomposition in Theorem~\ref{thm:main-matrix}, 
we can establish the following non-asymptotic theoretical guarantees. 
\begin{theos} 
\label{thm:robust}
Suppose that the signal $\thetastar$ and the noise satisfy
%
\begin{subequations}
\begin{align}
\|\theta^{\star}\|_2 \asymp \|\varepsilon\|_2 \asymp 1
\end{align}
with probability at least $1 - O(n^{-10})$.  
Assume that
\begin{align}
n > p \asymp n, 
	\qquad \lambda \asymp 1/\sqrt{n} 
\qquad \text{and} \qquad 
\sigma^2 \asymp 1/n . 
\end{align}
\end{subequations}
Consider the denoising functions chosen as in \eqref{eqn:robust-F-G}. 
Then with probability exceeding $1 - O(n^{-10})$, 
 the AMP iterates \eqref{eqn:AMP} with denoising functions \eqref{eqn:Huber-func}  and threshold \eqref{eq:bt_choice} admit the decomposition
\begin{align}
\label{eqn:robust-residual}
\theta_{t+1} - \thetastar &= \sum_{k = 1}^{t} \alpha_{t}^k\psi_k + \zeta_{t}
	\qquad \text{and} \qquad
r_t - \varepsilon =  \sum_{k = 1}^{t} \gamma_{t}^k\phi_k + \xi_{t}
\end{align}
for every $t \lesssim \frac{n}{\log^4 n}$, 
where $\{\psi_j\}_{1\leq j\leq t}$ (resp.~$\{\phi_j\}_{1\leq j\leq t}$) are independent Gaussian vectors drawn from $\mathcal{N}\big(0, \frac{1}{n}I_p\big)$ (resp.~$\mathcal{N}\big(0, \frac{1}{n}I_n\big)$), 
the coefficient vectors $\alpha_t=[\alpha_t^j]_{1\leq j\leq t}$ and $\gamma_t=[\gamma_t^j]_{1\leq j\leq t}$ obey
\begin{subequations}
\begin{align}
\label{eqn:pretty-se-robust-bla}
\big|\ltwo{\alpha_t}^2 - \alpha^{\star 2}_t\big| \lesssim \Big(\frac{t\log^2 n}{n}\Big)^{1/3}
	\qquad \text{and} \qquad
\big|\ltwo{\gamma_t}^2 - \gamma^{\star 2}_{t}\big| \lesssim \Big(\frac{t\log^2 n}{n}\Big)^{1/3}, 
\end{align}
and the residuals $\{\xi_{t}\}$ and $\{\zeta_{t}\}$ satisfy
\begin{align}
\|\xi_{t}\|_2, \|\zeta_{t}\|_2 &\lesssim \Big(\frac{t\log^2 n}{n}\Big)^{\frac{1}{3}}.
\end{align}
\end{subequations}
\end{theos}
\noindent The proof of this result is provided in Section~\ref{sec:pf-thm-robust}.


In words, Theorem~\ref{thm:robust} characterizes the distribution of AMP updates with finite-sample guarantees. 
Akin to the sparse regression case, 
if the number of iterations is taken to be $t = O(\log n)$, 
then one can write 
\begin{align}
\theta_{t+1}=\thetastar+v_{t+1}+\zeta_{t}\qquad\text{with }W_{1}\left(\mu_{v_{t+1}},\mathcal{N}\left(0,\frac{(\alpha_{t}^{\star})^{2}}{n}I_{p}\right)\right)\lesssim\frac{\log n}{n^{1/2}}~\text{ and }~\|\zeta_{t}\|_{2}\lesssim\frac{\log n}{n^{1/3}}	
\end{align}
with high probability. 
%
%
%
In the meantime, with probability exceeding $1-O(n^{-10})$ one has
\begin{align}
	\ltwo{\theta_{t+1} - \thetastar} = \ltwo{F_t(\beta_t)} = \ltwo{\gamma_t}
	=
	\gamma_t^{\star } + O\Big(\frac{\log n}{n^{1/3}}\Big).
\end{align}
Evidently, these results recover \citet[Theorem 1.2]{donoho2016high} in the high-dimensional asymptotics with $n, p \to \infty$ for any fixed $t$,  
while at the same time improving upon \citet{donoho2016high} by offering non-asymptotic distributional guarantees that account for up to a polynomial number of iterations.

Before concluding this section, 
we note that the performance of the robust M-estimator in the high dimensional asymptotics has been 
studied in \citet{donoho2016high,donoho2015variance} with the aid of the AMP machinery. 
Certain regularized variants of the robust M-estimator have also been analyzed by means of the leave-one-out analysis \citep{karoui2013asymptotic,el2018impact} 
and convex Gaussian min-max (CGMT) theorem \citep{thrampoulidis2018precise,han2022universality}. 
The only result that offers explicit non-asymptotic guarantees is provided in \citep{han2022universality}, 
which leverages the CGMT technique to control the finite-sample error bound to be the order of $O(n^{-1/500})$. 
In contrast, the AMP analysis in Theorem~\ref{thm:robust} offers finite-sample error bound on the order of $O\big(\frac{\log n}{n^{1/3}}\big)$. 
It is worth noting, however, that \citet{han2022universality} is able to extend beyond i.i.d.~Gaussian design and unveil interesting universality phenomena.





\section{Key technical innovation: controlling the residuals}
\label{sec:technical}

It is worth noting that the decomposition \eqref{defi:dynamics} in Theorem~\ref{thm:main-matrix}, 
while being fully non-asymptotic and general, has not yet offered quantitative descriptions about the magnitudes of the residual terms $\xi_{t}$ and $\zeta_{t}$, 
making it insufficient to imply any distributional guarantees of the AMP iterates. 
The key innovation of the current paper thus lies in establishing effective control of $\xi_{t}$ and $\zeta_{t}$. 
In comparison, the approach adopted in our prior work \cite{li2022non} was insufficient to accommodate the most challenging SNR regime for sparse and robust regression; more discussions on this can be found at the end of this section.


\subsection{A fine-grained decomposition for the residuals}

A key ingredient of our analysis is to develop a finer representation of $\xi_{t}$ (resp.~$\zeta_{t}$) by analyzing its corresponding coefficient on the set $\{G_1(s_1),\ldots,G_{t-1}(s_{t-1})\}$ (resp. $\{F_1(\beta_1),\ldots,F_{t}(\beta_{t})\}$). 
Towards this end, we find it helpful to first construct a couple of auxiliary sequences, detailed next.



\paragraph{Auxiliary sequences.}
We now construct recursively a set of auxiliary iterates $\{\widehat{s}_t\}\subseteq \mathbb{R}^n$ and $\{\widehat{\beta}_t\}\subseteq \mathbb{R}^p$, 
as well as the coefficient vectors $\hatalpha_{t}=[\hatalpha_{t}^k]_{1\leq k\leq t}\in \real^t$ and $\hatgamma_{t}=[\hatgamma_{t}^k]_{1\leq k\leq t}\in \real^{t}$. 
Specifically,  
\begin{itemize}
	\item Let us start with $\hatalpha_{0} = 0$ and $\hats_{1} = u_1$ and $\hatbeta_{1} = v_{1}.$

	\item For each $t\geq 1$ and $\phi_{k} \overset{\text{i.i.d.}}{\sim} \mathcal{N}(0, \frac{1}{n} I_n)$ and $\psi_{k} \overset{\text{i.i.d.}}{\sim} \mathcal{N}(0, \frac{1}{n} I_p)$, 
	
	\begin{itemize}
		\item construct the vector $\hatgamma_{t}=[\hatgamma_{t}^k]_{1\leq k\leq t}\in \real^{t}$ 
\begin{subequations}
\label{eqn:sec-order-terms}
\begin{align} \label{eq:zeta-coeff}
\hatgamma_{t}^k \defn
\left\{ \begin{array}{lll}
	\langle G_{t}^{\prime}(\hats_{t}) - G_{t}^{\prime}(s_{t})\rangle + \frac{1}{\|\gamma_{t}\|_2^2}\lt\langle \sum_{k=1}^{t} \gamma_{t}^k\phi_k, G_{t}(s_{t}) - G_{t}(\hats_{t})\rt\rangle & \text{for} & k = t \\
\hatalpha_{t-1}^k \lt\langle G_{t}^{\prime}\lt(\hats_{t}\rt) \circ G_{k}^{\prime}\lt(u_{k}\rt)\rt\rangle & \text{for} & k < t.
\end{array}\right.
\end{align}

	\item compute
\begin{align}
	\hatbeta_{t+1} &\defn v_{t+1} + \sum_{k = 1}^{t} \hatgamma_{t}^k F_{k}(v_{k}), 
	\qquad v_{t+1} \defn \sum_{k = 1}^{t} \alpha_{t}^k\psi_k.	\label{eqn:hat-betat} 
\end{align}

	\item construct the vectors $\hatalpha_{t}=[\hatalpha_{t}^k]_{1\leq k\leq t}\in \real^t$ such that 
\begin{align} \label{eq:xi-coeff}
	\hatalpha_{t}^k \defn 
	\left\{ \begin{array}{lll}
	\langle F_{t+1}^{\prime}(\hatbeta_{t+1}) - F_{t+1}^{\prime}(\beta_{t+1})\rangle + \frac{1}{\|\alpha_t\|_2^2}\lt\langle \sum_{k=1}^t \alpha_t^k\psi_k, F_{t+1}(\beta_{t+1}) - F_{t+1}(\hatbeta_{t+1})\rt\rangle & \text{for} & k = t \\
	\hatgamma_{t}^{k+1} \lt\langle F_{t+1}^{\prime}(\hatbeta_{t+1}) \circ F_{k+1}^{\prime}(v_{k+1})\rt\rangle & \text{for} & k < t.
	\end{array}\right.
\end{align}

	\item compute
\begin{align}
	 \hats_{t+1} &\defn u_{t+1} + \sum_{k = 1}^{t} \hatalpha_{t}^k G_{k}(u_{k}), 
	 \qquad u_{t+1} \defn \sum_{k = 1}^{t+1} \gamma_{t+1}^k\phi_k . \label{eqn:hat-st}
\end{align}
\end{subequations}
\end{itemize}
\end{itemize} 
\noindent Crucially, the auxiliary sequence $(\hats_{t}, \hatbeta_{t+1})$ constructed above serves as a good proxy of $(s_t,\beta_{t+1}).$



\paragraph{Fine-grained representation.} 
Equipped with the above quantities, we are ready to state the following result, whose proof is deferred to Section~\ref{sec:pf-theorem-second}.


\begin{theos}
\label{thm:residuals}
Under the assumptions of Theorem~\ref{thm:main-matrix}, the residual terms in the decomposition~\eqref{defi:dynamics} can be written as  
\begin{subequations}
\label{eqn:residual-second}
\begin{align}
	\xi_{t+1} &= \sum_{k = 1}^{t} \hatalpha_{t}^k G_k(s_k) + \hatxi_{t+1}, \label{eq:xi_expansion} \\
	\zeta_{t+1} &= \sum_{k = 1}^{t+1} \hatgamma_{t+1}^k F_{k}(\beta_{k}) + \hatzeta_{t+1}, \label{eq:zeta_expansion}
\end{align}
\end{subequations}
where $\hatxi_{t+1}$ and $\hatzeta_{t+1}$ satisfy, with probability at least $1 - O(n^{-10})$,   that
\begin{subequations}
\begin{align}
\label{eq:hatxi}
	\notag \hatxi_{t+1} &= \sum_{k = 1}^{t} a_k\lt[\big\langle \psi_k, F_{t+1}(\hatbeta_{t+1})\big\rangle - \big\langle F_{t+1}^{\prime}(\hatbeta_{t+1}) \big\rangle \alpha_{t}^k - \sum_{j = k+1}^{t} \hatgamma_{t}^j \lt\langle F_{t+1}^{\prime}(\hatbeta_{t+1}) \circ F_{j}^{\prime}\lt(v_{j}\rt)\rt\rangle\alpha_{j-1}^k\rt] \notag\\
	&\quad+ \mathcal{P}_{G_t(s_t)}^{\perp}\sum_{k = 1}^{t} a_k\big\langle \psi_k, F_{t+1}(\beta_{t+1}) - F_{t+1}(\hatbeta_{t+1})\big\rangle + 
	\hatxi_{t+1,\mathsf{res}}
\end{align}
where 
\begin{align}
	\hatzeta_{t+1} &= \sum_{k = 1}^{t+1} b_k\lt[ \big \langle \phi_k, G_{t+1}(\hats_{t+1}) \big\rangle - \big\langle G_{t+1}^{\prime}(\hats_{t+1}) \big\rangle \gamma_{t+1}^k - \sum_{j = k}^{t} \hatalpha_{t}^j \lt\langle G_{t+1}^{\prime}(\hats_{t+1}) \circ G_{j}^{\prime}\lt(u_{j}\rt)\rt\rangle\gamma_{j}^k\rt] \notag\\
	&\quad+ \mathcal{P}_{F_{t+1}(\beta_{t+1})}^{\perp}\sum_{k = 1}^{t+1} b_k\big\langle \phi_k, G_{t+1}(s_{t+1}) - G_{t+1}(\hats_{t+1})\big\rangle +
	\hatzeta_{t+1,\mathsf{res}}
\end{align}
with 
\begin{equation}
	\big\| \hatxi_{t+1,\mathsf{res}} \big\|_2 \lesssim \sqrt{\frac{t\log n}{n}}\|\gamma_{t+1}\|_2
	\qquad \text{and} \qquad
	\big\| \hatzeta_{t+1,\mathsf{res}} \big\|_2 \lesssim \sqrt{\frac{t\log n}{n}}\|\gamma_{t+1}\|_2.
\end{equation}
\end{subequations}
Here, 
	$\mathcal{P}_{w}^{\perp}$ denotes the linear projection onto the subspace orthogonal to the vector $w$, 
	and $\{a_k\}_{k=1}^{t}$ (resp.~$\{b_k\}_{k=1}^{t}$) represents a set of orthogonal basis 
	(which are made precise in expression~$\eqref{eqn:defn-abx-a}$ (resp.~\eqref{eqn:defn-abx-b})). 
\end{theos}

Let us take a moment to provide some technical interpretations about the usefulness of Theorem~\ref{thm:residuals} in controlling the residual terms $\ltwo{\xi_t}$ and $\ltwo{\zeta_{t}}$.
The readers who are more interested in seeing direct consequences of this result can move directly to Section~\ref{sec:direct-bound}.
\begin{itemize}
\item 
Considering the first term in $\hatxi_{t+1}$.  
Given that $\psi_{k} \overset{\text{i.i.d.}}{\sim} \mathcal{N}(0, \frac{1}{n} I_p)$, 
a little algebra reveals that this term takes the following form
\begin{align*}
	X^\top f(X) - \textsf{div} f(X)
	\qquad
	\text{with} ~\textsf{div} f \defn \sum_{i} \frac{\partial f_{i}}{\partial x_i}
\end{align*}
for some function $f$ and some Gaussian random vector $X$. 
If we pretend that the function $f$ is statistically independent from $X$, 
then the celebrated Stein lemma tells us that this term has zero mean, 
which provides some intuition why one can expect it to be small. 
Similar messages continue to hold for the first term of $\hatzeta_{t+1}$.

\item Next, let us take a look at the second term in $\hatxi_{t+1}$. 
	One important component here is the projection operator $\mathcal{P}_{G_t(s_t)}^{\perp}$, 
	which plays a crucial role in achieving the desired bound. 
		To explain this, note that in order to bound the $\ell_2$-norm in on the right-hand side of \eqref{eq:hatxi}, one strategy is to look at every unit vector $w\perp G_{t}(s_{t})$ and bound
\begin{align*}
\bigg\langle w,\,\mathcal{P}_{G_{t}(s_{t})}^{\perp}\sum_{k=1}^{t}a_{k}\big\langle\psi_{k},F_{t+1}(\beta_{t+1})-F_{t+1}(\hatbeta_{t+1})\big\rangle\bigg\rangle & =\sum_{k=1}^{t}\langle w,a_{k}\rangle\big\langle\psi_{k},F_{t+1}(\beta_{t+1})-F_{t+1}(\hatbeta_{t+1})\big\rangle\\
 & \approx\bigg\langle\sum_{k=1}^{t}\langle w,a_{k}\rangle\psi_{k},\,F'_{t+1}(v_{t+1})\circ(\beta_{t+1}-\hatbeta_{t+1})\bigg\rangle\\
 & \leq\left\Vert \sum_{k=1}^{t}\langle w,a_{k}\rangle\psi_{k}\circ F'_{t+1}(v_{t+1})\right\Vert _{2}\big\|\beta_{t+1}-\hatbeta_{t+1}\big\|_{2},
\end{align*}
	where the second line makes use of the mean value theorem and the fact that $\beta_{t+1}\approx\hatbeta_{t+1}\approx v_{t+1}$ (rigorous derivations are given around inequality~\eqref{eqn:intuition}). 
	To see why the above bound is useful, we recall two important facts  
\begin{subequations}
\begin{align*}
	G_{t}(s_{t}) &= \sum_{k = 1}^{t} \alpha_{t}^ka_k \qquad \text{ for } 
	\alpha_{t}^k \defn \big\langle G_{t}(s_t), a_k\big\rangle \qquad (1\leq k \leq t), \label{eq:Gt-st-informal}\\
	v_{t+1} &= \sum_{k = 1}^{t} \alpha_t^k \psi_{k}.
\end{align*}
\end{subequations}
Recalling that $w\perp G_{t}(s_{t})$ and assuming that $w$ and $\alpha_t$ are all statistically independent from $\{\psi_k\}$, 
one can easily see that 
\begin{equation*}
	\sum_{k=1}^{t}\langle w,a_{k}\rangle\psi_{k} \quad \text{is independent from } v_{t+1}. 
\end{equation*}
This independence property plays a crucial role in improving the pre-constant on the bound of 
$\big\Vert \sum_{k=1}^{t}\langle w,a_{k}\rangle\psi_{k}\circ F'_{t+1}(v_{t+1})\big\Vert _{2}$ (compared to the case when no independence is assumed) thus controlling the speed for which $\ltwo{\widehat{\xi}_t}$ grows. 
All this is enabled by considering the projections to $G_{t}(s_{t})$ and its orthogonal space and treat them separately.

\item With the decomposition~\eqref{eqn:residual-second} in mind,  
	a natural strategy to bound $\ltwo{\xi_{t+1}}$ (resp.~$\ltwo{\zeta_{t+1}}$) is to control $\|\hatxi_{t+1}\|_2$ (resp.~$\|\hatzeta_{t+1}\|_2$) and
$\sum_{k = 1}^{t} \hatalpha_{t}^k G_k(s_k)$ (resp.~$\sum_{k = 1}^{t} \hatgamma_{t}^k F_{k}(\beta_{k})$) separately. 
Let us now take a moment to discuss the term $\sum_{k = 1}^{t} \hatalpha_{t}^k G_k(s_k)$ --- the message for $\sum_{k = 1}^{t} \hatgamma_{t}^k F_{k}(\beta_{k})$ is similar. 
As we shall justify in the analysis, the coefficient $|\hatalpha_{t}^{k}|$ decays exponentially in the sense that
\begin{equation}
	|\hatalpha_{t}^{k}| \lesssim (1-c)^{t-k} |\hatalpha_{t}^{t}|
\end{equation}
for some constant $c>0$ bounded away from 0. 
Taking this collectively with \eqref{eqn:norm-condition} and the property $\big\|\alpha_{t}\big\|_{2}\approx \alpha_t^{\star}\lesssim 1$ thus reveals that
\begin{align*}
\Big\|\sum_{k=1}^{t}\hatalpha_{t}^{k}G_{k}(s_{k})\Big\|_{2} & \leq\sum_{k=1}^{t}\big|\hatalpha_{t}^{k}\big|\big\| G_{k}(s_{k})\big\|_{2}\lesssim\big|\hatalpha_{t}^{t}\big|\sum_{k=1}^{t}(1-c)^{t-k}\big\|\alpha_{k}\big\|_{2}\\
 & \asymp\big|\hatalpha_{t}^{t}\big|\sum_{k=1}^{t}(1-c)^{t-k}\alpha_{k}^{\star}\asymp\big|\hatalpha_{t}^{t}\big|.
\end{align*}
Consequently, the analysis will focus on bounding the size of $\hatalpha_{t}^{t}$.

\end{itemize}

As it turns out, the above observations play an important role in obtaining an effective control of $\xi_{t}$ and $\zeta_{t+1}$ under some mild conditions, to be detailed in the next subsection.

\subsection{Bounding the residuals under some key assumptions}
\label{sec:direct-bound}

With the decomposition in Theorem~\ref{thm:residuals} in place, 
we further develop upper bounds on the sizes of $\xi_{t}$ and $\zeta_{t+1}$ under some conditions.

To do so, let us  begin with some notation.  
When the function $F: \real \to \real$ is Lipschitz continuous, 
we define $\rho_F\geq 1$ to be the smallest constant (larger than or equal to $1$) such that\footnote{Note that this definition of $\rho_F$ assumes $\rho_F \geq 1$ primarily for notational simplicity; our result would not change if we do not impose this restriction but simply replace $\rho_{F}$ with $\rho_{F}\vee 1.$} 
\begin{align*}
	|F(x_1) - F(x_2)| \leq \rho_F |x_1 - x_2|, \qquad \text{ for all } x_1, x_2.
\end{align*}
%
Analogously, we define $\rho_G\geq 1$ for the function $G.$ 
Additionally, suppose that the functions $F$ and $G$ are both differentiable except at a finite number of points.  
By defining their corresponding derivatives as $F'$ and $G'$ (except at the non-differentiable points), 
we can introduce the quantities $\rho_{1,F}$ and $\rho_{1,G}$ to represent respectively the maximum local Lipschitz constants of $F'$ and $G'$ over the set of differentiable points.  
Armed with the above notation, we can introduce the following assumptions regarding the denoising functions $F_{t}$ and $G_{t}$ and the AMP updates. 
\begin{assumption}
\label{ass:control-lipschitz}
	Suppose that for any $t \lesssim \frac{n}{\log^4 n}$, the aforementioned Lipschitz constants satisfy 
\begin{subequations} \label{eq:conditions-simple}
\begin{align}
\label{eq:conditions-simple-1}
	\rho_{F_t}, \rho_{G_t} \asymp 1,\qquad\text{and}\qquad \rho_{1,F_t}, \rho_{1,G_t} = 0,
\end{align}
In addition, for any $t \lesssim \frac{n}{\log^4 n}$,  suppose that conditional on  $\|\xi_{t}\|_2, \|\zeta_{t}\|_2 \lesssim 1$, one has  
	the coefficients $\gamma_{t}$ and $\alpha_{t}$ in decomposition~\eqref{defi:dynamics} satisfy 
\begin{align}
	\label{eq:conditions-simple-2}
	\|\gamma_t\|_2, \|\alpha_t\|_2 \asymp 1\qquad\text{and}\qquad \lt\|F_{t}(0)\rt\|_2, \lt\|G_{t}(0)\rt\|_2 \lesssim 1
\end{align}
\end{subequations}
with probability at least $1 - O(n^{-11})$. 
\end{assumption}

\begin{remark}
	For readers familiar with the AMP literature, $\|\gamma_t\|_2, \|\alpha_t\|_2$ can be viewed as finite-sample counterparts of the asymptotic state evolution \eqref{eq:SE-informal}. 
	Therefore, the assumption \eqref{eq:conditions-simple-2} requires the finite-sample state evolution for the corresponding problem to be somewhat regular, with extreme events 
	occuring only with very low probability. The result in our theorem below might not hold if the AMP path degenerates or explode at some point. 

\end{remark}

\begin{assumption} 
\label{ass:control-simple}
Let $\widetilde{u}_{t} = \|\gamma_t\|_2 g_{1}$ and $\widetilde{v}_{t+1} = \|\alpha_t\|_2 g_{2}$, where $g_{1} \sim \mathcal{N}(0, \frac{1}{n} I_n)$ and $g_{2} \sim \mathcal{N}(0, \frac{1}{n} I_p)$ are independent with $\|\gamma_t\|_2$ and $\|\alpha_t\|_2$.
Suppose that there exists some universal constant $0 < c < 1/2$ such that 
\begin{align}
\frac{1}{n^2}\mathbb{E}\big[\lt\|G_{t}^{\prime}(\widetilde{u}_{t})\rt\|_2^2 \mymid \|\gamma_t\|_2\big]\mathbb{E}\big[\lt\|F_{t+1}^{\prime}(\widetilde{v}_{t+1})\rt\|_2^2 \mymid \|\alpha_t\|_2\big] &< (1-2c)^2, \\
\frac{1}{n^2}\mathbb{E}\big[\lt\|F_{t+1}^{\prime}(\widetilde{v}_{t+1})\rt\|_2^2 \mymid \|\alpha_t\|_2\big]\mathbb{E}\big[\lt\|G_{t+1}^{\prime}(\widetilde{u}_{t+1})\rt\|_2^2 \mymid \|\gamma_{t+1}\|_2\big] &< (1-2c)^2.
\end{align} 
\end{assumption}


Under these two assumptions, 
we can obtain simple bound on the size of the residual terms in decomposition~\eqref{defi:dynamics} as follows;  
the proof is deferred to Section~\ref{sec:pf-prop-final}.

\begin{theos}
\label{prop:final}
Suppose that the assumptions of Theorem~\ref{thm:main-matrix} hold. 
Under Assumptions~\ref{ass:control-lipschitz} and \ref{ass:control-simple}, the residual terms in decomposition~\eqref{defi:dynamics} satisfy,  with probability at least $1 - O(n^{-10})$,  
\begin{align}
\label{eqn:rate-general}
	\|\xi_{t}\|_2, \|\zeta_{t}\|_2 \lesssim \lt(\frac{t\log^2 n}{n}\rt)^{\frac{1}{3}}
\end{align}
for every $1 < t \lesssim {n}/{\log^4 n}$. 
\end{theos}
\begin{remark}
	Note that in this theorem (and the proof), we allow $F_t$ and $G_t$ to be either deterministic functions, or some random functions. 
	When they are random functions, we assume the existence of a collection of functions $\{F(\cdot;\tau)\}$ (resp.~$\{G(\cdot;b)\}$) parameterized by some constant-dimensional $\tau$ (resp.~$b)$ with $\|\tau\|_2\lesssim 1$ (resp.~$\|b\|_2\lesssim 1$) such that $F(\cdot;\tau)$ (resp.~$G(\cdot;b)$) is $\mathsf{poly}(n)$-Lipschitz in $\tau$ (resp.~$b$). 
	It is then assumed that both $\{F(\cdot;\tau)\}$ and $\{G(\cdot;b)\}$ satisfy analogous assumptions as in Assumptions~\ref{ass:control-lipschitz} and \ref{ass:control-simple}, 
	and that $F_t=F(\cdot;\tau_t)$ (resp.~$G_t=G(\cdot;b_t)$) for some random quantity $\tau_t$ (resp.~$b_t$).   
The fact that the parameters $\tau$ and $b$ are constant-dimensional makes it feasible to apply a standard covering-based argument. 
\end{remark}

\begin{remark} 
	We would like to point out that the exponents in the probability $1 - O(n^{-10})$ and $1 - O(n^{-11})$ in our  assumption can be replaced with $1-O(n^{-c})$ with an arbitrarily large constant $c>0$. 
\end{remark}

Theorem~\ref{prop:final} delivers simple upper bounds on the sizes of $\|\xi_{t}\|_2$ and $\|\zeta_{t}\|_2$, 
as long as the required assumptions on $F_{t}$, $G_{t}$ and the AMP updates can be validated. 
If these assumptions were satisfied, then taking this result collectively with Theorem~\ref{thm:main-matrix} would ensure that both $s_{t}$ and $\beta_{t+1}$ are well approximated by Gaussian distributions with error terms bounded in size by $O(({t\log^2 n}/{n})^{\frac{1}{3}}).$


Consequently, in order to establish our results for sparse and robust regression in Sections~\ref{sec:sparse} and \ref{sec:robust-Huber}, 
everything boils down to verifying these assumptions in the two models of interest. 
It is worth noting that the applicability of Theorem~\ref{prop:final} can potentially extend beyond these two important regression problems.


\paragraph{Comparison with \citet{li2022non}.}
We now pause to emphasize the technical novelty of this paper compared to the prior work \cite{li2022non}. 
To begin with, while the general decomposition in Theorem~\ref{thm:main-matrix} shares similarity with the one adopted in \cite{li2022non} (although we now need to accommodate non-symmetric random matrices), 
the approach outlined in \citet[Theorem 2]{li2022non} falls short of obtaining effective control the most challenging sample-limited regime, particularly when the denoising functions lack smoothness. 
For instance, when addressing the case of a sparse $\vstar$, \citet[Theorem 5]{li2022non} requires the number of observations to exceed $n \gtrsim k\log n$, with $k$ the sparsity of the true signal. 
This requirement arises because in \citet{li2022non}, each direction of the residual terms is treated equivalently and its $\ell_{2}$ norm is then controlled directly.
However, if our goal is to handle the most challenging scenario where $k$ is of the same order of $n$ and $p$, a more fine-grained control over $\xi_{t}$ and $\zeta_{t}$ along different directions become imperative.  
More specifically, the residual term $\xi_{t+1}$ turns out to have a larger degree of growth along the direction $G_{t}(s_t)$, and therefore, it makes sense to single out this direction and control its corresponding size separately as in Theorem~\ref{thm:residuals}. 
In Theorem~\ref{prop:final}, 
we single out the quantities $\frac{1}{n^2}\mathbb{E}\big[\|F_t'\|_2^2\big]\mathbb{E}\big[\|G'_t\|_2^2\big]$ to help control how error terms propagate across iterations;  
in our specific examples, we have demonstrated that this new approach of controlling residuals allows for more effective bounding of this factor. 

\section{Discussion}

In this paper, we have established a general recipe for understanding the non-asymptotic distributions for the celebrated AMP algorithm, tailored to sparse and robust regression. 
Our framework decomposes the AMP iterates into Gaussian random vectors and residual terms with explicit expressions that are tractable under some mild conditions. 
For both sparse and robust regression, 
our results have provided the first finite-sample distributional guarantees for the AMP iterates that can accommodate up to a polynomial number of iterations, 
which is in sharp contrast to prior theory that cannot go beyond $o\big(\frac{\log n}{\log \log n}\big)$ iterations. 
Furthermore, our theory has led to to improved distributional guarantees (i.e., improved error rates) for 
the optimally-tuned Lasso and the robust M-estimator compared to other existing approaches.  
The insights offered by our non-asymptotic analysis framework have improved upon prior works based on asymptotic state-evolution-type analysis.


Before concluding this paper, let us highlight several possible directions worthy of future investigation. 

\begin{itemize}
	\item Recall that our results have provided improved bounds for the residual terms; 
		for instance, when $t \asymp \log n$, our theory is able to bound the size of the residual terms by ${O}(\log n/n^{1/3})$ for both sparse and robust regression. 
		A natural question is concerned with the tightness of this error bound. 
		Our current conjecture is that the sharp bound on the residual terms should be  
	${O}(\poly (\log n)/n^{1/2})$; establishing or disproving this conjecture require more delicate analyses that go beyond the present analyses.  

	\item 
	Thus far, our framework confirms the validity of Gaussian approximation of AMP up to $O(n/\poly(\log n))$ iterations. 
	It remains to understand the behavior of AMP as 
	the number of iterations further increases beyond this range. 
		Will the state evolution recursion continue to provide reliable predictions as $t$ continues to grow?


	\item Finally, there has been a recent surge of interest in understanding the performances of AMP beyond the i.i.d.~Gaussian design. 
		Certain {\em universality} phenomena have been empirically observed and theoretically investigated \citep{bayati2015universality,chen2021universality,wang2022universality,dudeja2022universality}.  
		For instance, the asymptotic theory for AMP has been extended to accommodate 
	a family of rotationally invariant designs \citep{fan2022approximate,mondelli2021pca,cademartori2023non,venkataramanan2021estimation}. 
	Whether our results can be further generalized beyond Gaussian designs remains an interesting open question for future studies.

\end{itemize}


\vspace{0.5cm}

\begin{center}
    {\large APPENDIX}
\end{center}

\appendix

\section{Proof for our general results}


We present the proofs of Theorem~\ref{thm:main-matrix}, \ref{thm:residuals} and \ref{prop:final} together in this section and defer other technical details and lemmas to the appendices. On the high level, the proof of Theorem~\ref{thm:main-matrix} resembles the proof of \cite[Theorem 1]{li2022non} and the proofs of Theorem~\ref{thm:residuals} and \ref{prop:final} are based on a crucial higher order decomposition and a fine-grained control the residual terms.

\subsection{Proof of Theorem~\ref{thm:main-matrix}}
\label{sec:pf-theorem-main}

\paragraph{Step 1: constructing a key set of auxiliary sequences.}
Let us first introduce a sequence of auxiliary vectors/matrices $\{a_k, b_k, X_k\}_{1 \leq t\leq \min\{n,p\}}$ in a recursive fashion as below:
\begin{itemize}
 	\item[(i)] With our design matrix $X$ and the initialization $\{s_{1}, \beta_{1}\}$ in place, we define 
 	\begin{align}
 	\label{eqn:a1-b1}
	a_1 \defn \frac{G_1(s_1)}{\lt\|G_1(s_1)\rt\|_2} \in \real^n,
	\qquad 
	b_1 \defn \frac{F_1(\beta_1)}{\lt\|F_1(\beta_1)\rt\|_2}\in \real^p, \qquad\text{and}\qquad X_1 \defn X \in \real^{n\times p}; 
	\end{align}

	\item[(ii)] For every $2 \leq t < \min\{p, n\}$, concatenating the $a_{k}$'s and $b_{k}$'s into matrices 
	$U_{t-1} = [a_k]_{1 \le k \leq t-1} \in \real^{n\times (t-1)}$, $V_{t-1} = [b_k]_{1 \le k \leq t-1} \in \real^{p\times (t-1)}$, we can further define 
	\begin{subequations}
	\label{eqn:defn-abx}
		\begin{align}
		a_t &\defn \frac{\lt(I - U_{t-1}U_{t-1}^{\top}\rt)G_{t}(s_{t})}{\lt\|\lt(I - U_{t-1}U_{t-1}^{\top}\rt)G_{t}(s_{t})\rt\|_2}, \label{eqn:defn-abx-a}\\
		b_t &\defn \frac{\lt(I - V_{t-1}V_{t-1}^{\top}\rt)F_{t}(\beta_{t})}{\lt\|\lt(I - V_{t-1}V_{t-1}^{\top}\rt)F_{t}(\beta_{t})\rt\|_2}, \label{eqn:defn-abx-b}\\
		X_t &\defn \lt(I_n - a_{t-1}a_{t-1}^{\top}\rt)X_{t-1}\lt(I_p - b_{t-1}b_{t-1}^{\top}\rt),
		\end{align}
	\end{subequations}
	where the pair $(s_t, \beta_t)$ is generated by iteration~\eqref{eqn:AMP-general}.
 \end{itemize} 
By virtual of these definitions above, it is easily seen that vectors $\{a_k\}_{1\leq k\leq \min\{n,p\}}$ form an orthonormal basis and so are $\{b_k\}_{1\leq k\leq \min\{n,p\}}.$
By construction, $G_{t}(s_{t})$ lies in the span of $\{a_1, \ldots, a_t\}$ and similarly, 
$F_{t}(\beta_{t}) \in \textsf{span}\{b_1, \ldots, b_t\}$.
It is therefore legitimate to write 
\begin{subequations}
\label{eqn:orth-expa}
\begin{align}
	G_{t}(s_{t}) &= \sum_{k = 1}^{t} \alpha_{t}^ka_k, \qquad\text{for }\alpha_{t}^k \defn \lt\langle G_{t}(s_t), a_k\rt\rangle \qquad (1\leq k \leq t), \\
	F_{t}(\beta_{t}) &= \sum_{k = 1}^{t} \gamma_{t}^kb_k, \qquad\text{for }\gamma_{t}^k \defn \lt\langle F_{t}(\beta_{t}), b_k\rt\rangle\qquad (1\leq k \leq t),
\end{align}
which satisfies 
\begin{align*}
		\|\gamma_{t}\|_2 = \|F_t(\beta_t)\|_2\qquad\text{and}\qquad\|\alpha_{t}\|_2 = \|G_t(s_t)\|_2.
\end{align*}
\end{subequations}

\paragraph{Step 2: deriving distributional properties of $X_{k}b_{k}$ and $X_{k}^\top a_{k}$.} Next, we aim to establish some distributional characterizations of $X_{k}b_{k}$ and $X_{k}^\top a_{k}$. Towards this end, let us first consider another set of auxiliary vectors defined as below
\begin{subequations}
\label{eqn:def-phi-psi}
\begin{align}
\label{eqn:def-phi-k}
	\phi_k &= X_kb_k + \sum_{i = 1}^{k - 1} g^i_k a_i, \\
\label{eqn:def-psi-k}
	\psi_k &= \lt(I - b_{k}b_{k}^{\top}\rt)X_k^{\top}a_k + \sum_{i = 1}^{k} q^i_k b_i,
\end{align}
\end{subequations}
where each $g_i^{k}$ is i.i.d. generated from $\mathcal{N}(0,\frac{1}{n}).$
It turns out that $\phi_{k}$ and $\psi_{k}$ admit clean distributional guarantees summarized in the following lemma. 
\begin{lems}
\label{lem:phi-psi-distr}
	With $\{a_k,b_k,X_k\}_{1\leq k\leq \min{\{n,p\}}}$ defined in~\eqref{eqn:defn-abx}, it obeys 
	\begin{align*}
	 	\phi_k \stackrel{\text{i.i.d.}}{\sim} \mathcal{N}\left(0,\frac{1}{n}I_n\right),
	 	\qquad
	 	\text{and }~\psi_k \stackrel{\text{i.i.d.}}{\sim} \mathcal{N}\left(0,\frac{1}{n}I_p\right),
	 \end{align*} 
	 for all $1\leq k\leq \min{\{n,p\}}.$
\end{lems}
\noindent The proof of this result is postponed to Section~\ref{pf:lem:phi-psi-distr}.
We make note here that the covariance matrices of both $\phi_{k}$ and $\psi_{k}$ are identity matrices with normalized constant $1/n.$

\paragraph{Step 3: establishing two key decompositions as in \eqref{defi:dynamics}.} 
Let us start by showing relation~\eqref{defi:dynamics-s}.
First, we find it helpful to express $X_{1}$ as 
\begin{align}
\label{eqn:X-expa1}
	X_1 = X_t + \sum_{k = 1}^{t-1}(X_k - X_{k+1})
	= X_t + \sum_{k = 1}^{t-1} \Big[X_kb_kb_k^{\top} + a_ka_k^{\top}X_k\lt(I - b_{k}b_{k}^{\top}\rt)\Big].
\end{align}
For every $t\geq 1$, plugging the expansions (as in~\eqref{eqn:orth-expa}) that $F_{t}(\beta_{t}) = \sum_{k = 1}^{t} \gamma_{t}^kb_k$ and $G_{t-1}(s_{t-1}) = \sum_{k = 1}^{t-1} \alpha_{t-1}^ka_k$ leads to
\begin{align}
\label{eqn:st-xk-bk-ak}
\notag	s_{t} &= X_1F_{t}(\beta_t) - \langle F_t^{\prime}\rangle\sum_{k = 1}^{t-1} \alpha_{t-1}^ka_k \\
	&= \sum_{k = 1}^{t} \gamma_{t}^kX_kb_k + \sum_{k = 1}^{t-1} a_k\lt[\lt\langle \lt(I - b_{k}b_{k}^{\top}\rt)X_k^{\top}a_k, F_{t}(\beta_t)\rt\rangle - \langle F_t^{\prime}\rangle \alpha_{t-1}^k\rt],
\end{align}
where the last relation invokes the decomposition~\eqref{eqn:X-expa1}.
Substitution of the definition for $\phi_{k}$ and reorganizing terms further yield 
\begin{align}
\label{eqn:s-xi}
	\notag s_{t} &= \sum_{k = 1}^{t} \gamma_{t}^k\Big(\phi_k - \sum_{i = 1}^{k - 1} g_i^ka_i\Big) 
	+ \sum_{k = 1}^{t-1} a_k \Bigg[\Big\langle \psi_k - \sum_{i = 1}^{k} q_i^kb_i, F_{t}(\beta_t)\Big\rangle - \langle F_t^{\prime}\rangle \alpha_{t-1}^k\Bigg] \\
	&= \sum_{k = 1}^{t} \gamma_{t}^k\phi_k + 
	\underbrace{\sum_{k = 1}^{t-1} a_k\lt[\lt\langle \psi_k, F_{t}(\beta_t)\rt\rangle - \langle F_t^{\prime}\rangle \alpha_{t-1}^k - \sum_{i = 1}^{k} \gamma_t^iq_i^k - \sum_{i = k+1}^{t} \gamma_t^ig_k^i\rt]}_{=:\xi_t}.
\end{align}
As a consequence, we have established~\eqref{defi:dynamics-s} with $\xi_t \in \textsf{span}\{a_1,\ldots,a_{t-1}\}.$

As for property~\eqref{defi:dynamics-beta}, it is useful to write 
\begin{align}
\label{eqn:X-expa2}
	X_1 = X_t\lt(I - b_{t}b_{t}^{\top}\rt) + X_tb_{t}b_{t}^{\top} + \sum_{k = 1}^{t-1} \lt[X_kb_kb_k^{\top} + a_ka_k^{\top}X_k\lt(I - b_{k}b_{k}^{\top}\rt)\rt].
\end{align}
Again, invoking the expansions $F_{t}(\beta_{t}) = \sum_{k = 1}^{t} \gamma_{t}^kb_k$ and $G_{t}(s_{t}) = \sum_{k = 1}^{t} \alpha_{t}^ka_k$ gives 
\begin{align}
\notag \beta_{t+1} &= X_1^{\top}G_{t}(s_{t}) - \langle G_t^{\prime}\rangle\sum_{k = 1}^{t} \gamma_{t}^kb_k \\
\label{eqn:betat-xk-bk-ak}	&= \sum_{k = 1}^{t} \alpha_{t}^k\lt(I - b_{k}b_{k}^{\top}\rt)X_k^{\top}a_k + \sum_{k = 1}^{t} b_k\lt[\langle X_k^{\top}b_k, G_{t}(s_{t})\rangle - \langle G_t^{\prime}\rangle \gamma_{t}^k\rt] \\
\notag	&= \sum_{k = 1}^{t} \alpha_{t}^k\lt(\psi_k - \sum_{i = 1}^{k} q_i^kb_i\rt) + \sum_{k = 1}^{t} b_k\lt[\lt\langle \phi_k - \sum_{i = 1}^{k - 1} g_i^ka_i, G_{t}(s_t)\rt\rangle - \langle G_t^{\prime}\rangle \gamma_{t}^k\rt] \\
	&= \sum_{k = 1}^{t} \alpha_{t}^k\psi_k 
	+ 
	\underbrace{\sum_{k = 1}^{t} b_k\lt[\lt\langle \phi_k, G_{t}(s_t)\rt\rangle - \langle G_t^{\prime}\rangle \gamma_{t}^k - \sum_{i = 1}^{k - 1} \alpha_t^ig_i^k - \sum_{i = k}^{t} \alpha_t^iq_k^i\rt]}_{=: \zeta_t},\label{eqn:beta-zeta}
\end{align}
where the penultimate line uses the definitions of $\phi_{k}$ and $\psi_{k}$ (as of \eqref{eqn:def-phi-psi}).
Therefore, inequality \eqref{defi:dynamics-beta} holds with $\zeta_t \in \textsf{span}\{b_1,\ldots,b_{t-1}\}.$

\subsection{Proof of Theorem~\ref{thm:residuals}}
\label{sec:pf-theorem-second}

We move on to the proof of Theorem~\ref{thm:residuals}.

\paragraph{Controlling the residual terms $\xi_t$ and $\zeta_t$.}
In view of the definition of $\xi_{t}$ (cf.~\eqref{eqn:s-xi}), let us write
\begin{align*}
	\xi_t - \sum_{k = 1}^{t-1} a_k\lt[\lt\langle \psi_k, F_{t}(\beta_t)\rt\rangle - \langle F_t^{\prime}(\beta_t)\rangle \alpha_{t-1}^k\rt] 
	=
	\sum_{k = 1}^{t-1} a_k\Big[\sum_{i = 1}^{k} \gamma_t^iq_i^k - \sum_{i = k+1}^{t} \gamma_t^ig_k^i\Big].
\end{align*}
We aim to control the magnitude of the right hand side from above. First, we recall that the $q_i^{k}$'s and $g_k^{i}$'s are independently drawn from $\mathcal{N}(0,\frac{1}{n})$ --- independently of the randomness in the system, as a means to ensure the distributional characterization in Lemma~\ref{lem:phi-psi-distr}.
Towards bounding this quantity, recognizing that $\{a_1,\ldots,a_{t-1}\}$ forms an orthonormal basis, there exist a unit vector $\mu_t = [\mu_t^k]_{1\leq k \leq t} \in \mathbb{R}^{t}$ where 
\begin{align*}
	\Big\|\sum_{k = 1}^{t-1} a_k\Big[\sum_{i = 1}^{k} \gamma_t^iq_i^k - \sum_{i = k+1}^{t} \gamma_t^ig_k^i\Big]\Big\|_2
	=
	\sum_{k = 1}^{t-1} \mu_t^k \Big(\sum_{i = 1}^{k} \gamma_t^iq_i^k - \sum_{i = k+1}^{t} \gamma_t^ig_k^i\Big).
\end{align*}
Therefore, this quantity can be handled via standard Gaussian concentration inequalities as detailed in \cite[Lemma 3]{li2022non}. 
We now state the result directly without repeating its proof. With probability at least $1 - O(n^{-11})$, it satisfies 
\begin{align}
	\xi_t &= \sum_{k = 1}^{t-1} a_k\lt[\lt\langle \psi_k, F_{t}(\beta_t)\rt\rangle - \langle F_t^{\prime}(\beta_t)\rangle \alpha_{t-1}^k\rt]  + O\lt(\sqrt{\frac{t\log n}{n}}\|\gamma_t\|_2\rt).
\end{align}
By similar analysis, we are also ensured that with probability at least $1 - O(p^{-11})$, 
\begin{align}
    \zeta_t &= \sum_{k = 1}^{t} b_k\lt[\lt\langle \phi_k, G_{t}(s_t)\rt\rangle - \langle G_t^{\prime}(s_t)\rangle \gamma_{t}^k \rt] + O\lt(\sqrt{\frac{t\log n}{n}}\|\alpha_t\|_2\rt)
\end{align}
holds true.

\paragraph{Establishing the expansions~\eqref{eqn:residual-second}.}
Let us start with the term $\xi_{t}$. For every $t\geq 0$, first recall that
\begin{align*}
\xi_{t+1} 
&= \sum_{k = 1}^{t} a_k\lt[\lt\langle \psi_k, F_{t+1}(\beta_{t+1})\rt\rangle - \langle F_{t+1}^{\prime}(\beta_{t+1})\rangle \alpha_{t}^k\rt] + O\lt(\sqrt{\frac{t\log n}{n}}\|\gamma_{t+1}\|_2\rt) \\
&= \underbrace{\sum_{k = 1}^{t} a_k \lt[\langle \psi_k, F_{t+1}(\beta_{t+1}) - F_{t+1}(\hatbeta_{t+1})\rangle \rt]}_{=: \mathcal{R}_1} \\
&\qquad+
\sum_{k = 1}^{t} a_k \lt[\langle \psi_k, F_{t+1}(\hatbeta_{t+1})\rangle - \langle F_{t+1}^{\prime}(\beta_{t+1})\rangle \alpha_{t}^k\rt] + O\lt(\sqrt{\frac{t\log n}{n}}\|\gamma_{t+1}\|_2\rt).
%
\end{align*}
Intuitively, the magnitude of $\mathcal{R}_{1}$ is determined by the difference between $\beta_{t+1}$ and $\hatbeta_{t+1}$. If the expansion \eqref{eq:zeta_expansion} were true, the difference between $\beta_{t+1}$ and $\hatbeta_{t+1}$ arises from $\hatzeta_{t}$ as well as the difference between $\sum_{k = 1}^{t} \hatgamma_{t}^k F_{k}(\beta_{k})$ versus $\sum_{k = 1}^{t} \hatgamma_{t}^k F_{k}(v_{k})$. 

Next, if we project the term $\mathcal{R}_{1}$ to the direction that aligns with vector $G_t(s_t)$ and its orthogonal linear space, we end up with decomposition 
\begin{align*}
&\qquad \mathcal{R}_1 \\
&= \frac{1}{\|\alpha_t\|_2^2}\cdot G_t(s_t)^\top \mathcal{R}_1 \cdot G_t(s_t) + \mathcal{P}_{G_t(s_t)}^{\perp}\mathcal{R}_1\\
&= G_t(s_t) \cdot \frac{1}{\|\alpha_t\|_2^2} \sum_{k = 1}^{t} \alpha_t^k \lt[\langle \psi_k, F_{t+1}(\beta_{t+1}) - F_{t+1}(\hatbeta_{t+1})\rangle \rt]
+ \mathcal{P}_{G_t(s_t)}^{\perp} \sum_{k = 1}^{t} a_k\lt[\lt\langle \psi_k, F_{t+1}(\beta_{t+1})-  F_{t+1}(\hatbeta_{t+1})\rt\rangle\rt],
\end{align*}
where the first equality uses property $\|G_t(s_t)\|_2 = \ltwo{\alpha_t}$ and the second equality follows from the expansion~\eqref{eqn:orth-expa}.
In addition, due to the property of the orthogonal projection, we also find 
\begin{align}
\label{eqn:orthogonal}
\notag &\lt\langle G_{t}(s_{t}), \mathcal{P}_{G_t(s_t)}^{\perp} \sum_{k = 1}^{t} a_k\lt[\lt\langle \psi_k, F_{t+1}(\beta_{t+1})-  F_{t+1}(\hatbeta_{t+1})\rt\rangle\rt] \rt\rangle \\
	&= \lt\langle \sum_{k = 1}^{t} \alpha_{t}^ka_k, \mathcal{P}_{G_t(s_t)}^{\perp} \sum_{k = 1}^{t} a_k\lt[\lt\langle \psi_k, F_{t+1}(\beta_{t+1})-  F_{t+1}(\hatbeta_{t+1})\rt\rangle\rt] \rt\rangle
	= 0.
\end{align}

Putting the pieces above together, $\xi_{t+1}$ admits the following expression 
\begin{align}
\label{eqn:beethoven}
\notag \xi_{t+1} &= G_t(s_t)\cdot \frac{1}{\|\alpha_t\|_2^2} \Big\langle\sum_{k = 1}^{t} \alpha_t^k \psi_k, F_{t+1}(\beta_{t+1}) - F_{t+1}(\hatbeta_{t+1})\Big\rangle 
+ \sum_{k = 1}^{t} a_k\lt[\lt\langle \psi_k, F_{t+1}(\hatbeta_{t+1})\rt\rangle - \langle F_{t+1}^{\prime}(\beta_{t+1})\rangle \alpha_{t}^k\rt] \\
&\qquad + \mathcal{P}_{G_t(s_t)}^{\perp}\sum_{k = 1}^{t} a_k\lt\langle \psi_k, F_{t+1}(\beta_{t+1}) - F_{t+1}(\hatbeta_{t+1})\rt\rangle + O\lt(\sqrt{\frac{t\log n}{n}}\|\gamma_{t+1}\|_2\rt),
\end{align}
Contrasting what we have shown above with our target, it is sufficient to consider the second term above, which shall be done as follows. 


In the following, we establish the claim by decomposing the second term into two parts, corresponding to the influence of $\hatgamma_{t}^k F_{k}(\beta_{k})$'s and the randomness of Onsager term, respectively.

To simply the notation, let us define  
\begin{align*}
	A_t^k &\defn \lt\langle \psi_k, F_{t+1}(\hatbeta_{t+1})\rt\rangle - \langle F_{t+1}^{\prime}(\hatbeta_{t+1}) \rangle \alpha_{t}^k - \sum_{j = k+1}^{t} \hatgamma_{t}^j \lt\langle F_{t+1}^{\prime}(\hatbeta_{t+1}) \circ F_{j}^{\prime}\lt(v_{j}\rt)\rt\rangle\alpha_{j-1}^k.
\end{align*}
In view of this piece of notation and for each $i\geq 1$, $\sum_{k = 1}^{i} a_k\alpha_{i}^k = G_{i}(s_{i})$, a little algebra leads to 
\begin{align}
\label{eqn:violin-sonata}
\notag	&\sum_{k = 1}^{t} a_k\lt[\lt\langle \psi_k, F_{t+1}(\hatbeta_{t+1})\rt\rangle - \langle F_{t+1}^{\prime}(\beta_{t+1})\rangle \alpha_{t}^k\rt] \\
\notag	&= \sum_{k = 1}^{t} a_k\lt[\langle F_{t+1}^{\prime}(\hatbeta_{t+1}) - F_{t+1}^{\prime}(\beta_{t+1})\rangle \alpha_{t}^k + \lt\langle \psi_k, F_{t+1}(\hatbeta_{t+1})\rt\rangle - \langle F_{t+1}^{\prime}(\hatbeta_{t+1})\rangle \alpha_{t}^k\rt] \\
\notag	&= \sum_{k = 1}^{t} a_k\lt[\langle F_{t+1}^{\prime}(\hatbeta_{t+1}) - F_{t+1}^{\prime}(\beta_{t+1})\rangle \alpha_{t}^k + \sum_{j = k+1}^{t} \hatgamma_{t}^j \lt\langle F_{t+1}^{\prime}(\hatbeta_{t+1}) \circ F_{j}^{\prime}\lt(v_{j}\rt)\rt\rangle\alpha_{j-1}^k + A_t^k\rt] \\
	&=  G_{t}(s_{t}) \cdot  \langle F_{t+1}^{\prime}(\hatbeta_{t+1}) - F_{t+1}^{\prime}(\beta_{t+1})\rangle + \sum_{k = 1}^{t-1} \hatalpha_{t}^k G_{k}(s_{k}) + \sum_{k = 1}^{t} a_kA_t^k.
\end{align}
Here, the last relation follows from 
\begin{align*}
	\sum_{k = 1}^{t} a_k \sum_{j = k+1}^{t} \hatgamma_{t}^j \lt\langle F_{t+1}^{\prime}(\hatbeta_{t+1}) \circ F_{j}^{\prime}\lt(v_{j}\rt)\rt\rangle\alpha_{j-1}^k 
	&=
	\sum_{k = 1}^{t}\sum_{j = k+1}^{t}
	\widehat{\alpha}_t^{j-1} \alpha_{j-1}^k a_k \\
	&= 
	\sum_{j = 2}^{t} 
	\widehat{\alpha}_t^{j-1} \sum_{k =1}^{j-1} \alpha_{j-1}^k a_k
	=
	\sum_{k = 1}^{t-1} \hatalpha_{t}^k G_{k}(s_{k}),
\end{align*}
where we remind the readers that $\hatalpha_{t}^k$ is defined as of expression~\eqref{eq:xi-coeff}. 
Taking \eqref{eqn:violin-sonata} collectively with \eqref{eqn:beethoven} and recognizing the definition of $\hatalpha_{t}^t$, we end up with
\begin{align*}
	\xi_{t+1} &= \hatalpha_{t}^t G_{t}(s_{t}) + \sum_{k = 1}^{t-1} \hatalpha_{t}^k G_{k}(s_{k}) \\
	&\qquad + \sum_{k = 1}^{t} a_kA_t^k + \mathcal{P}_{G_t(s_t)}^{\perp}\sum_{k = 1}^{t} a_k\lt\langle \psi_k, F_{t+1}(\beta_{t+1}) - F_{t+1}(\hatbeta_{t+1})\rt\rangle + O\lt(\sqrt{\frac{t\log n}{n}}\|\gamma_{t+1}\|_2\rt),
\end{align*}
which validates the expansion~\eqref{eq:xi_expansion} for $t+1$. It is also worth noting that 
$\hatxi_{t+1}$ is defined exactly as the sum of the last three terms above.

When it comes to the expansion~\eqref{eq:zeta_expansion} at $t+1$, repeating a symmetric argument above leads to the required result. We omit its proof for brevity.

\subsection{Proof of Theorem~\ref{prop:final}}
\label{sec:pf-prop-final}

In order to prove this result, let us first state a key auxiliary lemma.  


\begin{lems}
\label{lem:magic-flute-simple}
Under the decomposition~\eqref{defi:dynamics} with~\eqref{eqn:residual-second} and Assumption~\ref{ass:control-lipschitz}, the Claim~\ref{claim:main}, stated below, holds for $t = 1$ with probability at least $1 - O(n^{-10})$. 
In addition, with probability at least $1 - O(n^{-10})$, for every 
\begin{align}
1 < t \lesssim \frac{n}{\log^4 n}, 
\end{align}
if the Claim~\ref{claim:main}, Assumption~\ref{ass:control-lipschitz} and \ref{ass:control-simple} all hold true for $t$, then Claim~\ref{claim:main} holds for $t+1$.
\end{lems}
\noindent The proof of this result can be found in Section~\ref{sec:pf-thm-magic-flute}.

\begin{claim}
\label{claim:main}
There exists universal constant $0 < c < 1$, such that 
the following set of inequalities hold true 
\begin{subequations} 
\label{eqn:main-simple}
\begin{align}
&\|\hatxi_{t}\|_2 \lesssim \sqrt{\frac{t\log^2 n}{n}}
\qquad\text{and}\qquad
\|\hatzeta_{t}\|_2 \lesssim \sqrt{\frac{t\log^2 n}{n}},\label{eqn:main-simple-A}\\
&\hatalpha_{t-1}^{t-1} \lesssim \lt(\frac{t\log^2 n}{n}\rt)^{\frac{1}{3}} 
\qquad\text{and}\qquad
\hatgamma_t^t \lesssim \lt(\frac{t\log^2 n}{n}\rt)^{\frac{1}{3}},
\label{eqn:main-simple-B}\\
&|\hatalpha_{t-1}^k| 
\le \left\{
\begin{array}{lll}
(1-c)^{t-k-1} \Big|\hatalpha_{\frac{t+k-1}{2}}^{\frac{t+k-1}{2}}\Big| & \text{if} & t-1-k = 2m, \\[0.4cm]
(1-c)^{t-k-2}\rho_F^2 \Big|\hatgamma_{\frac{t+k}{2}}^{\frac{t+k}{2}}\Big| & \text{if} & t-1-k = 2m+1,
\end{array}
\right. \label{eqn:main-simple-C}\\
&
|\hatgamma_{t}^k| \le \left\{
\begin{array}{lll}
(1-c)^{\frac{t-k}{2}}\Big|\hatgamma_{\frac{t+k}{2}}^{\frac{t+k}{2}}\Big| & \text{if} & t-k = 2m, \\
(1-c)^{\frac{t-k}{2}}\rho_G^2 \Big|\hatalpha_{\frac{t+k-1}{2}}^{\frac{t+k-1}{2}}\Big|& \text{if} & t-1-k = 2m.
\end{array} \label{eqn:main-simple-D}
\right.
\end{align}
\end{subequations}
\end{claim}

Based on this lemma, we first make the observation that if Claim~\ref{claim:main} holds true at iteration $t$, we arrive at 
\begin{subequations}
\begin{align}
\|\xi_{t}\|_2 &\le \sum_{k = 1}^{t-1} |\hatalpha_{t-1}^k| \ltwo{G_k(s_k)} + \ltwo{\hatxi_t} 
\lesssim \sum_{k=1}^{t-1} |\hatalpha_{t-1}^k|\cdot\|\alpha_k\|_2 + \sqrt{\frac{t\log^2 n}{n}} 
\lesssim \lt(\frac{t\log^2 n}{n}\rt)^{\frac{1}{3}}.
\end{align}
Here the first line invokes the relation that $\ltwo{G_t(s_t)} = \ltwo{\alpha_t}$, and the second line uses the inductive assumption~\eqref{eqn:main-simple-A} and the geometric decay of $\hatalpha_{t}^k$ in expression~\eqref{eqn:main-simple-C}. 
Similarly, we can deduce 
\begin{align}
\|\zeta_{t}\|_2 \lesssim \lt(\frac{t\log^2 n}{n}\rt)^{\frac{1}{3}}. 
\end{align}
\end{subequations}

Now if Assumptions~\ref{ass:control-lipschitz} and \ref{ass:control-simple} hold true over the execution of the AMP iterations, 
Claim~\ref{claim:main} is established by induction, since Lemma~\ref{lem:magic-flute-simple} validates both the initial condition and the inductive argument for Claim~\ref{claim:main}.

\section{Proof for sparse and robust regression}

\subsection{Proof of Theorem~\ref{thm:lasso}}
\label{sec:pf-thm-lasso}

The proof of this result is built upon Theorem~\ref{prop:final}. 
To show the residuals satisfy relation~\eqref{eqn:lasso-residual}, it is sufficient to validate Assumptions~\ref{ass:control-lipschitz} and \ref{ass:control-simple} over the execution of the AMP iterations.
We leave the arguments about the state evolution to Section~\ref{sec:lasso-se-pf}.

\subsubsection{Validating Assumption~\ref{ass:control-lipschitz}}
\label{sec:coro1-assmp1}

First, we make the direct observations that $F_{t}$ and $G_{t}$ defined in \eqref{eqn:lasso-F-G} satisfy $\rho_{F_{t}}, \rho_{G_{t}} = 1$ and $\rho_{1,F_{t}}, \rho_{1,G_{t}} = 0$ 
and $\lt\|F_{t}(0)\rt\|_2, \lt\|G_{t}(0)\rt\|_2 \lesssim 1$ given $\ltwo{\thetastar}, \ltwo{\varepsilon} \asymp 1.$
Then it is sufficient to verify that with high probability, 
\begin{align}
\label{eqn:norm-bound}
\|\gamma_t\|_2 \asymp \|\alpha_t\|_2 \asymp 1.
\end{align}

Towards this goal, recalling the initial choice where
\begin{align*}
	\beta_1 = -\thetastar 
	\qquad \text{ and } \qquad s_1 = Y - \varepsilon,
\end{align*}
and the norm property~\eqref{eqn:norm-condition}, we find $\|\gamma_1\|_2 = \ltwo{F_1(\beta_1)} = \|\theta^{\star}\|_2$.
In addition, notice that if $\tau_{t}$ is selected to be $\infty$, we observe 
\begin{align*}
\|r_t\|_2 = \|\varepsilon + X^{\top}\theta^{\star}\|_2 \lesssim 1.
\end{align*}
As $\tau_{t}$ is selected as the one that minimizes $\ltwo{r_t}$, it implies that
\begin{align} 
\label{eqn:gamma-norm-upper}
\Big\|\varepsilon + \sum_{j = 1}^t \gamma_t^j\phi_j+ \xi_t\Big\|_2 \lesssim 1.
\end{align}
Consequently, regarding $\alpha_t$, we bound 
\begin{align}
\label{eqn:trick}
\|\alpha_t\|_2 &= \|G_t(s_t)\|_2 = \Big\|\varepsilon + \sum_{j = 1}^{t} \gamma_{t}^j\phi_j + \xi_{t}\Big\|_2 \lesssim 1.
\end{align}
To further control the right hand side above, it is helpful to notice that 
\begin{align*}
\Big\|\varepsilon + \sum_{j = 1}^{t} \gamma_{t}^j\phi_j + \xi_{t}\Big\|_2 = \Big\|\varepsilon + \sum_{j = 1}^{t} \gamma_{t}^j\phi_j\Big\|_2 + O(\|\xi_{t}\|_2),
\end{align*}
and with probability at least $1-O(n^{-10})$,
\begin{align}
\label{eqn:dvorak-romantic}
\notag \lt\|\varepsilon + \sum_{j = 1}^{t} \gamma_{t}^j\phi_j\rt\|_2^2 &= \lt\|\varepsilon\rt\|_2^2 + \lt\|\sum_{j = 1}^{t} \gamma_{t}^j\phi_j\rt\|_2^2 + 2\varepsilon^{\top}\sum_{j = 1}^{t} \gamma_{t}^j\phi_j \\
&= \lt\|\varepsilon\rt\|_2^2 + \lt(1 + O\Big(\sqrt{\frac{t\log n}{n}}\Big)\rt)\lt\|\gamma_{t}\rt\|_2^2 + O\Big(\sqrt{\frac{t\log n}{n}}\lt\|\varepsilon\rt\|_2\lt\|\gamma_{t}\rt\|_2\Big),
\end{align}
where the last inequality invokes the spectral property as in \eqref{eqn:simple-rm-phi}.
Putting everything together, we arrive at
\begin{align}
\label{eqn:alpha-norm-bound}
\|\alpha_t\|_2 &= \sqrt{\|\gamma_t\|_2^2 + \|\varepsilon\|_2^2} + O\lt(\sqrt{\frac{t\log n}{n}}(\|\gamma_t\|_2 + \|\varepsilon\|_2) + \|\xi_t\|_2\rt).
%
\end{align}
Together with the relation~\eqref{eqn:trick} and $\ltwo{\varepsilon} \asymp 1$, the relation above 
leads to 
\begin{align}
\label{eqn:part-schubert}
	\|\alpha_t\|_2 \asymp 1, \qquad \ltwo{\gamma_t} \lesssim 1.
\end{align}

Finally, let us establish a proper lower bound for $\|\gamma_{t+1}\|_2$. 
Again, as a consequence of the norm relation~\eqref{eqn:norm-condition} and the Lipschitz property of soft-thresholding function, we write  
\begin{align}
\notag \|\gamma_{t+1}\|_2 = \|F_{t+1}(\beta_{t+1})\|_2  
&= \Big\|\theta^{\star} - \mathsf{ST}_{\tau_{t+1}}\Big(\theta^{\star} + \sum_{k = 1}^{t} \alpha_{t}^k\psi_k + \zeta_{t}\Big)\Big\|_2 \\
\notag &= \Big\|\theta^{\star} - \mathsf{ST}_{\tau_{t+1}}\Big(\theta^{\star} + \sum_{k = 1}^{t} \alpha_{t}^k\psi_k\Big)\Big\|_2 + O\Big(\|\zeta_{t}\|_2\Big) \\
%
&= \mathbb{E}\Big[\big\|\theta^{\star} - \mathsf{ST}_{\tau_{t+1}}(\theta^{\star} + \|\alpha_{t}\|_2g)\big\|_2 \mid \ltwo{\alpha_t} \Big] + O\lt(\sqrt{\frac{t\log n}{n}}+ \|\zeta_{t}\|_2\rt). \label{eqn:hippo}
%
\end{align}
Here, we invoke standard concentration inequality for Lipschitz function of Gaussian random variables \citep{borell1975brunn}.
To accommodate the randomness in $\alpha_{t}\in \real^{t}$, we take a union bound over a covering set of $\mathcal{S}^{t-1}$ of accuracy $\frac{1}{n}$.
Putting these ideas together, with probability at least $1 - O(n^{-10})$, it is ensured that 
\begin{align}
\Big\|\theta^{\star} - \mathsf{ST}_{\tau_{t}}\Big(\theta^{\star} + \sum_{j = 1}^{t} \alpha_{t}^j\psi_j\Big)\Big\|_2 
- \mathbb{E} \Bigg[\Big\|\theta^{\star} - \mathsf{ST}_{\tau_{t+1}}(\theta^{\star} + \ltwo{\alpha_t} g)\Big\|_2 \mid \ltwo{\alpha_t} \Bigg] = O\Big(\sqrt{\frac{t\log n}{n}}\Big).
\end{align}

Hence, it suffices to bound the right hand side of \eqref{eqn:hippo} from below which shall be down as follows. 
Towards this goal, for $\mu \asymp 1$, independent of $g$, if we define event  
\begin{align}
\label{eqn:set-e-hp}
\mathcal{E} := \lt\{g\sim \mathcal{N}(0, \frac{1}{n} I_p) ~|~ \big\|\theta^{\star} - \mathsf{ST}_{\tau_{t+1}}(\theta^{\star} + \mu g)\big\|_2 
= \mathbb{E}\lt[\big\|\theta^{\star} - \mathsf{ST}_{\tau_{t+1}}(\theta^{\star} + \mu g)\big\|_2 \mid \mu\rt] + O\lt(\sqrt{\frac{\log n}{n}}\rt)\rt\},	
\end{align}
as discussed above, the event $\mathcal{E}$ happens with probability at least $1 - O(n^{-10})$.
In view of this set, we write 
\begin{align}
\notag \mathbb{E}\Big[\lt\|\theta^{\star} - \mathsf{ST}_{\tau_{t+1}}(\theta^{\star} + \mu g)\rt\|_2^2\Big] 
&\stackrel{(\text{i})}= \mathbb{P}(\mathcal{E})\lt(\mathbb{E}\lt[\big\|\theta^{\star} - \mathsf{ST}_{\tau_{t+1}}(\theta^{\star} + \mu g)\big\|_2\rt] + O\lt(\sqrt{\frac{\log n}{n}}\rt)\rt)^2 \\
\notag &\qquad+ O\lt(\mathbb{E}\lt[(1 + \|g\|_2^2) \ind(\mathcal{E}^{\mathrm{c}})\rt]\rt)\\
\notag &= \lt(1-O\lt(n^{-10}\rt)\rt) \lt(\mathbb{E}\lt[\big\|\theta^{\star} - \mathsf{ST}_{\tau_{t+1}}(\theta^{\star} + \mu g)\big\|_2\rt] + O\lt(\sqrt{\frac{\log n}{n}}\rt)\rt)^2 \\
&\notag \qquad+ O\lt(\mathbb{P}(\mathcal{E}^{\mathrm{c}})\log n + \mathbb{E}\lt[(1 + \|g\|_2^2) \ind(\|g\|_2^2 \gtrsim \log n)\rt]\rt), \\
&\stackrel{(\text{ii})}{=} \lt(\mathbb{E}\lt[\big\|\theta^{\star} - \mathsf{ST}_{\tau_{t+1}}(\theta^{\star} + \mu g)\big\|_2\rt]\rt)^2  + O\lt(\sqrt{\frac{\log n}{n}}\rt), \label{eqn:penguin}
\end{align}
Here (i) uses the fact that
\begin{align*}
	\lt\|\theta^{\star} - \mathsf{ST}_{\tau_{t+1}}(\theta^{\star} + \mu g)\rt\|_2^2 
	& \lesssim \|\thetastar\|^2 +  \|\theta^{\star} + \mu g\|^2
	\lesssim 1 + \|g\|_2^2,
\end{align*}
by recognizing $\ltwo{\thetastar} \asymp$ and $\mu \asymp 1;$
(ii) invokes the basic relation for Gaussian random variable where $\mathbb{E}[\|g\|_2^2 \ind(\|g\|_2^2 \gtrsim \log n)] \lesssim \sqrt{\frac{\log n}{n}}$ for $g\sim \mathcal{N}(0,\frac{1}{n}I_p).$

Let us proceed to controlling the size of $\mathbb{E}[\lt\|\theta^{\star} - \mathsf{ST}_{\tau_{t+1}}(\theta^{\star} + \ltwo{\alpha_t} g)\rt\|_2^2]$ which in turn, 
provides the control of quantity $\mathbb{E}\big\|\theta^{\star} - \mathsf{ST}_{\tau_{t+1}}(\theta^{\star} + \|\alpha_{t}\|_2g)\big\|_2$. 
We claim that 
\begin{align}
\label{eqn:omelett}
	\mathbb{E}\Big[\lt\|\theta^{\star} - \mathsf{ST}_{\tau_{t+1}}(\theta^{\star} + \|\alpha_{t}\|_2g)\rt\|_2^2\Big] \gtrsim 1.
\end{align}
In the following, we prove the above claim by diving into two different cases and considering them separately. 
\begin{itemize}
	\item First consider the case when $\tau_{t+1}$ satisfies
\begin{align}
\label{eqn:tau-t-choice}
	\tau_{t+1} < \frac{\|\theta^{\star}\|_1}{4k} \leq \frac{\sqrt{k}\ltwo{\thetastar}}{k} \asymp \frac{1}{\sqrt{n}},
\end{align} 
for $\|\thetastar\|_1$ obeying \eqref{eqn:lasso-con2} and $n > 2k\log(p/k)$.
In this case, we find 
\begin{align*}
	\mathbb{E}\Big[\lt\|\theta^{\star} - \mathsf{ST}_{\tau_{t+1}}(\theta^{\star} + \|\alpha_{t}\|_2g)\rt\|_2^2 \Big]
	&\ge \mathbb{E}\Big[\lt\|\mathsf{ST}_{\tau_{t+1}}(\|\alpha_{t}\|_2g) \circ \ind(\theta^{\star} = 0)\rt\|_2^2\Big] \\
	&\gtrsim \mathbb{E} \Big[\Big\| |g| \ind(|g|\gtrsim 2\tau_{t+1}) \Big\|_2^2\Big]
	\gtrsim 1.
\end{align*}
\item On the other hand, when the relation~\eqref{eqn:tau-t-choice} is violated, we make the observation that 
\begin{align*}
	\|\theta^{\star}\|_1 \le 2\tau_{t+1}k + \|\theta^{\star}\|_2\sqrt{\|\ind(|\theta^{\star}| \ge 2\tau_{t+1})\|_0}, 
\end{align*}
which, together with \eqref{eqn:lasso-con2}, gives
\begin{align*}
	\|\ind(|\theta^{\star}| \ge 2\tau_{t+1})\|_0 \gtrsim k.
\end{align*}
Based on this property, we write 
\begin{align*}
\mathbb{E}\Big[\lt\|\theta^{\star} - \mathsf{ST}_{\tau_{t+1}}(\theta^{\star} + \|\alpha_{t}\|_2g)\rt\|_2^2\Big]
 &\ge \tau_{t+1}^2\mathbb{E}\Big[\lt\|\ind(|\theta^{\star}| \ge 2\tau_{t+1}) \circ \ind(\|\alpha_{t}\|_2|g| < \tau_{t+1})\rt\|_2^2 \Big] \gtrsim 1,
\end{align*}
where in the last inequality, we make the observation that $\tau_{t+1} > \lone{\thetastar}/4k \gtrsim \frac{1}{\sqrt{k}}$ and $\mprob(\ind(\|\alpha_{t}\|_2|g_i| < \tau_{t+1})) = O(1).$
\end{itemize}
Combining \eqref{eqn:hippo}, \eqref{eqn:penguin} and \eqref{eqn:omelett} leads to 
\begin{align}
	\|\gamma_{t+1}\|_2 \gtrsim 1, \qquad \text{for } t \lesssim \frac{n}{\log n}.
\end{align}
Taking this together with \eqref{eqn:part-schubert}, we have completed the proof of \eqref{eqn:norm-bound}
and thus justified Assumption~\ref{ass:control-lipschitz}.

\subsubsection{Validating Assumption~\ref{ass:control-simple}}
It is easily seen that $\frac{1}{n}\mathbb{E}\lt\|G_{t}^{\prime}(u_{t})\rt\|_2^2 = 1$.
Therefore validating Assumption~\ref{ass:control-simple} is equivalent to validating 
\begin{align} 
\label{eq:taustar}
\frac{1}{n}\mathbb{E}\lt\|F_{t+1}^{\prime}(v_{t+1})\rt\|_2^2 = \frac{1}{n}\mathbb{E}\lt\|\ind\lt(|\theta^{\star} + \|\alpha_{t}\|_2g| \ge \tau_{t+1}\rt)\rt\|_0 < (1-2c)^2,
\end{align}
for some constant $0 < c < 1/2.$
To establish this result, we find it helpful to first make the observation that for some small constant $c' > 0$
\begin{align}
\label{eqn:tau-set-re}
\mathcal{S}' \defn \lt\{\tau : -c'\sqrt{k} < \nabla_{\tau} \mathbb{E}\lt\|\theta^{\star} - \mathsf{ST}_{\tau}(\theta^{\star} + \|\alpha_{t}\|_2g)\rt\|_2^2 < c'\sqrt{k}\rt\} \subset \mathcal{S}
\end{align}
with set 
\begin{align*}
\mathcal{S} \defn \Big\{\tau : \mathbb{E}\Big[\lt\|\ind\lt(|(\theta^{\star} + \|\alpha_{t}\|_2g)| \ge \tau\rt)\rt\|_0\Big] < (1-2c)^2n\Big\}.
\end{align*}
Let us take the relation~\eqref{eqn:tau-set-re} as given for the moment, and come back to its proof at the end of this section. 
Based on this result, we shall prove Claim~\eqref{eq:taustar} by showing that
\begin{align}
\label{eqn:tau-s-prime}
	\tau_{t+1} \in \mathcal{S}'.
\end{align}

\paragraph{Proof of Claim~\eqref{eqn:tau-s-prime}.} 
We first prove that, if $\tau_{t+1} \notin \mathcal{S}'$, it must satisfy 
\begin{align}
\label{eqn:gap-tau}
	\inf_{\tau \in \mathcal{S}''} |\tau_{t+1} - \tau| > c'' \frac{1}{\sqrt{n}}, 
\end{align}
for some constant $c'' > 0.$
Here let us define an auxiliary set  
\begin{align}
	\mathcal{S}'' \defn \lt\{\tau : -\frac{c'}{2}\sqrt{k} < \nabla_{\tau} \mathbb{E}\lt\|\theta^{\star} - \mathsf{ST}_{\tau}(\theta^{\star} + \|\alpha_{t}\|_2g)\rt\|_2^2 < \frac{c'}{2}\sqrt{k}\rt\}.
\end{align}
In order to see this, observe that 
\begin{align}
\notag \nabla_{\tau} \mathbb{E}\lt\|\theta^{\star} - \mathsf{ST}_{\tau}(\theta^{\star} + \|\alpha_{t}\|_2g)\rt\|_2^2 &= \mathbb{E}\lt\langle\theta^{\star} - \mathsf{ST}_{\tau}(\theta^{\star} + \|\alpha_{t}\|_2g), \mathrm{sign}(\theta^{\star} + \|\alpha_{t}\|_2g)\ind(|\theta^{\star} + \|\alpha_{t}\|_2g| > \tau)\rt\rangle \\
&= \mathbb{E} \Big[\lt\langle\tau - \|\alpha_{t}\|_2g\mathrm{sign}(\theta^{\star} + \|\alpha_{t}\|_2g), \ind(|\theta^{\star} + \|\alpha_{t}\|_2g| > \tau)\rt\rangle\Big] \notag\\
&= \sum_{i:\theta_i^{\star}\ne 0}\mathbb{E}\Big[\lt(\tau - \|\alpha_{t}\|_2g\mathrm{sign}(\theta_i^{\star} + \|\alpha_{t}\|_2g)\rt) \ind(|\theta_i^{\star} + \|\alpha_{t}\|_2g| > \tau)\Big] \notag\\
&\qquad- (p-k)\mathbb{E}\lt[\mathsf{ST}_{\tau}(|\|\alpha_{t}\|_2g|)\rt].
\label{eqn:tau-tart-2}
\end{align}
The density function $|p_{\theta^{\star}_i + \|\alpha_{t}\|_2g_i}| \lesssim \sqrt{n}$ for $\ltwo{\alpha_t}\asymp 1$ and $g_{i} \sim \mathcal{N}(0,1/n)$, and therefore the right hand is a $O(n)$-Lipschitz function of $\tau$. Consequently, for $\tau\in \mathcal{S}''$, we deduce 
\begin{align*}
\frac{c'}{2} \sqrt{k}
\leq 
\Big|\nabla_{\tau_{t+1}} \mathbb{E}\lt\|\theta^{\star} - \mathsf{ST}_{\tau_{t+1}}(\theta^{\star} + \|\alpha_{t}\|_2g)\rt\|_2^2
-
\nabla_{\tau} \mathbb{E}\lt\|\theta^{\star} - \mathsf{ST}_{\tau}(\theta^{\star} + \|\alpha_{t}\|_2g)\rt\|_2^2\Big|
\lesssim n \inf_{\tau \in \mathcal{S}''} |\tau_{t+1} - \tau|,
\end{align*}
which proves the claimed gap~\eqref{eqn:gap-tau} between $\tau_{t+1}$ to $\mathcal{S}''$.


Given the relation~\eqref{eqn:gap-tau}, for $\widehat{\tau} = \inf\{\tau \in \mathcal{S}'', \tau > \tau_{t+1}\}$, this further implies
\begin{align}
\label{eqn:tau-t-conflict}
\notag &\mathbb{E}\Big\|\theta^{\star} - \mathsf{ST}_{\tau_{t+1}}(\theta^{\star} + \|\alpha_{t}\|_2g)\Big\|_2^2 - \inf_{\tau} \mathbb{E}\Big\|\theta^{\star} - \mathsf{ST}_{\tau}(\theta^{\star} + \|\alpha_{t}\|_2g)\Big\|_2^2 \\
\notag &\ge \mathbb{E}\Big\|\theta^{\star} - \mathsf{ST}_{\tau_{t+1}}(\theta^{\star} + \|\alpha_{t}\|_2g)\Big\|_2^2 - \mathbb{E}\Big\|\theta^{\star} - \mathsf{ST}_{\widehat{\tau}}(\theta^{\star} + \|\alpha_{t}\|_2g)\Big\|_2^2 \\
&\stackrel{(\text{i})}{\gtrsim} c'\sqrt{n}(\widehat{\tau} - \tau_{t+1}) \ge c'\sqrt{n} \cdot \inf_{\tau \in \mathcal{S}''} |\tau_{t+1} - \tau| \gtrsim 1,
\end{align}
where (i) is a consequence of the mean value theorem. 
Next we show that since $\tau_{t+1}$ is selected to minimize $\|\varepsilon + \sum_{k = 1}^{t+1} \gamma_{t+1}^k\phi_k + \xi_{t}\|_2$, the above relation contradicts with the choice of $\tau_{t+1}$. 
More specifically, as is shown in \eqref{eqn:alpha-norm-bound}, \eqref{eqn:hippo} and \eqref{eqn:penguin}, we can write 
\begin{align}
\label{eqn:dragon}
\notag &\Big\|\varepsilon + \sum_{k = 1}^{t+1} \gamma_{t+1}^k\phi_k + \xi_{t}\Big\|_2 \\
\notag &= 
\sqrt{\|\gamma_{t+1}\|_2^2 + \|\varepsilon\|_2^2} + O\lt(\sqrt{\frac{t\log n}{n}} + \|\xi_{t+1}\|_2\rt) \\
&= \sqrt{\mathbb{E}\Big[\lt\|\theta^{\star} - \mathsf{ST}_{\tau_{t+1}}(\theta^{\star} + \|\alpha_{t}\|_2g)\rt\|_2^2\Big] + \|\varepsilon\|_2^2 + O\lt(\sqrt{\frac{t\log n}{n}}+\|\xi_{t+1}\|_2\rt)} + O\lt(\sqrt{\frac{t\log n}{n}}+\|\xi_{t+1}\|_2\rt).
\end{align}
The threshold $\tau_{t+1}$ therefore cannot satisfy \eqref{eqn:tau-t-conflict}, as otherwise it does not minimize $\|\varepsilon + \sum_{k = 1}^{t+1} \gamma_{t+1}^k\phi_k + \xi_{t}\|_2$,  which in turn, validates the claimed relation~\eqref{eq:taustar}.

\paragraph{Proof of Property~\eqref{eqn:tau-set-re}.} 
It can be seen from numerical calculations (see Figure~\ref{fig:lasso} and the discussions around inequality~\eqref{eq:lasso-numerical1}) that for $G \sim \mathcal{N}(0, 1)$, if 
\begin{align} \label{eq:tau-tart-inv}
\sup_{\theta}\mathbb{E}\Big[\lt(\omega - G\mathrm{sign}(\theta + G)\rt) \circ \ind(|\theta + G| > \tau)\Big] - (p/k-1)\mathbb{E}\lt[\mathsf{ST}_{\omega}(|G|)\rt] \in (-c, c).
\end{align}
we have
\begin{align} 
\label{eq:tau-prob-inv}
1 - \frac{1 + (\frac{p}{k}-1)\mathbb{P}\lt(|G| \ge \omega\rt)}{2\log \frac{p}{k}} > 0.
\end{align}
The above result tells us that, for all $\tau$ satisfying
\begin{align}
\notag \nabla_{\tau} \mathbb{E}\lt\|\theta^{\star} - \mathsf{ST}_{\tau}(\theta^{\star} + \|\alpha_{t}\|_2g)\rt\|_2^2 
&= \sum_{i:\theta_i^{\star}\ne 0}\mathbb{E}\lt[\lt(\tau - \|\alpha_{t}\|_2g\mathrm{sign}(\theta_i^{\star} + \|\alpha_{t}\|_2g)\rt) \ind(|\theta_i^{\star} + \|\alpha_{t}\|_2g| > \tau)\rt] \notag\\
&\qquad- (p-k)\mathbb{E}\lt[\mathsf{ST}_{\tau}(|\|\alpha_{t}\|_2g|)\rt] \in (-c'\sqrt{k}, c'\sqrt{k}),
\label{eqn:tau-tart}
\end{align}
which implies that for $\omega := \frac{\sqrt{n}\tau}{\|\alpha_{t}\|_2}$,
\begin{align*}
k\sup_{\theta}\mathbb{E}\lt[\lt(\omega - G\mathrm{sign}(\theta + G)\rt) \ind(|\theta + G| > \tau)\rt] - (p-k)\mathbb{E}\lt[\mathsf{ST}_{\omega}(|G|)\rt] \in (-c'\sqrt{nk}/\|\alpha_{t}\|_2, c'\sqrt{nk}/\|\alpha_{t}\|_2),
\end{align*}
then we have 
\begin{align} \label{eq:tau-prob}
\mathbb{P}\lt(|\|\alpha_{t}\|_2g| \ge \tau\rt) = \mathbb{P}\lt(|G| \ge \omega\rt) \le \frac{2(1-4c)^2\log \frac{p}{k}-1}{\frac{p}{k}-1},
\end{align}
which establishes the property~\eqref{eqn:tau-set-re} immediately.
In order to see this, plugging in $n > 2k\log \frac{p}{k}$, inequality~\eqref{eq:tau-prob} ensures 
\begin{align*}
	\mathbb{P}\lt(|\|\alpha_{t}\|_2g| \ge \tau\rt) < \frac{(1-4c)^2n-k}{p-k},
\end{align*}
and hence, 
\begin{align*}
	\mathbb{E}\Big[\lt\|\ind\lt(|(\theta^{\star} + \|\alpha_{t}\|_2g)| \ge \tau\rt)\rt\|_0\Big] \le k + (p-k)\mathbb{P}\lt(|\|\alpha_{t}\|_2g| \ge \tau\rt) < (1-4c)^2n.
\end{align*}

\paragraph{Upper bound for $\tau_{t+1}$.}
Finally, let us establish a property of $\tau_{t+1}$. Specifically, we shall prove that $\tau_{t+1} \lesssim 1/\sqrt{n}$ when $p/k \lesssim 1$.
Before proceeding, let us make the following two observations both of which result from direct Gaussian integral. 
For $G\sim \mathcal{N}(0,1)$, it satisfies  
\begin{align*}
\mathbb{E}\lt[\mathsf{ST}_{\omega}(|G|)\rt] = \sqrt{\frac{2}{\pi}}\int_{\omega}^{\infty} x\exp\Big(-\frac{x^2}{2}\Big)\mathrm{d}x \lesssim \exp\Big(-\frac{\omega^2}{2}\Big),
\end{align*}
and
\begin{align*}
&\mathbb{E}\lt[\lt(\omega - G\mathrm{sign}(\theta_i + G)\rt) \circ \ind(|\theta_i + G| > \omega)\rt] \\
&= \mathbb{E}\lt[\lt(\omega - G\rt) \circ \ind(\omega - |\theta_i| < G < \omega + |\theta_i|)\rt] + 2\mathbb{E}\lt[\lt(\omega - G\rt) \circ \ind(G > \omega + |\theta_i|)\rt] \\
&\ge \frac{1}{\sqrt{2\pi}}\int_0^{|\theta_i|} x\exp\Big(-\frac{(x - \omega)^2}{2}\Big)\mathrm{d}x - \sqrt{\frac{2}{\pi}}\int_0^{\infty} x\exp\Big(-\frac{(x + \omega)^2}{2}\Big)\mathrm{d}x \\
&= \frac{1}{\sqrt{2\pi}}\exp\Big(-\frac{\omega^2}{2}\Big)\bigg\{\int_0^{|\theta_i|} x\exp(\omega x)\exp\Big(-\frac{x^2}{2}\Big)\mathrm{d}x - 2\int_0^{\infty} x\exp(-\omega x)\exp\Big(-\frac{x^2}{2}\Big)\mathrm{d}x\bigg\} \\
&\ge \frac{1}{\sqrt{2\pi}}\exp\Big(-\frac{\omega^2}{2}\Big)\bigg\{\int_0^{|\theta_i|} \frac{\omega^2x^3}{2}\exp\Big(-\frac{x^2}{2}\Big)\mathrm{d}x - 2\int_0^{\infty} x\exp\Big(-\frac{x^2}{2}\Big)\mathrm{d}x\bigg\} \\
&\ge \frac{1}{\sqrt{2\pi}}\exp\Big(-\frac{\omega^2}{2}\Big)\Big[\omega^2\Big(1 - \Big(\frac{\theta_i^2}{2} + 1\Big)\exp\Big(- \frac{\theta_i^2}{2}\Big)\Big) - 2\Big].
\end{align*}
Based on these two observations, for $\|\theta\|_1 \gtrsim \sqrt{k}$ and $\omega$ large enough, it obeys 
\begin{align}
\label{eqn:lj}
	\frac{1}{k}\sum_{i:\theta_i\ne 0}\mathbb{E}\lt[\lt(\omega - G\mathrm{sign}(\theta_i + G)\rt) \circ \ind(|\theta_i + G| > \omega)\rt] \gtrsim \omega^2\exp\Big(-\frac{\omega^2}{2}\Big).
\end{align}
Now, considering the transformation $\omega = \sqrt{n}\tau/\|\alpha_t\|_2$ and $\theta_i = \theta_i^{\star}/\|\alpha_t\|_2$, we obtain 
\begin{align*}
&\nabla_{\tau} \mathbb{E}\lt\|\theta^{\star} - \mathsf{ST}_{\tau}(\theta^{\star} + \|\alpha_{t}\|_2g)\rt\|_2^2 \\
&= \sum_{i:\theta_i^{\star}\ne 0}\mathbb{E}\lt[\lt(\tau - \|\alpha_{t}\|_2g\mathrm{sign}(\theta_i^{\star} + \|\alpha_{t}\|_2g)\rt) \ind(|\theta_i^{\star} + \|\alpha_{t}\|_2g| > \tau)\rt] - (p-k)\mathbb{E}\lt[\mathsf{ST}_{\tau}(|\|\alpha_{t}\|_2g|)\rt] \\
&= 
\frac{\ltwo{\alpha_t}}{\sqrt{n}} \sum_{i:\theta_i\ne 0}\mathbb{E}\lt[\lt(\omega - G\mathrm{sign}(\theta_i + G)\rt) \circ \ind(|\theta_i + G| > \omega)\rt] 
- \ltwo{\alpha_t} \frac{p-k}{\sqrt{n}}\mathbb{E}\lt[\mathsf{ST}_{\omega}(|G|)\rt].
\end{align*}
Taking the above together with \eqref{eqn:lj}, we have that given $\ltwo{\alpha_t} \asymp 1$ 
and $\tau > C/\sqrt{n}$ for some constant $C$ large enough, 
\begin{align*}
\nabla_{\tau} \mathbb{E}\lt\|\theta^{\star} - \mathsf{ST}_{\tau}(\theta^{\star} + \|\alpha_{t}\|_2g)\rt\|_2^2 > c\sqrt{k},
\end{align*}
for some constant $c > 0$. 
In view of the property~\eqref{eqn:tau-s-prime}, this property ensures that $\tau_{t+1} \lesssim 1/\sqrt{n}$.
As a result, it also leads to 
\begin{align}
\label{eqn:useful}
	\omega \defn \sqrt{n}\tau_{t+1}/\ltwo{\alpha_t} \lesssim 1,
	\qquad
	\mathbb{E}\lt[\mathsf{ST}_{\omega}(|G|)\rt] \gtrsim 1.
\end{align}

\subsubsection{State evolution}
\label{sec:lasso-se-pf}

Our final goal is to bound the difference between the non-asymptotic SE $(\alpha_t, \gamma_{t+1})$ to the deterministic SE defined in expression~\eqref{eqn:se-sparse}. 
In the following, we shall use the induction method to achieve this goal. 
In particular, it is easily validated that the set of relation~\eqref{eqn:pretty-se-lasso} holds true for $t=1$. Assuming that for some $t \geq 1$, 
\begin{align}
\big|\ltwo{\alpha_{t-1}}^2 - \alpha^{\star 2}_{t-1}\big| \lesssim \Big(\frac{t\log^2 n}{n}\Big)^{1/3}
\qquad \text{and} \qquad
	\big|\ltwo{\gamma_t}^2 - \gamma^{\star 2}_{t}\big| \lesssim \Big(\frac{t\log^2 n}{n}\Big)^{1/3},
\end{align}
it is thus sufficient to verify them for $t+1$.

Towards this end, let us first recall the expression~\eqref{eqn:alpha-norm-bound} that
\begin{subequations}
\label{eqn:se-lasso-finale}
\begin{align}
\label{eqn:se-alpha-lasso}
\|\alpha_t\|_2^2 &= \|\gamma_t\|_2^2 + \|\varepsilon\|_2^2 + O\lt(\Big(\frac{t\log^2 n}{n}\Big)^{\frac{1}{3}}\rt), 
\end{align}
where $\|\gamma_t\|_2, \|\varepsilon\|_2 \asymp 1$. 
In addition, combining expressions~\eqref{eqn:hippo} and~\eqref{eqn:penguin} yields  
\begin{align}
\label{eqn:se-gamma-lasso}
\|\gamma_{t+1}\|_2 
\notag & = \mathbb{E} \Bigg[\Big\|\theta^{\star} - \mathsf{ST}_{\tau_{t+1}}(\theta^{\star} + \ltwo{\alpha_t} g)\Big\|_2 \mid \ltwo{\alpha_t} \Bigg] 
+O\lt(\Big(\frac{t\log^2 n}{n}\Big)^{\frac{1}{3}}\rt)\\
\notag & \stackrel{(\text{i})}{=} \mathbb{E} \Bigg[\Big\|\theta^{\star} - \mathsf{ST}_{\tau_{t+1}}(\theta^{\star} + \sqrt{\|\gamma_t\|_2^2 + \|\varepsilon\|_2^2} g)\Big\|_2 \mid \ltwo{\alpha_t} \Bigg] +O\lt(\Big(\frac{t\log^2 n}{n}\Big)^{\frac{1}{3}}\rt)\\
& \stackrel{(\text{ii})}{=} \Big(\mathbb{E}\Big[\big\|\theta^{\star} - \mathsf{ST}_{\tau_{t+1}}\big(\theta^{\star} + \sqrt{\|\gamma_t\|_2^2 + \|\varepsilon\|_2^2}g\big)\big\|_2^2 \mid \ltwo{\alpha_t} \Big]\Big)^{1/2}
+ O\lt(\Big(\frac{t\log^2 n}{n}\Big)^{\frac{1}{3}}\rt).
\end{align}
\end{subequations}
Here, for inequality~$(\text{i})$, we have plugged in the relationship~\eqref{eqn:se-alpha-lasso} and invoked the Lipschitz property 
\begin{align*}
	\ltwo{\thetastar - \mathsf{ST}_{\tau_{t+1}}(\theta^{\star} + (\omega+\Delta) g}
	\leq 
	\ltwo{\thetastar - \mathsf{ST}_{\tau_{t+1}}(\theta^{\star} + \omega g)} + \ltwo{\Delta g}.
\end{align*}
We remind the readers that vector $g \in \mathcal{N}(0,\frac{1}{n}I_p)$, and is independent with the $(\alpha_t,\gamma_{t+1})$ sequence. 
For inequality~$(\text{ii})$ to hold, we recall the relation~\eqref{eqn:penguin}. 
According to the optimality of $\tau_{t+1}$ (in \eqref{eqn:se-sparse}), we also find
\begin{align}
\|\gamma_{t+1}\|_2 
&= \|F_{t+1}(\beta_{t+1})\|_2  
= \Big\|\theta^{\star} - \mathsf{ST}_{\tau_{t+1}}(\theta_{t+1})\Big\|_2  \\
&\le \big\|\theta^{\star} - \mathsf{ST}_{\tau_{t+1}^{\star}}\big(\theta_{t+1}\big)\big\|_2 + O\lt(\Big(\frac{t\log^2 n}{n}\Big)^{\frac{1}{3}}\rt) \notag\\
& \le \Big(\mathbb{E}\Big[\big\|\theta^{\star} - \mathsf{ST}_{\tau_{t+1}^{\star}}\big(\theta^{\star} + \sqrt{\|\gamma_t\|_2^2 + \|\varepsilon\|_2^2}g\big)\big\|_2^2 \mid \ltwo{\alpha_t} \Big]\Big)^{1/2}
+ O\lt(\Big(\frac{t\log^2 n}{n}\Big)^{\frac{1}{3}}\rt),
\label{eqn:se-gamma-lasso-star}
\end{align}
where the last line is derived again by uniform concentration inequalities, similar to relation~\eqref{eqn:se-gamma-lasso}.

Armed with the recursive formula~\eqref{eqn:se-lasso-finale}, controlling the difference between the non-asymptotic state evolution and its asymptotic analogue boils down to considering the growth of function  
\begin{align}
h_{\tau}(\mu) := \mathbb{E}\Big[\big\|\theta^{\star} - \mathsf{ST}_{\tau}(\theta^{\star} + \sqrt{\mu} g)\big\|_2^2\Big],
\end{align}
for every value of $\tau > 0$.
Direct computations yield 
\begin{align}
	\big|h'_{\tau}(\mu)\big| &= \frac{1}{\sqrt{\mu}}\Big|\mathbb{E}\Big[\big\langle\theta^{\star} - \mathsf{ST}_{\tau}(\theta^{\star} + \sqrt{\mu} g),  - \ind(|\theta^{\star} + \sqrt{\mu} g| > \tau) \circ g\big\rangle\Big]\Big|.
\end{align}
Considering the new rescaling
\begin{align}
\label{eqn:rescaling}
	\theta := \sqrt{n}\theta^{\star}/\sqrt{\mu}, 
	\qquad \text{and} \qquad \omega := \sqrt{n}\tau/\sqrt{\mu}, 
\end{align}
we can rewrite 
\begin{align*}
	h_{\tau}(\mu) := \frac{\mu}{n} 
	\mathbb{E}\Big[\big\|\theta - \mathsf{ST}_{\omega}(\theta + G)\big\|_2^2\Big].
\end{align*}
where $G\sim \mathcal{N}(0,1)$. 
In terms of the new scaling, for $k$-sparse $\thetastar$, some direct calculations lead to 
\begin{align} 
\big|h'_{\tau}(\mu)\big| 
&=
\frac{1}{n}\Big|\mathbb{E}\Big[\big\langle\theta - \mathsf{ST}_{\omega}(\theta + G),  \ind(|\theta + G| > \omega) \circ G\big\rangle\Big]\Big| \notag \\
&=
\Big|\frac{k}{n}\mathbb{E}\Big[\big(\mathsf{ST}_{\omega}(\theta + G)  - \theta\ind(|\theta + G| > \omega)\big) G\Big] + \frac{p-k}{n}\mathbb{E}\big[\mathsf{ST}_{\omega}(G)G\big]\Big|.
\end{align}
We claim that there exists constant $c \in (0,1)$ that only depends on the ratio $p/n$ and $k/p$ such that 
\begin{align}
\big|h'_{\tau}(\mu)\big|
&\le \Big|\frac{1}{n}\sum_{i:\theta_i\ne 0} \mathbb{E}\Big[\big(\mathsf{ST}_{\omega}(\theta_i + G)  - \theta\ind(|\theta + G| > \omega)\big) G\Big] + \frac{p-k}{n}\mathbb{E}\big[\mathsf{ST}_{\omega}(G)G\big]\Big| \le 1 - c, \label{eq:SE-lasso}
\end{align}
for both $\tau = \tau_{t+1}$ and $\tau = \tau_{t+1}^{\star}$.
Let us take this result as given for the moment and leave the proof of this claim to the end of this section. 

Given this result, we can bound the difference 
\begin{align}
\label{eqn:emp-recur}
\big|\|\gamma_{t+1}\|_2^2 - \gamma_{t+1}^{\star 2}\big| 
&\le \max_{\tau = \tau_{t+1}, \tau_{t+1}^{\star}}\big|h_{\tau}\Big(\|\gamma_t\|_2^2 + \|\varepsilon\|_2^2\big) - h_{\tau}\big(\gamma_{t}^{\star 2} + \|\varepsilon\|_2^2\big)\Big| + O\lt(\Big(\frac{t\log^2 n}{n}\Big)^{\frac{1}{3}}\rt) \notag\\
&\le (1 - c)\big|\|\gamma_{t}\|_2^2 - (\gamma_{t}^{\star})^2\big| + O\lt(\Big(\frac{t\log^2 n}{n}\Big)^{\frac{1}{3}}\rt),
\end{align}
as well as 
\begin{align}
\label{eqn:pop-recur}
	|\gamma_{t+1}^{\star 2} - \gamma_{t}^{\star 2}| &\le \max_{\tau = \tau_{t}^{\star}, \tau_{t+1}^{\star}}\Big|h_{\tau}\big(\gamma_{t}^{\star 2} + \|\varepsilon\|_2^2\big)
	- h_{\tau}\big(\gamma_{t-1}^{\star 2} + \|\varepsilon\|_2^2\big) \Big|
	\leq (1-c) \big|\gamma_{t}^{\star 2} - \gamma_{t-1}^{\star 2}\big|.
\end{align}
The last inequality ensures that sequence $\gamma_t^{\star}$ converges to some fixed point $\gamma^{\star}$. 
Recalling the initialization $\gamma_{1} = \gamma^\star_{1} = \ltwo{\theta^\star}$, invoking the relation~\eqref{eqn:emp-recur} recursively leads to 
\begin{align*}
	|\|\gamma_{t+1}\|_2^2 - (\gamma_{t+1}^{\star})^2| \lesssim \Big(\frac{t\log^2 n}{n}\Big)^{\frac{1}{3}}.
\end{align*}
Taking this together with \eqref{eqn:se-alpha-lasso} ensures that 
\begin{align*}
\big|\|\alpha_{t}\|_2^2 - (\alpha_{t}^{\star})^2\big| \lesssim \Big(\frac{t\log^2 n}{n}\Big)^{\frac{1}{3}}.
\end{align*}
In addition, the relation~\eqref{eqn:pop-recur} ensures $\gamma_{t}^{\star}$ converges to some $\gamma^{\star}$ exponentially with $t.$

\paragraph{Proof of Claim~\eqref{eq:SE-lasso}.}
We start by proving the claim~\eqref{eq:SE-lasso} for $\tau_{t+1}$. Let us recall that $\tau_{t+1}$ satisfies inequalities \eqref{eqn:tau-s-prime}. It thus yields 
\begin{align*}
	\bigg| \sum_{i:\theta_i^{\star}\ne 0}\mathbb{E}\Big[\lt(\tau_{t+1} - \|\alpha_{t}\|_2g\mathrm{sign}(\theta_i^{\star} + \|\alpha_{t}\|_2g)\rt) \ind(|\theta_i^{\star} + \|\alpha_{t}\|_2g| > \tau_{t+1})\Big] 
	- (p-k)\mathbb{E}\lt[\mathsf{ST}_{\tau_{t+1}}(|\|\alpha_{t}\|_2g|)\rt]\bigg| < c'\sqrt{k},
\end{align*}
for some small constant $c' > 0.$
It is thus sufficient to consider $\tau_{t+1}$ such that the above inequality holds true. 
Letting $\omega := \frac{\sqrt{n}\tau_{t+1}}{\|\alpha_{t}\|_2}$, the above relation further leads to  
\begin{align} 
\label{eq:tau-cond}
\bigg| \frac{1}{k}\sum_{i:\theta_i\ne 0}\mathbb{E}
\Big[\lt(\omega - G\mathrm{sign}(\theta_i + G)\rt) \circ \ind(|\theta_i + G| > \omega)\Big] - (p/k-1)\mathbb{E}\lt[\mathsf{ST}_{\omega}(|G|)\rt]\bigg| < c'', 
\end{align} 
for some small constant $c'' > 0$. 
Here, we remind the readers that in Section~\ref{sec:coro1-assmp1}, we have shown $\ltwo{\alpha_t} \asymp 1.$
Recall that we assume $p > 2.3k$ and $n > 2k\log \frac{p}{k}$.
To prove Claim~\eqref{eq:SE-lasso}, it is thus sufficient for us to show that for $e_3, c_3$ small enough,
\begin{align}
\label{eqn:pk-range}
2.3 < \frac{p}{k} = 1 + e_3 + \frac{\frac{1}{k}\sum_{i:\theta_i\ne 0}\mathbb{E}\lt[\lt(\omega - G\mathrm{sign}(\theta_i + G)\rt) \circ \ind(|\theta_i + G| > \omega)\rt]}{\mathbb{E}\lt[\mathsf{ST}_{\omega}(|G|)\rt]},
\end{align}
it satisfies  
\begin{align}
\label{eqn:lhs-tmp} 
\text{LHS} \defn 
&\Big|\frac{1}{2\log \frac{p}{k}}\frac{1}{k}\sum_{i:\theta_i\ne 0}\mathbb{E}\Big[\big(\mathsf{ST}_{\omega}(\theta_i + G)  - \theta\ind(|\theta_i + G| > \omega)\big) G\Big] + \frac{\frac{p}{k}-1}{2\log \frac{p}{k}}\mathbb{E}\big[\mathsf{ST}_{\omega}(G)G\big]\Big| < 1 - c_3.
\end{align}
It is easily seen that the above relation leads to the advertised bound~\eqref{eq:SE-lasso} by recognizing $n > 2k\log \frac{p}{k}$.

In order to prove the required inequality~\eqref{eqn:lhs-tmp}, 
plugging the expression for $\frac{p}{k}$ as in \eqref{eqn:pk-range} and in view of the the concavity of $\log(\cdot)$, 
we obtain 
\begin{align}
\text{LHS}
&\le \Bigg|\frac{\sum_{i:\theta_i\ne 0}\mathbb{E}\Big[\big(\mathsf{ST}_{\omega}(\theta_i + G)  - \theta\ind(|\theta_i + G| > \omega)\big) G + \frac{\mathbb{E}\big[\mathsf{ST}_{\omega}(G)G\big]}{\mathbb{E}\lt[\mathsf{ST}_{\omega}(|G|)\rt]}\lt[\lt(\omega - G\mathrm{sign}(\theta_i + G)\rt) \circ \ind(|\theta_i + G| > \omega)\rt]\Big]}
{2\sum_{i:\theta_i\ne 0}\log \Big(1 + \frac{\mathbb{E}\lt[\lt(\omega - G\mathrm{sign}(\theta_i + G)\rt) \circ \ind(|\theta_i + G| > \omega)\rt]}{\mathbb{E}\lt[\mathsf{ST}_{\omega}(|G|)\rt]}\Big)}
\Bigg| + e_4 \notag\\
&\le \sup_{\theta}\Bigg|\frac{\mathbb{E}\Big[\big(\mathsf{ST}_{\omega}(\theta + G)  - \theta\ind(|\theta + G| > \omega)\big) G + \frac{\mathbb{E}\big[\mathsf{ST}_{\omega}(G)G\big]}{\mathbb{E}\lt[\mathsf{ST}_{\omega}(|G|)\rt]}\lt[\lt(\omega - G\mathrm{sign}(\theta + G)\rt) \circ \ind(|\theta + G| > \omega)\rt]\Big]}
{2\log \Big(1 + \frac{\mathbb{E}\lt[\lt(\omega - G\mathrm{sign}(\theta + G)\rt) \circ \ind(|\theta + G| > \omega)\rt]}{\mathbb{E}\lt[\mathsf{ST}_{\omega}(|G|)\rt]}\Big)}
\Bigg| + e_4. \label{eq:lasso-numerical}
\end{align}

\begin{figure}[t]
\centering
\includegraphics[width=0.85\textwidth]{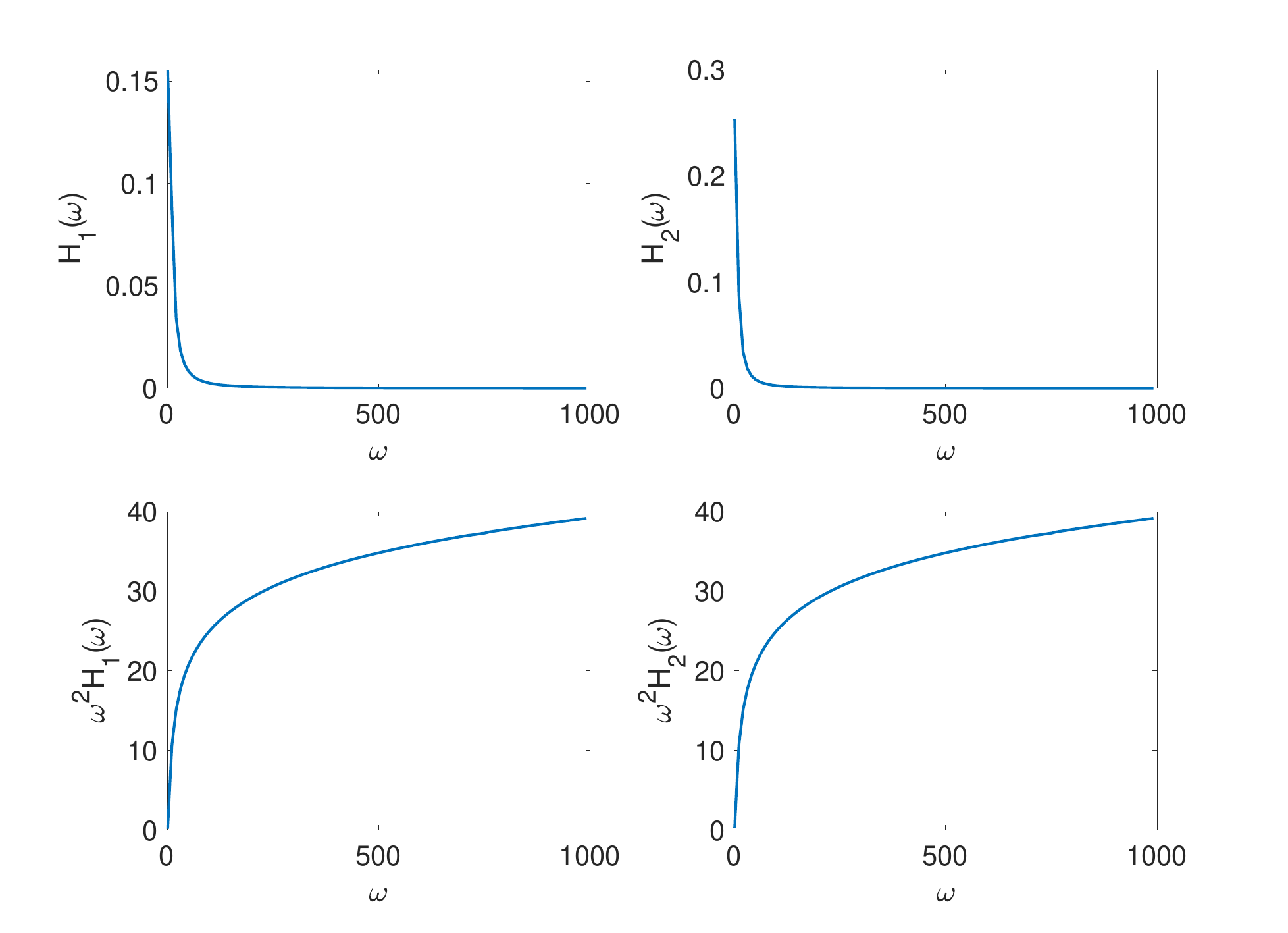}
	\caption{
	Numerical calculations for $H_1(\omega)$ and $H_2(\omega)$ of~\eqref{eq:defi-H-1} and~\eqref{eq:defi-H-2} such that $p/k \ge 2.3$. 
	}
\label{fig:lasso}
\end{figure}

Let us define 
\begin{align} 
\label{eq:defi-H-2}
H_2(\omega) := 1 - \sup_{\theta} \Bigg|\frac{\mathbb{E}\Big[\big(\mathsf{ST}_{\omega}(\theta + G)  - \theta\ind(|\theta + G| > \omega)\big) G + \frac{\mathbb{E}\big[\mathsf{ST}_{\omega}(G)G\big]}{\mathbb{E}\lt[\mathsf{ST}_{\omega}(|G|)\rt]}\lt[\lt(\omega - G\mathrm{sign}(\theta + G)\rt) \circ \ind(|\theta + G| > \omega)\rt]\Big]}
{2\log \Big(1 + \frac{\mathbb{E}\lt[\lt(\omega - G\mathrm{sign}(\theta + G)\rt) \circ \ind(|\theta + G| > \omega)\rt]}{\mathbb{E}\lt[\mathsf{ST}_{\omega}(|G|)\rt]}\Big)}
\Bigg|.
\end{align}
Figure~\ref{fig:lasso} demonstrates that for some constant $c_4 > 0$,
\begin{align*}
H_2(\omega) > \min\{c_4, ~\frac{c_4}{\omega^2}\}.
\end{align*}
Putting this together with the upper bound for $\omega$ in \eqref{eqn:useful} establishes \eqref{eqn:lhs-tmp}, and thus the Claim~\eqref{eq:SE-lasso} for $\tau_{t+1}$.

In addition, we make observations that $\tau_{t+1}^{\star}$ satisfies
\begin{align*}
\nabla_{\tau} \mathbb{E}\lt\|\theta^{\star} - \mathsf{ST}_{\tau_{t+1}^{\star}}(\theta^{\star} + \alpha_{t}^{\star}g)\rt\|_2^2 = 0,
\end{align*}
and $\alpha_{t}^{\star} = \|\alpha_{t}^{\star}\|_2 + o(1)$, which follows from~\eqref{eqn:se-alpha-lasso} and the induction condition.
As a result, $\tau_{t+1}^{\star}$ also satisfies the relation \eqref{eqn:tau-s-prime}, which by similar argument as above validates the Claim~\eqref{eq:SE-lasso} for $\tau_{t+1}^{\star}$.

\paragraph{Proof of inequality~\eqref{eq:tau-prob-inv}.}
Similarly, we make the observation that 
\begin{align}
\label{eq:lasso-numerical1}
	\frac{1 + \big(\frac{p}{k}-1\big)\mathbb{P}\lt(|G| \ge \omega\rt)}{2\log \frac{p}{k}} 
	&\le \sup_{\theta}\frac{1 + \frac{\mathbb{P}\lt(|G| \ge \omega\rt)\mathbb{E}\lt[\lt(\omega - G\mathrm{sign}(\theta + G)\rt) \circ \ind(|\theta + G| > \omega)\rt]}{\mathbb{E}\lt[\mathsf{ST}_{\omega}(|G|)\rt]}}{2\log \Big(1 + \frac{\mathbb{E}\lt[\lt(\omega - G\mathrm{sign}(\theta + G)\rt) \circ \ind(|\theta + G| > \omega)\rt]}{\mathbb{E}\lt[\mathsf{ST}_{\omega}(|G|)\rt]}\Big)} + c_4,
\end{align}
and define
\begin{align}
\label{eq:defi-H-1}
H_1(\omega) := 1 - \sup_{\theta}\frac{1 + \frac{\mathbb{P}\lt(|G| \ge \omega\rt)\mathbb{E}\lt[\lt(\omega - G\mathrm{sign}(\theta + G)\rt) \circ \ind(|\theta + G| > \omega)\rt]}{\mathbb{E}\lt[\mathsf{ST}_{\omega}(|G|)\rt]}}{2\log \Big(1 + \frac{\mathbb{E}\lt[\lt(\omega - G\mathrm{sign}(\theta + G)\rt) \circ \ind(|\theta + G| > \omega)\rt]}{\mathbb{E}\lt[\mathsf{ST}_{\omega}(|G|)\rt]}\Big)}.
\end{align}
As a consequence, Figure~\ref{fig:lasso} demonstrates the relation~\eqref{eq:tau-prob-inv} through a similar argument as above.

\subsection{Proof of Theorem~\ref{thm:robust}}
\label{sec:pf-thm-robust}


This result is again a consequence of Theorem~\ref{prop:final}. 
To establish the relation~\eqref{eqn:robust-residual}, we proceed by validating Assumptions~\ref{ass:control-lipschitz} and \ref{ass:control-simple} over the execution of the AMP iterations.
We derive the state evolution results in Section~\ref{sec:robust-se-pf}.

\subsubsection{Validating Assumption~\ref{ass:control-lipschitz}}

Before proceeding, we make a remark that $b_{t}$ is chosen according to \eqref{eq:bt_choice} such that
\begin{align}
\label{eqn:defn-kt}
	p(1 + 1/b_t) = \big\||r_t| < \lambda(1+b_t) \big\|_0 =: k_t. 
\end{align}
Such $b_{t}$ exists since since when $b_t \to \infty$, $k_t = n > p = \lim_{b_t \to \infty} p(1 + 1/b_t)$, and when $b_t \to 0$, $k_t \le n < \infty = \lim_{b_t \to 0} p(1 + 1/b_t)$.

Now we proceed to validate Assumption~\ref{ass:control-lipschitz} in this case. 
For denoising functions defined in~\eqref{eqn:robust-F-G}, requirement~\eqref{eq:conditions-simple-1} satisfies straightforwardly. 
In addition, we note that $\ltwo{F_t(0)} = 0$ and 
\begin{align*}
	\ltwo{G_t(0)} = \ltwo{g_t(\varepsilon)} 
	&= \frac{n b_t}{p(1 + b_t)}\lt\|\psi\lt(\varepsilon; \lambda(1+b_t)\rt)\rt\|_2 \\
	&\leq \frac{n b_t}{p(1 + b_t)}\ltwo{\varepsilon} \lesssim 1,
\end{align*}
Therefore, we only need to verify that $\|\gamma_t\|_2 \asymp \|\alpha_t\|_2 \asymp 1$.
We claim that this indeed the case.

\paragraph{Proof of $\|\gamma_t\|_2, \ltwo{\alpha_t} \lesssim 1$.} Specifically, if $\|\xi_t\|_2, \|\zeta_{t-1}\|_2 = o(1)$, we shall prove that there exists some large enough constant $\overline{\gamma}$ with $\|\varepsilon\|_2/\overline{\gamma} < c$ for some constant $c > 0$ small enough, such that
\begin{align}
\label{eqn:gamma-uni}
	\|\gamma_t\|_2 < \overline{\gamma}, \qquad \text{for every } t \geq 1.
\end{align}
It therefore implies $\|\gamma_t\|_2 \lesssim 1$.
We begin by noticing that $\|\gamma_1\|_2 = \|\theta^{\star}\|_2 \lesssim \overline{\gamma}$. 
In addition, with probability at least $1 - O(n^{-10})$, one has 
\begin{align}
\label{eqn:robust-alpha-sim}
\notag \|\alpha_t\|_2 &= \|G_t(s_t)\|_2 = \frac{n b_t}{p(1 + b_t)}\lt\|\psi\lt(s_t + \varepsilon; \lambda(1+b_t)\rt)\rt\|_2 \\
\notag &\le \frac{n}{p}\|s_t+ \varepsilon\|_2 = \frac{n}{p} \Big \| \sum_{k = 1}^{t} \gamma_{t}^k\phi_k + \xi_{t} + \varepsilon \Big\|_2 \\
&= \frac{n}{p} \sqrt{\|\gamma_t\|^2_2 + \|\varepsilon\|_2^2} +  O\Big(\sqrt{\frac{t\log n}{n}}(\|\gamma_t\|_2+\ltwo{\varepsilon}) + \ltwo{\xi_t}\Big),
\end{align}
where the last line follows from similar argument to \eqref{eqn:alpha-norm-bound}.
In addition, invoking the spectral properties as in \eqref{eqn:simple-rm} again gives 
\begin{align}
\label{eqn:robust-gamma-sim}
\|\gamma_{t+1}\|_2 &= \|F_{t+1}(\beta_{t+1})\|_2 
= \lt\|\sum_{k = 1}^{t} \alpha_{t}^k\psi_k + \zeta_{t}\rt\|_2 = \sqrt{\frac{p}{n}}\|\alpha_t\|_2 + O\lt(\sqrt{\frac{t\log n}{p}}\|\alpha_t\|_2 + \|\zeta_{t}\|_2\rt).
\end{align}
If relation~\eqref{eqn:gamma-uni} does not hold, then there exists some $t$ such that 
\begin{align*}
	\|\gamma_{t+1}\|_2 \ge \overline{\gamma} > \|\gamma_t\|_2.
\end{align*}
then as a consequence of \eqref{eqn:robust-alpha-sim} and \eqref{eqn:robust-gamma-sim}, one has $\|\gamma_t\|_2 \gtrsim \overline{\gamma}$.
By definition of $\overline{\gamma}$, this means $\|\varepsilon\|_2/\|\gamma_t\|_2 = o(1)$. 
In the following, we show this is impossible.

In order to see this, let us first consider quantity $\|\alpha_t\|_2$. 
Recognizing the Lipschitz property of the function $\ltwo{\psi(\cdot \mid \lambda(1+b_t))}$, we write 
\begin{align}
\label{eqn:alphat-allegro}
\notag \|\alpha_t\|_2 
&= \|G_t(s_t)\|_2 = \frac{n b_t}{p(1 + b_t)}\lt\|\psi\lt(s_t + \varepsilon; \lambda(1+b_t)\rt)\rt\|_2 \\
&= \frac{n b_t}{p(1 + b_t)}\Big\|\psi\Big(\sum_{k = 1}^{t} \gamma_{t}^k\phi_k; \lambda(1+b_t)\Big)\Big\|_2 + O(\|\varepsilon + \xi_t\|_2).
\end{align}

Conditional on $\ltwo{\gamma_t}$, we invoke the standard concentration result for Lipschitz function of Gaussian random variables \citep{borell1975brunn} to bound $\|\psi(\sum_{k = 1}^{t} \gamma_{t}^k\phi_k; \lambda(1+b_t))\|_2.$
To accommodate the randomness in $\gamma_{t}\in \real^{t}$, we take a union bound over a covering set of $\mathcal{S}^{t-1}$ of accuracy $\frac{1}{n}$.
Combining these ideas together, with probability at least $1 - O(n^{-10})$, it satisfies 
\begin{align}
\label{eqn:alphat-adagio}
	\Big\|\psi\Big(\sum_{k = 1}^{t} \gamma_{t}^k\phi_k; \lambda(1+b_t)\Big)\Big\|_2 
	- \mathbb{E}\Big[\|\psi(\|\gamma_t\|_2g; \lambda(1+b_t))\|_2 \mid \ltwo{\gamma_t}\Big] = O\Big(\|\varepsilon\|_2 + \sqrt{\frac{t\log n}{n}}\|\gamma_{t}\|_2\Big).
\end{align}
It therefore leads to 
\begin{align}
\label{eqn:alphat-rondo}
\|\alpha_t\|_2 
=
\mathbb{E}\Big[\|\psi(\|\gamma_t\|_2g; \lambda(1+b_t))\|_2 \mid \ltwo{\gamma_t}\Big] + O\Big(\|\varepsilon\|_2 + \sqrt{\frac{t\log n}{n}}\|\gamma_{t}\|_2\Big).
\end{align}

Recall the definition of $\psi(z; \lambda) = \min\{\max\{z, -\lambda\}, \lambda\}$ to obtain 
\begin{align}
	\notag \mathbb{E}\Big[\|\psi(\|\gamma_t\|_2g; \lambda(1+b_t))\|_2 \mid \ltwo{\gamma_t} \Big] 
	&= \|\gamma_t\|_2 \mathbb{E}\Big[\Big\|\psi(g; \frac{\lambda(1+b_t)}{\|\gamma_t\|_2})\Big\|_2\Big]\\
	&\leq \|\gamma_t\|_2\sqrt{n\mathbb{E}\Big[\min\lt\{\tilde{g}^2, \lambda^2(1+b_t)^2/\|\gamma_t\|_2^2\rt\}\Big]}.
\end{align}
Here, in the last equality, we denote $\tilde{g} \sim \mathcal{N}(0,\frac{1}{n})$. 
Therefore in order to control $\|\alpha_t\|_2$, it suffices to bound the quantity $\mathbb{E}[\min\{\tilde{g}^2, \lambda^2(1+b_t)^2/\|\gamma_t\|_2^2\}].$
We claim that it satisfies 
\begin{align}
\label{eqn:phenix}
\mathbb{E}\Big[\min\lt\{\tilde{g}^2, \lambda^2(1+b_t)^2/\|\gamma_t\|_2^2\rt\}\Big] \le (1-2c)^2\frac{p}{n^2}.
\end{align}

Putting everything together, we obtain 
\begin{align}
\label{eqn:robust-alpha-star}
	\ltwo{\alpha_t}\leq (1-2c)\ltwo{\gamma_t}\sqrt{\frac{n}{p}} 
	+ O\Big(\ltwo{\varepsilon + \xi_t} + \sqrt{\frac{t\log n}{n}}\|\gamma_{t}\|_2\Big).
\end{align}
Combining the relation~\eqref{eqn:robust-alpha-star} with \eqref{eqn:robust-gamma-sim}, we end up with 
\begin{align}
	\|\gamma_{t+1}\|_2 \le (1 - c)\|\gamma_{t}\|_2 + o(1) < \|\gamma_{t}\|_2
\end{align}
which is contradicted with $\|\gamma_{t+1}\|_2 \geq \overline{\gamma} > \|\gamma_{t}\|_2$. 
Hence, we conclude $\|\gamma_{t}\|_2 \lesssim 1$ for every $t\ge 1$,
which in turn leads to $\ltwo{\alpha_t}\lesssim 1$ by virtue of \eqref{eqn:robust-alpha-sim}.

\paragraph{Proof of inequality~\eqref{eqn:phenix}.}
First, we claim that the following relation holds true for $\tilde{g} \sim \mathcal{N}(0,\frac{1}{n})$, 
\begin{align}
\label{eqn:cat}
	\mprob\Big(|\tilde{g}| < \lambda(1+b_t)/\|\gamma_t\|_2 \Big) = \frac{p}{n} + o(1).
\end{align}
Let us take this relation as give for the moment and come back to its proof later. 
Based on this, we obtain
\begin{align*}
	\frac{n\mathbb{E}\Big[\min\lt\{\tilde{g}^2, \lambda^2(1+b_t)^2/\|\gamma_t\|_2^2\rt\}\Big]}{\mathbb{P}\Big(|\tilde{g}| < \lambda(1+b_t)/\|\gamma_t\|_2\Big)} \le 1-5c,
\end{align*} 
for some constant $c$ depending on $\frac{p}{n}$.
This can be seen from the numerical simulation in Figure~\ref{fig:robust}.
Here, we let $\tau := \sqrt{n}\lambda(1+b_t)/\|\gamma_t\|_2$ and 
\begin{align}\label{eq:defi-H1} 
H_1(\tau) := \frac{1}{\mathbb{P}\big(|G| > \tau\big)}\Big(1 - \frac{\mathbb{E}\big[\min\lt\{G^2, \tau^2\rt\}\big]}{\mathbb{P}\big(|G| < \tau\big)}\Big).
\end{align}
Putting these two things together finishes the proof of inequality~\eqref{eqn:phenix}.


Finally, we conclude by proving relation~\eqref{eqn:cat}. 
In view of the expression $r_t = \sum_{k = 1}^{t} \gamma_{t}^k\phi_k + \xi_{t} + \varepsilon$, first we make the observation that 
\begin{align*}
\Bigg|\Big\|\ind\lt(|r_t| < \lambda(1+b_t)\rt)\Big\|_0 
- \Big\|\ind\Big(|\sum_{k = 1}^{t} \gamma_{t}^k\phi_k| < \lambda(1+b_t)\Big)\Big\|_0\Bigg|
\leq 
\Bigg\|\ind\lt(|r_t| < \lambda(1+b_t)\rt) - \ind\Big(|\sum_{k = 1}^{t} \gamma_{t}^k\phi_k| < \lambda(1+b_t)\Big)\Bigg\|_1.
\end{align*}
Lemma~\ref{lem:cover}~\eqref{eqn:Fprime-o} bounds the $\ell_{1}$ discrepancy of $F_{t}'$ after pertrubing the input. By similar argument, we can derive perturbation results for $G_t'$. In particular, 
we bound 
\begin{align}
\label{eqn:piazzola}
	\Big\|\ind\lt(|r_t| < \lambda(1+b_t)\rt) - \ind\Big(|\sum_{k = 1}^{t} \gamma_{t}^k\phi_k| < \lambda(1+b_t)\Big)\Big\|_1
	\leq 
	t \log n + n \Big(\frac{\ltwo{\xi_t + \varepsilon}}{\ltwo{\gamma_t}}\Big)^{2/3}.
\end{align}
Combining the above two inequalities yields 
\begin{align}
\label{eqn:tosca}
\Bigg|\frac{1}{n}\Big\|\ind\lt(|r_t| < \lambda(1+b_t)\rt)\Big\|_0 - \frac{1}{n}\Big\|\ind\Big(|\sum_{k = 1}^{t} \gamma_{t}^k\phi_k|  < \lambda(1+b_t)\Big)\Big\|_0 \Bigg|
\leq 
\frac{t\log n}{n} + o(1),
\end{align}
where the last inequality uses the assumption $\ltwo{\xi_t} = o(1)$, $\ltwo{\varepsilon} \asymp 1$, and the relation $\ltwo{\gamma_t} \gtrsim \overline{\gamma}.$
Now we move on to control the term $\frac{1}{n}\|\ind(|\sum_{k = 1}^{t} \gamma_{t}^k\phi_k| < \lambda(1+b_t))\|_0.$ Some direct algebra leads to 
\begin{align}
\label{eqn:sibelius}
	\notag \frac{1}{n}\Big\|\ind\Big(|\sum_{k = 1}^{t} \gamma_{t}^k\phi_k| < \lambda(1+b_t)\Big)\Big\|_0
	&= \frac{1}{n}\lt\|\ind\lt(\frac{|\sum_{k = 1}^{t} \gamma_{t}^k\phi_k|}{\|\gamma_t\|_2} < \frac{\lambda(1+b_t)}{\|\gamma_t\|_2}\rt)\rt\|_0 \\
	\notag &\leq 
	\sup_{\zeta \in \mathcal{S}^{t-1}}\frac{1}{n}\Big\|\ind\Big(|\sum_{k = 1}^{t} \zeta^k\phi_k| < \frac{\lambda(1+b_t)}{\|\gamma_t\|_2}\Big)\Big\|_0 \\
	&\leq  
	\sup_{\zeta \in \mathcal{N}_\epsilon(\mathcal{S}^{t-1})}\frac{1}{n}\Big\|\ind\Big(|\sum_{k = 1}^{t} \zeta^k\phi_k| < \frac{\lambda(1+b_t)}{\|\gamma_t\|_2} + \epsilon \Big)\Big\|_0,
	%
\end{align}
where $\mathcal{N}_\epsilon(\mathcal{S}^{t-1})$ forms an $\epsilon$-cover of $\mathcal{S}^{t-1}$. 
Here in the last inequality, we make the observation that for every $\zeta' \in \mathcal{S}^{t-1}$, there exists $\zeta \in \mathcal{N}_\epsilon(\mathcal{S}^{t-1})$ such that $\ltwo{\zeta - \zeta'} \leq \epsilon$ and hence 
\begin{align}
\Bigg||\sum_{k = 1}^{t} \zeta'^{k}\phi_{kj}| - |\sum_{k = 1}^{t} \zeta^k\phi_{kj}|\Bigg|
\leq 
\Big|\sum_{k = 1}^{t} (\zeta^k-\zeta'^{k}) \phi_{kj}\Big| 	
\leq 
\ltwo{\zeta-\zeta'} \Big(\sum_{j=1}^t \phi^2_{kj}\Big)^{1/2} \leq \epsilon,
\end{align}
with probability at least $1 - O(n^{-11})$. 

Fix each $\zeta$ independent of $\{\phi_{k}\}$, $\|\ind(|\sum_{k = 1}^{t} \zeta^k\phi_k| < \frac{\lambda(1+b_t)}{\|\gamma_t\|_2} + \epsilon )\|_0$ is the summation of $n$ independent Bernoulli distribution with parameter $\mprob(|\tilde{g}| < \frac{\lambda(1+b_t)}{\|\gamma_t\|_2} + \epsilon),$ for $\tilde{g} \sim\mathcal{N}(0,1/n)$. 
According to standard concentration result for summation of independent Bernoulli's, we obtain 
\begin{align*}
	\frac{1}{n}\Big\|\ind\Big(|\sum_{k = 1}^{t} \zeta^k\phi_k| < \frac{\lambda(1+b_t)}{\|\gamma_t\|_2} + \epsilon \Big)\Big\|_0
	&\leq 
	\mprob\Big(|\tilde{g}| < \frac{\lambda(1+b_t)}{\|\gamma_t\|_2} + \epsilon\Big) + O\Big(\sqrt{\frac{\log \frac{1}{\delta}}{n}}\Big),
\end{align*}
with probability at least $1 - \delta$.
If we set $\epsilon = \frac{1}{n}$ and take a union bound over elements in $\mathcal{N}_\epsilon(\mathcal{S}^{t-1})$,   it holds that 
\begin{align}
\notag	\sup_{\zeta \in \mathcal{N}_\epsilon(\mathcal{S}^{t-1})}\frac{1}{n}\Big\|\ind\Big(|\sum_{k = 1}^{t} \zeta^k\phi_k| < \frac{\lambda(1+b_t)}{\|\gamma_t\|_2} + \epsilon \Big)\Big\|_0
	&\leq 
	\mprob\Big(|\tilde{g}| < \frac{\lambda(1+b_t)}{\|\gamma_t\|_2} + \frac{1}{n}\Big) + O\Big(\sqrt{\frac{t\log n}{n}}\Big)\\
	&\leq  
	\mprob\Big(|\tilde{g}| < \frac{\lambda(1+b_t)}{\|\gamma_t\|_2} \Big) + O\Big(\sqrt{\frac{t\log n}{n}}\Big),
\end{align}
with probability at least $1 - O(n^{-10})$. 
Here the last inequality uses $P\Big(\lambda(1+b_t)/\|\gamma_t\|_2  < |\tilde{g}| < \lambda(1+b_t)/\|\gamma_t\|_2 + 1/n \Big) < \sqrt{2/(\pi n)}$ since the density function of $|\tilde{g}|$ is bounded by $\sqrt{2n/\pi}.$ 
Combining with \eqref{eqn:sibelius}, the above relation leads to 
\begin{align}
\label{eqn:sibelius-part1}
	\frac{1}{n}\Big\|\ind\Big(|\sum_{k = 1}^{t} \gamma_{t}^k\phi_k| < \lambda(1+b_t)\Big)\Big\|_0
	\leq  
	\mprob\Big(|\tilde{g}| < \frac{\lambda(1+b_t)}{\|\gamma_t\|_2} \Big) + O\Big(\sqrt{\frac{t\log n}{n}}\Big). 
\end{align}
Similarly, one can deduce 
\begin{align}
\label{eqn:sibelius-part2}
	\notag \frac{1}{n}\Big\|\ind\Big(|\sum_{k = 1}^{t} \gamma_{t}^k\phi_k| < \lambda(1+b_t)\Big)\Big\|_0
	%
	&\geq 
	\inf_{\zeta \in \mathcal{S}^{t-1}}\frac{1}{n}\Big\|\ind\Big(|\sum_{k = 1}^{t} \zeta^k\phi_k| < \frac{\lambda(1+b_t)}{\|\gamma_t\|_2}\Big)\Big\|_0 \\
	\notag &\geq  
	\inf_{\zeta \in \mathcal{N}_\epsilon(\mathcal{S}^{t-1})}\frac{1}{n}\Big\|\ind\Big(|\sum_{k = 1}^{t} \zeta^k\phi_k| < \frac{\lambda(1+b_t)}{\|\gamma_t\|_2} - \epsilon \Big)\Big\|_0\\
	&\geq \mprob\Big(|\tilde{g}| < \frac{\lambda(1+b_t)}{\|\gamma_t\|_2} \Big) + O\Big(\sqrt{\frac{t\log n}{n}}\Big). 
\end{align} 
Putting these two parts with \eqref{eqn:tosca}, we conclude 
\begin{align}
\label{eq:ratio}
	\frac{k_t}{n} = \frac{1}{n}\Big\|\ind\lt(|r_t| < \lambda(1+b_t)\rt)\Big\|_0 
	= 
	\mprob\Big(|\tilde{g}| < \frac{\lambda(1+b_t)}{\|\gamma_t\|_2} \Big) + o(1).
\end{align}

To establish \eqref{eqn:cat}, it is sufficient to notice that $k_{t} \asymp p$. 
By definition of $k_{t}$ (cf.~\eqref{eqn:defn-kt}), it satisfies straightforwardly that $k_{t} > p.$
According to \eqref{eq:ratio}, it then implies
\begin{align}
\label{eqn:sunset}
	\mprob\Big(|\tilde{g}| < \frac{\lambda(1+b_t)}{\|\gamma_t\|_2} \Big) \geq \frac{p}{n} + o(1).
\end{align}
Given $\tilde{g} \sim \mathcal{N}(0,1/n)$, the above relation ensures 
$b_t \asymp \overline{\gamma}$ for $\lambda \asymp 1/\sqrt{n}$ and $\lambda(1+b_t)/\|\gamma_t\|_2 \asymp \frac{1}{\sqrt{n}}$, which in turns gives
\begin{align*}
	k_{t} = p(1 + 1/b_t) = p(1 + o(1)).
\end{align*}
We thus finish the proof of inequality~\eqref{eqn:cat}.

\paragraph{Proof of $\|\gamma_t\|_2, \ltwo{\alpha_t} \gtrsim 1$.} 
We are only left to show $\ltwo{\alpha_t}\gtrsim 1$. 
Similar to \eqref{eqn:alphat-allegro}, invoking the Lipschitz property of the function $\ltwo{\psi(\cdot \mid \lambda(1+b_t))}$ gives 
\begin{align}
\label{eqn:robust-alpha-lb}
\notag \|\alpha_t\|_2 
&= \|G_t(s_t)\|_2 = \frac{n b_t}{p(1 + b_t)}\lt\|\psi\lt(s_t + \varepsilon; \lambda(1+b_t)\rt)\rt\|_2 \\
\notag &= \frac{n b_t}{p(1 + b_t)}\Big\|\psi\Big(\varepsilon + \sum_{k = 1}^{t} \gamma_{t}^k\phi_k; \lambda(1+b_t)\Big)\Big\|_2 + O(\|\xi_t\|_2)\\
&\geq \frac{n b_t}{p(1 + b_t)}\Big\|\psi\Big(\varepsilon + \sum_{k = 1}^{t} \gamma_{t}^k\phi_k; \lambda(1+b_t)\Big)\circ \ind_{\mathcal{I}} \Big\|_2 + O(\|\xi_t\|_2)
\end{align}
where we define $\mathcal{I} := \{i:\varepsilon_i \sim \mathcal{N}(0, \sigma^2)\}$.
Recall that $\varepsilon_i$ is drawn from the mixture distribution of $\mathcal{N}(0, \sigma^2)$ and some other distribution $H$ as in \eqref{eqn:robust-error}.

In order to show $\ltwo{\alpha_t}\gtrsim 1$, it  suffices to lower bound $\|\psi(\varepsilon + \sum_{k = 1}^{t} \gamma_{t}^k\phi_k; \lambda(1+b_t))\circ \ind_{\mathcal{I}}\|_2$. 
Before proceeding, we make two key observations.
\begin{itemize}
\item We first make the remark that $\{\phi_{k}\}$ are independent of $\varepsilon$. 
In fact, when constructing $\{\phi_{k}\}$ (cf.~\eqref{eqn:def-phi-psi}), we have viewed $\varepsilon$ as a deterministic vector and each $\phi_{k}$ admits a fixed distribution $\mathcal{N}(0,\frac{1}{n}I_n)$ no matter what value $\varepsilon$ takes. In other words, $\{\phi_{k}\}$ are independent of $\varepsilon$.

\item From our discussions around \eqref{eqn:sunset}, it satisfies $\lambda(1+b_t) \asymp \frac{\ltwo{\gamma_t}}{\sqrt{n}}$, and hence, 
\begin{align}
\label{eqn:robust-boundness}
	\Big\|\psi\Big(\varepsilon + \sum_{k = 1}^{t} \gamma_{t}^k\phi_{k}; \lambda(1+b_t)\Big) \circ \ind_{\mathcal{I}}\Big\|_{2}
	\leq 
	\|\lambda(1+b_t)\circ \ind_{\mathcal{I}}\|_2 \lesssim 1.
\end{align}
\end{itemize}

With these two facts in mind, similar to \eqref{eqn:sibelius-part2}, we obtain 
\begin{align}
\label{eqn:apple}
\notag \Big\|\psi\Big(\varepsilon + \sum_{k = 1}^{t} \gamma_{t}^k\phi_{k}; \lambda(1+b_t)\Big) \circ \ind_{\mathcal{I}}\Big\|_{2} 
&\geq \inf_{\zeta\in \mathcal{N}_{\epsilon}(S^{t-1})}
\Big\|\psi\Big(\varepsilon + \ltwo{\gamma_t}\sum_{k = 1}^{t} \zeta^k \phi_{k} - \epsilon; \lambda(1+b_t)\Big) \circ \ind_{\mathcal{I}}\Big\|_{2}  \\
&\geq
\mathbb{E}\Big[\|\psi(\varepsilon+\|\gamma_t\|_2g - 1/n; \lambda(1+b_t))\circ \ind_{\mathcal{I}}\|_2 \mid \ltwo{\gamma_t}\Big] + O\lt(\sqrt{\frac{t\log n}{n}}\ltwo{\gamma_t}\rt),
\end{align}
with probability at least $1 - O(n^{-11})$. Here we take $\epsilon = 1/n.$
Following the exact same argument for deriving~\eqref{eqn:penguin}, since $\Big\|\psi\Big(\varepsilon + \sum_{k = 1}^{t} \gamma_{t}^k\phi_{k}; \lambda(1+b_t)\Big) \circ \ind_{\mathcal{I}}\Big\|_{2}$ concentrates tightly around its mean, and is bounded from above, one can deduce 
\begin{align}
	\notag &\mathbb{E}\Big[\big\|\psi(\varepsilon+\|\gamma_t\|_2g - 1/n; \lambda(1+b_t))\circ \ind_{\mathcal{I}}\big\|^2_2 \mid \ltwo{\gamma_t}\Big]\\
	&=
	\Big(\mathbb{E}\Big[\|\psi(\varepsilon+\|\gamma_t\|_2g - 1/n; \lambda(1+b_t))\circ \ind_{\mathcal{I}}\|_2 \mid \ltwo{\gamma_t}\Big]\Big)^2 + O\lt(\sqrt{\frac{t\log n}{n}}\rt).
\end{align}
When it comes to further bounding quantity $\mathbb{E}[\|\psi(\varepsilon+\|\gamma_t\|_2g - 1/n; \lambda(1+b_t))\circ \ind_{\mathcal{I}}\|^2_2 \mid \ltwo{\gamma_t}]$, 
some direct algebra gives 
\begin{align}
\label{eqn:robust-lb}
\notag \mathbb{E}\Big[\big\|\psi(\varepsilon+\|\gamma_t\|_2g - 1/n; \lambda(1+b_t))\circ \ind_{\mathcal{I}}\big\|^2_2 \mid \ltwo{\gamma_t}\Big]
&=   
\mathbb{E} \Big[\sum_{k\in \mathcal{I}} \min \big\{(\varepsilon_k + \|\gamma_t\|_2g_k - 1/n)^2,~\lambda^2(1+b_t)^2\big\}\Big]\\
\notag &\gtrsim   
\sum_{k\in \mathcal{I}} \frac{1}{2}\lambda^2(1+b_t)^2 \mprob\Big(\varepsilon_k > 0, \|\gamma_t\|_2 g_k - \frac{1}{n} > \frac{1}{2}\lambda(1+b_t), k\in\mathcal{I}\Big)\\
&\gtrsim 1,
\end{align}
where in the last line, we recall the independence between $\varepsilon$ and $\{\phi_k\}$ and conclude 
\begin{align*}
	\mprob\Big(\varepsilon_k > 0, \|\gamma_t\|_2 g_k - \frac{1}{n} > \frac{1}{2}\lambda(1+b_t), k\in\mathcal{I}\Big)
	=
	\mprob\Big(\varepsilon_k > 0\Big) \mprob\Big(\|\gamma_t\|_2 g_k - \frac{1}{n} > \frac{1}{2}\lambda(1+b_t), k\in\mathcal{I}\Big)
	\gtrsim 1.
\end{align*}
Putting together inequalities~\eqref{eqn:robust-alpha-lb}, \eqref{eqn:apple} and \eqref{eqn:robust-lb}
concludes $\|\alpha_t\|_2 \gtrsim 1$ and thus $\|\gamma_t\|_2 \gtrsim 1.$

\medskip

Combining these two parts together, we have shown that $\ltwo{\alpha_t}, \ltwo{\gamma_t} \asymp 1$, thus validating the Assumption~\ref{ass:control-lipschitz}.

\subsubsection{Validating Assumption~\ref{ass:control-simple}}
In view of the definition of \eqref{eqn:robust-F-G}, it is easily seen that $\frac{1}{n}\mathbb{E}\lt\|F_{t+1}^{\prime}(v_{t+1})\rt\|_2^2 = \frac{p}{n}$. 
Therefore, in order to justify Assumption~\ref{ass:control-simple}, we are only left with computing 
$\frac{1}{n}\mathbb{E}[\|G_{t+1}^{\prime}(u_{t+1})\|_2^2]$. 
Towards this, according to the choice of $b_t$ (cf.~\eqref{eqn:defn-kt}), we first make the observation that  
\begin{align*}
	\frac{1}{n}\lt\|G_{t}^{\prime}(s_{t})\rt\|_2^2 
	&= \frac{1}{n}\lt\|g_t^{\prime}(r_t)\rt\|_2^2 = \Big(\frac{n b_t}{p(1 + b_t)}\Big)^2 \cdot \frac{k_t}{n} = \frac{n b_t}{p(1 + b_t)}.
\end{align*}
%
Next we shall compute the difference between $\|G_{t}^{\prime}(u_{t})\|_2^2$ and $\|G_{t}^{\prime}(s_{t})\|_2^2$. Some direct algebra gives 
\begin{align*}
\lt|\lt\|G_{t}^{\prime}(u_{t})\rt\|_2^2 - \lt\|G_{t}^{\prime}(s_{t})\rt\|_2^2\rt| 
&= \Big(\frac{n b_t}{p(1 + b_t)}\Big)^2 \cdot  \Big|\big\||\varepsilon + u_t| < \lambda(1+b_t)\big\|_0 
- \big\||r_t| < \lambda(1+b_t)\big\|_0\Big| \\
&\leq \Big(\frac{n b_t}{p(1 + b_t)}\Big)^2 \cdot  \Big|\ind(\varepsilon + u_t| < \lambda(1+b_t)) 
- \ind(|r_t| < \lambda(1+b_t))\Big| \\
&\lesssim t\log n + n\Big(\frac{\|\xi_t\|_2}{\ltwo{\varepsilon + u_t}}\Big)^{\frac{2}{3}} \ll n.
\end{align*}
Here the penultimate inequality results similarly from the relation~\eqref{eqn:piazzola} as a consequence of concentration lemma~\ref{lem:ind-noise}; the last inequality invokes the assumption $\ltwo{\xi_t} = o(1)$ and 
inequality~\eqref{eqn:dvorak-romantic} that 
\begin{align*}
	\ltwo{\varepsilon + u_t} = \Big\|\varepsilon + \sum_{k=1}^t \gamma_t^k\phi_k\Big\|_2 =
	\sqrt{\lt\|\varepsilon\rt\|_2^2 + \lt(1 + O\Big(\sqrt{\frac{t\log n}{n}}\Big)\rt)\lt\|\gamma_{t}\rt\|_2^2 + O\Big(\sqrt{\frac{t\log n}{n}}\lt\|\varepsilon\rt\|_2\lt\|\gamma_{t}\rt\|_2\Big)} \gtrsim 1.
\end{align*}

Combining the derivations above, we arrive at 
\begin{align}
\frac{1}{n}\mathbb{E}\lt\|G_{t}^{\prime}(u_{t})\rt\|_2^2 
&\leq \frac{n}{p}\Big(1 - \frac{1}{1+b_t}\Big) + o(1), 
\end{align}
provided that $t \lesssim \frac{n}{\log^4 n}$.
Putting everything together, we have verified that 
\begin{align}
	\frac{1}{n^2}\mathbb{E}\lt\|F_{t+1}^{\prime}(v_{t+1})\rt\|_2^2\mathbb{E}\lt\|G_{t+1}^{\prime}(u_{t+1})\rt\|_2^2 &< 1 - \frac{1}{1+ c' \overline{\gamma}},
\end{align}
for some constant $c'$. Here we recall $b_{t} \asymp \overline{\gamma}$ as discussed around inequality~\eqref{eqn:sunset}.
We have thus validated Assumption~\ref{ass:control-simple}.

\subsubsection{State evolution}
\label{sec:robust-se-pf}

Again, our final goal is to bound the difference between the non-asymptotic SE $(\alpha_t, \gamma_{t+1})$ to the deterministic SE defined in expression~\eqref{eqn:robust-SE}. 
We proceed by using the induction method to achieve this goal. 
Firstly, it is easily seen that the set of relation~\eqref{eqn:pretty-se-robust-bla} holds true for $t=1$. Assuming that for some $t \geq 1$, 
\begin{align}
\big|\ltwo{\alpha_{t-1}}^2 - \alpha^{\star 2}_{t-1}\big| \lesssim \Big(\frac{t\log^2 n}{n}\Big)^{1/3}
\qquad \text{and} \qquad
	\big|\ltwo{\gamma_t}^2 - \gamma^{\star 2}_{t}\big| \lesssim \Big(\frac{t\log^2 n}{n}\Big)^{1/3},
\end{align}
it is thus sufficient to verify them for $t+1$.

In the derivations above, we have shown that $\ltwo{\alpha_t}, \ltwo{\gamma_t} \asymp 1$. Based on these relations, we claim that 
\begin{align}
\label{eqn:logan}
\|\gamma_{t+1}\|_2^2 
&= \frac{n}{p}\Big(\frac{b_t}{1+b_t}\Big)^2\mathbb{E}\Big[\big\|\psi\big(\varepsilon + \|\gamma_{t}\|_2g; \lambda(1+b_t)\big)\big\|_2^2 \mid \ltwo{\gamma_t}, \varepsilon \Big] + O\lt(\Big(\frac{t\log^2 n}{n}\Big)^{\frac{1}{3}}\rt).
\end{align}

\paragraph{Proof of relation~\eqref{eqn:logan}.}
First we recall inequality~\eqref{eqn:robust-gamma-sim} to obtain that 
\begin{align}
\label{eqn:gamma-t+1-doupi}
\|\gamma_{t+1}\|_2^2 &= \frac{p}{n}\|\alpha_t\|_2^2 + O\lt(\Big(\frac{t\log^2 n}{n}\Big)^{\frac{1}{3}}\rt). 
\end{align}
The definition of $\alpha_{t}$ directly yields 
\begin{align}
\label{eqn:mozart-17}
	\ltwo{\alpha_t} = \ltwo{G_t(s_t)} = \frac{n b_t}{p(1 + b_t)}\lt\|\psi\lt(s_t + \varepsilon; \lambda(1+b_t)\rt)\rt\|_2.
\end{align}
Next, in view of the Lipschitz property of the function $\psi$ and the decomposition~\eqref{defi:dynamics-s} of $s_{t}$, we can further conclude 
\begin{align*}
\|\alpha_t\|_2 
&= \frac{n b_t}{p(1 + b_t)}\Big\|\psi\Big(\varepsilon + \sum_{k = 1}^{t} \gamma_{t}^k\phi_k; \lambda(1+b_t)\Big)\Big\|_2 + O\lt(\Big(\frac{t\log^2 n}{n}\Big)^{\frac{1}{3}}\rt) \notag\\
&= \frac{n b_t}{p(1 + b_t)}
\mathbb{E}\Big[\big\|\psi\big(\varepsilon + \|\gamma_{t}\|_2g; \lambda(1+b_t)\big)\big\|_2\mid \ltwo{\gamma_t}, \varepsilon \Big] + O\lt(\Big(\frac{t\log^2 n}{n}\Big)^{\frac{1}{3}}\rt),
\end{align*}
with probability at least $1 - O(n^{-10}).$
Here we invoke the relation $\ltwo{\xi_t} \lesssim O((\frac{t\log^2 n}{n})^{\frac{1}{3}})$ in the first line, and the concentration property of Gaussian vectors as in \eqref{eqn:alphat-rondo} in the second line.
In addition, since $\psi\big(\varepsilon + \|\gamma_{t}\|_2g; \lambda(1+b_t)\big)$ concentrates well around its expectation, with similar argument as in \eqref{eqn:penguin}, one can derive 
\begin{align*}
\|\alpha_t\|_2 
&= \frac{n b_t}{p(1 + b_t)}\sqrt{\mathbb{E}\Big[\big\|\psi\big(\varepsilon + \|\gamma_{t}\|_2g; \lambda(1+b_t)\big)\big\|_2^2 \mid \ltwo{\gamma_t}, \varepsilon\Big]} + O\lt(\Big(\frac{t\log^2 n}{n}\Big)^{\frac{1}{3}}\rt),
\end{align*}
with probability at least $1 - O(n^{-11}).$ Putting the above together with \eqref{eqn:gamma-t+1-doupi}, we establish the advertised relation~\eqref{eqn:logan}.

Now based on the recursion~\eqref{eqn:logan}, in order to study how $\ltwo{\gamma_{t}}$ evolves with $t$, it is sufficient to study the following function for any value $b > 0$,    
\begin{align}
h_b(\mu) := \mathbb{E}\Big[\big\|\psi\big(\varepsilon + \sqrt{\mu}g; \lambda(1+b)\big)\big\|_2^2 \mid \varepsilon \Big].
\end{align}
With this definition, the limiting state-evolution~\eqref{eqn:robust-SE} satisfies 
\begin{align}
\label{eqn:se-iterate-robust}
	\gamma_{t+1}^{\star 2} = \frac{n}{p}\Big(\frac{b_t^{\star}}{1+b_t^{\star}}\Big)^2 h_{b_t^{\star}}(\gamma_{t}^{\star 2}),
\end{align}
while its non-asymptotic analogue satisfies 
\begin{align}
\label{eqn:se-non-iterate}
	\|\gamma_{t+1}\|_2^2 = \frac{n}{p}\Big(\frac{b_t}{1+b_t}\Big)^2 h_{b_t}(\|\gamma_{t}\|_2^2) + O\lt(\Big(\frac{t\log^2 n}{n}\Big)^{\frac{1}{3}}\rt).
\end{align}

Now our goal is to control the difference between $\|\gamma_{t+1}\|_2^2$ and $\gamma_{t+1}^{\star 2}$. Towards this end, we first make note of the following two properties regarding the function $h_{b}$.
\begin{itemize}
\item First, we claim that
\begin{align} 
\big|h_b'(\mu)\big| 
&= \frac{1}{\sqrt{\mu}}\Big|\mathbb{E}\Big[\big\langle\psi\big(\varepsilon + \sqrt{\mu}g; \lambda(1+b)\big), \ind(|\varepsilon + \sqrt{\mu}g| < \lambda(1+b)) \circ g\big\rangle \mid \varepsilon \Big]\Big| \notag\\
&\le \frac{1-2c}{\mu}\mathbb{E}\Big[\big\|\psi\big(\varepsilon + \sqrt{\mu}g; \lambda(1+b)\big)\big\|_2^2 \mid \varepsilon \Big]
\label{eq:SE-robust}\\
&= \frac{1-2c}{\mu}h_b(\mu). \notag 
\end{align}
Let us take inequality~\eqref{eq:SE-robust} as given for the moment, and come back to its proof at the end of this section. 
Note that for $\mu \asymp 1$, $h_b(\mu)/{\mu} = O(1)$ and the above property guarantees that $\frac{1}{\mu}h(\mu)$ is $O(1)$-Lipschitz continuous function of $\mu.$ As a result, one has 
\begin{align*}
	(1-2c)\frac{h_b(\mu)}{\mu} \le (1-2c)\frac{h_b(\mu')}{\mu'} + O\left(\Big(\frac{t\log^2 n}{n}\Big)^{\frac{1}{3}}\right)
	\leq (1-c)\frac{h_b(\mu')}{\mu'},
\end{align*}
for $|\mu - \mu'| = O((\frac{t\log^2 n}{n})^{\frac{1}{3}})$.
Putting the things above together yields
\begin{align}
\label{eqn:simple-prop}
\big|h_b'(\mu)\big| &\leq (1-c)\frac{h_b(\mu')}{\mu'},
\end{align}
for some constant $c \in (0,1)$, $b = b_t, b_t^{\star}$ and $|\mu - \mu'| = O((\frac{t\log^2 n}{n})^{\frac{1}{3}})$.

\item 
For any value of $\mu$, regarding $b^2h_b(\mu)/(1+b)^{2}$ as a function of $b$, notice that
\begin{align}
\frac{\partial \big(b^2\psi(u/(1+b); \lambda)^2\big)}{\partial b} = \frac{2b\psi(u/(1+b); \lambda)^2}{1+b} \ge 0.
\end{align}
Therefore, $b^2h_b(\mu)/(1+b)^{2}$ is a non-decreasing function of $b$. 
\end{itemize}

With these properties in place, putting relations \eqref{eqn:se-iterate-robust} and \eqref{eqn:se-non-iterate} together gives 
\begin{align}
&\big|\|\gamma_{t+1}\|_2^2 - (\gamma_{t+1}^{\star})^2\big| \notag\\
&\le \frac{n}{p}\max\bigg\{\Big(\frac{b_t^{\star}}{1+b_t^{\star}}\Big)^2\big|h_{b_t^{\star}}\big(\|\gamma_t\|_2^2\big) - h_{b_t^{\star}}\big((\gamma_{t}^{\star})^2\big)\big|, \Big(\frac{b_t}{1+b_t}\Big)^2\big|h_{b_t}\big(\|\gamma_t\|_2^2\big) - h_{b_t}\big((\gamma_{t}^{\star})^2\big)\big|\bigg\} + O\lt(\Big(\frac{t\log^2 n}{n}\Big)^{\frac{1}{3}}\rt) \notag\\
&\lesssim \big|\|\gamma_{t}\|_2^2 - (\gamma_{t}^{\star})^2\big| + O\lt(\Big(\frac{t\log^2 n}{n}\Big)^{\frac{1}{3}}\rt) = o(1),
\end{align}
where the last relation invokes our inductive assumption. Furthermore, we can write 
\begin{align}
\label{eqn:se-robust-key}
&(\gamma_{t+1}^{\star})^{-2}\big|\|\gamma_{t+1}\|_2^2 - (\gamma_{t+1}^{\star})^2\big| \notag\\
&\le \max\big\{\big(h_{b_t^{\star}}(\gamma_{t}^{\star 2})\big)^{-1}\big|h_{b_t^{\star}}\big(\|\gamma_t\|_2^2\big) - h_{b_t^{\star}}\big((\gamma_{t}^{\star})^2\big)\big|, (1+o(1))\big(h_{b_t}(\|\gamma_{t}\|_2^2)\big)^{-1}\big|h_{b_t}\big(\|\gamma_t\|_2^2\big) - h_{b_t}\big((\gamma_{t}^{\star})^2\big)\big|\big\} + O\lt(\Big(\frac{t\log^2 n}{n}\Big)^{\frac{1}{3}}\rt) \notag\\
&\lesssim (1 - c)(\gamma_{t}^{\star})^{-2}\big|\|\gamma_{t}\|_2^2 - (\gamma_{t}^{\star})^2\big| + O\lt(\Big(\frac{t\log^2 n}{n}\Big)^{\frac{1}{3}}\rt), 
\end{align}
where the second line holds since for $b_t^{\star} \ge b_t$, 
\begin{align*}
h_{b_t}\big(\|\gamma_t\|_2^2\big) - h_{b_t}\big((\gamma_{t}^{\star})^2\big) 
\le h_{b_t}\big(\|\gamma_t\|_2^2\big) - h_{b_t^{\star}}\big((\gamma_{t}^{\star})^2\big) 
\le h_{b_t^{\star}}\big(\|\gamma_t\|_2^2\big) - h_{b_t^{\star}}\big((\gamma_{t}^{\star})^2\big)
\end{align*}
and the last line makes use of~\eqref{eqn:simple-prop} with $\mu = c'\|\gamma_{t}\|_2^2 + (1-c')(\gamma_{t}^{\star})^2$ and $\mu' = (\gamma_{t}^{\star})^2$ for some $0 \le c' \le 1$.


In view of the initialization $\gamma_{1} = \gamma^\star_{1} = \ltwo{\theta^\star}$ and $\|\gamma_{t}\|_2 \asymp 1$, invoking the above relation~\eqref{eqn:se-robust-key} recursively leads to 
\begin{align*}
\big|\|\gamma_{t+1}\|_2^2 - (\gamma_{t+1}^{\star})^2\big| &\lesssim \Big(\frac{t\log^2 n}{n}\Big)^{\frac{1}{3}}.
\end{align*}
which in turns gives that
\begin{align*}
\big|\|\alpha_{t}\|_2^2 - (\alpha_{t}^{\star})^2\big| \lesssim \Big(\frac{t\log^2 n}{n}\Big)^{\frac{1}{3}}.
\end{align*}

\paragraph{Proof of Claim~\eqref{eq:SE-robust}.}
To establish the expression~\eqref{eq:SE-robust} for $b = b_{t} \text{ and } b^\star_{t}$,  
consider a change of variable $\tau \defn \frac{\sqrt{n}\lambda(1+b)}{\sqrt{\mu}}$. 
It is then sufficient to prove that for any $\tau > 0$,
\begin{align} 
\label{eq:defi-H2} 
H_2(\tau) := 1 - \sup_{\varepsilon} \frac{\Big|\mathbb{E}\Big[G(\varepsilon + G)\ind(|\varepsilon + G| < \tau)\Big]\Big|}{\mathbb{E}\big[(\varepsilon + G)^2 \wedge \tau^2\big]} > c, 
\end{align}
for some constant $c \in (0,1)$.

\begin{figure}[t]
\centering
\includegraphics[width=0.85\textwidth]{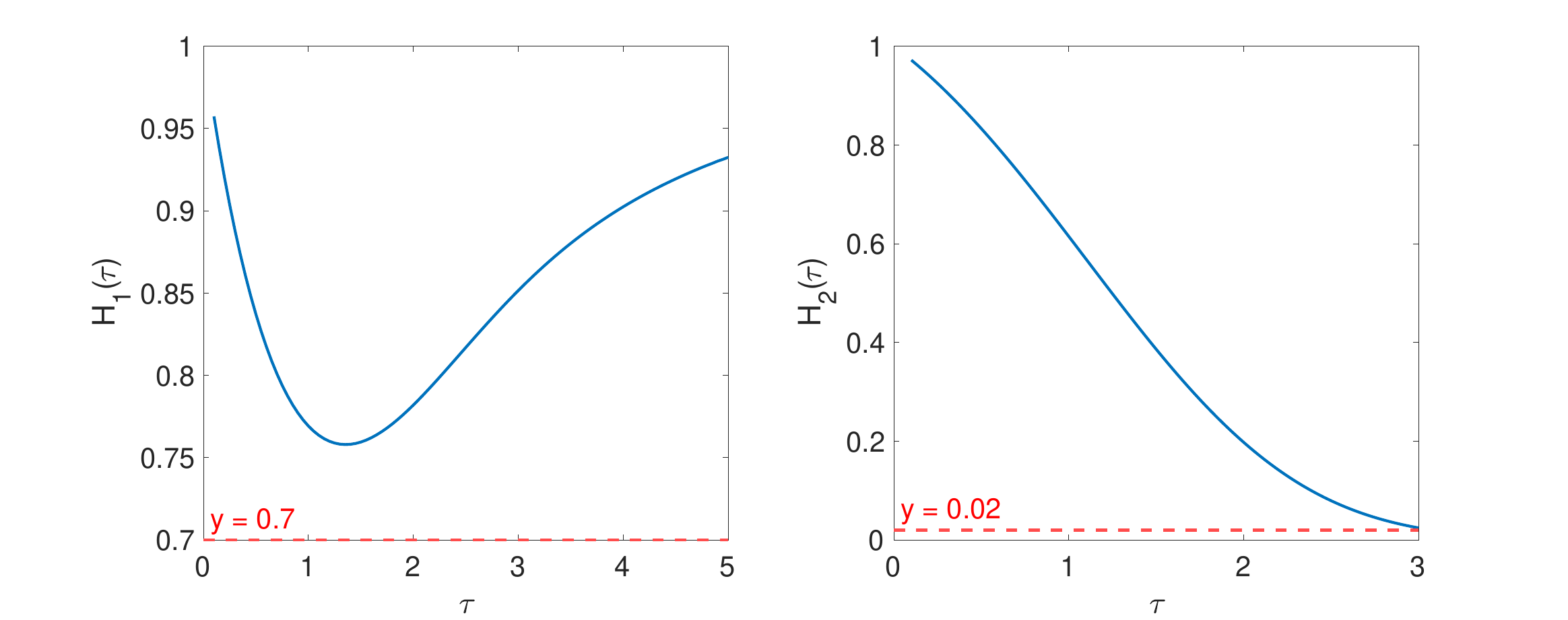}
	\caption{
	Numerical calculations for $H_1(\tau)$ of~\eqref{eq:defi-H1} with $\tau \in (0,5)$
	and $H_2(\tau)$ of~\eqref{eq:defi-H2} with $\tau \in (0,3)$. 
	}
\label{fig:robust}
\end{figure}

\begin{itemize}
\item 
For $\tau \in (0,3)$, the required inequality~\eqref{eq:defi-H2} can be directly obtained from the numerical simulation in Figure~\ref{fig:robust}. Here the value $H_2(\tau) \geq 0.02$ for $\tau \in (0,3)$. 

\item For $\tau \ge 3$, due to symmetry, it is enough to consider the case $\varepsilon \ge 0$.
When $\varepsilon \ge \tau$, some direct calculations yield 
\begin{align*}
\Big|\mathbb{E}\Big[G(\varepsilon + G)\ind(|\varepsilon + G| < \tau)\Big]\Big| &\le \max\Big\{\mathbb{E}\Big[G(\varepsilon + G)\ind(-\varepsilon - \tau < G < -\varepsilon)\Big], -\mathbb{E}\Big[G(\varepsilon + G)\ind(-\varepsilon < G < \tau - \varepsilon)\Big]\Big\} \\
&\le \max\Big\{\mathbb{E}\Big[G^2\ind(G < -\tau)\Big], -\tau\mathbb{E}\Big[G\ind(G < 0)\Big]\Big\}
\le \tau,
\end{align*}
and
\begin{align*}
\mathbb{E}\big[(\varepsilon + G)^2 \wedge \tau^2\big] &> \tau^2\mathbb{E}\big[\ind(G > 0)\big] = \frac{\tau^2}{2}.
\end{align*}
In this case, $H_{2}(\tau) \geq 1 - \frac{2}{\tau} \geq 1/3.$
In the other case, for $0 \le \varepsilon \le \tau$, we obtain 
\begin{align*}
\Big|\mathbb{E}\Big[G(\varepsilon + G)\ind(|\varepsilon + G| < \tau)\Big]\Big| &\le \Big|\mathbb{E}\Big[G(\varepsilon + G)\ind(-\varepsilon - \tau < G < \varepsilon - \tau)\Big]\Big| + 2\mathbb{E}\Big[G^2\ind(0 < G < \tau - \varepsilon)\Big] \\
&< \varepsilon\Big|\mathbb{E}\Big[(\varepsilon + G)\ind(-\varepsilon - \tau < G < \varepsilon - \tau)\Big]\Big| + \mathbb{E}\Big[(\varepsilon + G)^2\ind(-\varepsilon - \tau < G < \varepsilon - \tau)\Big] \\
&\qquad+ 2\mathbb{E}\Big[G^2\ind(0 < G < \tau - \varepsilon)\Big] \\
&< \tau^2\mathbb{E}\Big[\ind(G < \varepsilon - \tau)\Big]+ \mathbb{E}\Big[(\varepsilon + G)^2\ind(-\varepsilon - \tau < G < \varepsilon - \tau)\Big] \\
&\qquad+ 2\mathbb{E}\Big[G^2\ind(0 < G < \tau - \varepsilon)\Big],
\end{align*}
and
\begin{align*}
\mathbb{E}\big[(\varepsilon + G)^2 \wedge \tau^2\big] &> \tau^2\mathbb{E}\Big[\ind(G > \tau - \varepsilon)\Big] + \mathbb{E}\Big[(\varepsilon + G)^2\ind(-\varepsilon - \tau < G < \varepsilon - \tau)\Big] \\
&\qquad + 2\mathbb{E}\Big[(\varepsilon^2 + G^2)\ind(0 < G < \tau - \varepsilon)\Big].
\end{align*}
As a result, $H_2(\tau) > 0.$ Since we only need to consider those $\tau$ such that $\tau \asymp 1$ which forms a compact set, therefore it implies there exists some constant $c$ such that 
 $H_2(\tau) > c.$
\end{itemize}

We have thus completed the proof of Claim~\eqref{eq:SE-robust}.

\section{Auxiliary concentration lemmas and their proofs}
\label{sec:concentration}

In this section, we collect a few concentration results for functions of random vectors that shall be used multiple times throughout this paper.

\subsection{Lemma statements}


The first result below considers the summation of independent sub-exponential random variables and develops a Bernstein-like concentration bound. 

\begin{lems}
\label{lem:concentration}
Suppose that $Z_i$'s are independent random variables satisfying
\begin{align}
	\mathbb{E}[Z_i] = 0 \qquad\text{and}\qquad 
	\mathbb{P}\Big(|Z_i| \geq B \log \frac{1}{\delta} \Big) \le \delta, 
	\qquad\text{for every } 0< \delta \leq \frac{1}{\poly(n)}, \label{eq:Z_condition}
\end{align}
for some $B \geq 0.$
Then with probability at least $1 - \delta$, one has 
\begin{align}
\label{eqn:sum-zi}
	\Big|\sum_{i = 1}^n Z_i\Big| \lesssim \sqrt{\sum_{i = 1}^n \Big(\mathsf{Var}(Z_i)+ \big(\frac{B\log n}{n}\big)^2 \Big)\log\frac{1}{\delta}} + B\log n\log\frac{1}{\delta}.
\end{align}
\end{lems}
\noindent The proof of this result can be found in Section~\ref{sec:pf-lem-concentration}.

It is worth pointing out that a direct application of Bernstein's inequality --- in view of the boundedness condition \eqref{eq:Z_condition} for each $Z_{i}$ --- adds an additional $\log \frac{1}{\delta}$ to the second term of \eqref{eqn:sum-zi} (see e.g.~\cite[Section 2.1.3]{wainwright2019high}). 
In that case, this additional term when combined with a covering argument shall result in an inferior dependence on the parameter $t$.

Next, we derive a useful concentration bound associated with indicator functions.
In particular, we count the number of times that a Lipschitz function crosses a certain threshold over multiple independent realizations. 
This concentration result turns out to be useful when dealing with discontinuous denoising functions with its proof postponed to Section~\ref{sec:pf-lem-ind-noise}. 

\begin{lems} \label{lem:ind-noise}
Consider independent random vectors $\{X_i\}_{i=1}^n$.
Suppose for each $i\in [n]$, $h_i(x; \theta)$ is a Lipschitz function w.r.t. $\theta \in \Theta$, with Lipschitz constant equals to $L$. 
Additionally, assume that for any fixed $\theta$, there exists some $\sigma > 0$ such that 
\begin{align}
\label{eqn:ind-h-bound}
	\mathbb{P}\lt(|h_i(X_i; \theta)| < \frac{s\sigma}{n}\rt) < \frac{s}{n},
	\qquad 
	\forall s \in [n]. 
\end{align}
Then for every $\varepsilon \in \real^n$, with probability at least $1 - O(n^{-10})$, it obeys 
\begin{align}
\label{eqn:ind-h-finale}
\sup_{\theta \in \Theta} \sum_{i = 1}^n \ind(|h_i(X_i; \theta)| < \varepsilon_i) \lesssim \log N\Big(\frac{\sigma}{n^2}, \Theta\Big) \log n + \lt(\frac{n \|\varepsilon\|_2}{\sigma}\rt)^{\frac{2}{3}}.
\end{align}
\end{lems}

Finally, we conclude this section by summarizing some standard concentration results of independent Gaussian random vectors. 

Again, denote a collection of independent Gaussian vectors $\{\phi_k\}_{1\leq k\leq t}$ and $\{\psi_k\}_{1\leq k\leq t}$, with $\phi_k \overset{\mathrm{i.i.d.}}{\sim} \mathcal{N}(0,\frac{1}{n}I_n)$ and $\psi_k \overset{\mathrm{i.i.d.}}{\sim} \mathcal{N}(0,\frac{1}{n}I_p)$. 
If we concatenate $\{\phi_{k}\}_{k=1}^{t}$ into a matrix $\Phi \in \real^{n\times t}$ as 
\begin{align}
	\Phi \defn \big[\phi_{1}, \ldots, \phi_{t}\big] \in \real^{n \times t}, 
	\qquad
	\text{where } \phi_k \stackrel{\mathrm{i.i.d.}}{\sim} 
	\mathcal{N}(0,\frac{1}{n} I_n),
\end{align}
its maximum singular value satisfies  
\begin{align}
\label{eqn:concent-op-norm}
	\mprob\left( \|\Phi\|_{\mathrm{op}} \geq 1 + \sqrt{\frac{t}{n}} +  \frac{\delta}{\sqrt{n}} \right) 
	\leq e^{-\delta^2/2}.
\end{align}
In addition, for Wishart matrices, invoking the above result together with a union bound tells us that  
\begin{subequations}
\label{eqn:simple-rm}
\begin{align}
\label{eqn:simple-rm-phi}
	\lt\|(\phi_1, \ldots, \phi_{t-1})^{\top}(\phi_1, \ldots, \phi_{t-1}) - I_{t-1}\rt\|_{\textsf{op}} \lesssim \sqrt{\frac{t\log \frac{n}{\delta}}{n}},\qquad\text{for every } 1 < t \leq n 
\end{align}
with probability at least $1 - \delta$ (see, also, \citet[Example 6.2]{wainwright2019high}). 
Similarly, for random vectors $\{\psi_k\}_{k=1}^{t-1}$, independently drawn from $\mathcal{N}(0,\frac{1}{n}I_p)$, one has 
\begin{align}
\label{eqn:simple-rm-psi}
	\lt\|\frac{n}{p}(\psi_1, \ldots, \psi_{t-1})^{\top}(\psi_1, \ldots, \psi_{t-1}) - I_{t-1}\rt\|_{\textsf{op}} \lesssim \sqrt{\frac{t\log \frac{p}{\delta}}{p}},\qquad\text{for every } 1 < t \leq n.
\end{align}
\end{subequations}

For our convenience, we also recall the following lemma from \cite{li2022non}. Here for every vector $x\in \real^{n}$, we follow the convention and write $|x|_{(i)}$ as its $i$-th largest entry in magnitude. 
\begin{lems}{\cite[Lemma 8]{li2022non}}
\label{lem:brahms-lemma}
With probability at least $1-\delta$, it holds that 
\begin{subequations}
\label{eqn:vive-brahms}
\begin{align}
	\label{eqn:brahms}
	\Big| \max_{1\leq k\leq t-1} \|\phi_k\|_2 -1 \Big| & \lesssim \sqrt{\frac{\log \frac{n}{\delta}}{n}}, \\
	\label{eqn:long}
	 \sup_{a = [a_k]_{1\leq k< t} \in \mathcal{S}^{t-2}} \bigg| \Big\|\sum_{k = 1}^{t-1} a_k\phi_k\Big\|_2 - 1 \bigg| & \lesssim \sqrt{\frac{t\log \frac{n}{\delta}}{n}}, \\
	\label{eqn:vive}
\sup_{a=[a_{k}]_{1\leq k<t}\in\mathcal{S}^{t-2}}\sum_{i=1}^{s}\Big|\sum_{k=1}^{t-1}a_{k}\phi_{k}\Big|_{(i)}^{2}
	&\lesssim\frac{(t+s)\log \frac{n}{\delta}}{n},\qquad \forall 1 \leq s\leq n.	
\end{align}
\end{subequations}
 \end{lems}

For a set of random vectors $\{\phi_k\}_{k=1}^{t}$ independently drawn from $\mathcal{N}(0,\frac{1}{n}I_n)$, 

%


\subsection{Proof of Lemma~\ref{lem:concentration}}
\label{sec:pf-lem-concentration}

For every integer $k\geq 2$, let us consider the $k$-moment of random variable $Z_{i}$. 
For notational convenience, define $Y_i \defn \frac{Z_i}{B \log n}$ and direct calculations yield 
\begin{align}
	\label{eqn:bernstein-condition}
	\notag \mathbb{E}\lt[\lt|Z_i\rt|^k\rt] &= \mathbb{E}\lt[\lt|Z_i\rt|^k\ind \lt(\lt|Z_i\rt| \le Bk\log n\rt)\rt] + \mathbb{E}\lt[\lt|Z_i\rt|^k\ind \lt(\lt|Z_i\rt| > Bk\log n\rt)\rt] \\
	\notag &\stackrel{(\mathrm{i})}{\lesssim} (B k\log n)^{k-2}\mathsf{Var}(Z_i) + (B k\log n)^k\exp(-k\log n) \\
	&\le (B k\log n)^k\lt(\mathsf{Var}(Y_i) + \frac{1}{n^2}\rt),
\end{align}
where the last step uses the definition of $Y_{i}$ and $k\geq 2.$
To verify the relation $(\mathrm{i})$, let us use $\mu_{\lt|Z_i\rt|}$ to denote the density of $\lt|Z_i\rt|$.
By direct calculations, one has 
\begin{align}
	\mathbb{E}\lt[\lt|Z_i\rt|^k\ind \lt(\lt|Z_i\rt| > Bk\log n\rt)\rt]
	&= \int_{Bk\log n}^{\infty} x^k\mu_{\lt|Z_i\rt|}(\dx) \notag\\
	&\le (B k\log n)^k\exp(-k\log n) + \int_{Bk\log n}^{\infty} \exp\Big(-\frac{x}{B}\Big)\dx^k \notag\\
	&= (B k\log n)^k\exp(-k\log n) + k B^k\int_{k\log n}^{\infty} x^{k-1}\exp(-x)\dx \label{eqn:exp-tail-tmp}
\end{align}
where the first inequality invokes the condition~\eqref{eq:Z_condition}. 
To continue, invoking the rule of integration by part, the right hand side of \eqref{eqn:exp-tail-tmp} equals to 
\begin{align}
\eqref{eqn:exp-tail-tmp}	
& = (B k\log n)^k\exp(-k\log n) + kB^k \exp(-k\log n)\sum_{m = 1}^{k-1}\frac{(k-1)!}{(k-m)!}(k\log n)^{k-m}  \notag\\
& \leq (B k\log n)^k\exp(-k\log n) + kB^k \exp(-k\log n)\sum_{m = 1}^{k-1} \frac{1}{\log^{m-1} n} (k\log n)^{k} \notag\\
&\lesssim (B k\log n)^k\exp(-k\log n), \label{eq:exp-tail}
\end{align}
where the first inequality is proved by upper bounding the ratio between consecutive terms by $\frac{1}{\log n}$.
Putting the pieces together establishes the relation $(\mathrm{i})$ and thus the Bernstein condition stated in \eqref{eqn:bernstein-condition}.

Given that each $Z_{i}$ satisfies the Bernstein-type condition~\eqref{eqn:bernstein-condition},
by the power series expansion, we obtain, 
\begin{align*}
\mathbb{E}\lt[\exp\lt(\lambda Z_i\rt)\rt] = \mathbb{E}\lt[\sum_{k=0}^\infty \frac{\lambda^k Z_i^k}{k!}\rt] = 1 + \mathbb{E}\lt[\sum_{k=2}^\infty \frac{\lambda^k Z_i^k}{k!}\rt]
&\le \exp\lt(\sum_{k = 2}^{\infty} \frac{\mathbb{E}\lt[\lambda^kZ_i^k\rt]}{k!}\rt) \\
&\le \exp\lt(\sum_{k = 2}^{\infty} \frac{\lambda^k\lt(Bk\log n\rt)^{k}}{k!}\lt(\mathsf{Var}(Y_i) + \frac{1}{n^2}\rt)\rt) \\
&\le \exp\lt(\frac{e^2\lambda^2}{2}\lt(B\log n\rt)^2\lt(\mathsf{Var}(Y_i) + \frac{1}{n^2}\rt)\rt),
\end{align*}
for any $0 < \lambda < \frac{1}{2eB\log n}$. 
Here in the last step, we use the fact that for $x = \lambda B\log n \le \frac{1}{2e}$, 
\begin{align*}
\sum_{k = 2}^{\infty} \frac{(kx)^k}{k!} 
~\stackrel{(\text{i})}{\le}~ \sum_{k = 2}^{\infty}  \frac{1}{\sqrt{2\pi k}}\lt(ex\rt)^k 
\leq e^2x^2 \sum_{k = 0}^{\infty} \frac{1}{2^k} = \frac{e^2x^2}{2}, 
\end{align*}
where (i) follows from Stirling's formula where $\sqrt{2\pi}k^{k+\frac{1}{2}}e^{-k}\leq k!$. 
The above bound of the moment generating function leads naturally to a high probability control where 
one can apply Markov’s inequality to arrive 
\begin{align*}
\mathbb{P}\lt(\sum_{i = 1}^n Z_i > t\rt) &\le 
\min_{0 < \lambda < \frac{1}{2eB\log n}}\Bigg\{\exp(-\lambda t)\cdot \mathbb{E}\Big[\exp\Big(\lambda \sum_{i = 1}^n Z_i\Big)\Big]\Bigg\} \\
&\le 
\min_{0 < \lambda < \frac{1}{2eB\log n}} \Bigg\{\exp(-\lambda t)\cdot 
\exp\Big(\frac{1}{2}\sum_{i = 1}^n e^2\lambda^2\lt(B\log n\rt)^2\big(\mathsf{Var}(Y_i) + \frac{1}{n^2}\big)\Big)\Bigg\}. 
\end{align*}
Selecting $\lambda$ according to   
\begin{align*}
	\lambda \asymp \min\lt\{\frac{1}{B\log n}, \frac{t}{\sum_{i = 1}^n \lt(B\log n\rt)^2\lt(\mathsf{Var}(Y_i) + \frac{1}{n^2}\rt)\Big)}\rt\},
\end{align*} 
it is ensured that with probability at least $1 - \delta/2$,
\begin{align}
	\sum_{i = 1}^n Z_i \le t \lesssim \max\lt\{B\log n\log\frac{1}{\delta}, ~\sqrt{\sum_{i = 1}^n \Big(\mathsf{Var}(Z_i) + \big(\frac{B\log n}{n}\big)^2\Big)\log\frac{1}{\delta}}\,\rt\}.
\end{align}
Repeating the same argument above for $-Z_{i}$ shows that the above inequality holds for $\sum_{i = 1}^n - Z_i $. 
Putting these together completes the proof of relation~\eqref{eqn:sum-zi}.


\subsection{Proof of Lemma~\ref{lem:ind-noise}}
\label{sec:pf-lem-ind-noise}

Given any fixed $\theta\in \Theta$, independent of $\{X_i\}$,  let us consider random variables $\ind(|h_i(X_i; \theta)| < \frac{s\sigma}{n})$ for $i\in [n]$. 
In view of inequality~\eqref{eqn:ind-h-bound}, $\ind(|h_i(X_i; \theta)| < \frac{s\sigma}{n})$ forms a set of independent Bernoulli random variables with parameter smaller than $\frac{s}{n}.$
Invoking Bernstein's inequality (see, e.g.~\cite[Chapter 2]{wainwright2019high}) ensures that 
\begin{align}
\label{eqn:h-bernstein}
	\sum_{i = 1}^n \ind\lt(|h_i(X_i; \theta)| < \frac{s\sigma}{n}\rt) \leq s + \sqrt{2 s\log \frac{1}{\delta}} + 2 \log \frac{1}{\delta},
\end{align}
with probability at least $1 - \delta.$ 
In order to deal with random $\thetahat$, consider an $\epsilon$-cover of $\Theta$ of $\ell_2$-norm and denote it by $\mathcal{N}_{\epsilon}$. 
By virtue of the Lipschitz property of $h$, there exists some $\theta \in \mathcal{N}_{\epsilon}$ such that 
\begin{align*}
	|h_i(X_i; \thetahat) - h_i(X_i; \theta)| \leq L\epsilon,
\end{align*}
which implies 
\begin{align*}
	\ind(|h_i(X_i; \thetahat)| < x_i) 
	\leq 
	\ind(|h_i(X_i; \theta)| < x_i + L\epsilon),
	\qquad
	\text{for every }x_i
\end{align*}
and hence, 
\begin{align}
	\sup_{\thetahat \in \Theta} \sum_{i = 1}^n \ind(|h_i(X_i; \thetahat)| < x_i) 
	\le 
	\sup_{\theta\in \mathcal{N}_{\epsilon}} \sum_{i = 1}^n \ind(|h_i(X_i; \theta)| < x_i + L\epsilon).
\end{align}

Before diving into our main proof, let us state a key property regarding $h_i(X_i; \thetahat)$'s based on the above observation. 
In particular, select parameters 
\begin{align*}
 	x_i = \frac{\sigma s}{4n}, \quad L\epsilon = \frac{\sigma}{100 n^2}, 
 	~\text{ and }~ \delta = \frac{1}{n^{11} N(\epsilon, \Theta)}.
\end{align*} 
Taking a union bound of \eqref{eqn:h-bernstein} over the $\epsilon$-cover $\mathcal{N}_{\epsilon}$, we arrive at 
\begin{align}
	\sup_{\thetahat \in \Theta} \sum_{i = 1}^n \ind\Big(|h_i(X_i; \thetahat)| < \frac{s \sigma}{4 n}\Big) 
	\le 
	\sup_{\theta\in \mathcal{N}_{\epsilon}} \sum_{i = 1}^n \ind\Big(|h_i(X_i; \theta)| < \frac{s \sigma}{4 n} + L\epsilon\Big)
	\le 
	\frac{s}{3} + O\Big(\log N(\frac{\sigma}{100 n^2}, \Theta) \log n\Big),
\end{align}
with probability at least $1 - O(n^{-11}).$
In words, for every $\thetahat \in \Theta$ and $s \gtrsim \log N(\frac{\sigma}{100 n^2}, \Theta)\log n$, the total number of index $i$, such that $|h_i(X_i; \thetahat)| < s \sigma/(4 n)$ is with high probability smaller than $s/2.$
It implies that, if we rank the magnitude of $|h_i(X_i; \thetahat)|$ from the smallest to the largest, then regardless of the value of $\thetahat$, it always holds true with probability at least $1 - O(n^{-11})$ that 
\begin{align}
\label{eqn:R}
	 \sum_{i=1}^s h^2_{(i)}(X_i, \thetahat) \geq \frac{s}{2} \cdot (\frac{s \sigma}{4 n})^2
	 =
	 \frac{s^3 \sigma^2}{n^2}.
\end{align}
We emphasize that \eqref{eqn:R} holds true for every $\thetahat$ and $s$, as long as $s \gtrsim \log N(\frac{\sigma}{100 n^2}, \Theta) \log n.$

We are now ready to establish the proof of \eqref{eqn:ind-h-finale}. 
For every $\varepsilon \in \real^n$ and $\theta \in \mathcal{N}_\epsilon$, let us define set 
\begin{align*}
 \mathcal{I}_\theta := \lt\{i : |h_i(X_i; \theta)| < \varepsilon_i + L\epsilon\rt\}.
\end{align*}
Then according to this definition, one naturally has 
\begin{align}
	\sum_{i = 1}^n \ind(|h_i(X_i; \theta)| < \varepsilon_i + L\epsilon)
	=
	|\mathcal{I}_\theta|,
\end{align}
and as well as 
\begin{align}
\label{eqn:J}
	\sum_{i = 1}^{|\mathcal{I}_\theta|} h_{(i)}^2 \le \sum_{i = 1}^{|\mathcal{I}_\theta|} (\varepsilon_i + L\epsilon)^2 \lesssim 
	\|\varepsilon\|^2_2 + \frac{\sigma^2}{n^4},
\end{align}
where the last steps invokes the choice of $L\epsilon = \frac{\sigma}{100 n^2}.$
If $|\mathcal{I}_\theta|\lesssim \log N(\frac{\sigma}{100 n^2}, \Theta)  \log n$ for every $\theta \in \mathcal{N}_\epsilon$, then it is straightforward to see that 
\begin{align*}
	\sup_{\thetahat \in \Ball(r)} \sum_{i = 1}^n \ind\Big(|h_i(X_i; \thetahat)| < \varepsilon_i\Big) 
	\le 
	\sup_{\theta\in \mathcal{N}_{\epsilon}} \sum_{i = 1}^n \ind\Big(|h_i(X_i; \theta)| < \varepsilon_i + L\epsilon\Big)
	=
	\sup_{\theta\in \mathcal{N}_{\epsilon}} |\mathcal{I}_\theta|
	\lesssim \log N(\frac{\sigma}{100 n^2}, \Theta) \log n.
\end{align*}
Otherwise, taking collectively inequality~\eqref{eqn:J} with inequality~\eqref{eqn:R}, one has 
\begin{align*}
	\forall \theta\in \mathcal{N}_\epsilon, \qquad \frac{|\mathcal{I}_\theta|^3\sigma^2}{n^2} \lesssim \|\varepsilon\|^2_2 + \frac{\sigma^2}{n^4},,
\end{align*}
which further implies 
\begin{align}
	\sum_{i = 1}^n \ind(|h_i(X_i; \theta)| < \varepsilon_i + L\epsilon)
	=
	|\mathcal{I}_\theta|
	\lesssim
	\lt(\frac{n \|\varepsilon\|_2}{\sigma}\rt)^{\frac{2}{3}}.
\end{align}
Thus, in this case, we are guaranteed that 
\begin{align}
	\sup_{\thetahat \in \Ball(r)} \sum_{i = 1}^n \ind\Big(|h_i(X_i; \thetahat)| < \varepsilon_i\Big) 
	\le 
	\sup_{\theta\in \mathcal{N}_{\epsilon}} \sum_{i = 1}^n \ind\Big(|h_i(X_i; \theta)| < \varepsilon_i + L\epsilon\Big)
	\lesssim \lt(\frac{n \|\varepsilon\|_2}{\sigma}\rt)^{\frac{2}{3}}.
\end{align}
Putting these two cases together completes the proof of property~\eqref{eqn:ind-h-finale}.

\section{Proof of Lemma~\ref{lem:magic-flute-simple}}
\label{sec:pf-thm-magic-flute}

\subsection{A general statement}

In this section, we prove a more general version of Lemma~\ref{lem:magic-flute-simple} without imposing the Assumption~\ref{ass:control-lipschitz}. This general result reduces to Lemma~\ref{lem:magic-flute-simple} in the special case.

To simplify our statement, we start by introducing some auxiliary notation. 
Specifically, let us define 
\begin{subequations}
\label{eqn:defn-alpha-gamma}
\begin{align}
\overalpha_t &:= \max_{k \le t} \Big\{\lt\|G_{k}(0)\rt\|_2, \|\alpha_k\|_2,\frac{1}{\poly(n)}\Big\}, \label{eqn:defn-alpha-t-bar}\\
\overgamma_{t} &:= \max_{k \le t} \Big\{\lt\|F_{k}(0)\rt\|_2, \|\gamma_{k}\|_2,\frac{1}{\poly(n)}\Big\}, \label{eqn:defn-gamma-t-bar}
\end{align}
and we write  
\begin{align}
\label{eqn:mu-nv-main}
\mu_t := \frac{\rho_F\rho_G\overalpha_{t-1} + \rho_G\overgamma_{t-1}}{\|\alpha_{t-1}\|_2}, 
\qquad
\nu_t := \frac{\rho_F\rho_G\overgamma_{t} + \rho_F\overalpha_{t-1}}{\|\gamma_t\|_2},
\end{align}
and
\begin{align}
\overline{\xi}_t &\defn \overalpha_{t-1}\Bigg(\hatalpha_{t-1}^{t-1} + \rho_F^2\hatgamma_{t-1}^{t-1} + \rho_F \sqrt{\frac{t\log^2 n}{n}}\Bigg) + \sqrt{\frac{t\log^2 n}{n}}\rho_F\rho_G\overgamma_{t}, \\
\overline{\zeta}_t &\defn \overgamma_{t}\Bigg(\hatgamma_{t}^{t} + \rho_G^2\hatalpha_{t-1}^{t-1}+  \rho_G\sqrt{\frac{t\log^2 n}{n}} \Bigg) + \sqrt{\frac{t\log^2 n}{n}}\rho_F\rho_G\overalpha_{t}.
\end{align}
\end{subequations}
Our goal is to establish the following claim in order to control the sizes of $\xi_{t}$ and $\zeta_{t}$. 
This claim is of the same form as in Claim~\ref{claim:main}. 

\begin{claim}
\label{claim:main-general}
There exists universal constant $0 < c < 1$, such that the following set of inequalities hold true 
\begin{subequations}
\label{eqn:main}
\begin{align}
&\|\hatxi_{t}\|_2 \lesssim \sqrt{\frac{t\log^2 n}{n}}\lt(\rho_F\rho_G\overgamma_{t} + \rho_F\overalpha_{t-1}\rt), \label{eq:xihat}\\
&\|\hatzeta_{t}\|_2 \lesssim \sqrt{\frac{t\log^2 n}{n}}\lt(\rho_F\rho_G\overalpha_{t} + \rho_G\overgamma_{t}\rt), \label{eq:zetahat}\\
&|\hatalpha_{t-1}^{t-1}| \lesssim \sqrt{\frac{t\log^2 n}{n^2}}\rho_{1, F}\lt(\rho_F\rho_G\overalpha_{t-1} + \rho_G\overgamma_{t-1}\rt) + \rho_F\lt(\frac{\mu_t^2t\log^2 n}{n}\rt)^{\frac{1}{3}}, \label{eq:alphahat}\\
&|\hatgamma_t^t| \lesssim \sqrt{\frac{t\log^2 n}{n^2}}\rho_{1, G}\lt(\rho_F\rho_G\overgamma_{t} + \rho_F\overalpha_{t-1}\rt) + \rho_G\lt(\frac{\nu_t^2t\log^2 n}{n}\rt)^{\frac{1}{3}}, \label{eq:gammahat}\\
&|\hatalpha_{t-1}^k| 
\le \left\{
\begin{array}{lll}
(1-c)^{t-k-1} \Big|\hatalpha_{\frac{t+k-1}{2}}^{\frac{t+k-1}{2}}\Big| & \text{if} & t-1-k = 2m, \\[0.4cm]
(1-c)^{t-k-2}\rho_F^2 \Big|\hatgamma_{\frac{t+k}{2}}^{\frac{t+k}{2}}\Big| & \text{if} & t-1-k = 2m+1,
\end{array}
\right. \label{eq:alphahat-k} \\
&|\hatgamma_{t}^k| \le \left\{
\begin{array}{lll}
(1-c)^{\frac{t-k}{2}} |\hatgamma_{\frac{t+k}{2}}^{\frac{t+k}{2}}| & \text{if} & t-k = 2m, \\
(1-c)^{\frac{t-k}{2}}\rho_G^2 |\hatalpha_{\frac{t+k-1}{2}}^{\frac{t+k-1}{2}}| & \text{if} & t-1-k = 2m.
\end{array}
\right. \label{eq:gammahat-k}
\end{align}
\end{subequations}
\end{claim}

Let us state an inductive results regarding the above Claim~\ref{claim:main-general}. The proof of this result is provided in Section~\ref{sec:pf-general-gem}.

\begin{lems}
\label{lem:magic-flute}
Under the decomposition~\eqref{defi:dynamics} with~\eqref{eqn:residual-second}, the bound~\eqref{eqn:main} holds for $t = 1$
with probability at least $1 - n^{-10}$. 
In addition, with probability at least $1 - O(n^{-10})$, for every $t$ satisfying  
\begin{align}
\label{eqn:assmp-t-range}
t \ll \frac{n}{\rho_F^4\rho_G^4\log^4 n},
\end{align}
if the bound~\eqref{eqn:main} and Assumption~\ref{ass:control} below hold for $t$, then the bound~\eqref{eqn:main} holds for $t+1$.
\end{lems}

\begin{assumption} \label{ass:control}
For the decomposition~\eqref{defi:dynamics} with~\eqref{eqn:residual-second}, assume the following conditions hold:
\begin{itemize}
\item there exists some universal constant $0 < c < 1/2$, such that 
\begin{align}
\label{eqn:assmp-prime}
\frac{1}{n^2}\mathbb{E}\big[\lt\|G_{t}^{\prime}(\widetilde{u}_{t})\rt\|_2^2 \mymid \|\gamma_t\|_2\big]\mathbb{E}\big[\lt\|F_{t+1}^{\prime}(\widetilde{v}_{t+1})\rt\|_2^2 \mymid \|\alpha_t\|_2\big] &< (1-2c)^2, \\
\frac{1}{n^2}\mathbb{E}\big[\lt\|F_{t+1}^{\prime}(\widetilde{v}_{t+1})\rt\|_2^2 \mymid \|\alpha_t\|_2\big]\mathbb{E}\big[\lt\|G_{t+1}^{\prime}(\widetilde{u}_{t+1})\rt\|_2^2 \mymid \|\gamma_{t+1}\|_2\big] &< (1-2c)^2.
\end{align}

\item for some universal constant $0 \leq c < 1$, such that 
\begin{subequations}
\begin{align}
\label{eqn:assmp-gamma-tail}
&|\hatalpha_{t-1}^k|, |\hatgamma_{t}^k| \le c^{t-k}\poly(n),\\
\label{eqn:assmp-alphat-gammat}
&\frac{1}{\poly(n)}\rho_{1, F}\overalpha_t \lesssim \rho_F \lesssim \poly(n) \qquad\text{and} \qquad \frac{1}{\poly(n)}\rho_{1, G}\overgamma_t \lesssim \rho_G \lesssim \poly(n), \\
\label{eqn:assmp-bach}
&\rho_G\sum_{k = 1}^{t} |\hatalpha_{t-1}^k| \overgamma_{t} \ll \|\gamma_t\|_2 
\qquad\text{and} \qquad 
\rho_F\sum_{k = 1}^{t} |\hatgamma_{t}^k| \overalpha_t \ll \|\alpha_t\|_2;
\end{align}
\end{subequations}

\item in addition, we assume 
\begin{subequations}
\label{ass:Onsager}
\begin{align}
\label{eqn:ballet}
&\rho_F\rho_G^2\sum_{k = 1}^{t-1} |\hatalpha_{t-1}^k| \overline{\xi}_{k - 1} \lesssim \sqrt{\frac{t\log^2 n}{n}}\lt(\rho_F\overalpha_{t-1} + \overgamma_{t}\rt), \\
\label{eqn:serenade}
&\rho_F^2\rho_G\sum_{k = 1}^{t} |\hatgamma_t^k| \overline{\zeta}_{k - 1} \lesssim \sqrt{\frac{t\log^2 n}{n}}\lt(\rho_G\overgamma_{t} + \overalpha_{t}\rt), \\
\label{eqn:assmp-rho-a}
&\rho_{1, G}\sqrt{\frac{t\log n}{n}}\overline{\xi}_t + \rho_G\sqrt{\log n}\lt(\frac{\overline{\xi}_t}{\|\gamma_t\|_2}\rt)^{\frac{1}{3}} \ll \frac{1}{\rho_F}, \\
\label{eqn:assmp-rho-b}
&\rho_{1, F}\sqrt{\frac{t\log n}{n}}\overline{\zeta}_t + \rho_F\sqrt{\log n}\lt(\frac{\overline{\zeta}_t}{\|\alpha_t\|_2}\rt)^{\frac{1}{3}} \ll \frac{1}{\rho_G}, \\
\label{eqn:assmp-rho-c}
&\frac{1}{n}\rho_G\lt[ \sqrt{n}\rho_{1, G}\overalpha_{t-1}(\hatalpha_{t-1}^{t-1} + \rho_F^2\hatgamma_{t-1}^{t-1}) + \rho_Gn\lt(\frac{\overalpha_{t-1}(\hatalpha_{t-1}^{t-1} + \rho_F^2\hatgamma_{t-1}^{t-1})}{\|\gamma_t\|_2}\rt)^{\frac{2}{3}}\rt] \ll \frac{1}{\rho_F^2}, \\
\label{eqn:assmp-rho-d}
&\frac{1}{n}\rho_F\lt[ \sqrt{n}\rho_{1, F}\overgamma_{t}(\hatgamma_{t}^{t} + \rho_G^2\hatalpha_{t-1}^{t-1}) + \rho_Fn\lt(\frac{\overgamma_{t}(\hatgamma_{t}^{t} + \rho_G^2\hatalpha_{t-1}^{t-1})}{\|\alpha_t\|_2}\rt)^{\frac{2}{3}}\rt] \ll \frac{1}{\rho_G^2}.
\end{align}
\end{subequations}
\end{itemize}
\end{assumption}

\begin{remark}
Under the Assumptions~\ref{ass:control-lipschitz} and \ref{ass:control-simple}, and further assuming that the Claim~\ref{claim:main} holds true at $t$, it is straightforward to verify Assumption~\ref{ass:control}. Crucially, taking $\rho_{1,F}, \rho_{1,G} = 0$, $\rho_F, \rho_G \asymp 1$ and $\|\gamma_t\|_2, \|\alpha_t\|_2 \asymp 1$ helps simplify the presentations of the above formulas to a large extent. 
\end{remark}

\subsection{Proof of Lemma~\ref{lem:magic-flute}}
\label{sec:pf-general-gem}

Before diving into the proof of Lemma~\ref{lem:magic-flute}, let us first state a few preliminaries. 
Throughout this proof, we condition on the event
\begin{equation}
\label{eqn:rm-event}
\begin{aligned}
	&\lt\|(\phi_1, \ldots, \phi_{t-1})^{\top}(\phi_1, \ldots, \phi_{t-1}) - I_{t-1}\rt\|_{\textsf{op}} \lesssim \sqrt{\frac{t\log \frac{n}{\delta}}{n}},\qquad\text{for every } 1 < t \leq n, \\
	\text{and }~
	&\lt\|\frac{n}{p}(\psi_1, \ldots, \psi_{t-1})^{\top}(\psi_1, \ldots, \psi_{t-1}) - I_{t-1}\rt\|_{\textsf{op}} \lesssim \sqrt{\frac{t\log \frac{p}{\delta}}{p}},\qquad\text{for every } 1 < t \leq n,
\end{aligned}
\end{equation}
both of which hold with probability at least $1-\delta$ according to inequalities \eqref{eqn:concent-op-norm} and \eqref{eqn:simple-rm-psi}; in this proof, we shall take $\delta = n^{-11}.$
In addition, the following results turn out to be essential for our analysis, whose proof is deferred to Section~\ref{sec:proof-lem:Onsager}.
\begin{lems} 
\label{lem:Onsager}
Under the assumptions \eqref{eqn:assmp-alphat-gammat} -- \eqref{eqn:assmp-t-range}, the following two inequalities hold true with probability at least $1 - O(n^{-11})$:
\begin{align}
\label{eqn:Onsager-brahms}
&\lt\|\sum_{k = 1}^{t} a_k\Big[\lt\langle \psi_k, F_{t+1}(\hatbeta_{t+1})\rt\rangle - \langle F_{t+1}^{\prime}(\hatbeta_{t+1}) \rangle \alpha_{t}^k - \sum_{j = k+1}^{t} \hatgamma_{t}^j \lt\langle F_{t+1}^{\prime}(\hatbeta_{t+1}) \circ F_{j}^{\prime}\lt(v_{j}\rt)\rt\rangle\alpha_{j-1}^k\Big]\rt\|_2 \notag\\
&\qquad \lesssim \sqrt{\frac{t\log^2 n}{n}}\lt(\overgamma_{t+1} + \rho_F\overalpha_{t}\rt),
\end{align}
and
\begin{align}
\label{eqn:Onsager-beethoven}
&\lt\|\sum_{k = 1}^{t+1} b_k
\Big[\lt\langle \phi_k, G_{t+1}(\hats_{t+1})\rt\rangle - \langle G_{t+1}^{\prime}(\hats_{t+1}) \rangle \gamma_{t+1}^k - \sum_{j = k}^{t} \hatalpha_{t}^j \lt\langle G_{t+1}^{\prime}(\hats_{t+1}) \circ G_{j}^{\prime}\lt(u_{j}\rt)\rt\rangle\gamma_{j}^k\Big]\rt\|_2 \notag\\
&\qquad \lesssim \sqrt{\frac{t\log^2 n}{n}}\lt(\overalpha_{t+1} + \rho_G\overgamma_{t+1}\rt). 
\end{align}
\end{lems}

\paragraph{Initial case for $t = 1$.} 
In order to prove inequalities~\eqref{eqn:main} for $t=1$, first recall the initialization condition that 
\begin{align*}
	s_{1} = r_{1} - \epsilon = X\thetastar.
\end{align*}
By construction of \eqref{eqn:a1-b1} and \eqref{eqn:def-phi-k}, $\phi_{1} = \thetastar/\ltwo{\thetastar}$ and $\gamma^1_{1} = \ltwo{\thetastar}$ and hence $\xi_{1} = \hatxi_1 = 0.$
In addition, since $\hats_{1} = s_{1}$, it is guaranteed that  
\begin{align*}
	\hatgamma_{1}^1 &\defn
	\langle G_{1}^{\prime}(\hats_{1}) - G_{1}^{\prime}(s_{1})\rangle + \frac{1}{\|\gamma_{1}\|_2}
	\langle \phi_1, G_{1}(s_{1}) - G_{1}(\hats_{1})\rangle = 0. 
\end{align*}
As a result, inequalities~\eqref{eq:xihat} and \eqref{eq:gammahat} hold for $t=1$. The requirements for coefficient $|\hatalpha_{t}|$ naturally holds as they equal to zero when $t=0.$
It suffices to establish the required result for $\ltwo{\hatzeta_1}$. Towards this, by construction~\eqref{eqn:orth-expa} and definition~\eqref{eqn:beta-zeta}, we observe that  
\begin{align*}
	\alpha^1_{1} &= \inprod{G_1(s_1)}{a_1} = \ltwo{G_1(s_1)} \\
	\zeta_1 &= \hatzeta_1 = b_1 \lt[\lt\langle \phi_1, G_{1}(s_1)\rt\rangle - \langle G_1^{\prime}\rangle \gamma_{1}^1 - \alpha_1^1q_1^1\rt].
\end{align*}
To obtain a control of the right hand side of $\hatzeta_{1}$, expression~\eqref{eqn:Onsager-beethoven} of Lemma~\ref{lem:Onsager} --- whose assumptions satisfy naturally as both $\hatalpha_{0}$ and $\hatgamma_{1}^1$ vanish --- ensures that 
\begin{align}
\ltwo{\hatzeta_1} = \lt|\lt\langle \phi_1, G_{1}(s_1)\rt\rangle - \langle G_1^{\prime}\rangle \gamma_{1}^1 - \alpha_1^1q_1^1\rt| \lesssim \sqrt{\frac{\log n}{n}}\lt(\overalpha_{1} + \rho_G\overgamma_{1}\rt),
\end{align}
from which, we complete the proof of \eqref{eq:zetahat}.


\paragraph{Inductive relation.} 
Suppose both Assumption~\ref{ass:control} and the target conclusion~\eqref{eqn:main} hold at the $t$-th iteration.
We shall prove the inequality set~\eqref{eqn:main} at the $t+1$-th iteration.
First, we remark that given~\eqref{eqn:main} holds at iteration $t$, decomposition~\eqref{eq:xi_expansion} leads to  
\begin{align}
\label{eqn:bruckner4}
\|\xi_{t}\|_2 &\le \sum_{k = 1}^{t-1} |\hatalpha_{t-1}^k| \ltwo{G_k(s_k)} + \ltwo{\hatxi_t} \notag \\
&\stackrel{(\mathrm{i})}{\lesssim} \sum_{k=1}^{t-1} |\hatalpha_{t-1}^k|\cdot\|\alpha_k\|_2 + \sqrt{\frac{t\log^2 n}{n}}\lt(\rho_F\rho_G\overgamma_{t} + \rho_F\overalpha_{t-1}\rt) \notag\\
&\stackrel{(\mathrm{ii})}{\lesssim} \overalpha_{t-1}\Bigg(\hatalpha_{t-1}^{t-1} + \rho_F^2\hatgamma_{t-1}^{t-1} + \rho_F \sqrt{\frac{t\log^2 n}{n}}\Bigg) + \sqrt{\frac{t\log^2 n}{n}}\rho_F\rho_G\overgamma_{t} =: \overline{\xi}_t.
\end{align}
Here $(\mathrm{i})$ invokes the relation that $\ltwo{G_t(s_t)} = \ltwo{\alpha_t}$ and the inductive assumption~\eqref{eq:xihat}; $(\mathrm{ii})$ follows from the geometric decay of $\hatalpha_{t}^k$ in expression~\eqref{eq:alphahat-k}. 
Similarly, in view of expression~\eqref{eq:zeta_expansion}, one has  
\begin{align}
\label{eqn:bruckner4-adagio}
\|\zeta_{t}\|_2  
&\le \sum_{k=1}^{t} |\hatgamma_{t}^k| \|F_k(\beta_k)\|_2 + \ltwo{\hatzeta_t}\notag \\
&\lesssim \sum_{k=1}^{t} |\hatgamma_{t}^k|\cdot\|\gamma_k\|_2 + \sqrt{\frac{t\log^2 n}{n}}\lt(\rho_F\rho_G\overalpha_{t} + \rho_G\overgamma_{t}\rt) \notag\\
&\lesssim \overgamma_{t}\Bigg(\hatgamma_{t}^{t} + \rho_G^2\hatalpha_{t-1}^{t-1}+  \rho_G\sqrt{\frac{t\log^2 n}{n}} \Bigg) + \sqrt{\frac{t\log^2 n}{n}}\rho_F\rho_G\overalpha_{t} =: \overline{\zeta}_t
\end{align}
In addition, we obtain in the similar fashion that 
\begin{align}
\lt\|\hats_{t} - u_{t}\rt\|_2 
&= \Big\| \sum_{k = 1}^{t-1} \hatalpha_{t-1}^k G_{k}(u_{k})\Big\|_2 
\lesssim \overalpha_{t-1}(\hatalpha_{t-1}^{t-1} + \rho_F^2\hatgamma_{t-1}^{t-1}) \le \overline{\xi}_t, \\
\lt\|\hatbeta_{t+1} - v_{t+1}\rt\|_2 
&= \Big\|\sum_{k = 1}^{t} \hatgamma_{t}^k F_{k}(v_{k})\Big\|_2 
\lesssim \overgamma_{t}(\hatgamma_{t}^{t} + \rho_G^2\hatalpha_{t-1}^{t-1}) \le \overline{\zeta}_t.\label{eqn:bruckner9}
\end{align}
With these control in place, let us verify the induction results for the next iteration based on Assumption~\ref{ass:control}. 


\subsubsection{Induction step for quantities $\ltwo{\hatxi_{t+1}}$ and $\ltwo{\hatzeta_{t+1}}$}

Let us start by showing that expression~\eqref{eq:xihat} holds at the $t+1$-th iteration. 
In view of expression~\eqref{eq:hatxi}, it directly satisfies that  
\begin{align}
	\notag \ltwo{\hatxi_{t+1}} & \leq 
	\Big\|\sum_{k = 1}^{t} a_k\Big[\lt\langle \psi_k, F_{t+1}(\hatbeta_{t+1})\rt\rangle - \langle F_{t+1}^{\prime}(\hatbeta_{t+1}) \rangle \alpha_{t}^k - \sum_{j = k+1}^{t} \hatgamma_{t}^j \lt\langle F_{t+1}^{\prime}(\hatbeta_{t+1}) \circ F_{j}^{\prime}\lt(v_{j}\rt)\rt\rangle\alpha_{j-1}^k\Big] \Big\|_2\\
	&\quad+ 
	\Big\|\mathcal{P}_{G_t(s_t)}^{\perp}\sum_{k = 1}^{t} a_k\lt\langle \psi_k, F_{t+1}(\beta_{t+1}) - F_{t+1}(\hatbeta_{t+1})\rt\rangle\Big\|_2 + O\lt(\sqrt{\frac{t\log n}{n}}\|\gamma_{t+1}\|_2\rt).
\end{align}
It is then sufficient to control the above two terms on the right accordingly. 
Recall that Lemma~\ref{lem:Onsager} ensures that 
\begin{align}
	\lt\|\sum_{k = 1}^{t} a_k
	\Big[\lt\langle \psi_k, F_{t+1}(\hatbeta_{t+1})\rt\rangle - \langle F_{t+1}^{\prime}(\hatbeta_{t+1}) \rangle \alpha_{t}^k - \sum_{j = k+1}^{t} \hatgamma_{t}^j \lt\langle F_{t+1}^{\prime}(\hatbeta_{t+1}) \circ F_{j}^{\prime}\lt(v_{j}\rt)\rt\rangle\alpha_{j-1}^k\Big]\rt\|_2
	&\lesssim \sqrt{\frac{t\log^2 n}{n}}\lt(\overgamma_{t+1} + \rho_F\overalpha_t\rt), \label{eq:Onsager-F} \\
	\lt\|\sum_{k = 1}^{t+1} b_k
	\Big[\lt\langle \phi_k, G_{t+1}(\hats_{t+1})\rt\rangle - \langle G_{t+1}^{\prime}(\hats_{t+1}) \rangle \gamma_{t+1}^k - \sum_{j = k}^{t} \hatalpha_{t}^j \lt\langle G_{t+1}^{\prime}(\hats_{t+1}) \circ G_{j}^{\prime}\lt(u_{j}\rt)\rt\rangle\gamma_{j}^k\Big]\rt\|_2
	&\lesssim \sqrt{\frac{t\log^2 n}{n}}\lt(\overalpha_{t+1} + \rho_G\overgamma_{t+1}\rt). \label{eq:Onsager-G}
\end{align}
which completes the control of the first term.

Regarding the second term, recall that $\{a_k\}$ (defined in expression~\eqref{eqn:defn-abx}) forms a set of orthonormal basis.
There exists a unit vector $w \in G_t(s_t)^{\perp} \cap \mathsf{span}\{a_k\}$ such that 
\begin{align}
\label{eqn:ortho-proj}
\notag \left\|\mathcal{P}_{G_t(s_t)}^{\perp}\sum_{k = 1}^{t} a_k\lt\langle \psi_k, F_{t+1}(\beta_{t+1}) - F_{t+1}(\hatbeta_{t+1})\rt\rangle  \right\|_2
	&= w^{\top}\mathcal{P}_{G_t(s_t)}^{\perp}\sum_{k = 1}^{t} a_k\lt\langle \psi_k, F_{t+1}(\beta_{t+1}) - F_{t+1}(\hatbeta_{t+1})\rt\rangle \\
	&= \lt\langle\sum_{k=1}^{t} \omega_k a_k,
	\sum_{k = 1}^{t} a_k\lt\langle \psi_k, F_{t+1}(\beta_{t+1}) - F_{t+1}(\hatbeta_{t+1})\rt\rangle 
	\rt\rangle,
\end{align}
where $\omega \in \mathcal{S}^{t-1}$ for $\omega_k = w^{\top}a_k$. 
In view of this decomposition, we remark that as stated in \eqref{eqn:orthogonal}, one has  
\begin{align}
\label{eqn:orthogonal-w-a}
 	0 = \lt\langle \sum_{k=1}^{t} \omega_k a_k, \sum_{k = 1}^{t} \alpha_{t}^ka_k \rt\rangle = 
 	\inprod{\omega}{\alpha_t}.
 \end{align} 
In addition, there exists some $\overline{\zeta}\in \real^p$ such that  
\begin{align}
\label{eqn:intuition}
\left\|\mathcal{P}_{G_t(s_t)}^{\perp}\sum_{k = 1}^{t} a_k\lt\langle \psi_k, F_{t+1}(\beta_{t+1}) - F_{t+1}(\hatbeta_{t+1})\rt\rangle \right\|_2
%
%
\notag &= \lt\langle \sum_{k = 1}^{t} \omega_k\psi_k, F_{t+1}(\beta_{t+1}) - F_{t+1}(\hatbeta_{t+1})\rt\rangle \\
\notag &= \lt\langle \sum_{k = 1}^{t} \omega_k\psi_k, F_{t+1}^{\prime}(v_{t+1} + \overline{\zeta}) \circ \lt(\beta_{t+1} - \hatbeta_{t+1}\rt)\rt\rangle \\
&\le \Big\| \sum_{k = 1}^{t} \omega_k\psi_k \circ F_{t+1}^{\prime}(v_{t+1} + \overline{\zeta})\Big\|_2 \lt\|\beta_{t+1} - \hatbeta_{t+1}\rt\|_2.
\end{align}
Therefore, it is enough to control the two terms on the right hand side above, which is what shall be done in the following.

\begin{itemize}
\item \textbf{Control of $\|\beta_{t+1} - \hatbeta_{t+1}\|_2$.} Here, the size of $\overline{\zeta}$ can be bounded by 
\begin{align}
\label{eqn:zeta-bar}
	\|\overline{\zeta}\|_2 \le \|\hatbeta_{t+1} - v_{t+1}\|_2 + \|\beta_{t+1} - v_{t+1}\|_2 
	\leq  \|\hatbeta_{t+1} - v_{t+1}\|_2 + \ltwo{\zeta_t}
	\lesssim \overline{\zeta}_t,
\end{align}
where the last relation follows from the derivations in \eqref{eqn:bruckner4-adagio} and \eqref{eqn:bruckner9}.
Putting the pieces above together, we arrive at 
\begin{align}
	\ltwo{\hatxi_{t+1}} &\le
	\Big\| \sum_{k = 1}^{t} \omega_k\psi_k \circ F_{t+1}^{\prime}(v_{t+1} + \overline{\zeta})\Big\|_2 \lt\|\beta_{t+1} - \hatbeta_{t+1}\rt\|_2
	+
	O\Big(\sqrt{\frac{t\log^2 n}{n}}\lt(\overgamma_{t+1} + \rho_F\overalpha_t\rt)\Big).
\end{align}
To further control the right hand side, recall the definitions that $\hatbeta_{t+1} = v_{t+1} + \sum_{k = 1}^{t} \hatgamma_{t}^k F_{k}(v_{k})$ and 
$\beta_{t+1} = v_{t+1} + \zeta_t = v_{t+1} + \sum_{k = 1}^{t} \hatgamma_{t}^k F_{k}(\beta_{k}) + \hatzeta_t$. 
Therefore, we have 
\begin{align}
\label{eqn:beta-betahat-1}
	\notag \|\beta_{t+1} - \hatbeta_{t+1}\|_2 
	&\leq \ltwo{\hatzeta_t} + \Big\|\sum_{k = 1}^{t} \hatgamma_{t}^k (F_{k}(\beta_{k}) - F_{k}(v_{k}))\Big\|_2\\
	&\leq \ltwo{\hatzeta_t}	+ O\Big(\frac{1}{\rho_F}\sqrt{\frac{t\log^2 n}{n}}\lt(\overgamma_{t+1} + \overalpha_t\rt)\Big).
\end{align}
Here, in the last inequality, we make the observation that 
\begin{align*}
	\Big\|\sum_{k = 1}^{t} \hatgamma_{t}^k (F_{k}(\beta_{k}) - F_{k}(v_{k}))\Big\|_2 
	\leq \sum_{k = 1}^{t} |\hatgamma_t^k| \rho_F \ltwo{\zeta_{k-1}} 
	\leq \rho_F \sum_{k = 1}^{t} |\hatgamma_t^k| \overline{\zeta}_{k-1} 
	\lesssim \frac{1}{\rho_F}\sqrt{\frac{t\log^2 n}{n}}\lt(\overgamma_{t} + \overalpha_t\rt),
\end{align*}
where the last inequality follows from the assumptions~\eqref{eqn:serenade} and $\rho_{G} \geq 1.$
Now invoking our inductive assumption~\eqref{eq:zetahat} and the assumption that $\rho_{F}, \rho_{G} \geq 1$, we arrive at 
\begin{align}
\label{eqn:beta-betahat-2}
\|\beta_{t+1} - \hatbeta_{t+1}\|_2 	
	&\lesssim \sqrt{\frac{t\log^2 n}{n}}\lt(\rho_F\rho_G\overalpha_{t} + \rho_G\overgamma_{t}\rt).
\end{align}

As a result, it can be concluded that 
\begin{align}
	\notag 
	&\ltwo{\hatxi_{t+1}} \\
	\notag &\le
	\Big\|\sum_{k = 1}^{t} \omega_k\psi_k \circ F_{t+1}^{\prime}(v_{t+1} + \overline{\zeta})\Big\|_2 \Big(\ltwo{\hatzeta_t} + O\Big(\frac{1}{\rho_F}\sqrt{\frac{t\log^2 n}{n}}\lt(\overgamma_{t+1} + \overalpha_t\rt)\Big)\Big)
	+ O\Big(\sqrt{\frac{t\log^2 n}{n}}\lt(\overgamma_{t+1} + \rho_F\overalpha_t\rt)\Big) \\
	\notag 
	&\leq \Big\| \sum_{k = 1}^{t} \omega_k\psi_k \circ F_{t+1}^{\prime}(v_{t+1} + \overline{\zeta})\Big\|_2 \ltwo{\hatzeta_t}
	+ \rho_F \Big\|\sum_{k = 1}^{t} \omega_k\psi_k\Big\|_2 O\Big(\frac{1}{\rho_F}\sqrt{\frac{t\log^2 n}{n}}\lt(\overgamma_{t} + \overalpha_t\rt)\Big) + O\Big(\sqrt{\frac{t\log^2 n}{n}}\lt(\overgamma_{t+1} + \rho_F\overalpha_t\rt)\Big)\\
	&\leq 
	\Big\| \sum_{k = 1}^{t} \omega_k \psi_k \circ F_{t+1}^{\prime}(v_{t+1} + \overline{\zeta})\Big\|_2 \ltwo{\hatzeta_t}
	+  O\Bigg(\sqrt{\frac{t\log^2 n}{n}}\lt(\overgamma_{t+1} + \rho_F\overalpha_t\rt)\Bigg),
	\label{eqn:xi-allegro}
\end{align}
where the last line follows since we condition on event~\eqref{eqn:rm-event}.

\item \textbf{Control of $\| \sum_{k = 1}^{t} \omega_k\psi_k \circ F_{t+1}^{\prime}(v_{t+1} + \overline{\zeta})\|_2$.}
Next, we shall focus our attention on quantity $\| \sum_{k = 1}^{t} \omega_k\psi_k \circ F_{t+1}^{\prime}(v_{t+1} + \overline{\zeta})\|_2$
with the size of $\overline{\zeta}$ bounded by \eqref{eqn:zeta-bar}.
To begin with, the property~\eqref{eqn:buble} summarized in Lemma~\ref{lem:cover} ensures 
\begin{align}
\notag &\lt\| \sum_{k = 1}^{t} \omega_k\psi_k \circ \lt[F_{t+1}^{\prime}(v_{t+1} + \overline{\zeta}) - F_{t+1}^{\prime}(v_{t+1})\rt]\rt\|_2 \\
&\le \rho_{1, F}\sqrt{\frac{t\log n}{n}}\|\overline{\zeta}\|_2 + \rho_F\lt(\sqrt{\frac{t\log^2 n}{n}} + \sqrt{\log n}\lt(\frac{\lt\|\overline{\zeta}\rt\|_2}{\|\alpha_t\|_2}\rt)^{\frac{1}{3}}\rt) \ll \frac{1}{\rho_G}. \label{eqn:bruckner7}
\end{align}
Here in the last inequality, we invoke the assumptions~\eqref{eqn:assmp-t-range} and \eqref{eqn:assmp-rho-b}.

Now in view of triangle's inequality, in order to bound the size of $\|\sum_{k = 1}^{t} \omega_k\psi_k \circ F_{t+1}^{\prime}(v_{t+1} + \overline{\zeta})\|$, it is sufficient to control the size of $\|\sum_{k = 1}^{t} \omega_k\psi_k \circ F_{t+1}^{\prime}(v_{t+1})\|_2.$
Towards this goal, we introduce the following lemma. 


\begin{lems} \label{lem:induction}
With probability at least $1 - O(n^{-10})$, for any $\omega \in \mathcal{S}^{t-1}$, it holds that 
\begin{subequations}
\begin{align}
\label{eqn:induction-con-a}
	\Bigg\| \sum_{k = 1}^{t} \omega_k\psi_k \circ F_{t+1}^{\prime}(v_{t+1})\Bigg\|_2^2 - \frac{1}{n}\mathbb{E}\big[\lt\|F_{t+1}^{\prime}(\widetilde{v}_{t+1})\rt\|_2^2 \mymid \|\alpha_t\|_2\big] 
	\lesssim \sqrt{\frac{t\log^2 n}{n}}\rho_F^2,
\end{align}
and
\begin{align}
\label{eqn:induction-con-b}
	\frac{1}{n}\lt\|F_{t+1}^{\prime}(v_{t+1})\rt\|_2^2 - 
	\frac{1}{n}\mathbb{E}\big[\lt\|F_{t+1}^{\prime}(\widetilde{v}_{t+1})\rt\|_2^2 \mymid \|\alpha_t\|_2\big]
	\lesssim \sqrt{\frac{t\log^2 n}{n}}\rho_F^2.
\end{align}
\end{subequations}
\end{lems}
\noindent The proof of this lemma is postponed to Section~\ref{sec:proof-lem:induction}.

Combining \eqref{eqn:induction-con-a} with the assumption~\eqref{eqn:assmp-t-range}, which suggests $\rho_F^2\sqrt{\frac{t\log^2 n}{n}} \ll 1/\rho_G^2$, we arrive at  
\begin{align*}
	\Bigg\| \sum_{k = 1}^{t} \omega_k\psi_k \circ F_{t+1}^{\prime}(v_{t+1})\Bigg\|_2^2 - \frac{1}{n}\mathbb{E}\big[\lt\|F_{t+1}^{\prime}(\widetilde{v}_{t+1})\rt\|_2^2 \mymid \|\alpha_t\|_2\big] 
	\ll \frac{1}{\rho_G^2}. 
\end{align*}
Taking this collectively with \eqref{eqn:bruckner7}, it is ensured that 
\begin{align}
\notag \Big\|\sum_{k = 1}^{t} \omega_k\psi_k \circ F_{t+1}^{\prime}(v_{t+1}+\overline{\zeta})\Big\|_2
&\leq
o\Big(\frac{1}{\rho_G}\Big) + \sqrt{\frac{1}{n}\mathbb{E}\big[\lt\|F_{t+1}^{\prime}(\widetilde{v}_{t+1})\rt\|_2^2 \mymid \|\alpha_t\|_2\big] 
	+ o\Big(\frac{1}{\rho_G^2}\Big)} \\
&\le  
 \sqrt{\frac{1}{n}\mathbb{E}\big[\lt\|F_{t+1}^{\prime}(\widetilde{v}_{t+1})\rt\|_2^2 \mymid \|\alpha_t\|_2\big]} + o\Big(\frac{1}{\rho_G}\Big) . \label{eqn:world}
\end{align}
Similarly, one can also conclude that 
\begin{align}
\label{eqn:cup}
	\Big\| \sum_{k = 1}^{t} \mu_k\phi_k \circ G_{t}^{\prime}(u_{t} + \overline{\xi})\Big\|_2 
	\le  
	\sqrt{\frac{1}{n}\mathbb{E}\big[\lt\|G_{t}^{\prime}(\widetilde{u}_{t})\rt\|_2^2 \mymid \|\gamma_t\|_2\big]}+o\Big(\frac{1}{\rho_F}\Big) .
\end{align}
\end{itemize}

\paragraph*{In summary.}
With these properties in place, we are ready to bound $\ltwo{\hatxi_{t+1}}.$ Recall expression~\eqref{eqn:xi-allegro} to obtain 
\begin{align}
\label{eqn:xi-t-recursion}
	\ltwo{\hatxi_{t+1}} 
	%
	&\le \Bigg(\sqrt{\frac{1}{n}\mathbb{E}\big[\lt\|F_{t+1}^{\prime}(\widetilde{v}_{t+1})\rt\|_2^2 \mymid \|\alpha_t\|_2\big]} + o\Big(\frac{1}{\rho_G}\Big)\Bigg) \|\hatzeta_{t}\|_2 + O\Big(\sqrt{\frac{t\log^2 n}{n}}\lt(\overgamma_{t+1} + \rho_F\overalpha_t\rt)\Big),
\end{align}
where we make use of the relation~\eqref{eqn:world}. Similar to \eqref{eqn:xi-t-recursion}, one can establish the recursive relation between $\|\hatzeta_{t}\|_2$ and $\|\hatxi_{t}\|_2$ as in 
\begin{align*}
	\|\hatzeta_{t}\|_2 \le \Bigg(\sqrt{\frac{1}{n}\mathbb{E}\big[\lt\|G_{t}^{\prime}(\widetilde{u}_{t})\rt\|_2^2 \mymid \|\gamma_t\|_2\big]} + o\Big(\frac{1}{\rho_F}\Big)\Bigg) \|\hatxi_{t}\|_2 + 
	O\Big(\sqrt{\frac{t\log^2 n}{n}}\lt(\overalpha_{t} + \rho_G\overgamma_t\rt)\Big).
\end{align*}
Consequently, the above two relations combined together yields 
\begin{align}
\label{eqn:christmas}
		\notag \ltwo{\hatxi_{t+1}} 
		&\le 
		\Bigg(\sqrt{\frac{1}{n}\mathbb{E}\big[\lt\|F_{t+1}^{\prime}(\widetilde{v}_{t+1})\rt\|_2^2 \mymid \|\alpha_t\|_2\big]} + o\Big(\frac{1}{\rho_G}\Big)\Bigg)
		\Bigg(\sqrt{\frac{1}{n}\mathbb{E}\big[\lt\|G_{t}^{\prime}(\widetilde{u}_{t})\rt\|_2^2 \mymid \|\gamma_t\|_2\big]} + o\Big(\frac{1}{\rho_F}\Big)\Bigg)
		\cdot\|\hatxi_{t}\|_2 \\
		&\qquad\qquad\qquad\qquad+ O\Big(\sqrt{\frac{t\log^2 n}{n}}\lt(\rho_F\rho_G\overgamma_{t+1} + \rho_F\overalpha_{t}\rt)\Big) \notag\\
		\notag &\le \Bigg(\sqrt{\frac{1}{n^2}\mathbb{E}\big[\lt\|G_{t}^{\prime}(\widetilde{u}_{t})\rt\|_2^2 \mymid \|\gamma_t\|_2\big]\mathbb{E}\big[\lt\|F_{t+1}^{\prime}(\widetilde{v}_{t+1})\rt\|_2^2 \mymid \|\alpha_t\|_2\big]} + o(1)\Bigg)\cdot\|\hatxi_{t}\|_2 + O\Big(\sqrt{\frac{t\log^2 n}{n}}\lt(\rho_F\rho_G\overgamma_{t+1} + \rho_F\overalpha_{t}\rt)\Big)\\
		&\le \Big(1 - c\Big)\|\hatxi_{t}\|_2 + O\Big(\sqrt{\frac{t\log^2 n}{n}}\lt(\rho_F\rho_G\overgamma_{t+1} + \rho_F\overalpha_{t}\rt)\Big),
\end{align}
where the last inequality follows from the assumption~\eqref{eqn:assmp-prime}.

Putting inequality~\eqref{eqn:christmas} and the inductive assumption~\eqref{eq:xihat} that $\|\hatxi_{t}\|_2 \lesssim \sqrt{\frac{t\log^2 n}{n}}\lt(\rho_F\rho_G\overgamma_{t} + \rho_F\overalpha_{t-1}\rt)$, 
we complete the proof of \eqref{eq:xihat} at $t+1.$
Additionally, the control of $\hatzeta_{t+1}$ can be derived in a similar way, which we omit here for brevity.

\subsubsection{Induction step for quantity $|\hatalpha_{t}^{t}|$ and $|\hatgamma_{t+1}^{t+1}|$}
\label{sec:lemma-magic-alpha}

Let us recall our definition of $\hatalpha_t^t$ in expression~\eqref{eq:xi-coeff} where 
\begin{align*} 
	\hatalpha_{t}^t \defn 
	\langle F_{t+1}^{\prime}(\hatbeta_{t+1}) - F_{t+1}^{\prime}(\beta_{t+1})\rangle + \frac{1}{\|\alpha_t\|_2^2}\lt\langle \sum_{k=1}^t \alpha_t^k\psi_k, F_{t+1}(\beta_{t+1}) - F_{t+1}(\hatbeta_{t+1})\rt\rangle.
\end{align*}
To establish~\eqref{eq:alphahat} at the $t+1$-th iteration, it is sufficient to bound the two terms on right of the above expression respectively.  

To begin with, according to inequality~\eqref{eqn:Fprime-a} in Lemma~\ref{lem:cover}, we are ensured that 
\begin{align*}
	\langle F_{t+1}^{\prime}(\hatbeta_{t+1}) - F_{t+1}^{\prime}(\beta_{t+1})\rangle &\le \frac{1}{\sqrt{n}}\rho_{1, F}\lt\|\beta_{t+1} - \hatbeta_{t+1}\rt\|_2 
	+ \rho_F\lt(\frac{t\log n}{n} + \lt(\frac{\|\beta_{t+1} - \hatbeta_{t+1}\|_2}{\|\alpha_t\|_2}\rt)^{\frac{2}{3}}\rt)\\
	&\lesssim \sqrt{\frac{t\log^2 n}{n^2}}\rho_{1, F}\lt(\rho_F\rho_G\overalpha_{t} + \rho_G\overgamma_{t}\rt) + \rho_F\lt(\frac{\mu_{t+1}^2t\log^2 n}{n}\rt)^{\frac{1}{3}}.
\end{align*} 
where we invoke the relation~\eqref{eqn:beta-betahat-2} and recall that 
$\mu_{t+1} := \frac{\rho_F\rho_G\overalpha_{t} + \rho_G\overgamma_{t}}{\|\alpha_{t}\|_2}$ as in expression~\eqref{eqn:mu-nv-main}. 
In addition, conditioning on event~\eqref{eqn:rm-event}, it is easily seen that 
\begin{align*}
\frac{1}{\|\alpha_t\|_2^2}\lt\langle \sum_k \alpha_t^k\psi_k, F_{t+1}(\beta_{t+1}) - F_{t+1}(\hatbeta_{t+1})\rt\rangle 
&\le \frac{1}{\|\alpha_t\|_2^2}\Big\|\sum_k \alpha_t^k\psi_k\Big\|_2 \cdot \rho_F\lt\|\beta_{t+1} - \hatbeta_{t+1}\rt\|_2 \\
&\lesssim \rho_F\frac{\lt\|\beta_{t+1} - \hatbeta_{t+1}\rt\|_2}{\|\alpha_t\|_2} \lesssim \rho_F\sqrt{\frac{\mu_{t+1}^2t\log^2 n}{n}}.
\end{align*}
Combining these two bounds establishes~\eqref{eq:alphahat} for $\hatalpha_t^t$. 
Moreover, the upper bound of $\hatgamma_{t+1}^{t+1}$ as in \eqref{eq:gammahat} can be derived in a similar manner, thus is omitted here.

\subsubsection{Induction step for quantities $|\hatalpha_{t-1}^k|$ }

For $k < t$, recall the definitions in \eqref{eq:xi-coeff} and \eqref{eq:zeta-coeff} that 
\begin{align*} 
	\hatalpha_{t}^k &\defn 
	\hatgamma_{t}^{k+1} \lt\langle F_{t+1}^{\prime}(\hatbeta_{t+1}) \circ F_{k+1}^{\prime}(v_{k+1})\rt\rangle \\
	\hatgamma_{t+1}^k &\defn 
	\hatalpha_{t}^k \lt\langle G_{t+1}^{\prime}(\hats_{t+1}) \circ G_{k}^{\prime}(u_{k})\rt\rangle
\end{align*}
When $k = t-1$, it is easily seen that $|\hatalpha_{t}^{t-1}| \defn |\hatgamma_{t}^{t} \langle F_{t+1}^{\prime}(\hatbeta_{t+1}) \circ F_{t}^{\prime}(v_{t})\rangle| \leq \rho^2_{F}  |\hatgamma_{t}^{t}|.$
So we only need to prove for $k \leq t-2$ where 
\begin{align}
\label{eqn:alpha-recursion} 
	\hatalpha_{t}^k &\defn 
	\hatgamma_{t}^{k+1} \lt\langle F_{t+1}^{\prime}(\hatbeta_{t+1}) \circ F_{k+1}^{\prime}(v_{k+1})\rt\rangle 
	=
	\hatalpha_{t-1}^{k+1} \lt \langle G_{t}^{\prime}(\hats_{t}) \circ G_{k+1}^{\prime}(u_{k+1})\rt\rangle
	\lt\langle F_{t+1}^{\prime}(\hatbeta_{t+1}) \circ F_{k+1}^{\prime}(v_{k+1})\rt\rangle.
\end{align}

Before proceeding, let us make note of the simple relation that 
\begin{align}
\label{eqn:golden}
	\lt\langle F_{t+1}^{\prime}(\hatbeta_{t+1}) \circ F_{k+1}^{\prime}(v_{k+1})\rt\rangle 
	&\le
	\lt\langle F_{t+1}^{\prime}(v_{t+1}) \circ F_{k+1}^{\prime}(v_{k+1})\rt\rangle 
	+
	\frac{1}{n}\rho_F\lt\| F_{t+1}^{\prime}(\hatbeta_{t+1}) - F_{t+1}^{\prime}(v_{t+1})\rt\|_1. 
\end{align}
We shall bound the two parts on the right respectively as below. 
\begin{itemize}
\item 
In view of Lemma~\ref{lem:cover} inequality~\eqref{eqn:Fprime-o}, it satisfies that 
\begin{align}
\label{eqn:messi}
\notag \frac{1}{n}\rho_F\lt\| F_{t+1}^{\prime}(\hatbeta_{t+1}) - F_{t+1}^{\prime}(v_{t+1})\rt\|_1 
	&\le 
	\frac{1}{n}\rho_F\Bigg[ \sqrt{n}\rho_{1, F}\|\hatbeta_{t+1} - v_{t+1}\|_2 
	+ 
	\rho_F\Bigg(t\log n + n\Bigg(\frac{\|\hatbeta_{t+1} - v_{t+1}\|_2}{\|\alpha_t\|_2}\Bigg)^{\frac{2}{3}}\Bigg)\Bigg] \\
	&
	\lesssim 
	\frac{1}{n}\rho_F\Bigg[ \sqrt{n}\rho_{1, F}\overgamma_{t}(\hatgamma_{t}^{t} + \rho_G^2\hatalpha_{t-1}^{t-1})
	+ 
	\rho_F n\Bigg(\frac{\overgamma_{t}(\hatgamma_{t}^{t} + \rho_G^2\hatalpha_{t-1}^{t-1})}{\|\alpha_t\|_2}\Bigg)^{\frac{2}{3}}\Bigg] + \frac{t\log n\rho^2_F}{n}
\end{align}
where the last step follows from the inequality~\eqref{eqn:bruckner9} where we proved $\|\hatbeta_{t+1} - v_{t+1}\|_2  \le  \overgamma_{t}(\hatgamma_{t}^{t} + \rho_G^2\hatalpha_{t-1}^{t-1})$.
Next, by virtue of assumptions \eqref{eqn:assmp-rho-d} and \eqref{eqn:assmp-t-range}, 
the right hand side of inequality~\eqref{eqn:messi} further satisfies 
\begin{align}
\label{eqn:ball}
\frac{1}{n}\rho_F\lt\| F_{t+1}^{\prime}(\hatbeta_{t+1}) - F_{t+1}^{\prime}(v_{t+1})\rt\|_1 
\ll
\frac{1}{\rho_G^2}.
\end{align}

\item It is then sufficient to consider the quantity $\langle F_{t+1}^{\prime}(v_{t+1}) \circ F_{k+1}^{\prime}(v_{k+1})\rangle.$
Towards this, it is easily seen that 
\begin{align*}
	\lt\langle F_{t+1}^{\prime}(v_{t+1}) \circ F_{k+1}^{\prime}(v_{k+1})\rt\rangle \le \frac{1}{n}\lt\| F_{t+1}^{\prime}(v_{t+1})\rt\|_2 \lt\| F_{k+1}^{\prime}(v_{k+1})\rt\|_2.
\end{align*}
To bound the right hand side, notice that according to Lemma~\ref{lem:induction}, for every $k$, it obeys 
\begin{align*}
	\frac{1}{n}\lt\| F_{k+1}^{\prime}(v_{k+1})\rt\|_2^2 - \frac{1}{n}\mathbb{E}\big[\lt\|F_{k+1}^{\prime}(\widetilde{v}_{k+1})\rt\|_2^2 \mymid \|\alpha_k\|_2\big] \lesssim \rho_F^2\sqrt{\frac{t\log^2 n}{n}} \ll \frac{1}{\rho_G^2}. 
\end{align*}
Therefore we arrive at 
\begin{align*}
	\lt\langle F_{t+1}^{\prime}(v_{t+1}) \circ F_{k+1}^{\prime}(v_{k+1})\rt\rangle 
	&\le \sqrt{\frac{1}{n^2} \lt\|F_{t+1}^{\prime}(v_{t+1})\rt\|^2_2 \lt\| F_{k+1}^{\prime}(v_{k+1})\rt\|^2_2} \\
	&\le \sqrt{ \frac{1}{n}\mathbb{E}\big[\lt\|F_{t+1}^{\prime}(\widetilde{v}_{t+1})\rt\|_2^2 \mymid \|\alpha_t\|_2\big] + o(\frac{1}{\rho_G^2})}
	\sqrt{ \frac{1}{n}\mathbb{E}\big[\lt\|F_{k+1}^{\prime}(\widetilde{v}_{k+1})\rt\|_2^2 \mymid \|\alpha_k\|_2\big] + o(\frac{1}{\rho_G^2})}
	%
\end{align*}
\end{itemize}
Taking the above inequality collectively with displays \eqref{eqn:golden} and \eqref{eqn:ball} gives us 
\begin{align}
\label{eqn:F-product}
	\lt\langle F_{t+1}^{\prime}(\hatbeta_{t+1}) \circ F_{k+1}^{\prime}(v_{k+1})\rt\rangle
	\le 
	\sqrt{ \frac{1}{n^2}\mathbb{E}\big[\lt\|F_{t+1}^{\prime}(\widetilde{v}_{t+1})\rt\|_2^2 \mymid \|\alpha_t\|_2\big]\mathbb{E}\big[\lt\|F_{k+1}^{\prime}(\widetilde{v}_{k+1})\rt\|_2^2 \mymid \|\alpha_k\|_2\big]} + o(\frac{1}{\rho^2_G}).
\end{align}
Similarly, one can develop a symmetric bound for $G^{\prime}$ as
\begin{align}
\label{eqn:G-product}
	\lt\langle G_{t}^{\prime}(\hats_{t}) \circ G_{k+1}^{\prime}(u_{k+1})\rt\rangle
	\le 
	\sqrt{ \frac{1}{n}\mathbb{E}\big[\lt\|G_{t}^{\prime}(\widetilde{u}_{t})\rt\|_2^2 \mymid \|\gamma_t\|_2\big]\mathbb{E}\big[\lt\|G_{k+1}^{\prime}(\widetilde{u}_{k+1})\rt\|_2^2 \mymid \|\gamma_{k+1}\|_2\big]} + o(\frac{1}{\rho^2_F}).
\end{align}
Therefore, under assumption~\eqref{eqn:assmp-prime}, it holds that  
\begin{align*}
	\lt\langle F_{t+1}^{\prime}(\hatbeta_{t+1}) \circ F_{k+1}^{\prime}(v_{k+1})\rt\rangle\lt\langle G_{t}^{\prime}\lt(\hats_{t}\rt) \circ G_{k+1}^{\prime}\lt(u_{k+1}\rt)\rt\rangle < (1-c)^2.
\end{align*}

As a consequence, we can establish the inductive relationship that 
\begin{align*} 
	|\hatalpha_{t}^{k}| &\leq 
	(1-c)^2 |\hatalpha_{t-1}^{k+1}|.
	%
\end{align*}
Invoking the above inequality recursively validates the relation~\eqref{eq:alphahat-k} at $t+1$-th step.
In addition, one can prove for inequality \eqref{eq:gammahat-k} in a similar manner.


\section{Proof of Lemma~\ref{lem:Onsager}}
\label{sec:proof-lem:Onsager}

Let us present the proof of inequality~\eqref{eqn:Onsager-brahms} and inequality~\eqref{eqn:Onsager-beethoven} can be established in the same fashion. 
Throughout this proof, let us condition on the event where both expressions~\eqref{eqn:simple-rm-phi} and \eqref{eqn:simple-rm-psi} hold true (with $\delta$ chosen as $\max(n,p)^{-11}$)  
\begin{subequations}
\label{eqn:lemma-onsager-event}
\begin{align}
	\lt\|(\phi_1, \ldots, \phi_{t-1})^{\top}(\phi_1, \ldots, \phi_{t-1}) - I_{t-1}\rt\|_{\text{op}} \lesssim \sqrt{\frac{t\log n}{n}},\qquad\text{for every } 1 < t \leq n, \\
	\lt\|\frac{n}{p}(\psi_1, \ldots, \psi_{t-1})^{\top}(\psi_1, \ldots, \psi_{t-1}) - I_{t-1}\rt\|_{\text{op}} \lesssim \sqrt{\frac{t\log p}{p}},\qquad\text{for every } 1 < t \leq n, 
\end{align}
\end{subequations}
with probability at least $1 - O(n^{-11})$.

We are ready to control the norm on the left hand side of inequality~\eqref{eqn:Onsager-brahms}. First, recalling that $\{a_k\}$ forms an orthogonal basis, 
therefore there exists a unit vector
\begin{align}
\label{eqn:def-omega}
	\omega \defn \sum_{k=1}^t \omega_k a_k \in \mathcal{S}^{t-1},
\end{align}
--- depending on the randomness of $\{\psi_{k}\}$ --- such that
\begin{align}
&\lt\|\sum_{k = 1}^{t} a_k\Big[\lt\langle \psi_k, F_{t+1}(\hatbeta_{t+1})\rt\rangle - \langle F_{t+1}^{\prime}(\hatbeta_{t+1}) \rangle \alpha_{t}^k - \sum_{j = k+1}^{t} \hatgamma_{t}^j \lt\langle F_{t+1}^{\prime}(\hatbeta_{t+1}) \circ F_{j}^{\prime}\lt(v_{j}\rt)\rt\rangle\alpha_{j-1}^k\Big]\rt\|_2 \notag\\
&= \omega^\top \sum_{k = 1}^{t} a_k\Big[\lt\langle \psi_k, F_{t+1}(\hatbeta_{t+1})\rt\rangle - \langle F_{t+1}^{\prime}(\hatbeta_{t+1}) \rangle \alpha_{t}^k - \sum_{j = k+1}^{t} \hatgamma_{t}^j \lt\langle F_{t+1}^{\prime}(\hatbeta_{t+1}) \circ F_{j}^{\prime}\lt(v_{j}\rt)\rt\rangle\alpha_{j-1}^k\Big] \notag \\
&= \Bigg\langle \sum_{k = 1}^{t} \omega_k\psi_k, F_{t+1}(\hatbeta_{t+1})\Bigg\rangle - \langle F_{t+1}^{\prime}(\hatbeta_{t+1}) \rangle \cdot \omega^\top \alpha_t - \sum_{j = 1}^{t-1} \hatgamma_{t}^{j+1} \lt\langle F_{t+1}^{\prime}(\hatbeta_{t+1}) \circ F_{j+1}^{\prime}\lt(v_{j+1}\rt)\rt\rangle\sum_{k = 1}^{j} \omega_k\alpha_{j}^k \notag\\
&=: \sum_{i = 1}^p X_i^0 + \sum_{i = 1}^p X_i - \sum_{i = 1}^p Y_i - \sum_{i = 1}^p Z_i, \label{eq:Onsager}
\end{align}
where to simplify our presentation, we introduce the following short-hand notation  
\begin{subequations}
\label{eqn:def-24601}
\begin{align}
X_i^0 &\defn \lt[\sum_{k = 1}^{t} \omega_k\psi_{k} \circ F_{t+1}\Big(\sum_{k = 1}^{t} \hatgamma_{t}^k F_{k}(0)\Big)\rt]_i \label{eqn:def-xi-bruckner}\\
X_i &\defn \lt[\sum_{k = 1}^{t} \omega_k\psi_{k} \circ \Bigg(F_{t+1}(\hatbeta_{t+1}) - F_{t+1}\Big(\sum_{k = 1}^{t} \hatgamma_{t}^k F_{k}(0)\Big)\Bigg)\rt]_i \label{eqn:def-xi-real-bruckner}\\
Y_i &\defn \frac{1}{n}\lt[F_{t+1}^{\prime}(\hatbeta_{t+1})\rt]_i \cdot \omega^\top \alpha_t \label{eqn:def-yi-bruckner} \\
Z_i &\defn \frac{1}{n}\sum_{j = 1}^{t-1} \hatgamma_{t}^{j+1} \lt[ F_{t+1}^{\prime}(\hatbeta_{t+1}) \circ F_{j+1}^{\prime}\lt(v_{j+1}\rt)\rt]_i\sum_{k = 1}^{j} \omega_k\alpha_{j}^k.
\label{eqn:def-zi-bruckner}
\end{align}
\end{subequations}
Here, in view of definitions in \eqref{eqn:sec-order-terms}, we remind the readers that, in these expressions above, the parameters concerned are
\begin{align} \label{defi:theta}
	\thetahat \defn \Big\{\omega, \{\alpha_k\}_{k\leq t}, \hatgamma_t, \{\tau_k\}_{k \le t+1} \Big\}, 
	\quad \text{where }  \omega \in \mathcal{S}^{t-1}, \alpha_k \in \real^k, \hatgamma_t \in \real^t, \text{ and }\tau_k \in \real^s.
\end{align}
Let us point out a few properties for the above parameters:
\begin{itemize}
	\item $\tau_k$ corresponds to the parameter used for defining function $F_k$. By assumption, $\tau_k$ is of finite and low dimension and $\ltwo{\tau_k}\leq c$ for some universal constant;

	\item according to the assumption~\eqref{eqn:assmp-alphat-gammat}, 
	$\ltwo{\alpha_k} \leq \overalpha_t \leq \poly(n)$, for every $k\leq t$;

	\item in view of the assumption~\eqref{eqn:assmp-gamma-tail}, $\ltwo{\hatgamma_t}^2 = \sum_{k=1}^t (\hatgamma_t^{k})^2 \leq \poly(n)/(1 - c^2).$

\end{itemize}
The value of $\thetahat$ depends on the randomness in $\{\psi_k\}$; we collect all the possible values of $\thetahat$ to be space $\Theta$, namely,
\begin{align}
\label{eqn:Theta}
	\Theta \defn \Big\{(\omega, \{\alpha_k\}_{k\leq t}, \hatgamma_t, \{\tau_k\}_{k \le t+1}) ~\big|~ \omega\in \mathcal{S}^{t-1}, \ltwo{\tau_k}\leq 1, \ltwo{\alpha_k} \leq \poly(n),  \ltwo{\hatgamma_t}\leq \frac{\poly(n)}{1 - c^2} \Big\}. 
\end{align}

Next, let us control the right hand side of inequality~\eqref{eq:Onsager}, which shall be done by bounding each term in the summation separately.

\subsection{Controlling~\eqref{eqn:def-xi-bruckner}}

Regarding the first term $\sum_{i = 1}^p X_i^0$, since $(\omega,\hatgamma)$ has complicated statistical dependence on the randomness of $\{\psi_k\}$, we find it helpful to construct an $\epsilon$-covering set $\mathcal{N}_{\epsilon}$ of the space $\mathcal{S}^{t-1}(1) \times \mathcal{S}^{t-1}(\poly(n))$ for $\epsilon = \frac{1}{\poly(n)}$, with its cardinality satisfying
\begin{align*}
	|\mathcal{N}_{\epsilon}| \lesssim \Big(\frac{\poly(n)}{\epsilon}\Big)^t = \poly(n)^t.
\end{align*}
Before diving into the main proof, we make note of the following property  
\begin{align} \label{eq:F0norm}
	\notag \Bigg\|F_{t+1}\Big(\sum_{k = 1}^{t} \hatgamma_{t}^k F_{k}(0)\Big)\Bigg\|_2 
	&\le \lt\|F_{t+1}(0)\rt\|_2 + \rho_F \Big\|\sum_{k = 1}^{t} \hatgamma_{t}^k F_{k}(0) \Big\|_2 \\
	&\le \lt\|F_{t+1}(0)\rt\|_2 + \rho_F\sum_{k = 1}^{t} |\hatgamma_{t}^k| \cdot \lt\|F_{k}(0)\rt\|_2 \lesssim \overgamma_{t+1},
\end{align}
where in the last inequality, we invoke assumption~\eqref{eqn:assmp-bach} which implies $\rho_F\sum_{k = 1}^{t} |\hatgamma_{t}^k| \ll 1.$

Given every $(\omega, \hatgamma_t) \in \mathcal{S}^{t-1}(1)\times \mathcal{S}^{t-1}(\poly(n))$, there exists $(\widetilde{\omega}, \widetilde{\gamma}_t) \in \mathcal{N}_{\epsilon}$ satisfying $\|\widetilde{\omega} - \omega\|_2 + \|\widetilde{\gamma}_t - \hatgamma_t\|_2 \le \epsilon  = \frac{1}{\poly(n)}$. This fact together with the Lipschitz property of function $F_t$ gives 
\begin{align}
\notag &\lt | \Big\langle \sum_{k = 1}^{t} \omega_k\psi_k, F_{t+1}\Big(\sum_{k = 1}^{t} \hatgamma_{t}^k F_{k}(0)\Big)\Big\rangle - \Big\langle \sum_{k = 1}^{t} \widetilde{\omega}_k\psi_k, F_{t+1}\Big(\sum_{k = 1}^{t} \widetilde{\gamma}_{t}^k F_{k}(0)\Big)\Big\rangle \rt|\\
\notag &\leq \Big\|\sum_{k = 1}^{t} \omega_k\psi_k\Big\|_2 \Big\|F_{t+1}\Big(\sum_{k = 1}^{t} \hatgamma_{t}^k F_{k}(0)\Big) - F_{t+1}\Big(\sum_{k = 1}^{t} \widetilde{\gamma}_{t}^k F_{k}(0)\Big)\Big\|_2
+ \Big\|F_{t+1}\Big(\sum_{k = 1}^{t} \hatgamma_{t}^k F_{k}(0)\Big)\Big\|_2 \Big\|\sum_{k = 1}^{t} \omega_k\psi_k - \sum_{k = 1}^{t} \widetilde{\omega}_k\psi_k \Big\|_2 \\
\notag 
&\stackrel{(\mathrm{i})}{\lesssim} \sqrt{\frac{p}{n}}\Big(1+\sqrt{\frac{t\log p}{p}}\Big) 
\cdot \Big(\rho_F \ltwo{\hatgamma_t - \widetilde{\gamma}_t} \overgamma_t +  
\ltwo{\omega - \widetilde{\omega}} \overgamma_{t+1} \Big)\\
&\lesssim \sqrt{\frac{p}{n}}\Big(1+\sqrt{\frac{t\log p}{p}}\Big) 
\cdot \epsilon \cdot \rho_F\overgamma_{t+1} 
\lesssim \frac{1}{\poly(n)}\rho_F\overgamma_{t+1}, \label{eq:cover-X0}
\end{align}
with probability at least $1 - O(n^{-11}).$
Here, $\mathrm{(i)}$ follows from the fact that we condition on the event \eqref{eqn:lemma-onsager-event}. 
Putting things together, we arrive at  
\begin{align*}
\notag \sum_{i = 1}^p X_i^0 &= \lt\langle \sum_{k = 1}^{t} \omega_k\psi_k, F_{t+1}\Big(\sum_{k = 1}^{t} \hatgamma_{t}^k F_{k}(0)\Big)\rt\rangle \\
\notag &\le \sup_{(\widetilde{\omega}, \widetilde{\gamma}) \in \mathcal{N}_{\epsilon}}\lt\langle \sum_{k = 1}^{t} \widetilde{\omega}_k\psi_k, F_{t+1}\Big(\sum_{k = 1}^{t} \widetilde{\gamma}_{t}^k F_{k}(0)\Big)\rt\rangle + \frac{1}{\poly(n)}\rho_F \overgamma_{t+1}.
\end{align*}
To further control the right hand side above, let us make note of the following two observations that 
(i) $\lt\langle \sum_{k = 1}^{t} \widetilde{\omega}_k\psi_k, F_{t+1}\Big(\sum_{k = 1}^{t} \widetilde{\gamma}_{t}^k F_{k}(0)\Big)\rt\rangle$
is stochastically dominated by $\mathcal{N}(0,\frac{\overgamma_{t+1}^2}{n})$ and (ii) the standard concentration result (see, e.g.~\cite[Exercise 2.12]{wainwright2019high}).
\begin{align}
\label{eqn:standard-gaussian}
	\mprob\left(\sup_{i \in [k]} X_i - \sqrt{2 \sigma^2 \log k} \geq t\right) 
	\leq 2e^{-\frac{t^2}{2\sigma^2}},
\end{align}
for $X_i \stackrel{\textrm{i.i.d}}{\sim} \mathcal{N}(0,\sigma^2).$
Consequently, we conclude that 
\begin{align}
\label{eqn:du-pre}
\sum_{i = 1}^n X_i^0 
	&\lesssim \sqrt{\frac{t\log (n)}{n}}\overgamma_{t+1} + \frac{1}{\poly(n)}\rho_F\overgamma_{t+1},
\end{align}
with probability at least $1 - O(n^{-11}).$

\subsection{Controlling the remaining terms}

For notational simplicity, if we concatenate $\{\psi_{k}\}_{k=1}^{t}$ into matrix $\Psi$, namely, 
\begin{align}
	\Psi \defn \big[\psi_{1}, \ldots, \psi_{t}\big] \in \real^{p \times t}, 
	\qquad
	\text{where } \psi_k \stackrel{\mathrm{i.i.d.}}{\sim} 
	\mathcal{N}(0,\frac{1}{n} I_p),
\end{align}
it suffices to control the summation of the remaining three terms as 
\begin{align}
\label{eqn:H-thetahat}
	H(\Psi; \thetahat) \defn \sum_{i = 1}^p X_i - \sum_{i = 1}^p Y_i - \sum_{i = 1}^p Z_i. 
\end{align} 
For any fixed parameter $\theta\in \Theta$ and fixed $i \in [p]$, it is easily seen from definitions~\eqref{eqn:def-24601}, random vector $(X_i, Y_i, Z_i)$ only depends on the $i$-th row of $\Psi$ matrix, namely, $[\psi_{k, i}]_{k = 1}^t$, which implies that $\{(X_i, Y_i, Z_i)\}_{i = 1}^p$ are independent for different $i$. 
In addition, in view of Stein's lemma of Gaussian random vectors --- which ensures $\mathbb{E}_{X \sim \mathcal{N}(0, 1)}\big[Xf(X)\big] = \mathbb{E}_{X \sim \mathcal{N}(0, 1)}\big[f^{\prime}(X)\big]$ --- one can verify  
\begin{align}
\label{eqn:mean-zero}
	\Exs \big[H(\Psi; \theta)\big] = 0. 
\end{align}
However, the above property holds only when $\theta$ is held as a fixed vector, independent of $\{\psi_k\}.$
When a random $\thetahat$ is concerned, due to the statistical dependence between $\thetahat$ and $(X_i, Y_i, Z_i)$'s, the mean zero property does not hold anymore. 

On the high level, to control $H(\Psi; \thetahat)$, the idea is to invoke Lemma~\ref{lem:concentration} to bound $H(\Psi; \theta)$ for each fixed $\theta \in \Theta$, which is achieved via Step 1-3 below. 
In Step 4, we develop a uniform control of $H(\Psi; \theta)$ over the space $\Theta$ in order to deal with the statistical dependence involved in $\thetahat.$

In order to apply Lemma~\ref{lem:concentration}, it boils down to computing the variance of $H(\Psi; \theta)$ and validating the property~\eqref{eq:Z_condition}, as shall be done as follows. 

\paragraph{Step 1: variances control.}
We claim that for every fixed $\theta \in \Theta$ (independent of $\{\psi_k\}$), the variance terms satisfy the following relations respectively 
\begin{subequations}
	\begin{align}
		\sum_{i = 1}^p \mathsf{Var}(X_i) &\lesssim \frac{\log n}{n}(\overgamma_{t+1}^2 + \rho_F^2\overalpha_{t}^2) \label{eqn:var-x}\\
		\sum_{i = 1}^p \mathsf{Var}(Y_i) &\lesssim \frac{\rho_F^2}{n^2}\lt\|\alpha_t\rt\|_2^2 \label{eqn:var-y}\\
		\sum_{i = 1}^p \mathsf{Var}(Z_i) &\lesssim \frac{\rho_F^2}{n^2}\lt\|\alpha_t\rt\|_2^2. \label{eqn:var-z}
	\end{align}
\end{subequations}
Given the above relations, the variance of $H(\Psi; \theta)$ satisfies 
\begin{align*}
	\var(H(\Psi; \theta)) \lesssim \frac{\log n}{n}(\overgamma_{t+1}^2 + \rho_F^2\overalpha_{t}^2).
\end{align*}
Let us establish these three claims respectively.  
We remark that throughout this step, $\theta$ should always be viewed as a fixed constant that does not dependent on any randomness of the problem. With a slight abuse of notation, we still write parameters such as $\alpha_t$, $\hatbeta_t$, $\hatgamma_t$, but here they should be understood as fixed parameters. 

\begin{itemize}
\item 
First, by noticing $|F_{t+1}^{\prime}(\hatbeta_{t+1})| \leq \rho_{F}$ and $|\omega^{\top}\alpha_t| \leq \ltwo{\alpha_t}$, we obtain  
\begin{align*}
	\sum_{i = 1}^p \mathsf{Var}(Y_i) \le \frac{1}{n^2}\sum_{i=1}^p \mathbb{E} [F_{t+1}^{\prime}(\hatbeta_{t+1})]_i^2 \lt\|\alpha_t\rt\|_2^2 \lesssim \frac{\rho_F^2}{n}\lt\|\alpha_t\rt\|_2^2,
\end{align*}
which establishes relation~\eqref{eqn:var-y}. 

\item In terms of the relation~\eqref{eqn:var-z}, for each fixed $\theta \in \Theta$, we remind the readers that $[F_{t+1}^{\prime}(\hatbeta_{t+1}) \circ F_{j+1}^{\prime}\lt(v_{j+1}\rt)]_i$ are independent of each other. This fact leads to 
\begin{align*}
\sqrt{\sum_{i = 1}^p \mathsf{Var}(Z_i)} 
&= 
\frac{1}{n}\sqrt{
\sum_{i = 1}^p  \mathsf{Var} \Bigg(\sum_{j = 1}^{t-1} \Big(\sum_{k = 1}^{j} \omega_k\alpha_{j}^k\Big)\cdot\hatgamma_{t}^{j+1} \lt[ F_{t+1}^{\prime}(\hatbeta_{t+1}) \circ F_{j+1}^{\prime}\lt(v_{j+1}\rt)\rt]_i\Bigg)} \\
&= 
\frac{1}{n}\sqrt{
\mathsf{Var} \Bigg(\sum_{j = 1}^{t-1} \Big(\sum_{k = 1}^{j} \omega_k\alpha_{j}^k\Big)\cdot\hatgamma_{t}^{j+1} \sum_{i = 1}^p \lt[ F_{t+1}^{\prime}(\hatbeta_{t+1}) \circ F_{j+1}^{\prime}\lt(v_{j+1}\rt)\rt]_i\Bigg)} \\
&= 
\frac{1}{n}\sqrt{
\Exs\Bigg[\sum_{j = 1}^{t-1} \Big(\sum_{k = 1}^{j} \omega_k\alpha_{j}^k\Big)\cdot\hatgamma_{t}^{j+1} \sum_{i = 1}^p \lt[F_{t+1}^{\prime}(\hatbeta_{t+1}) \circ F_{j+1}^{\prime}\lt(v_{j+1}\rt)\rt]_i
- \Exs \lt[ F_{t+1}^{\prime}(\hatbeta_{t+1}) \circ F_{j+1}^{\prime}\lt(v_{j+1}\rt)\rt]_i \Bigg]^2} \\
&\leq 
\frac{1}{n}\sum_{j = 1}^{t-1} \Big|\sum_{k = 1}^{j} \omega_k\alpha_{j}^k\cdot\hatgamma_{t}^{j+1}\Big| \sqrt{
\Exs\Bigg[ \sum_{i = 1}^p \lt[F_{t+1}^{\prime}(\hatbeta_{t+1}) \circ F_{j+1}^{\prime}\lt(v_{j+1}\rt)\rt]_i
- \Exs \lt[ F_{t+1}^{\prime}(\hatbeta_{t+1}) \circ F_{j+1}^{\prime}\lt(v_{j+1}\rt)\rt]_i \Bigg]^2}, 
\end{align*}
where the last used the basic property that $\sqrt{\Exs [\sum_i X_i]^2} \leq \sum_i\sqrt{\Exs [X_i^2]}.$
To avoid confusion, we remind the readers that the parameters here are treated as fixed and independent of the randomness in the problem. 
Now, in view of the basic inequality $|\sum_{k=1}^j \omega_k \alpha_j^k| \leq \ltwo{\alpha_j}$, we can further bound 
\begin{align}
\notag
\sqrt{\sum_{i = 1}^p \mathsf{Var}(Z_i)} 
&\leq \frac{1}{n}\sum_{j = 1}^{t-1} |\hatgamma_{t}^{j+1}|\cdot\|\alpha_j\|_2 \sqrt{\mathbb{E}\lt\|F_{t+1}^{\prime}(\hatbeta_{t+1}) \circ F_{j+1}^{\prime}\lt(v_{j+1}\rt)\rt\|_2^2} \\
&\le \frac{1}{n}\sum_{j = 1}^{t-1} |\hatgamma_{t}^{j+1}| \|\alpha_j\|_2 \rho_F^2\lesssim \frac{\rho_F}{n}\lt\|\alpha_t\rt\|_2,
\end{align}
where the last inequality invokes the assumption~\eqref{eqn:assmp-bach}.
This completes the proof of inequality \eqref{eqn:var-z}. 

\item When it comes to the inequality~\eqref{eqn:var-x}, basic inequality ensures that 
\begin{align*}
\sum_{i = 1}^p \mathsf{Var}(X_i) \lesssim \sum_{i = 1}^p \mathsf{Var}(X_i^0 + X_i)  + \sum_{i = 1}^p \mathsf{Var}(X_i^0).
\end{align*}
We shall compute these two terms on the right respectively. 
Note that for fixed $\omega \in \mathcal{S}^{t-1}$ independent of $\{\psi_k\}$, $\sum_k \omega_k\psi_k \sim \mathcal{N}(0,\frac{1}{n} I_p)$, therefore 
\begin{align*}
\sum_{i = 1}^p \mathsf{Var}(X_i^0) 
&\le \mathbb{E}\lt\|\sum_{k = 1}^{t} \omega_k\psi_k \circ F_{t+1}\Big(\sum_{k = 1}^{t} \hatgamma_{t}^k F_{k}(0)\Big)\rt\|_2^2 
= \frac{1}{n}\lt\|F_{t+1}\Big(\sum_{k = 1}^{t} \hatgamma_{t}^k F_{k}(0)\Big)\rt\|_2^2 \lesssim \frac{1}{n}\overgamma_{t+1}^2,
\end{align*}
where the last inequality makes use of~\eqref{eq:F0norm}.

\begin{lems} \label{lem:var-Xi}
	Under the assumptions of Lemma~\ref{lem:Onsager}, it satisfies 
	\begin{align}
	\sum_{i = 1}^p \mathsf{Var}(X_i^0 + X_i)
	\lesssim
	\frac{\log n}{n} (\overgamma_{t+1}^2 + \rho_F^2 \overalpha_t^2).
	\end{align}
\end{lems}
\noindent The proof of this result is postponed to Section~\ref{sec:Pf-lem-var-Xi}.

Putting these pieces together proves the claimed result in inequality~\eqref{eqn:var-x}. 

\end{itemize}

\paragraph{Step 2: sizes control.}
To apply Lemma~\ref{lem:concentration}, one needs to check condition \eqref{eq:Z_condition}.
Since the zero mean condition holds as in \eqref{eqn:mean-zero}, it only requires us to provide a high probability control on the size of $H(\Psi,\theta)$ for each $\theta \in \Theta.$
In the following, we bound the sizes of $X_{i}, Y_{i}$ and $Z_{i}$ respectively. 
\begin{itemize}
\item 
By definition of $X_{i}$ in expression~\eqref{eqn:def-xi-real-bruckner}, we claim that 
\begin{align}
\label{eqn:Xi-unif-bound}
\notag |X_i| &\le \Big|\sum_{k = 1}^{t} \omega_k\psi_{k, i}\Big| \cdot \lt|F_{t+1}(\hatbeta_{t+1}) - F_{t+1}\lt(\sum_{k = 1}^{t} \hatgamma_{t}^k F_{k}(0)\rt)\rt|_{i} \\
&\lesssim \sqrt{\frac{\log \frac{1}{\delta}}{n}} \cdot \rho_F\sqrt{\frac{\log\frac{1}{\delta}}{n}}\lt\|\alpha_t\rt\|_2 \lesssim \frac{\log \frac{1}{\delta}}{n}\rho_F\overalpha_t,
\end{align}
with probability at least $1-\delta$, for every $\delta \lesssim \frac{1}{t}$. 

In order to see this, first notice that for every $k\geq 1$ and $t\leq n$, $v_{k} \sim \mathcal{N}(0,\frac{\ltwo{\alpha_{k-1}}^2}{n})$ and $\sum_{k = 1}^{t} \omega_k\psi_{k, i} \sim \mathcal{N}(0,\frac{1}{n})$. Standard concentration inequality ensures that
\begin{align*}
\lt|\sum_{k = 1}^{t} \omega_k\psi_{k, i}\rt| \lesssim \sqrt{\frac{\log \frac{1}{\delta}}{n}} 
\qquad \text{ and }
 	\|v_{k}\|_{\infty} \lesssim \ltwo{\alpha_{k-1}} \sqrt{\frac{\log\frac{1}{\delta}}{n}},
 \end{align*} 
with probability at least $1 - \delta$.
In addition, taking the Lipschitz property of $F_{t+1}$ together with the above concentration results ensures 
\begin{align}
\notag \lt|F_{t+1}(\hatbeta_{t+1}) - F_{t+1}\lt(\sum_{k = 1}^{t} \hatgamma_{t}^k F_{k}(0)\rt)\rt|_{i} &\le \rho_F\lt|\hatbeta_{t+1} - \sum_{k = 1}^{t} \hatgamma_{t}^k F_{k}(0)\rt|_{i} \\
\notag &= \rho_F \lt|v_{t+1} + \sum_{k = 1}^{t} \hatgamma_{t}^k F_{k}(v_{k}) - \sum_{k = 1}^{t} \hatgamma_{t}^k F_{k}(0)\rt|_{i}\\
\notag &\lesssim \rho_F\sqrt{\frac{\log\frac{1}{\delta}}{n}}\lt[\|\alpha_t\|_2 + \rho_F\sum_{k = 1}^{t} \hatgamma_{t}^k \ltwo{\alpha_{k-1}}\rt] \\
&\lesssim \rho_F\sqrt{\frac{\log\frac{1}{\delta}}{n}}\overalpha_t, \label{eqn:pineapple}
\end{align}
with probability at least $1-\delta$.
Here the last inequality uses the definition~\eqref{eqn:defn-alpha-t-bar} and the condition~\eqref{eqn:assmp-bach}, 
Putting everything together completes the proof of \eqref{eqn:Xi-unif-bound}.

\item As for $Y_{i}$, in view of the definition \eqref{eqn:def-yi-bruckner}, some direct algebra leads to 
\begin{align*}
|Y_i| = \lt|\frac{1}{n}\lt[F_{t+1}^{\prime}(\hatbeta_{t+1})\rt]_i \cdot \omega^\top \alpha_t\rt| \leq 
 \frac{1}{n}\lt\|F_{t+1}^{\prime}(\hatbeta_{t+1})\rt\|_{\infty}\|\alpha_t\|_2 \le \frac{\rho_F}{n}\|\alpha_t\|_2.
\end{align*}

\item Regarding $Z_{i}$, it is easily seen that 
\begin{align*}
|Z_i| 
&=
\Big|\frac{1}{n}\sum_{j = 1}^{t-1} \hatgamma_{t}^{j+1} \lt[ F_{t+1}^{\prime}(\hatbeta_{t+1}) \circ F_{j+1}^{\prime}\lt(v_{j+1}\rt)\rt]_i\sum_{k = 1}^{j} \omega_k\alpha_{j}^k \Big|\\
&\le \frac{1}{n}\sum_{j = 1}^{t-1} |\hatgamma_{t}^{j+1}|\cdot \lt\|F_{t+1}^{\prime}(\hatbeta_{t+1}) \circ F_{j+1}^{\prime}\lt(v_{j+1}\rt)\rt\|_{\infty}\|\alpha_j\|_2 \le \frac{\rho_F^2}{n}\sum_{j = 1}^{t-1} |\hatgamma_{t}^{j+1}| \|\alpha_j\|_2 \lesssim \frac{\rho_F}{n}\lt\|\alpha_t\rt\|_2,
\end{align*}
where the last inequality follows from the condition $\rho_F\sum_{k = 1}^{t} |\hatgamma_{t}^k| \ll 1.$
\end{itemize}

Combining the pieces above, we are ensured that 
\begin{align}
\label{eqn:dante-Paradiso}
	\mprob\lt(|X_i|+|Y_i|+|Z_i| \lesssim \frac{\rho_F\overalpha_t}{n}\log \frac{1}{\delta} \rt)
	\geq 
	1 - \delta.
\end{align}

\paragraph{Step 3: Putting everything together.}
Equipped with the above variances control and sizes control, a direct application of Lemma~\ref{lem:concentration} yields that for each $\theta\in \Theta$ independent of the problem randomnesses, 
\begin{align}
\label{eqn:dante-peom}
	\lt| H(\Psi; \theta)\rt| 
	= 
	\lt|\sum_{i = 1}^p X_i - \sum_{i = 1}^p Y_i - \sum_{i = 1}^p Z_i\rt| 
	\lesssim \sqrt{\frac{\log n\log \frac{1}{\delta}}{n}}(\overgamma_{t+1} + \rho_F\overalpha_t)
\end{align}
with probability at least $1-\delta$.

\paragraph{Step 4: a covering argument.}

Thus far, we have established an upper bound for $\lt| H(\Psi; \theta)\rt|$ regarding any fixed $\theta \in \Theta.$ In the following, we aim to develop a control of $\lt| H(\Psi; \theta)\rt|$ uniformly over the space $\Theta$ via a standard covering argument in addition to a uniform bound. 
A direct covering of space $\Theta$ requires a covering number of order
$O\Big((\frac{1}{\epsilon})^{t^2}\Big),$
which leads to a squared dependence of $t$ in expression~\eqref{eqn:Onsager-brahms}. 
Next, we show that it is sufficient to construct a set of cardinality 
	$O\Big((\frac{1}{\epsilon})^{t\log n}\Big),$
where function $H(\Psi; \cdot)$ lies at most $\epsilon$ apart from each other. Implementing this idea leads to the right dependence of $t$ in expression~\eqref{eqn:Onsager-brahms}.

Before proceeding, let us denote $\mathcal{M}_{\epsilon}$ as an $\epsilon$-net for a subset of $\Theta$ (referred to as $\Theta_0$) where
\begin{align}
\label{eqn:Theta-0}
	\Theta_0 \defn
	 \Bigg\{ \Big(\omega, \{\alpha_k\}_{k\leq t}, \widetilde{\gamma}_t, \{\tau_k\}_{k \le t+1}\Big)
	 ~\big |
	 \text{ for every } \Big(\omega, \{\alpha_k\}_{k\leq t}, \hatgamma_t, \{\tau_k\}_{k \le t+1}\Big) \in \Theta \Bigg\}, 
\end{align}
where for each $k\leq t$, 
\begin{align}
\label{eqn:gamma-tilde}
\widetilde{\gamma}_t^k = 
	\begin{cases}
	0, & \text{for }k \le t - O(\log n),\\
	\hatgamma_t^k & \text{o.w.}
	\end{cases}
\end{align}
In words, $\Theta_0$ stands for the subset of $\Theta$ where the corresponding $\hatgamma_t$ is restricted to have zero entries except for the last $O(\log n)$ coordinates, namely, $\hatgamma_t^k = 0$ for $k \le t - O(\log n)$.
For every $\widehat{\theta} = \Big(\omega, \{\alpha_k\}_{k\leq t}, \hatgamma_t, \{\tau_k\}_{k \le t+1}\Big) \in \Theta$, assumption~\eqref{eqn:assmp-gamma-tail} ensures that $|\hatgamma_t^k| \le \epsilon$ for $k \le t - O(\log n)$.
Therefore, there exists some $\thetatilde \in \mathcal{M}_{\epsilon}$ with
\begin{align*} 
	\thetatilde \defn \Big(\widetilde{\omega}, \{\widetilde{\alpha}_k\}_{k\leq t}, \widetilde{\gamma}_t, \{\widetilde{\tau}_k\}_{k \le t+1} \Big), 
\end{align*}
such that $\|\omega - \widetilde{\omega}\|_2 \leq \epsilon$, and for every $k$, $\|\alpha_k - \widetilde{\alpha}_k\|_2 \leq \epsilon$, $\|\tau_k - \widetilde{\tau}_k\|_2 \leq \epsilon$ and $|\hatgamma_t^k - \widetilde{\gamma}_t^k| \leq \epsilon$. 
We claim that 
\begin{align}
\label{eqn:stability-H}
	|H(\Psi; \thetahat) - H(\Psi; \thetatilde)| \lesssim \frac{t\log^2 n}{n}(\overgamma_{t+1} + \rho_F\overalpha_t).
\end{align}

Let us take inequality~\eqref{eqn:stability-H} as given for the moment and come back to its proof after establishing our main result. 
It turns out that given the structure of $\Theta_0$, it is sufficient to consider a subset $\Theta^\prime_0 \subset \Theta_0$ where the corresponding $\alpha_k = 0 \in \real^k \text{ for }k \le t - O(\log n).$
More specifically, define set 
\begin{align}
	\Theta_0^{\prime} \defn \Bigg\{
	 \thetatilde^{\prime} = (\widetilde{\omega}^{\prime}, \{\widetilde{\alpha}_k^{\prime}\}_{k \le t} , \widetilde{\gamma}_t^{\prime}, \{\widetilde{\tau}_k^{\prime}\}_{k \le t+1}) \in \Theta_0
	 ~\big |
	 \text{ for every } \widetilde{\theta} = \Big(\widetilde{\omega}, \{\widetilde{\alpha}_k\}_{k\leq t}, \widetilde{\gamma}_t, \{\widetilde{\tau}_k\}_{k \le t+1} \Big) \in \mathcal{M}_\epsilon \Bigg\}, 
\end{align}
where 
\begin{gather*}
	\widetilde{\omega}^{\prime} = \widetilde{\omega}, 
	~\widetilde{\gamma}_t^{\prime} = \widetilde{\gamma}_t, 
	~\widetilde{\tau}_k^{\prime} = \widetilde{\tau}_k, ~\forall k\leq t+1,  \\
	\widetilde{\alpha}_k^{\prime} = 0 \in \real^k \text{ for }k \le t - O(\log n), 
	\quad
	\widetilde{\alpha}_k^{\prime} = \widetilde{\alpha}_k \text{ for } t - O(\log n)\leq k \leq t.
\end{gather*}
Given a $\thetatilde$ and $\thetatilde^{\prime}$ pair, since the corresponding $\widetilde{\gamma}_{t}^k = \widetilde{\gamma}_{t}^{\prime k} = 0$, we are guaranteed that   
\begin{align}
\label{eqn:dim-red-beta}
	\widetilde{\beta}_{t+1} \defn v_{t+1} + \sum_{k = 1}^{t} \widetilde{\gamma}_{t}^k F_{k}(v_{k})
	=
	v_{t+1} + \sum_{k = t - O(\log n)}^{t} \widetilde{\gamma}_{t}^{\prime k} F_{k}(v_{k})
	=
	\widetilde{\beta}_{t+1}^{\prime},
\end{align}
as $v_{k+1}$ is determined by $\alpha_{k}\in \real^k$. In addition, in view of the definition of $H$ function (cf.~\eqref{eqn:H-thetahat}), the corresponding $H$ functions have the same value as 
\begin{align}
\label{eqn:H-function-equal}
	H(\Psi; \thetatilde) = H(\Psi; \thetatilde^{\prime}).
\end{align}
These observations imply that it is enough to restrict to set $\Theta_0^{\prime}$ when consider a covering for function $H(\Psi; \theta)$.
By construction, the cardinality of $\Theta_0^{\prime}$ equals to 
\begin{align}
	|\Theta_0^{\prime}| \lesssim \Big(\frac{1}{\epsilon}\Big)^{t\log n},
\end{align}
which yields a much small size compared to $|\mathcal{M}_\epsilon|$.

Armed with the properties above, we are ready to control $\sup_{\theta \in \Theta} H(\Psi; \theta)$. Specifically, as a consequence of relation~\eqref{eqn:stability-H}, we find that 
\begin{align*}
	\sup_{\theta \in \Theta} H(\Psi; \theta) &\leq \sup_{\thetatilde \in \mathcal{M}_\epsilon} H(\Psi; \thetatilde) + C \frac{t\log^2 n}{n}(\overgamma_{t+1} + \rho_F\overalpha_t) \\
	&= \sup_{\thetatilde^{\prime} \in \Theta_0^{\prime}} H(\Psi; \thetatilde^{\prime}) + C \frac{t\log^2 n}{n}(\overgamma_{t+1} + \rho_F\overalpha_t), 
\end{align*}
for some universal constant $C$. 
Here the last equality follows from property \eqref{eqn:H-function-equal}. 
Now in order to control quantity $\sup_{\thetatilde^{\prime} \in \Theta_0^{\prime}} H(\Psi; \thetatilde^{\prime})$, recall that we have shown that for every fixed $\theta \in \Theta$, \eqref{eqn:dante-peom} holds true with probability at least $1 - \delta$.
Now setting $\delta = n^{-11} \epsilon^{t\log n}$ and taking a uniform bound over $|\Theta_0^{\prime}|$ ensure 
\begin{align}
	\sup_{\thetatilde^{\prime} \in \Theta_0^{\prime}} H(\Psi; \thetatilde^{\prime}) &\lesssim \sqrt{\frac{t\log^2 n}{n}}(\overgamma_{t+1} + \rho_F\overalpha_t),
\end{align}
and hence, 
\begin{align}
\label{eqn:elgar}
\sup_{\theta \in \Theta} H(\Psi; \theta) &\leq \sup_{\thetatilde \in \mathcal{M}_\epsilon} H(\Psi; \thetatilde) + C \frac{t\log^2 n}{n}(\overgamma_{t+1} + \rho_F\overalpha_t) \lesssim \sqrt{\frac{t\log^2 n}{n}}(\overgamma_{t+1} + \rho_F\overalpha_t),
\end{align}
with probability at least $1 - O(n^{-11}).$

\paragraph*{In summary.} 
Putting together inequality~\eqref{eq:Onsager} with \eqref{eqn:du-pre} and \eqref{eqn:elgar}, we conclude that 
\begin{align*}
&\sup_{\theta \in \Theta}~ \Bigg\|\sum_{k = 1}^{t} a_k\Big[\lt\langle \psi_k, F_{t+1}(\hatbeta_{t+1})\rt\rangle - \langle F_{t+1}^{\prime}(\hatbeta_{t+1}) \rangle \alpha_{t}^k - \sum_{j = k+1}^{t} \hatgamma_{t}^j \lt\langle F_{t+1}^{\prime}(\hatbeta_{t+1}) \circ F_{j}^{\prime}\lt(v_{j}\rt)\rt\rangle\alpha_{j-1}^k\Big]\Bigg\|_2 \notag\\
&\leq \sup_{\theta \in \Theta} ~\sum_{i = 1}^p X_i^0 + \sum_{i = 1}^p X_i - \sum_{i = 1}^p Y_i - \sum_{i = 1}^p Z_i \\
&\leq \sup_{\theta \in \Theta} ~\sum_{i = 1}^p X_i^0 + \sup_{\theta \in \Theta} ~H(\Psi; \theta) \\
&\lesssim \sqrt{\frac{t\log^2 n}{n}}(\overgamma_{t+1} + \rho_F\overalpha_t),
\end{align*}
with probability at least $1 - O(n^{-11}).$
Thus, we complete the proof of the targeted bound~\eqref{eqn:Onsager-brahms}.

\subsection{Other auxiliary details for Lemma~\ref{lem:Onsager}}

\paragraph{Proof of inequality~\eqref{eqn:stability-H}.}
Let us begin by considering quantity $\sum_{i = 1}^p (X_i - \widetilde{X}_i)$ with $X_i$ and $\widetilde{X}_i$ associated with $\thetahat$ and $\thetatilde$ respectively. 
Here $\thetahat = \Big(\omega, \{\alpha_k\}_{k\leq t}, \hatgamma_t, \{\tau_k\}_{k \le t+1}\Big)$ and $\thetatilde = \Big(\widetilde{\omega}, \{\widetilde{\alpha}_k\}_{k\leq t}, \widetilde{\gamma}_t, \{\widetilde{\tau}_k\}_{k \le t+1} \Big)$ satisfy
\begin{align*}
	\|\omega - \widetilde{\omega}\|_2 \leq \epsilon,
\end{align*}
and for every $k$,
\begin{align*}
	\|\alpha_k - \widetilde{\alpha}_k\|_2 \leq \epsilon, \quad
	\|\tau_k - \widetilde{\tau}_k\|_2 \leq \epsilon, \quad
	|\hatgamma_t^k - \widetilde{\gamma}_t^k| \leq \epsilon.
\end{align*}

We aim to show that 
\begin{align}
\label{eqn:x-diff-in-H}
\sum_{i = 1}^p (X_i - \widetilde{X}_i) \lesssim \frac{1}{\poly(n)}\rho_F(\overalpha_{t} + \overgamma_{t+1}) \asymp \frac{1}{\poly(n)}(\rho_F\overalpha_{t} + \overgamma_{t+1}).
\end{align}
As already shown by inequality \eqref{eq:cover-X0}, for every 
$(\omega, \hatgamma_t)$ and $(\widetilde{\omega}, \widetilde{\gamma}_t)$ pairs satisfying $\|\widetilde{\omega} - \omega\|_2 + \|\widetilde{\gamma}_t - \hatgamma_t\|_2 \le \epsilon  = \frac{1}{\poly(n)}$,
one has 
\begin{align}
\lt| \Big\langle \sum_{k = 1}^{t} \omega_k\psi_k, F_{t+1}\Big(\sum_{k = 1}^{t} \hatgamma_{t}^k F_{k}(0)\Big)\Big\rangle - \Big\langle \sum_{k = 1}^{t} \widetilde{\omega}_k\psi_k, F_{t+1}\Big(\sum_{k = 1}^{t} \widetilde{\gamma}_{t}^k F_{k}(0)\Big) \Big\rangle\rt|
&\lesssim \frac{1}{\poly(n)}\rho_F\overgamma_{t+1}.
\end{align}
Similarly, the Lipschitz property of function $F_t$ ensures that 
\begin{align}
\notag &\lt| \Big\langle \sum_{k = 1}^{t} \omega_k\psi_k, F_{t+1}\Big(\hatbeta_{t+1}\Big)\Big\rangle - \Big\langle \sum_{k = 1}^{t} \widetilde{\omega}_k\psi_k, F_{t+1}\Big(\widetilde{\beta}_{t+1}\Big)\Big\rangle\rt| \\
\notag &\leq \Big\|\sum_{k = 1}^{t} \omega_k\psi_k\Big\|_2 \Big\|F_{t+1}\Big(\hatbeta_{t+1}\Big) - F_{t+1}\Big(\widetilde{\beta}_{t+1}\Big)\Big\|_2
+ \Big\|F_{t+1}\Big(\widetilde{\beta}_{t+1}\Big)\Big\|_2 \Big\|\sum_{k = 1}^{t} \omega_k\psi_k - \sum_{k = 1}^{t} \widetilde{\omega}_k\psi_k \Big\|_2 \\
&\stackrel{(\mathrm{i})}{\lesssim} \frac{p}{n}\Big(1+\sqrt{\frac{t\log p}{p}}\Big) \left(\rho_F \ltwo{\hatbeta_{t+1} - \widetilde{\beta}_{t+1}} 
+ \Big\|F_{t+1}\Big(\widetilde{\beta}_{t+1}\Big)\Big\|_2
 \ltwo{\omega - \widetilde{\omega}} \right) \label{eqn:hk-egg-tart-pre} \\
&\lesssim \frac{1}{\poly(n)}\rho_F\overalpha_{t} \label{eqn:hk-egg-tart}.
\end{align}
Here $(\mathrm{i})$ follows from the event~\eqref{eqn:lemma-onsager-event}
and we leave the proof of \eqref{eqn:hk-egg-tart} to the end of this step. 
Putting these two relations together yields the claimed bound in \eqref{eqn:x-diff-in-H}.

Regarding quantity $\sum_{i = 1}^p (Y_i- \widetilde{Y}_i)$, some direct algebra shows that 
\begin{align*}
\Big|\sum_{i = 1}^p (Y_i- \widetilde{Y}_i)\Big| 
& = 
\Big|\langle F_{t+1}^{\prime}(\hatbeta_{t+1}) \rangle \cdot \widehat{\omega}^\top \widehat{\alpha}_t
- \langle F_{t+1}^{\prime}(\widetilde{\beta}_{t+1}) \rangle \cdot \widetilde{\omega}^\top \widetilde{\alpha}_t\Big| \\
&\leq 
\Big|\langle F_{t+1}^{\prime}(\hatbeta_{t+1}) \rangle\Big| \cdot |\widehat{\omega}^\top \widehat{\alpha}_t - \widetilde{\omega}^\top \widetilde{\alpha}_t|
+ \Big|\langle F_{t+1}^{\prime}(\widehat{\beta}_{t+1}) \rangle - \langle F_{t+1}^{\prime}(\widetilde{\beta}_{t+1}) \rangle\Big| \cdot|\widetilde{\omega}^\top \widetilde{\alpha}_t\Big| \\
&\leq 
\rho_F \cdot (\ltwo{\widehat{\alpha}_t - \widetilde{\alpha}_t}+\ltwo{\widehat{\omega} - \widetilde{\omega}}) 
+ \Big|\langle F_{t+1}^{\prime}(\widehat{\beta}_{t+1}) \rangle - \langle F_{t+1}^{\prime}(\widetilde{\beta}_{t+1}) \rangle\Big| \cdot\ltwo{\widetilde{\alpha}_t}.
\end{align*}
In order to bound the right hand side above, Lemma~\ref{lem:cover} provides an upper bound of quantity $\Big|\langle F_{t+1}^{\prime}(\widehat{\beta}_{t+1}) \rangle - \langle F_{t+1}^{\prime}(\widetilde{\beta}_{t+1}) \rangle\Big|$ in expression \eqref{eq:cover-Y}.
Taking this upper bound together with the fact that $\ltwo{\widehat{\alpha}_t - \widetilde{\alpha}_t}, \ltwo{\widehat{\omega} - \widetilde{\omega}} \leq \epsilon$, we obtain 
\begin{align}
\label{eqn:stab-yi}
\notag \Big|\sum_{i = 1}^p (Y_i- \widetilde{Y}_i)\Big| 
& \lesssim \rho_F \epsilon + \Big(\frac{t\log^2 n}{n}\rho_F + \frac{1}{\poly(n)}\rho_{1, F}\|\alpha_t\|_2\Big)\cdot\ltwo{\widetilde{\alpha}_t} \\
& \lesssim \frac{t\log^2 n}{n}\rho_F\overalpha_t,
\end{align}
where the last step invokes assumption~\eqref{eqn:assmp-alphat-gammat} and $\epsilon = 1/\poly(n).$

It remains to consider quantity $\sum_{i=1}^p (Z_i-\widetilde{Z}_i).$
Recalling the definition of $Z_i$ (cf.~\eqref{eqn:def-zi-bruckner}), we can write  
\begin{align}
\label{eqn:Z-decomp}
\notag &\Big|\sum_{i=1}^p (Z_i-\widetilde{Z}_i)\Big|\\
\notag =& \Big|\sum_{j = 1}^{t-1} \hatgamma_{t}^{j+1} \lt\langle F_{t+1}^{\prime}(\hatbeta_{t+1}) \circ F_{j+1}^{\prime}\lt(v_{j+1}\rt)\rt\rangle\sum_{k = 1}^{j} \widehat{\omega}_k\widehat{\alpha}_{j}^k 
-
\sum_{j = 1}^{t-1} \widetilde{\gamma}_{t}^{j+1} \lt\langle F_{t+1}^{\prime}(\widetilde{\beta}_{t+1}) \circ F_{j+1}^{\prime}\lt(v_{j+1}\rt)\rt\rangle\sum_{k = 1}^{j} (\widetilde{\omega}_k \widetilde{\alpha}_{j}^k\Big| \\
\notag = & \Big|\sum_{j = 1}^{t-1} \hatgamma_{t}^{j+1} \sum_{k = 1}^{j} \widehat{\omega}_k\widehat{\alpha}_{j}^k\Big(\lt\langle F_{t+1}^{\prime}(\hatbeta_{t+1}) \circ F_{j+1}^{\prime}\lt(v_{j+1}\rt)\rt\rangle - \lt\langle F_{t+1}^{\prime}(\widetilde{\beta}_{t+1}) \circ F_{j+1}^{\prime}\lt(v_{j+1}\rt)\rt\rangle \Big)
\\ 
\notag &\quad + 
\sum_{j = 1}^{t-1} \Big(
\widehat{\gamma}_{t}^{j+1}\sum_{k = 1}^{j} \widehat{\omega}_k \widehat{\alpha}_{j}^k - \widetilde{\gamma}_{t}^{j+1}\sum_{k = 1}^{j} \widetilde{\omega}_k \widetilde{\alpha}_{j}^k \Big)
\lt\langle F_{t+1}^{\prime}(\widetilde{\beta}_{t+1}) \circ F_{j+1}^{\prime}\lt(v_{j+1}\rt)\rt\rangle\Big| \\
\notag \leq &
\sum_{j = 1}^{t-1} |\hatgamma_{t}^{j+1}| \ltwo{\widehat{\alpha}_{j}}
\cdot \Big|\lt\langle F_{t+1}^{\prime}(\hatbeta_{t+1}) \circ F_{j+1}^{\prime}\lt(v_{j+1}\rt)\rt\rangle - \lt\langle F_{t+1}^{\prime}(\widetilde{\beta}_{t+1}) \circ F_{j+1}^{\prime}\lt(v_{j+1}\rt)\rt\rangle \Big|\\
&\quad +
\sum_{j = 1}^{t-1} \Big|
\widehat{\gamma}_{t}^{j+1}\sum_{k = 1}^{j} \widehat{\omega}_k \widehat{\alpha}_{j}^k - \widetilde{\gamma}_{t}^{j+1}\sum_{k = 1}^{j} \widetilde{\omega}_k \widetilde{\alpha}_{j}^k \Big| \cdot
\Big|\lt\langle F_{t+1}^{\prime}(\widetilde{\beta}_{t+1}) \circ F_{j+1}^{\prime}\lt(v_{j+1}\rt)\rt\rangle\Big|.
\end{align}
It is then sufficient to bound the two terms above respectively. 
First, note that 
\begin{align*}
	\Big|\widehat{\gamma}_{t}^{j+1}\sum_{k = 1}^{j} \widehat{\omega}_k \widehat{\alpha}_{j}^k - \widetilde{\gamma}_{t}^{j+1}\sum_{k = 1}^{j} \widetilde{\omega}_k \widetilde{\alpha}_{j}^k\Big|
	&=
	|\widehat{\gamma}_{t}^{j+1}| \cdot \Big|\sum_{k = 1}^{j} \widehat{\omega}_k \widehat{\alpha}_{j}^k - \sum_{k = 1}^{j} \widetilde{\omega}_k \widetilde{\alpha}_{j}^k\Big| + 
	|\sum_{k = 1}^{j} \widetilde{\omega}_k \widetilde{\alpha}_{j}^k|\cdot \Big|\widehat{\gamma}_{t}^{j+1} - \widetilde{\gamma}_{t}^{j+1}\Big| \\
	&\leq |\widehat{\gamma}_{t}^{j+1}| (\ltwo{\widehat{\alpha}_{j} - \widetilde{\alpha}_{j}}+\ltwo{\widehat{\omega} - \widetilde{\omega}}) 
	+
	\ltwo{\widetilde{\alpha}_{j}}\cdot \Big|\widehat{\gamma}_{t}^{j+1} - \widetilde{\gamma}_{t}^{j+1}\Big| \\
	&\lesssim \epsilon \cdot (|\widehat{\gamma}_{t}^{j+1}|+ \ltwo{\widetilde{\alpha}_{j}}),
\end{align*}
where the last step uses the relation between $\widehat{\theta}$ and $\widetilde{\theta}.$
Therefore, the second term of \eqref{eqn:Z-decomp} satisfies 
\begin{align}
\label{eqn:Zi-part-1}
\notag \sum_{j = 1}^{t-1} \Big|
\widehat{\gamma}_{t}^{j+1}\sum_{k = 1}^{j} \widehat{\omega}_k \widehat{\alpha}_{j}^k - \widetilde{\gamma}_{t}^{j+1}\sum_{k = 1}^{j} \widetilde{\omega}_k \widetilde{\alpha}_{j}^k \Big| \cdot
\Big|\lt\langle F_{t+1}^{\prime}(\widetilde{\beta}_{t+1}) \circ F_{j+1}^{\prime}\lt(v_{j+1}\rt)\rt\rangle\Big|
&\leq 
\frac{\rho_F^2}{\poly(n)}\ \sum_{j = 1}^{t-1}(|\widehat{\gamma}_{t}^{j+1}|+ \ltwo{\widetilde{\alpha}_{j}}) \\
&\leq \frac{\rho_F}{\poly(n)}\overalpha_t,
\end{align}
recognizing the uniform bound $|F'_k|\leq \rho_F.$

When it comes to the first term in expression \eqref{eqn:Z-decomp}, Lemma~\ref{lem:cover} controls the difference $|\langle F_{t+1}^{\prime}(\hatbeta_{t+1}) \circ F_{j+1}^{\prime}\lt(v_{j+1}\rt)\rangle - \langle F_{t+1}^{\prime}(\widetilde{\beta}_{t+1}) \circ F_{j+1}^{\prime}\lt(v_{j+1}\rt)\rangle|$ in expression~\eqref{eq:cover-Z}; 
invoking this bound, we arrive at 
\begin{align*}
&\sum_{j = 1}^{t-1} |\hatgamma_{t}^{j+1}| \ltwo{\widehat{\alpha}_{j}}
\cdot \Big|\lt\langle F_{t+1}^{\prime}(\hatbeta_{t+1}) \circ F_{j+1}^{\prime}\lt(v_{j+1}\rt)\rt\rangle - \lt\langle F_{t+1}^{\prime}(\widetilde{\beta}_{t+1}) \circ F_{j+1}^{\prime}\lt(v_{j+1}\rt)\rt\rangle \Big| \\
\lesssim &
\sum_{j = 1}^{t-1} |\hatgamma_{t}^{j+1}| \ltwo{\widehat{\alpha}_{j}}
\Big(\frac{t\log^2 n}{n}\rho_F^2 + \frac{1}{\poly(n)}\rho_F\rho_{1, F}\|\alpha_t\|_2 \Big)\\
\lesssim &
\frac{t\log^2 n}{n}\rho_F \overalpha_t
+
\frac{1}{\poly(n)}\rho_{1, F}\overalpha_t^2 \lesssim \frac{t\log^2 n}{n}\rho_F \overalpha_t,
\end{align*}
where we plug in the assumptions~\eqref{eqn:assmp-alphat-gammat} and \eqref{eqn:assmp-bach}. 
Combining pieces together, we conclude that 
\begin{align}
\label{eqn:stab-zi}
	\Big|\sum_{i=1}^p (Z_i-\widetilde{Z}_i)\Big| \lesssim \frac{t\log^2 n}{n}\rho_F \overalpha_t.
\end{align}

To conclude, by combining the three parts above in expressions~\eqref{eqn:x-diff-in-H}, \eqref{eqn:stab-yi} and \eqref{eqn:stab-zi} together, we end up with 
\begin{align*}
	|H(\Psi; \thetahat) - H(\Psi; \thetatilde)| 
	\leq 
	\Big|\sum_{i=1}^p (X_i-\widetilde{X}_i)\Big| +
	\Big|\sum_{i=1}^p (Y_i-\widetilde{Y}_i)\Big| +
	\Big|\sum_{i=1}^p (Z_i-\widetilde{Z}_i)\Big| 
	\lesssim \frac{t\log^2 n}{n}(\overgamma_{t+1} + \rho_F\overalpha_t),
\end{align*}
thus completing the inequality~\eqref{eqn:stability-H}.



\paragraph*{Proof of inequality~\eqref{eqn:hk-egg-tart}.} First, conditioning on event~\eqref{eqn:lemma-onsager-event}, we make note of the following two relations where 
\begin{align}
\label{eqn:Fk-vk-bound}
	\ltwo{F_{k}(v_{k})} \leq \ltwo{F_{k}(0)} + \rho_F \ltwo{v_k} \leq \overline{\alpha}_k + \rho_F \frac{p}{n}\Big(1+\sqrt{\frac{t\log p}{p}}\Big) \ltwo{\alpha_{k-1}} \lesssim \rho_F \overline{\alpha}_k,
\end{align}
and 
\begin{align}
\label{eqn:betahat-bound}
\notag \ltwo{\hatbeta_{t+1}} \le \|v_{t+1}\|_2 + \sum_{k = 1}^{t} \|\hatgamma_{t}^k F_{k}(v_{k})\|_2 
&\lesssim  \frac{p}{n}\Big(1+\sqrt{\frac{t\log p}{p}}\Big) \|\alpha_t\| + \sum_{k = 1}^{t} |\hatgamma_{t}^k|\|F_{k}(v_{k})\|_2 \\
&\lesssim  \frac{p}{n}\Big(1+\sqrt{\frac{t\log p}{p}}\Big)\|\alpha_t\| + \sum_{k = 1}^{t} |\hatgamma_{t}^k| \rho_F\overline{\alpha}_t 
\lesssim \overline{\alpha}_t.
\end{align}
Here inequality~\eqref{eqn:betahat-bound} holds for $\ltwo{\widetilde{\beta}_{t+1}}$ similarly. 
Combining the above two relations, we are guaranteed that 
\begin{align*}
	\lt\|F_{t+1}\Big(\widetilde{\beta}_{t+1}\Big)\rt\|_2
	\leq 
	\ltwo{F_{t+1}(0)} + \rho_F\ltwo{\widetilde{\beta}_{t+1}} 
	\leq 
	\rho_F\overline{\alpha}_k.
\end{align*}
In addition, some direct algebra together with the Lipschitz property of function $F_k$ leads to 
\begin{align}
\notag \ltwo{\hatbeta_{t+1} - \widetilde{\beta}_{t+1}} &\le \|v_{t+1} - \widetilde{v}_{t+1}\|_2 + \sum_{k = 1}^{t} \|\hatgamma_{t}^k F_{k}(v_{k}) - \widetilde{\gamma}_{t}^k F_{k}(\widetilde{v}_{k})\|_2 \\
\notag &\lesssim \frac{p}{n}\Big(1+\sqrt{\frac{t\log p}{p}}\Big)\|\alpha_t - \widetilde{\alpha}_t\|_2 + \sum_{k = 1}^{t} |\hatgamma_{t}^k - \widetilde{\gamma}_{t}^k|\|F_{k}(\widetilde{v}_{k})\|_2 + \sum_{k = 1}^{t} |\hatgamma_{t}^k| \|F_{k}(v_{k}) - F_{k}(\widetilde{v}_{k})\|_2 \\
&\lesssim \epsilon + \epsilon\sum_{k = 1}^{t} \rho_F\overline{\alpha}_k 
+ \sum_{k = 1}^{t} |\hatgamma_{t}^k| \rho_F \ltwo{v_k - \widetilde{v}_{k}}.
\label{eqn:tmp-beta-diff}
\end{align}
Conditioning on event~\eqref{eqn:lemma-onsager-event}, for each $1 \leq k \leq t$, one has 
\begin{align}
\label{eqn:v-diff}
\ltwo{v_k - \widetilde{v}_{k}}
\lesssim 
\frac{p}{n}\Big(1+\sqrt{\frac{t\log p}{p}}\Big) \|\alpha_{k-1} - \widetilde{\alpha}_{k-1}\|_2
\lesssim \epsilon,
\end{align}
which implies that \eqref{eqn:tmp-beta-diff} can be further bounded as  
\begin{align}
\ltwo{\hatbeta_{t+1} - \widetilde{\beta}_{t+1}} &\le 
\epsilon + \epsilon\sum_{k = 1}^{t} \rho_F\overline{\alpha}_k + 
\epsilon \sum_{k = 1}^{t} |\hatgamma_{t}^k| \rho_F 
\lesssim \frac{1}{\poly(n)}\rho_F\overalpha_{t}. \label{eqn:beta-diff}
\end{align}
Substituting the above relations into \eqref{eqn:hk-egg-tart-pre} establishes relation~\eqref{eqn:hk-egg-tart}.

\section{Proof of auxiliary lemmas}

\subsection{Proof of Lemma~\ref{lem:phi-psi-distr}}
\label{pf:lem:phi-psi-distr}
	
To facilitate our analysis, let us first introduce some definitions and basic properties. 
Recall that we define two sets of orthonormal basis $\{a_k\}_{1\leq k\leq \min\{n,p\}}$ and $\{b_k\}_{1\leq k\leq \min\{n,p\}}$ and for each $t$, concatenate them into orthonormal matrices
\begin{align*}
	U_{t} = [a_k]_{1 \le k \leq t} \in \real^{n\times t},
	\qquad V_{t} = [b_k]_{1 \le k \leq t} \in \real^{p\times t}.
\end{align*}
For every $1 \leq k \leq \min\{n,p\}$, we write the orthogonal complement of $U_{k}$ as $U_{k}^{\perp} \in \real^{n\times (n-k)}$ which satisfies $U_{k}^\top U_{k}^{\perp} = 0$ and $U_{k}^{\perp \top} U_{k}^{\perp} = I_{n-k}.$
Similarly, we write the orthogonal complement of $V_{k}$ as $V_{k}^{\perp} \in \real^{p \times (p-k)}.$
Additionally, we find it helpful to consider the projection where the rows of $X$ are projected to the $p-k$-dimensional space $V_{k}^{\perp}$, and the columns to the $n-k$-dimensional space $U_{k}^{\perp}$, which is denoted by 
\begin{align}
	\widetilde{X}_{k+1} \defn U_{k}^{\perp\top} X V_{k}^{\perp} \in \real^{(n-k)\times(p-k)}.
\end{align}
With these notation in place, $X_{k+1}$ (defined as in \eqref{eqn:defn-abx}) obeys 
\begin{align*}
	X_{k+1} &= \lt(I_n - a_{k}a_{k}^{\top}\rt)X_{k}\lt(I_p - b_{k}b_{k}^{\top}\rt)
	    = \cdots = \lt(I_n - U_{k}U_{k}^{\top}\rt) X \lt(I_p - V_{k}V_{k}^{\top}\rt)\\
	    &= U_{k}^{\perp}U_{k}^{\perp\top} X V_{k}^{\perp} V_{k}^{\perp\top}
	    = U_{k}^{\perp} \widetilde{X}_{k+1}  V_{k}^{\perp\top}.
\end{align*}

\begin{claim}
\label{claim:beethoven}
For every $1\leq k\leq \min\{n,p\}$, conditional on $\{a_i, b_i\}_{1\leq i\leq k}$ and $(s_0, \beta_1)$,
the following properties hold true:
\begin{itemize}
	\item $\widetilde{X}_{k+1}$ is a rescaled Wigner matrix in $\real^{(n-k)\times(p-k)}$, with $(\widetilde{X}_{k+1})_{ij} \stackrel{\text{i.i.d.}}{\sim} \mathcal{N}(0, \frac{1}{n})$; 

	\item $X_{k+1}$ is conditional independent of $\{X_{i}b_i, X_i^\top a_i\}_{1\leq i\leq k};$

	\item the randomness of $(s_{k}, \beta_{k+1})$ and $\{a_{k+1}, b_{k+1}\}$ comes purely from $\{X_{i}b_i,  X_i^\top a_i\}_{1\leq i\leq k}$, 
	and hence $(s_{k}, \beta_{k+1})$ and $\{a_{k+1}, b_{k+1}\}$ are conditionally independent of $W_{k+1}$. 
\end{itemize}	
\end{claim}
\noindent 
In view of relations~\eqref{eqn:st-xk-bk-ak} and \eqref{eqn:betat-xk-bk-ak}, the proof of Claim~\ref{claim:beethoven} proceeds by the same induction method as in the proof of \cite[Claim 1]{li2022non}. We thus omit its details here. 

Equipped with the results above, let us characterize the distribution of $W_{k}b_{k}$ for $2\leq k \leq \min\{n,p\}.$ To begin with, conditional on $\{a_i, b_i\}_{1\leq i\leq k-1}$ and $(s_0, \beta_1)$, we have 
\begin{align*}
	a_i^\top X_{k}b_{k} &= a_i^\top U_{k-1}^{\perp} \widetilde{X}_{k}  V_{k-1}^{\perp\top} b_k = 0
	\qquad \text{for } i\leq k-1; \\
	U_{k-1}^{\perp \top} X_{k}b_{k} 
	&= (U_{k-1}^{\perp \top} U_{k-1}^{\perp}) \widetilde{X}_{k}  (V_{k-1}^{\perp\top}b_{k})
	\sim 
	\mathcal{N}\left(0, \frac{1}{n}I_{n-k+1}\right),
\end{align*}
where we make use of the fact that $\widetilde{X}_{k}$ is a rescaled Wigner matrix in $\real^{(n-k+1)\times(p-k+1)}$ conditionally independent of $b_{k}$.
As a consequence, if one generates i.i.d.~Gaussian random variables $g_k^{i} \sim \mathcal{N}(0,\frac{1}{n})$ for all $1 \leq i < k$, then conditional on  $\{a_i, b_i\}_{1\leq i\leq k-1}$ and $(s_0, \beta_1)$,
it follows that 
\begin{align}
\label{eqn:phi-k-brahms}
	\phi_k &= X_kb_k + \sum_{i = 1}^{k - 1} g^i_k a_i \sim \mathcal{N}\left(0, \frac{1}{n}I_{n}\right).
\end{align}
Similarly, one can characterize the distribution of $a_k^\top X_{k}(I_p - b_kb_k^\top)$ by noticing 
\begin{align*}
	a_k^\top X_{k} (I_n - b_kb_k^\top) b_k &= 0;  \\
	a_k^\top X_{k}(I_p - b_kb_k^\top) b_{i} &= a_k^\top X_{k} b_{i} = a_k^\top U_{k-1}^{\perp} \widetilde{X}_{k}  V_{k-1}^{\perp\top} b_i  = 0 \qquad \text{for } i\leq k-1; \\
	a_k^\top X_{k}(I_p - b_kb_k^\top) V^\perp_{k} &= a_k^\top U_{k-1}^{\perp} \widetilde{X}_{k}  V_{k-1}^{\perp\top}(I_p - b_kb_k^\top) V^\perp_{k}
	=
	(a_k^\top U_{k-1}^{\perp}) \widetilde{X}_{k}  (V_{k-1}^{\perp\top} V^\perp_{k})
	\sim 
	\mathcal{N}\left(0, \frac{1}{n}I_{n-k}\right).
\end{align*}
Here the last relation follows since conditioning on $\{a_i, b_i\}_{1\leq i\leq k-1}$ and $(s_0, \beta_1)$,  $\widetilde{X}_{k}$ is a rescaled Wigner matrix in $\real^{(n-k+1)\times(p-k+1)}$ independent of $(a_{k}, b_{k})$.
It thus obeys that 
\begin{align*}
	\psi_k = \lt(I - b_{k}b_{k}^{\top}\rt)X_k^{\top}a_k + \sum_{i = 1}^{k} q^i_k b_i \sim \mathcal{N}\left(0,\frac{1}{n}I_p\right).
\end{align*}

Finally, we make the observation that $\{\phi_i\}_{1\leq i\leq k}$ are independent, so as $\{\psi_i\}_{1\leq i\leq k}$. In order to see this, first note that each $\phi_{k}$ is independent of $\{a_i, b_i\}_{1\leq i\leq k-1}$ and $(s_0, \beta_1)$ which follows immediately from the conditional distributional guarantee established in~\eqref{eqn:phi-k-brahms}. Next, putting together Claim~\ref{claim:beethoven} with the definition of $\phi_{k}$ implies that conditional on $\{a_i, b_i\}_{1\leq i\leq k-1}$ and $(s_0, \beta_1)$, $\phi_{k}$ --- whose randomness comes purely from $X_{k}b_{k}$ and $g^i_{k}$ --- is statistically independent of $\phi_{1},\ldots,\phi_{k-1}$. Again, as the distribution of $\phi_{k}$ does not relies on $\{a_i, b_i\}_{1\leq i\leq k-1}$ and $(s_0, \beta_1)$, therefore, we conclude $\{\phi_i\}_{1\leq i\leq k}$ are statistically independent. 
Similarly, one can also validate $\{\psi_i\}_{1\leq i\leq k}$ are statistically independent. 

\subsection{Proof of Lemma~\ref{lem:induction}}
\label{sec:proof-lem:induction}

To control quantity $\| \sum_{k = 1}^{t} \omega_k\psi_k \circ F_{t+1}^{\prime}(v_{t+1})\|_2^2$ with $\omega \perp \alpha_t$ (see~\eqref{eqn:orthogonal-w-a}), the idea is to invoke Lemma~\ref{lem:concentration} for any fixed $\omega \in \mathcal{S}^{t-1}$ --- independent of $v_{t+1}$
and apply a standard covering argument. 
Given any fixed $\alpha_t$ and $\omega \perp \alpha_t \in \mathcal{S}^{t-1}$, $\sum_{k = 1}^{t} \omega_k\psi_k$ follows $ \mathcal{N}(0,\frac{1}{n}I_n)$, which is independent with $v_{t+1} \defn \sum_{k = 1}^{t} \alpha_{t}^k\psi_k$.
This implies that
\begin{align*}
\mathbb{E}\Bigg[\Bigg\| \sum_{k = 1}^{t} \omega_k\psi_k \circ F_{t+1}^{\prime}(v_{t+1})\Bigg\|_2^2\Bigg] = \frac{1}{n}\mathbb{E}\Big[\lt\|F_{t+1}^{\prime}(v_{t+1})\rt\|_2^2\Big].
\end{align*}
Recognizing that $|F_{t+1}^{\prime}| \leq \rho_F$, it can be easily verified via properties for Gaussian distribution that 
\begin{align*}
	\mathbb{E}\lt[ \Big(\sum_{k = 1}^{t} \omega_k\psi_k \circ F_{t+1}^{\prime}(v_{t+1})\Big)_i^4 \rt] \lesssim \frac{1}{n^2}\rho_F^4,
\end{align*}
and 
\begin{align*}
	\mprob\left(\max_{i \in [n]}\Big(\sum_{k = 1}^{t} \omega_k\psi_k \circ F_{t+1}^{\prime}\lt(v_{t+1}\rt)\Big)_i^2 \lesssim \frac{\rho_F^2}{n}\log \frac{n}{\delta}\right) \geq 1 - \delta. 
\end{align*}
In view of Lemma~\ref{lem:concentration}, we obtain 
\begin{align}
\Bigg\| \sum_{k = 1}^{t} \omega_k\psi_k \circ F_{t+1}^{\prime}(v_{t+1})\Bigg\|_2^2 - \frac{1}{n}\mathbb{E}\Big[\lt\|F_{t+1}^{\prime}(v_{t+1})\rt\|_2^2\Big] 
\lesssim \sqrt{\frac{\log \frac{1}{\delta}}{n}}\rho_F^2 + \frac{\rho_F^2}{n}\log^2 n\log \frac{1}{\delta},
\end{align}
which holds with probability at least $1 - \delta.$
To take care of the statistical dependence between $\omega, \alpha_t$ and $\psi_{k}$, let us consider an $\epsilon$-cover of $\mathcal{S}^{t-1}$ in terms of the $\ell_2$-norm, denoted by $\mathcal{N}_\epsilon$. With this definition, we can write 
\begin{align*}
&\sup_{\omega \perp \alpha \in \mathcal{S}^{t-1}}~ \Bigg\{\Big\| \sum_{k = 1}^{t} \omega_k\psi_k \circ F_{t+1}^{\prime}(v_{t+1})\Bigg\|_2^2 - \frac{1}{n}\mathbb{E}\Big[\lt\|F_{t+1}^{\prime}(v_{t+1})\rt\|_2^2\Big] \Bigg\} \\
&\leq \sup_{\omega \perp \alpha \in \mathcal{N}_\epsilon}~ \Bigg\{\Big\| \sum_{k = 1}^{t} \omega_k\psi_k \circ F_{t+1}^{\prime}(v_{t+1})\Bigg\|_2^2 - \frac{1}{n}\mathbb{E}\Big[\lt\|F_{t+1}^{\prime}(v_{t+1})\rt\|_2^2\Big] \Bigg\} + \poly(n)\cdot\epsilon \\
%
&\leq \sup_{\omega \perp \alpha \in \mathcal{N}_\epsilon}~ \Bigg\{\Big\| \sum_{k = 1}^{t} \omega_k\psi_k \circ F_{t+1}^{\prime}(v_{t+1})\Bigg\|_2^2 - \frac{1}{n}\mathbb{E}\Big[\lt\|F_{t+1}^{\prime}(v_{t+1})\rt\|_2^2\Big] \Bigg\} + \poly(n)\cdot\epsilon\\
%
&\lesssim \sqrt{\frac{\log \Big(\frac{N(\epsilon, \mathcal{S}^{t-1})}{\delta}\Big)}{n}}\rho_F^2 + \frac{\rho_F^2}{n}\log^2 n\log \Big(\frac{N(\epsilon, \mathcal{S}^{t-1})}{\delta}\Big) + \poly(n)\cdot\epsilon,
\end{align*}
where the last inequality holds with probability $1 - \delta$ and the second inequality follows from that conditioning on the event in \eqref{eqn:simple-rm-psi}, $\|v_{t+1} - \widetilde{v}_{t+1}\|_2 \le \poly(n)\cdot\epsilon$ and $\big\|\sum_{k=1}^t (w_k - \widetilde{w}_k) \psi_k\big\|_2 \le \poly(n)\cdot\epsilon$.
Selecting parameters  
\begin{align*}
	\delta = \frac{1}{n^{10}}
	\qquad
	\text{and }~
	\epsilon = \frac{1}{\poly(n)},
\end{align*}
gives 
\begin{align*}
&\sup_{\omega \in \mathcal{S}^{t-1}}~ \Bigg\{\Big\| \sum_{k = 1}^{t} \omega_k\psi_k \circ F_{t+1}^{\prime}(v_{t+1})\Bigg\|_2^2 - \frac{1}{n}\mathbb{E}\Big[\lt\|F_{t+1}^{\prime}(v_{t+1})\rt\|_2^2\Big] \Bigg\} \lesssim \sqrt{\frac{t\log^2 n}{n}}\rho_F^2,
\end{align*}
which completes the proof of the targeted bound~\eqref{eqn:induction-con-a}.
Following similar argument, one can also derive inequality~\eqref{eqn:induction-con-b}.



\subsection{Proof of Lemma~\ref{lem:var-Xi}}
\label{sec:Pf-lem-var-Xi}

Recalling the definition in expression~\eqref{eqn:def-24601}, we begin by directly decomposing the quantity of interest as 
\begin{align}
	\notag &\sum_{i = 1}^n \mathsf{Var}(X_i^0 + X_i)  \\
	\notag &\le \mathbb{E}\lt\|\sum_{k = 1}^{t} \omega_k\psi_k \circ F_{t+1}(\hatbeta_{t+1})\rt\|_2^2 \\
	\notag &\lesssim \mathbb{E}\lt\|\sum_{k = 1}^{t} \omega_k\psi_k \circ \ind \Big(\sum_{k = 1}^{t} \omega_k\psi_k \lesssim \sqrt{\frac{\log n}{n}}\Big) \circ F_{t+1}(\hatbeta_{t+1})\rt\|_2^2 
	+ \mathbb{E}\lt\|\sum_{k = 1}^{t} \omega_k\psi_k \circ \ind \Big(\sum_{k = 1}^{t} \omega_k\psi_k \gtrsim \sqrt{\frac{\log n}{n}}\Big) \circ F_{t+1}(\hatbeta_{t+1})\rt\|_2^2 \\
	&\lesssim \frac{\log n}{n}\mathbb{E}\lt\|F_{t+1}(\hatbeta_{t+1})\rt\|_2^2 + \mathbb{E}\lt\|\sum_{k = 1}^{t} \omega_k\psi_k \circ \ind\Big(\sum_{k = 1}^{t} \omega_k\psi_k \gtrsim \sqrt{\frac{\log n}{n}}\Big) \circ F_{t+1}(\hatbeta_{t+1})\rt\|_2^2.  
	\label{eqn:bound-x-allegro}
\end{align}
Next, we control these two parts above separately. 

\begin{itemize}

\item  Regarding the first part, following by the Lipschitz property of $F_{t+1}$ and relation \eqref{eq:F0norm}, it satisfies 
\begin{align}
\label{eqn:dante}
\notag \lt\|F_{t+1}(\hatbeta_{t+1})\rt\|_2^2  
&\lesssim \Big\|F_{t+1}\Big(\sum_{k = 1}^{t} \hatgamma_{t}^k F_{k}(0)\Big)\Big\|_2^2
+ \rho_F^2 \Big\|\hatbeta_{t+1} - \sum_{k = 1}^{t} \hatgamma_{t}^k F_{k}(0)\Big\|_2^2 \\
&\lesssim \overgamma_{t+1}^2 + \rho_F^2 \Big\|\hatbeta_{t+1} - \sum_{k = 1}^{t} \hatgamma_{t}^k F_{k}(0)\Big\|_2^2.
\end{align}
Here, recall $\hatbeta_{t+1} \defn v_{t+1} + \sum_{k = 1}^{t} \hatgamma_{t}^k F_{k}(v_{k})$ and $v_{k+1} = \sum_{k = 1}^{t} \alpha_{t}^k\psi_k$ to obtain 
\begin{align}
\label{eqn:d-comedy}
\Big\|\hatbeta_{t+1} - \sum_{k = 1}^{t} \hatgamma_{t}^k F_{k}(0)\Big\|_2^2 
= \Big\|v_{t+1} + \sum_{k = 1}^{t} \hatgamma_{t}^k F_{k}(v_{k}) - \sum_{k = 1}^{t} \hatgamma_{t}^k F_{k}(0)\Big\|_2^2 
&\lesssim \lt\|v_{t+1}\rt\|_2^2 + \rho_F^2 \Big\|\sum_{k = 1}^{t} |\hatgamma_{t}^k v_{k}|\Big\|_2^2,
\end{align}
where the last inequality again invokes the Lipschitz property of $F_{k}$. 
Taking the expectation on both sides, we arrive at  
\begin{align}
	\Exs \Big\|\hatbeta_{t+1} - \sum_{k = 1}^{t} \hatgamma_{t}^k F_{k}(0)\Big\|_2^2 
	&\lesssim \|\alpha_t\|_2^2 + \rho_F^2 \sum_{i, j=1}^t|\hatgamma_{t}^i\hatgamma_{t}^j| \cdot \mathbb{E}[\lt\|v_i\rt\|_2\lt\|v_j\rt\|_2]. \label{eqn:cellist}
\end{align}
Here, again, we remind the readers that $\hatgamma_{t}$ is regarded as a fixed parameter. 
Now in order to bound the right hand side of \eqref{eqn:cellist}, since $v_{i+1} \sim \mathcal{N}(0,\frac{\ltwo{\alpha_i}^2}{n} I_p)$ for every fixed $\alpha_{i}$, it obeys that 
\begin{align*}
	\Exs \Bigg[\frac{\ltwo{v_{i+1}}}{\ltwo{\alpha_i}}\frac{\ltwo{v_{j+1}}}{\ltwo{\alpha_j}}\Bigg]
	&\leq \Exs \Big[\max\{\ltwo{X}^2, \ltwo{Y}^2\}\Big]
	\qquad \text{where } X, Y \sim \mathcal{N}\Big(0,\frac{1}{n} I_p\Big)\\
	&\leq \Exs \ltwo{X}^2 +  \Exs \ltwo{Y}^2 \lesssim \frac{p}{n}.
\end{align*}
Therefore, the right hand side of \eqref{eqn:cellist} further satisfies 
\begin{align*}
\mathbb{E}\Big\|\hatbeta_{t+1} - \sum_{k = 1}^{t} \hatgamma_{t}^k F_{k}(0)\Big\|_2^2 
&\lesssim \|\alpha_t\|_2^2 +  \frac{p\rho_F^2}{n} \sum_{i, j=0}^{t-1}|\hatgamma_{t}^i\hatgamma_{t}^j| 
\|\alpha_i\|_2\|\alpha_j\|_2  \\
&\lesssim \lt(\|\alpha_t\|_2 +  \sqrt{\frac{p}{n}}\rho_F\sum_{k = 1}^{t} \hatgamma_{t}^k \|\alpha_{k-1}\|_2\rt)^2 \lesssim \overalpha_{t}^2.
\end{align*}
Combining with \eqref{eqn:dante} ensures 
\begin{align}
\label{eqn:dante-purgatorio}
\Exs \lt\|F_{t+1}(\hatbeta_{t+1})\rt\|_2^2 \lesssim \overgamma_{t+1}^2 + \rho_F^2\overalpha_{t}^2.
\end{align}
Thus, we complete the control of the first term in \eqref{eqn:bound-x-allegro}.

\item It then suffices to control the second term, which shall again be done by means of concentration of measure. 
We claim that 
\begin{align}
\label{eqn:dante-infeno}
	\mathbb{E}\lt\|\sum_{k = 1}^{t} \omega_k\psi_k \circ \ind\lt(\sum_{k = 1}^{t} \omega_k\psi_k \gtrsim \sqrt{\frac{\log n}{n}}\rt) \circ F_{t+1}(\hatbeta_{t+1})\rt\|_2^2 \lesssim \frac{1}{\poly(n)}
	(\overgamma_{t+1}^2+\overalpha_t^2).
\end{align}
In order to see this, first, by putting together inequalities \eqref{eqn:dante}, \eqref{eqn:d-comedy} and \eqref{eq:F0norm}, we have  
\begin{align*}
	\Big\|F_{t+1}(\hatbeta_{t+1}) \Big\|_2
	&\lesssim \overgamma_{t+1} + \ltwo{v_{t+1}} + \rho_F \Big\|\sum_{k = 1}^{t} |\hatgamma_{t}^k v_{k}|\Big\|_2 \\
	&\leq \overgamma_{t+1} + \ltwo{v_{t+1}} + \rho_F \sum_{k = 1}^{t} |\hatgamma_{t}^k| \|v_{k}\|_2 
	\le \overgamma_{t+1} +  2\|\Psi\|_{\mathrm{op}} \overalpha_t,
\end{align*}
where the last inequality uses the fact that for each $k$, $v_{k+1} = \sum_{k = 1}^{t} \alpha_{t}^k\psi_k$ and $\rho_F \sum_{k = 1}^{t} |\hatgamma_{t}^k| \ll 1.$
In view of this relation, we can decompose the quantity of interest as 
\begin{align}
	\notag &\mathbb{E} \lt[\Big\|\sum_{k = 1}^{t} \omega_k\psi_k \circ \ind\Big(\sum_{k = 1}^{t} \omega_k\psi_k \gtrsim \sqrt{\frac{\log n}{n}}\Big) \circ F_{t+1}(\hatbeta_{t+1})\Big\|_2^2 \rt] \\
	\notag &\lesssim
	\overgamma_{t+1}^2 \mathbb{E}\lt[\Big\|\sum_{k = 1}^{t} \omega_k\psi_k \circ \ind\Big(\sum_{k = 1}^{t} \omega_k\psi_k \gtrsim \sqrt{\frac{\log n}{n}}\Big) \Big\|_2^2\rt]
	+
	\overalpha_t^2 
	\mathbb{E} \lt[\Big \|\sum_{k = 1}^{t} \omega_k\psi_k \circ \ind\Big(\sum_{k = 1}^{t} \omega_k\psi_k \gtrsim \sqrt{\frac{\log n}{n}}\Big)\Big\|_2^2 \cdot
	\|\Psi\|_{\mathrm{op}}^2 \rt] \\
	&\lesssim 
	\frac{\overgamma_{t+1}^2}{\poly(n)} + \overalpha_t^2 
	\mathbb{E} \lt[\Big \|\sum_{k = 1}^{t} \omega_k\psi_k \circ \ind\Big(\sum_{k = 1}^{t} \omega_k\psi_k \gtrsim \sqrt{\frac{\log n}{n}}\Big)\Big\|_2^2 \cdot
	\|\Psi\|_{\mathrm{op}}^2 \rt]. \label{eqn:cats-lots-of-them}
\end{align}
Here, the last inequality uses the property that 
\begin{align}
\label{eqn:simple-tail}
	\Exs \Bigg[X^2_i \ind\Big(X_i \gtrsim \sqrt{\frac{\log n}{n}}\Big)\Bigg]\leq \frac{1}{\poly(n)}, 
	\qquad \text{for } X_i \sim \mathcal{N}\Big(0,\frac{1}{n}\Big),
\end{align}
and given a fixed vector $\omega \in \mathcal{S}^t$, $\sum_{k = 1}^{t} \omega_k\psi_k \sim \mathcal{N}(0, \frac{1}{n}I_p)$.
We now turn to the upper bound of the second term on the right of expression~\eqref{eqn:cats-lots-of-them}.
\begin{align*}
	&\mathbb{E} \lt[\Big \|\sum_{k = 1}^{t} \omega_k\psi_k \circ \ind\Big(\sum_{k = 1}^{t} \omega_k\psi_k \gtrsim \sqrt{\frac{\log n}{n}}\Big)\Big\|_2^2 \cdot
	\|\Psi\|_{\mathrm{op}}^2 \rt] \\
	&= 
	\mathbb{E} \lt[\Big \|\sum_{k = 1}^{t} \omega_k\psi_k \circ \ind\Big(\sum_{k = 1}^{t} \omega_k\psi_k \gtrsim \sqrt{\frac{\log n}{n}}\Big)\Big\|_2^2 \cdot
	\|\Psi\|_{\mathrm{op}}^2 \ind (\|\Psi\|_{\mathrm{op}} - 1 \lesssim \sqrt{\frac{\log n}{n}})\rt] \\
	&\qquad +
	\mathbb{E} \lt[\Big \|\sum_{k = 1}^{t} \omega_k\psi_k \circ \ind\Big(\sum_{k = 1}^{t} \omega_k\psi_k \gtrsim \sqrt{\frac{t\log n}{n}}\Big)\Big\|_2^2 \cdot
	\|\Psi\|_{\mathrm{op}}^2 \ind (\|\Psi\|_{\mathrm{op}} - 1 \gtrsim \sqrt{\frac{t\log n}{n}})\rt] \\
	&\stackrel{(\mathrm{i})}{\lesssim} \frac{1}{\poly(n)} + \mathbb{E}\lt[ \|\Psi\|_{\mathrm{op}}^4 \ind (\|\Psi\|_{\mathrm{op}} - 1 \gtrsim \sqrt{\frac{t\log n}{n}})\rt] \\
	&\stackrel{(\mathrm{ii})}{\lesssim} \frac{1}{\poly(n)},
\end{align*}
where (i) results from relation~\eqref{eqn:simple-tail} and (ii) follows from the concentration result for $\|\Psi\|_{\mathrm{op}}$ (cf.~\eqref{eqn:concent-op-norm}).
\end{itemize}
Finally, combining relations \eqref{eqn:dante-purgatorio} and \eqref{eqn:dante-infeno} leads to our target bound.

\subsection{Covering lemmas}

As defined around display~\eqref{eqn:Theta-0}, for every $\thetahat = \Big(\omega, \{\alpha_k\}_{k\leq t}, \hatgamma_t, \{\tau_k\}_{k \le t+1}\Big) \in \Theta$,  $\thetatilde \defn \Big(\widetilde{\omega}, \{\widetilde{\alpha}_k\}_{k\leq t}, \widetilde{\gamma}_t, \{\widetilde{\tau}_k\}_{k \le t+1} \Big)$ is a point that lies in the $\epsilon$-cover $\mathcal{M}_\epsilon$ of $\Theta_0$ which satisfies 
\begin{align*}
	\|\omega - \widetilde{\omega}\|_2 \leq \epsilon,
	\quad
	\|\alpha_k - \widetilde{\alpha}_k\|_2 \leq \epsilon, \quad
	\|\tau_k - \widetilde{\tau}_k\|_2 \leq \epsilon, \quad
	|\hatgamma_t^k - \widetilde{\gamma}_t^k| \leq \epsilon,
\end{align*}
for every $k$ and $\epsilon = 1/\poly(n).$ We also record that 
\begin{align*}
\widetilde{\gamma}_t^k = 
	\begin{cases}
	0, & \text{for }k \le t - O(\log n),\\
	\hatgamma_t^k & \text{o.w.}
	\end{cases}
\end{align*}

\begin{lems} 
\label{lem:cover}

Under the assumptions \eqref{eqn:assmp-alphat-gammat} -- \eqref{eqn:assmp-t-range}, the following set of relations holds with probability at least $1 - O(n^{-10})$ 
\begin{itemize}
\item
\begin{subequations}
\begin{align}
\label{eq:cover-Z}
&\lt|\lt\langle F_{t+1}^{\prime}(\hatbeta_{t+1}) \circ F_{j}^{\prime}\lt(v_{j}\rt)\rt\rangle - \lt\langle F_{t+1}^{\prime}(\widetilde{\beta}_{t+1}) \circ F_{j}^{\prime}\lt(\widetilde{v}_{j}\rt)\rt\rangle\rt| 
\,\lesssim\, 
\frac{t\log^3 n}{n}\rho_F^2 + \frac{1}{\poly(n)}\rho_F\rho_{1, F}\overalpha_t, \\
&\lt|\langle F_{t+1}^{\prime}(\hatbeta_{t+1}) \rangle - \langle F_{t+1}^{\prime}(\widetilde{\beta}_{t+1}) \rangle\rt| 
\,\lesssim\, \frac{t\log^3 n}{n}\rho_F + \frac{1}{\poly(n)}\rho_{1, F}\overalpha_t; \label{eq:cover-Y}
\end{align}
\end{subequations}

\item
\begin{subequations}
\begin{align}
\label{eqn:buble}
&\lt\| \sum_{k = 1}^{t} \omega_k\psi_k \circ \lt[F_{t+1}^{\prime}(v_{t+1} + \varepsilon) - F_{t+1}^{\prime}(v_{t+1})\rt]\rt\|_2 
\,\lesssim\, \rho_{1, F}\sqrt{\frac{t\log n}{n}}\|\varepsilon\|_2 + \rho_F\Bigg(\sqrt{\frac{t\log^3 n}{n}} + \sqrt{\log n}\Big(\frac{\lt\|\varepsilon\rt\|_2}{\|\alpha_t\|_2}\Big)^{\frac{1}{3}}\Bigg), \\
\label{eqn:Fprime-o}
&\lt\|F_{t+1}^{\prime}(\hatbeta_{t+1}) - F_{t+1}^{\prime}(v_{t+1})\rt\|_1 
\,\lesssim\, \sqrt{n}\rho_{1, F}\lt\|\hatbeta_{t+1} - v_{t+1}\rt\|_2 
+ \rho_F\Bigg(t\log n + n\Big(\frac{\|\hatbeta_{t+1} - v_{t+1}\|_2}{\|\alpha_t\|_2}\Big)^{\frac{2}{3}}\Bigg), \\
\label{eqn:Fprime-a}
&\lt|\langle F_{t+1}^{\prime}(\hatbeta_{t+1}) - F_{t+1}^{\prime}(\beta_{t+1})\rangle\rt| 
\,\lesssim\, \frac{1}{\sqrt{n}}\rho_{1, F}\lt\|\hatbeta_{t+1} - \beta_{t+1}\rt\|_2 + \frac{1}{n}\rho_F\Bigg(t\log^2 n + n\Big(\frac{\|\hatbeta_{t+1} - \beta_{t+1}\|_2}{\|\alpha_t\|_2}\Big)^{\frac{2}{3}}\Bigg);
\end{align} 
\end{subequations}
\item

\begin{subequations}
\begin{align}
\label{eqn:latte}
&\lt|\lt\| \sum_{k = 1}^{t} \omega_k\psi_k \circ F_{t+1}^{\prime}(v_{t+1})\rt\|_2^2 - \lt\| \sum_{k = 1}^{t} \widetilde{\omega}_k\psi_k \circ F_{t+1}^{\prime}(\widetilde{v}_{t+1})\rt\|_2^2\rt| 
\,\lesssim\, \frac{t\log^3 n}{n}\rho_F^2 + \frac{1}{\poly(n)}\rho_{1, F}\rho_F, \\
\label{eqn:expresso}
&\lt|\lt\| F_{t+1}^{\prime}(v_{t+1})\rt\|_2^2 - \lt\| F_{t+1}^{\prime}(\widetilde{v}_{t+1})\rt\|_2^2\rt| \,\lesssim\,  \rho_F^2t\log^2 n + \frac{1}{\poly(n)}\rho_{1, F}\rho_F. 
\end{align}
\end{subequations}
\end{itemize}
\end{lems}
\noindent The proof of this lemma is provided in Section~\ref{sec:pf-lem-cov}.

\subsection{Proof of Lemma~\ref{lem:cover}}
\label{sec:pf-lem-cov}

Before diving into details, let us first describe a general framework for bounding the fluctuation of a function when its input is perturbed slightly. 
Validating each inequality of Lemma~\ref{lem:cover} then boils down to computing specific parameters in the general framework.  
Throughout this proof, we condition on the event where both \eqref{eqn:simple-rm-phi} and \eqref{eqn:simple-rm-psi} satisfy with $\delta$ selected as $\max(n,p)^{-11}$.

\subsubsection{A general framework}
Let us first set up the stage. 
The multivariate mapping and its perturbation that we are interested in are of the form 
\begin{align*}
	H(x, \theta) = \Big[c_i(x_i, \theta)h_i(u_i(x_i, \theta))\Big]_{i=1}^n
\qquad
\text{and}
\qquad
	H_{\varepsilon}(x, \theta) = \Big[c_i(x_i, \theta)h_i(u_i(x_i, \theta) + \varepsilon_i)\Big]_{i=1}^n,
\end{align*}
for perturbation vector $\varepsilon \in \real^n$ and parameter $\theta \in \real^d.$
Here $c_i$ and $u_i$ denote $\poly(n)$-Lipschitz continuous functions of $\theta$ and $h_i$ stands for functions with finite jump points. 
Specifically, consider functions $h_i$ that can be decomposed into a continuous component and a discontinuous component  
\begin{align}
	h_i(u) =  h_i^{\mathsf{cont}}(u) + h_i^{\mathsf{dis}}(u).
\end{align}
We assume the continuous part of $h_i$ is $L$-Lipschitz and the discontinuous component takes the form  
\begin{align*}
	h_i^{\mathsf{dis}}(u) := \sum_{k=1}^{M_i} s_i^k \ind(u > \tau_i^k).
\end{align*}
Here for each $i\in [n]$, we denote the discontinuous points of $h_i$ as $\{\tau_i^k\}_{k=1}^{M_i}$, and the size of their jumps as $\{s_i^k\}_{k=1}^{M_i}$.

Given every $x$ and $\varepsilon$, in order to compute the difference between $H(x, \theta)$ and $H_\varepsilon(x, \theta)$, it is critical to track where $h_i^{\mathsf{dis}}(u)$ and $h_i^{\mathsf{dis}}(u+\varepsilon)$ differ. 
For this purpose, let us define the index set  
\begin{align*}
	\mathcal{I} := \Big\{i : \ind(u_i(x_i, \theta) > \tau_i^k) \ne \ind(u_i(x_i, \theta) + \varepsilon_i > \tau_i^k)\text{ for some }k \Big\}.
\end{align*}
In words, $h_i^{\mathsf{dis}}(u_i) = h_i^{\mathsf{dis}}(u_i + \varepsilon_i)$ for all $i \in [n]$, on set $\mathcal{I}.$
In terms of this notation, the Lipschitz property of $h^{\mathsf{cont}}_i$ ensures that 
\begin{subequations}
\label{eqn:H-bound-123}
\begin{align}
	\|H(x, \theta) - H_{\varepsilon}(X, \theta)\|_1 &\lesssim \sum_{i \in \mathcal{I}} B|c_i(x_i, \theta)| + \sum_{i \notin \mathcal{I}} L|c_i(x_i, \theta)\varepsilon_i|, \label{eq:H1}\\
	\|H(x, \theta) - H_{\varepsilon}(X, \theta)\|_2^2 &\lesssim \sum_{i \in \mathcal{I}} B^2|c_i(x_i, \theta)|^2 + \sum_{i \notin \mathcal{I}} L^2|c_i(x_i, \theta)\varepsilon_i|^2, \label{eq:H2}
\end{align}
provided that $|h_i(x_i, \theta)| \lesssim B$ for every $i$.

Additionally, consider mappings  
\begin{align*}
	H^j(x, \theta) = \Big[h_i^j(u_i^j(x_i, \theta))\Big]_{i=1}^n
	\qquad \text{and} \qquad
	H_{\varepsilon}^j(x, \theta) = \Big[h_i^j(u_i^j(x_i, \theta) + \varepsilon_i^j)\Big]_{i=1}^n,
\end{align*}
for $j = 1$ or $2.$
Under the assumption $|h_i^j| \lesssim B$, in view of the Lipschitz property for the continuous part of $h_i^j$, we can conclude similarly that 
\begin{align} \label{eq:H12}
	\|H^1(x, \theta) \circ H^2(x, \theta) - H_{\varepsilon^1}^1(x, \theta) \circ H_{\varepsilon^2}^2(x, \theta)\|_1 
	\lesssim 
	\sum_{i \in \widetilde{\mathcal{I}}} B^2 + \sum_{i \notin \widetilde{\mathcal{I}}} LB(|\varepsilon_i^1| + |\varepsilon_i^2|),
\end{align}
for the index set  
\begin{align*}
	\widetilde{\mathcal{I}} := \Big\{i : \ind(u_i^j(x_i, \theta) > \tau_i^{j, k}) \ne \ind(u_i^j(x_i, \theta) + \varepsilon_i^j > \tau_i^{j, k})\text{ for some }j, k\Big\}.
\end{align*}
\end{subequations}

We shall employ these three relations above to establish Lemma~\ref{lem:cover}, which boils down to compute the right hand side of each inequality in \eqref{eqn:H-bound-123}. 
Towards this goal, the idea is to apply the concentration results developed in Section~\ref{sec:concentration}. 
Below, we state two key observations and then turn to the calculations of each inequality individually. 

Consider a random vector $X \in \real^n$. For any fixed $\theta \in \Theta$, suppose there exists some $\sigma > 0$ such that 
\begin{align}
	\mathbb{P}\lt(|u^j_i(X_i; \theta)| < \frac{s\sigma}{n}\rt) < \frac{s}{n}, \label{eq:con-jump}
\end{align}
for every $s \in [n]$ and $j = 1,2$ if there are two sets of $u_i^j$ concerned. 
In view of Lemma~\ref{lem:ind-noise} and~\eqref{eq:con-jump}, we have
\begin{align}
\label{eqn:set-size}
|\mathcal{I}| \lesssim \log N\Big(\frac{\sigma}{100n^2},\Theta\Big)\log n + \lt(\frac{n \|\varepsilon\|_2}{\sigma}\rt)^{\frac{2}{3}},
\end{align}
where $N(\frac{\sigma}{100n^2},\Theta)$ denotes the covering number of $\Theta$.

In addition, suppose that for every fixed $\theta$ and $i\in [n]$, $c_i(x_i; \theta)$ is sub-exponential with 
\begin{align}
\label{eqn:c_i-size}
	\mathbb{P}\Big(|c_i(x_i; \theta)| \le M\log\frac{1}{\delta}\Big) \ge 1 - \delta,
\end{align}
for some $M > 0$. It is easily seen that $\mathbb{E}\big[|c_i(x_i, \theta)|\big] \lesssim M$ and $\mathsf{Var}\big(|c_i(x_i, \theta)|\big) \lesssim M^2.$
Conditioning on the cardinality of $\mathcal{I}$, let us consider the quantity $\sum_{i = 1}^n w_ic_i(x_i, \theta)$ where $w \in \{0, \pm 1\}^n$ and $\|w\|_1 = |\mathcal{I}|$.
For every fixed $w$ and $\theta \in \Theta$, by virtue of Lemma~\ref{lem:concentration}, it holds true that 
\begin{align}
\label{eqn:fix-s-c}
\sum_{i = 1}^n w_ic_i(x_i, \theta) \leq M|\mathcal{I}| + \sum_{i = 1}^n w_i\big(c_i(x_i, \theta) - \mathbb{E}\big[w_i(x_i, \theta)\big]\big) \lesssim M|\mathcal{I}| + M\sqrt{|\mathcal{I}|\log\frac{1}{\delta}} +  M\log n\log\frac{1}{\delta}
\end{align}
with probability at least $1 - \delta$. 
Here we make use of the following relations  
\begin{align*}
\mathbb{E}\big[w_i(x_i, \theta)\big] \lesssim M \qquad \text{and }~
\sum_{i = 1}^n \mathsf{Var}\big(s_ic_i(x_i, \theta)\big) \lesssim M^2|\mathcal{I}|.
\end{align*}
Note that \eqref{eqn:fix-s-c} holds true for every fixed $w \in \{0, \pm 1\}^n$ and $\theta\in \Theta$. In order to accommodate the possible statistical dependences, we consider an $\epsilon$-cover of $\Theta$.
Selecting parameters  
\begin{align*}
	\delta = \frac{1}{n^{11} N(\epsilon, \Theta) 2^{|\mathcal{I}|}{n \choose |\mathcal{I}|}},
	\quad
	\epsilon = \frac{1}{\poly(n)}
\end{align*}
and taking union bound of \eqref{eqn:fix-s-c} over possible choices of $w$ and $\theta \in \Theta$ give 
\begin{align}
	\notag \sup_{\thetahat \in \Theta}~\sum_{i \in \mathcal{I}} |c_i(x_i, \thetahat)| 
	& \le \sup_{\substack{\thetahat \in \Theta, w\in \{0,\pm 1\}^m \\ \|w\|_1 = |\mathcal{I}|}} ~\sum_{i = 1}^{n} s_ic_{i}(x_i, \theta) \\
	\notag &\stackrel{(\mathrm{i})}{\le} \sup_{\substack{\theta \in \mathcal{N}_\epsilon, w\in \{0,\pm 1\}^m \\ \|w\|_1 = |\mathcal{I}|}} ~\sum_{i = 1}^{n} w_ic_{i}(x_i, \theta) + \frac{1}{\poly(n)}\\
	&\lesssim  M|\mathcal{I}|\log n + M \log N\Big(\frac{1}{\poly(n)},\Theta\Big) \log^2 n, \label{eqn:sum-ci}
\end{align}
where (i) follows from the choice of $\epsilon$ and the Lipschitz property of each $c_i$.







\subsubsection{Validating inequalities of Lemma~\ref{lem:cover}}

To validate Lemma~\ref{lem:cover}, we follow the general recipe provided above for specific choices of functions $h_i$, $u_i$ and $c_i$. 
In particular, we shall select $h_i$ as either $F_{t+1,i}^{\prime}$ or $(F_{t+1, i}^{\prime})^2$, $u_i$ as either $\hatbeta_{t+1, i}$ or $v_{t+1, i}$, and $c_i(x_i, \theta)$ as $1$, $\sum_{k = 1}^{t} \omega_k\psi_{k, i}$, or $(\sum_{k = 1}^{t} \omega_k\psi_{k, i})^2$.

As discussed above, we make note of the following observations:
\begin{itemize}
	\item Inequality set \eqref{eqn:H-bound-123} requires a uniform bound $B$ for function $h_i$, in which case, we can take $B = \rho_F$ when 
	$h_i = F_{t+1,i}^{\prime}$ and $B = \rho^2_F$ when $(F_{t+1, i}^{\prime})^2$. 

	\item For assumption~\eqref{eqn:c_i-size}, $M$ can be set as $1$, $\frac{1}{\sqrt{n}}$, and $\frac{1}{n}$, respectively;

	\item Regarding assumption~\eqref{eq:con-jump}, it suffices to select $\sigma$ parameter as $\frac{\overline{\alpha}_t}{\sqrt{n}}$ for both 
	$\hatbeta_{t+1}$ and $v_{t+1}$. We leave the proof of this fact to the end of this section. 
\end{itemize}

\paragraph*{Proof of inequality~\eqref{eq:cover-Z}.}
The idea is to apply inequality~\eqref{eq:H12} for proper choices of $H^1$ and $H^2$. 
Specifically, set
\begin{align*}
	H^1(\Psi, \theta) \defn F_{t+1}^{\prime}(\widetilde{\beta}_{t+1})
	\quad
	\text{and }
	H^2(\Psi, \theta) \defn F_{j}^{\prime}\lt(\widetilde{v}_{j}\rt),
\end{align*}
and
\begin{align*}
	H_{\varepsilon}^1(\Psi, \theta) \defn  F_{t+1}^{\prime}(\hatbeta_{t+1})
	\quad
	\text{and }
	H_{\varepsilon}^2(\Psi, \theta) \defn F_{j}^{\prime}\lt(v_{j}\rt).
\end{align*}
With these choices in mind, $\varepsilon^1 = \hatbeta_{t+1} - \widetilde{\beta}_{t+1}$, $\varepsilon^2 = v_{j} - \widetilde{v}_{j}$, and they satisfy  
\begin{align}
\label{eqn:epsilon-norm}
	\ltwo{\varepsilon^1} 
	& \lesssim \frac{1}{\poly(n)}\rho_F\overalpha_{t},
	\qquad
	\ltwo{\varepsilon^2} \lesssim \frac{1}{\poly(n)}, 
\end{align}
in view of relations \eqref{eqn:v-diff} and \eqref{eqn:beta-diff}. 
Then according to~\eqref{eq:H12}, we have
\begin{align*}
\lt|\lt\langle F_{t+1}^{\prime}(\hatbeta_{t+1}) \circ F_{j}^{\prime}\lt(v_{j}\rt)\rt\rangle - \lt\langle F_{t+1}^{\prime}(\widetilde{\beta}_{t+1}) \circ F_{j}^{\prime}\lt(\widetilde{v}_{j}\rt)\rt\rangle\rt| &\le \frac{1}{n}\lt\|F_{t+1}^{\prime}(\hatbeta_{t+1}) \circ F_{j}^{\prime}\lt(v_{j}\rt) - F_{t+1}^{\prime}(\widetilde{\beta}_{t+1}) \circ F_{j}^{\prime}\lt(\widetilde{v}_{j}\rt)\rt\|_1\\
&\lesssim \frac{1}{n}\rho_F^2 |\widetilde{\mathcal{I}}| + \frac{1}{\poly(n)}\rho_F\rho_{1, F}\overalpha_t.
\end{align*}
To establish inequality~\eqref{eq:cover-Z}, it suffices to bound the cardinality of $\widetilde{\mathcal{I}}$ which shall be done via inequality~\eqref{eqn:set-size}. 
Specifically, inequality~\eqref{eqn:set-size} requires bounding the covering number of the space $\{(\widetilde{\beta}_{t+1}, \widetilde{v}_j)\}$.

Recall the definitions 
\begin{align*}
	\widetilde{v}_j = \sum_{k=1}^{j-1} \widetilde{a}_{j-1}^i \psi_k
	\quad \text{and }~\widetilde{\beta}_{t+1} \defn \sum_{k = 1}^{t} \widetilde{\alpha}_{t}^k\psi_k + \sum_{k = 1}^{t} \widetilde{\gamma}_{t}^k F_{k}(\widetilde{v}_{k}).
\end{align*}
For every $\widetilde{\theta} = (\widetilde{\omega}, \{\widetilde{\alpha}_k\}_{k = 1}^t, \widetilde{\gamma}_t, \{\widetilde{\tau}_j\}_{j = 1}^{t+1})$, let
\begin{align*}
	\widetilde{\theta}^{\prime} = (\widetilde{\omega}^{\prime}, \{\widetilde{\alpha}_k^{\prime}\}_{k = 1}^t, \widetilde{\gamma}_t^{\prime}, \{\widetilde{\tau}_j^{\prime}\}_{j = 1}^{t+1}),
\end{align*}
where $\widetilde{\alpha}_k^{\prime} = \widetilde{\alpha}_k$ for $k > t - O(\log n)$, $\widetilde{\alpha}_k^{\prime} = 0$ for $k \le t - O(\log n)$.
It is proved in \eqref{eqn:dim-red-beta} that $\widetilde{\beta}_{t+1} = \widetilde{\beta}_{t+1}^{\prime}$. Therefore to construct a $\epsilon$-cover for space $\{(\widetilde{\beta}_{t+1}, \widetilde{v}_j)\}$, it is sufficient to consider a $\epsilon$ cover for $\widetilde{a}_{j-1}$ together with $\widetilde{\theta}^{\prime}$. The total dimension is of order $t \log n$. 
As a result, inequality~\eqref{eqn:set-size} gives 
\begin{align}
	|\widetilde{\mathcal{I}}| \lesssim t\log^3 n + \lt(\frac{n \|\varepsilon\|_2}{\sigma}\rt)^{\frac{2}{3}} \lesssim t\log^3 n. 
\end{align}
Here $\varepsilon = (\varepsilon_1, \varepsilon_2)$ satisfies inequality \eqref{eqn:epsilon-norm}, and $\sigma = \sqrt{\frac{2}{n\pi}}\overline{\alpha}_t.$
The last relation follows from assumption \eqref{eqn:assmp-alphat-gammat} that 
$\overalpha_t \leq \poly(n)$. 

Putting things together completes the proof of inequality~\eqref{eq:cover-Z}.
Inequality~\eqref{eq:cover-Y} is a direct consequence of inequality~\eqref{eq:cover-Z} by directly setting $H^2(\Psi, \theta) = H_\varepsilon^2(\Psi, \theta) = 1.$

\paragraph*{Proof of inequality~\eqref{eqn:buble}.}
In order to prove inequality~\eqref{eqn:buble}, let us take $c_i = (\sum_{k = 1}^{t} \omega_k\psi_k)_i$, 
\begin{align*}
	H(\Psi, \theta) &\defn \sum_{k = 1}^{t} \omega_k\psi_k \circ F_{t+1}^{\prime}(v_{t+1})\\
	\text{ and } H_{\varepsilon}(\Psi, \theta) &\defn \sum_{k = 1}^{t} \omega_k\psi_k \circ F_{t+1}^{\prime}(v_{t+1} + \varepsilon).
\end{align*}
Then according to~\eqref{eq:H2}, we have
\begin{align}
\label{eqn:buble-adagio}
\lt\| \sum_{k = 1}^{t} \omega_k\psi_k \circ \lt[F_{t+1}^{\prime}(v_{t+1} + \varepsilon) - F_{t+1}^{\prime}(v_{t+1})\rt]\rt\|_2^2 &\lesssim \sum_{i \in \mathcal{I}} \rho_F^2\Big|\sum_{k = 1}^{t} \omega_k\psi_{k, i}\Big|^2 + \sum_{i \notin \mathcal{I}} \rho_{1, F}^2\Big|\sum_{k = 1}^{t} \omega_k\psi_{k, i}\varepsilon_i\Big|^2.
\end{align}
We control each term of the right hand side of \eqref{eqn:buble-adagio} respectively. 
To begin with, the parameter that we shall build a $\epsilon$-cover with is $v_{t+1}$ which is determined by $\hatalpha_t \in \real^t$. 
In view of inequality~\eqref{eqn:set-size}, we have 
\begin{align}
\label{eqn:I4buble}
	|\mathcal{I}| \lesssim \log N\Big(\frac{\sigma}{100n^2},\Theta\Big) \log n + \lt(\frac{n \|\varepsilon\|_2}{\sigma}\rt)^{\frac{2}{3}} \lesssim t\log^2 n + \lt(\frac{n \|\varepsilon\|_2}{\|\hatalpha_t\|/\sqrt{n}}\rt)^{\frac{2}{3}},
\end{align}
where we recall the $\sigma$ parameter for $v_{t+1}$ equals to $\sqrt{\frac{2}{n\pi}}\overline{\alpha}_t.$
In view of the relation~\eqref{eqn:vive} in Lemma~\ref{lem:brahms-lemma}, one has 
\begin{align}
\label{eqn:basics4buble}
	\Big\|\sum_{k = 1}^{t} \omega_k\psi_k\Big\|_{\infty} \lesssim \frac{t\log n}{n},
~\text{ and }~
\sum_{i \in \mathcal{I}} \Big|\sum_{k = 1}^{t} \omega_k\psi_{k, i}\Big|^2 \lesssim \frac{(t + |\mathcal{I}|)\log n}{n},
\end{align}
with probability at least $1 - O(n^{-10}).$
Taking everything collectively, we arrive at
\begin{align*}
\lt\| \sum_{k = 1}^{t} \omega_k\psi_k \circ \lt[F_{t+1}^{\prime}(v_{t+1} + \varepsilon) - F_{t+1}^{\prime}(v_{t+1})\rt]\rt\|_2^2 
&\lesssim \rho_F^2 \lt(\frac{t\log^3 n}{n} + \log n\lt(\frac{\lt\|\varepsilon\rt\|_2}{\|\hatalpha_t\|_2}\rt)^{\frac{2}{3}}\rt) + \rho_{1, F}^2\frac{t\log n}{n}\|\varepsilon\|_2^2,
\end{align*}
from which the advertised claim in \eqref{eqn:buble} follows.  
The proofs of \eqref{eqn:Fprime-o} and \eqref{eqn:Fprime-a} can be established in the same manner, by invoking relation~\eqref{eq:H1} with $c_i = 1$ and $h_i = F'_{t+1,i}$. Here $\hatbeta_{t+1} - v_{t+1}$ and $\hatbeta_{t+1} - \beta_{t+1}$ play the role of $\varepsilon$ in these cases.

\paragraph*{Proof of inequality~\eqref{eqn:latte}.}
To establish inequality~\eqref{eqn:latte}, consider $c_i = (\sum_{k = 1}^{t} \omega_k\psi_k)^2_i$ and 
\begin{align*}
	H(\Psi, \theta) &\defn \lt( \sum_{k = 1}^{t} \omega_k\psi_k \circ F_{t+1}^{\prime}(v_{t+1})\rt)^2 \in \real^p,\\
	H_{\varepsilon}(\Psi, \theta) &\defn \lt( \sum_{k = 1}^{t} \omega_k\psi_k \circ F_{t+1}^{\prime}(\widetilde{v}_{t+1})\rt)^2 \in \real^p.
\end{align*}
By virtue of \eqref{eqn:v-diff}, $\|\varepsilon\|_2 = \|\widetilde{v}_{t+1} - v_{t+1}\|_2 \lesssim \epsilon = \frac{1}{\poly(n)}$. 
Some basic algebra leads to 
\begin{align}
	\notag \lt|\lt\| \sum_{k = 1}^{t} \omega_k\psi_k \circ F_{t+1}^{\prime}(v_{t+1})\rt\|_2^2 - \lt\| \sum_{k = 1}^{t} \omega_k\psi_k \circ F_{t+1}^{\prime}(\widetilde{v}_{t+1})\rt\|_2^2\rt| 
	&\le \lt\|\lt( \sum_{k = 1}^{t} \omega_k\psi_k \circ F_{t+1}^{\prime}(v_{t+1})\rt)^2 - \lt( \sum_{k = 1}^{t} \omega_k\psi_k \circ F_{t+1}^{\prime}(\widetilde{v}_{t+1})\rt)^2\rt\|_1 \\
	\notag &\lesssim \sum_{i \in \mathcal{I}} 
	\rho_F^2\Big|\sum_{k = 1}^{t} \omega_k\psi_{k, i}\Big|^2 
	+ \sum_{i \notin \mathcal{I}} \rho_{1, F}\rho_F\Big|\sum_{k = 1}^{t} \omega_k\psi_{k, i}\Big|^2|\varepsilon_i| 
\end{align}
where the last line follows from relation~\eqref{eq:H1}. 
Similar to the discussions around display~\eqref{eqn:I4buble}, inequality~\eqref{eqn:set-size} gives 
\begin{align*}
	|\mathcal{I}| \lesssim \log N\Big(\frac{\sigma}{100n^2},\Theta\Big) \log n + \lt(\frac{n \|\varepsilon\|_2}{\sigma}\rt)^{\frac{2}{3}} \lesssim t\log^2 n + \lt(\frac{n \|\varepsilon\|_2}{\|\hatalpha_t\|/\sqrt{n}}\rt)^{\frac{2}{3}}
	\lesssim t\log^2 n,
\end{align*}
where the last inequality invokes $\|\varepsilon\| \lesssim \epsilon = \frac{1}{\poly(n)}$. 
Taking this together with concentration bounds in~\eqref{eqn:basics4buble} further leads to 
\begin{align}
\label{eqn:latte1}
	\notag \lt|\lt\| \sum_{k = 1}^{t} \omega_k\psi_k \circ F_{t+1}^{\prime}(v_{t+1})\rt\|_2^2 - \lt\| \sum_{k = 1}^{t} \omega_k\psi_k \circ F_{t+1}^{\prime}(\widetilde{v}_{t+1})\rt\|_2^2\rt| 
	&\lesssim \frac{t\log^3 n}{n}\rho_F^2 + \frac{t^2\log^2 n}{n^2}\rho_{1, F}\rho_F\|\varepsilon\|_1 \\
	&\lesssim \frac{t\log^3 n}{n}\rho_F^2 + \frac{1}{\poly(n)}\rho_{1, F}\rho_F.
\end{align}
Contrasting the above to our target bound~\eqref{eqn:latte}, we are only left to consider replacing $\omega_k$ to $\widetilde{\omega}_k$ in the second term of the left hand side. 
Specifically, note that 
\begin{align}
\label{eqn:latte2}
\notag &\lt|\lt\| \sum_{k = 1}^{t} \omega_k\psi_k \circ F_{t+1}^{\prime}(\widetilde{v}_{t+1})\rt\|_2^2 - \lt\| \sum_{k = 1}^{t} \widetilde{\omega}_k\psi_k \circ F_{t+1}^{\prime}(\widetilde{v}_{t+1})\rt\|_2^2\rt| \\
\notag &\lt|
\Big(\sum_{k = 1}^{t} \omega_k\psi_k \circ F_{t+1}^{\prime}(\widetilde{v}_{t+1}) - \sum_{k = 1}^{t} \widetilde{\omega}_k\psi_k \circ F_{t+1}^{\prime}(\widetilde{v}_{t+1})\Big)^\top
\Big(\sum_{k = 1}^{t} \omega_k\psi_k \circ F_{t+1}^{\prime}(\widetilde{v}_{t+1}) + \sum_{k = 1}^{t} \widetilde{\omega}_k\psi_k \circ F_{t+1}^{\prime}(\widetilde{v}_{t+1})\Big)\rt|\\
\notag &\le \sqrt{p}\rho_F \Big\|\sum_{k = 1}^{t} \omega_k\psi_k\Big\|_\infty \Big\| \sum_{k = 1}^{t} (\omega_k- \widetilde{\omega}_k)\psi_k  \circ F_{t+1}^{\prime}(\widetilde{v}_{t+1})\Big\|_2\\
\notag &\lesssim \sqrt{p} \rho_F \frac{t\log n}{n} \cdot \rho_F \frac{p}{n}\Big(1+\sqrt{\frac{t\log n}{p}}\Big)\ltwo{\omega- \widetilde{\omega}} \\
&\lesssim \frac{1}{\poly(n)}\rho_F^2,
\end{align}
where for the penultimate line, recall that we condition on the event~\eqref{eqn:simple-rm-psi} with probability at least $1 - O(n^{-11})$; the last line invokes the assumption that $\ltwo{\omega - \widetilde{\omega}} \leq \epsilon.$
Combined with inequality~\eqref{eqn:latte1}, the above relation implies that 
\begin{align*}
&\lt|\lt\| \sum_{k = 1}^{t} \omega_k\psi_k \circ F_{t+1}^{\prime}(v_{t+1})\rt\|_2^2 - \lt\| \sum_{k = 1}^{t} \widetilde{\omega}_k\psi_k \circ F_{t+1}^{\prime}(\widetilde{v}_{t+1})\rt\|_2^2\rt| \\
&\le \lt|\lt\| \sum_{k = 1}^{t} \omega_k\psi_k \circ F_{t+1}^{\prime}(v_{t+1})\rt\|_2^2 - \lt\| \sum_{k = 1}^{t} \omega_k\psi_k \circ F_{t+1}^{\prime}(\widetilde{v}_{t+1})\rt\|_2^2\rt| + \frac{1}{\poly(n)}\rho_F^2 \\
&\lesssim \frac{1}{\poly(n)}\rho_{1, F}\rho_F + \frac{t\log^2 n}{n}\rho_F^2.
\end{align*}
We thus complete the proof of inequality~\eqref{eqn:latte}.
Similarly, by taking $c_i = 1$ for every $i \in [n]$, inequality~\eqref{eqn:expresso} follows by the same argument above immediately.

The remaining terms can be proved in a similar way, which is omitted here for simplicity.

\subsubsection{Other auxiliary details}

\paragraph{$\sigma$-parameter for $\hatbeta_{t+1}$ and $v_{j+1}$.}
Recall that ${v}_{j+1} = \sum_{k=1}^{j} \widehat{\alpha}_{j}^i \psi_k$ and $\hatbeta_{t+1} \defn \sum_{k = 1}^{t} \widehat{\alpha}_{t}^k\psi_k + \sum_{k = 1}^{t} \widehat{\gamma}_{t}^k F_{k}({v}_{k}).$ 
Given every fixed $\theta\in \Theta$ and $i \in [n]$, by definition, each $v_{j+1, i}$ follows $\mathcal{N}(0, \frac{\|\widehat{\alpha}_j\|_2^2}{n})$, and hence, the density function of $|v_{j+1, i}|$ is uniformly bounded by $\sqrt{\frac{2}{n\pi}}\|\widehat{\alpha}_j\|_2.$
Therefore, it is sufficient to set $\sigma = \sqrt{\frac{2}{n\pi}}\|\widehat{\alpha}_j\|_2$ for assumption~\eqref{eq:con-jump}.

Additionally, the quantity of interest $\hatbeta_{t+1}$ yields the following decomposition 
\begin{align}
\label{eqn:hatbeta-decomp}
	\hatbeta_{t+1} &= v_{t+1} + \sum_{k = 1}^{t} \hatgamma_{t}^k F_{k}(v_{k}) 
	= v_{t+1} + \sum_{k = 1}^{t} \hatgamma_{t}^k F_{k}(v_{k}^{\parallel} + v_{k}^{\perp}) \\
	\notag \text{where }\qquad 
	v_{k}^{\parallel} &= \frac{v_k^\top v_{t+1}}{\ltwo{v_{t+1}}^2}v_{t+1},
	\quad v_{k}^{\perp} = v_k - v_{k}^{\parallel},
\end{align}
where $v_{k}^{\parallel}$ denotes the component that aligns with $v_{t+1}$ while 
$v_{k}^{\perp}$ denotes the component that is orthogonal to $v_{t+1}.$
As discussed previously, given every fixed $\theta$, each ${v}_k$ follows a Gaussian distribution $\mathcal{N}(0, \frac{\|\widehat{\alpha}_j\|_2^2}{n}I_p)$, therefore $\hatbeta_{t+1}$ is a function of Gaussian vectors. 
In addition, we make the following observation that 
\begin{align*}
	\frac{|v_k^\top v_{t+1}|}{\ltwo{v_{t+1}}^2}
	&=
	\Bigg|\Big(\sum_{j=1}^{k-1} \widehat{\alpha}_{k-1}^j \psi_j\Big)^\top \sum_{j=1}^{t} \widehat{\alpha}_{t}^j \psi_j\Bigg|
	\cdot
	\Big\|\sum_{j=1}^{t} \widehat{\alpha}_{t}^j \psi_j\Big\|_2^{-2}\\
	&= 
	\Big|(\widehat{\alpha}_{k-1}^1,\ldots,\widehat{\alpha}_{k-1}^{k-1},0,\ldots,0) \Psi^\top \Psi \hatalpha_t\Big|
	\cdot 
	(\hatalpha_t^\top \Psi^\top \Psi \hatalpha_t)^{-1}\\
	&\leq 
	\frac{p}{n}\Big(1+\sqrt{\frac{t\log\frac{p}{\delta}}{p}}\Big)\ltwo{\hatalpha_{k-1}}\ltwo{\hatalpha_t}
	\cdot
	\frac{n}{p}\Big(1-\sqrt{\frac{t\log\frac{p}{\delta}}{p}}\Big)^{-1}\ltwo{\hatalpha_t}^{-2} 
	\lesssim 
	\frac{\ltwo{\hatalpha_{k-1}}}{\ltwo{\hatalpha_t}}
	\le 
	\frac{\overline{\alpha}_{k-1}}{\ltwo{\hatalpha_t}}.
\end{align*}
Recalling the assumption~\eqref{eqn:assmp-bach} that $\rho_F\sum_{k = 1}^{t} \hatgamma_{t}^k \overline{\alpha}_t \ll \|\hatalpha_t\|_2 $, 
it therefore implies that conditioning on any value of $v_{t+1}^\perp$, $\hatbeta_{t+1}$ is a Lipschitz function of $v_{t+1}$ with Lipschitz constant of order $1$. 
As a result, for every $i \in [n]$ and interval $\mathcal{I}_i$ of length $\varepsilon$, it holds that 
\begin{align*}
	\mathbb{P}_{v_{t+1}}\Big(\frac{1}{\|\alpha_{t}\|_2/\sqrt{n}} \hatbeta_{t+1, i} \in \mathcal{I}_i ~\big|~ v_{t+1}^\perp\Big) \lesssim \varepsilon,
\end{align*}
by noticing that $\hatbeta_{t+1, i}$ is a $\Theta(1)$-Lipschitz function of $v_{t+1,i}$. 
This implies that assumption~\eqref{eq:con-jump} holds with $\sigma = \|\alpha_{t}\|_2/\sqrt{n}.$

\section*{Acknowledgment}

This work was partially supported by the NSF grant DMS 2147546/2015447, the NSF CAREER award DMS-2143215, and the Google Research Scholar Award. Part of this work was done while G.~Li and Y.~Wei were visiting the Simons Institute for the Theory of Computing.

\bibliographystyle{apalike}
\bibliography{reference-amp}


\end{document}